\documentclass[english,11pt]{article}
\usepackage{setup}
\usepackage{newcommands}
\usepackage{xr}
\usepackage{amsthm,amsmath,amssymb}
\usepackage{hyperref}
\hypersetup{
	colorlinks,%
	citecolor=blue,%
	filecolor=blue,%
	linkcolor=blue,%
	urlcolor=blue
}
\usepackage{hypernat}
\usepackage[utf8x]{inputenc}
\usepackage{amsmath}
\usepackage{graphicx}
\usepackage{verbatim}
\usepackage[colorinlistoftodos]{todonotes}
\usepackage[T1]{fontenc}
\usepackage{amssymb}
\usepackage{mathtools}
\usepackage{bbm}
\usepackage{multirow}
\usepackage{float}

\title{A Cross Validation Framework for Signal Denoising with Applications to Trend Filtering, Dyadic CART and Beyond}

\author{Anamitra Chaudhuri%\thanks{Department of Statistics, UIUC} 
\and 
  Sabyasachi Chatterjee\thanks{Supported by NSF Grant DMS-1916375} \\
  \and
  Department of Statistics, University of Illinois, Urbana Champaign
  }

\date{}

\begin{document}

\maketitle

\begin{abstract}
			
			This paper formulates a general cross validation framework for signal denoising. The general framework is then applied to nonparametric regression methods such as Trend Filtering and Dyadic CART. The resulting cross validated versions are then shown to attain nearly the same rates of convergence as are known for the optimally tuned analogues. There did not exist any previous theoretical analyses of cross validated versions of Trend Filtering or Dyadic CART. To illustrate the generality of the framework we also propose and study cross validated versions of two fundamental estimators; lasso for high dimensional linear regression and singular value thresholding for matrix estimation. Our general framework is inspired by the ideas in~\cite{chatterjee2015prediction} and is potentially applicable to a wide range of estimation methods which use tuning parameters.
			
	
\end{abstract}
%\begin{keyword}
%	Cross Validation, Trend Filtering Dyadic CART, Singular Value Threshholding, Lasso, Adaptive Risk Bounds.
%\end{keyword}

%\begin{keyword}
	
%\end{keyword}
		
		%\end{frontmatter}
	
%Cross Validation, Trend Filtering Dyadic CART, Singular Value Threshholding, Lasso, Adaptive Risk Bounds.

%\input{intro}
\section{Introduction}

%This paper presents a general cross validation framework for estimators in a sequence model. The main motivation for this article comes from the fact that Trend Filtering, a relatively recently proposed (see )nonparametric regression method, yet does not appear to have a cross validated version with provable theoretical guarantees. This is in spite of several works (see )in understanding the risk behaviour of Trend Filtering, all of which assumes a certain theoretical choice (which depends on unknown problem parameters) of the tuning parameter. 

Cross Validation (CV) is a general statistical technique for choosing tuning parameters in a data driven way and is heavily used in practice for a wide variety of statistical methods. In spite of this, there is very little theoretical understanding of most CV algorithms used in practice. Within the nonparametric regression literature, rigorous theoretical guarantees for cross validated methods are limited to kernel smoothers, local linear regression methods or ridge regression (see~\cite{wong1983consistency},~\cite{shao1993linear},~\cite{li2004cross},~\cite{golub1979generalized}) which are all linear functions of the data $y$. There appears to be a need for theoretically backed general framework for building cross validation procedures for modern nonlinear regression methods. In this paper we attempt to start filling this gap in the literature by providing a general recipe to build provably adaptive and rate optimal CV estimators for some nonparametric estimation methods of current interest.

As an illustrative modern nonparametric regression method, we consider Trend Filtering (TF), proposed by~\cite{kim2009ell_1}; see~\cite{tibshirani2020divided} for a comprehensive overview. %The trend filtering estimate is defined as the minimizer of a penalized least squares criterion, in which the penalty term sums the absolute $r \geq 1$th order discrete derivatives over the design points. 
TF estimators, of order $r \geq 1$, fit $r$th degree (discrete) splines (piecewise polynomials with certain regularity). In contrast to classical nonparametric regression methods such as local polynomials, splines, kernels etc., TF is a \textit{spatially adaptive}  method as the knots of the piecewise polynomials are chosen in a data driven fashion. The last few years have seen a flurry of research (e.g., ~\cite{tibshirani2014adaptive},~\cite{guntuboyina2020adaptive},~\cite{ortelli2019prediction}) in trying to understand the theoretical properties of TF. However, all the existing guarantees hold when the tuning parameter is chosen in an optimal way depending on problem parameters which are typically unknown. On the other hand, the practical applications of TF almost always involves cross validating the tuning parameter. This motivates the following natural question. \textit{\textbf{Is it possible to define a cross validated version of Trend Filtering which provably maintains all the risk guarantees known for optimally tuned Trend Filtering?}} This is an important open question which motivates the study in this paper.

Our main focus is on developing theoretically tractable CV versions for modern fixed design nonparametric regression/signal denoising methods such as Trend Filtering, Dyadic CART, other image/matrix denoising methods, etc. Inspired by the idea underlying the cross validation method for Lasso, proposed by~\cite{chatterjee2015prediction}, we formalize a general cross validation framework for estimators in the so called \textit{sequence model}. %Our CV framework is a variant of $K$ fold cross validation. 
This framework, a variant of $K$ fold CV, provides a unified, theoretically principled and computationally efficient way to design CV versions for a variety of estimation methods. In particular, we establish a general result about any CV estimator (which fits in our framework) in Theorem~\ref{thm:main} which can then be used to obtain rate optimal guarantees for different estimators of interest.

We use this framework to propose and study a cross validated version of Trend Filtering with nearly matching theoretical guarantees known for the corresponding optimally tuned version; thereby answering our main question (in bold) posed above in the affirmative. To the best of our knowledge, before our work there has been no study done on the theoretical properties of a cross validated version of Trend Filtering. In practice, a particular CV version, implemented in the Rpackage Genlasso~\cite{arnold2020package}, is commonly used. However, no theoretical guarantees are available for this particular version. We outline the differences and similarities of our CV version with this one and present simulations which suggest that our CV version exhibits competitive finite sample performance as compared to this version.

%Our main focus is on developing theoretically tractable CV versions for modern fixed design nonparametric regression/signal denoising methods such as Trend Filtering, Dyadic CART, other image/matrix denoising methods, et cetera. In this paper we propose a general CV framework, which is a variant of $K$ fold CV. This framework provides a unified, theoretically principled and computationally efficient way to design CV versions for a variety of nonparametric estimation methods. 

%Now, generally speaking, there are several possible variations of CV (like $K$ fold, leave one out etc) that can be done for a given estimation method. For some of these methods (like Trend Filtering) there appears to be a particular CV version that is used heavily in practice (implemented in the Rpackage Genlasso~\cite{arnold2020package}). For other methods (like image denoising methods) there appears to be no one particular version that dominates but instead several different adhoc versions of CV are used in practice. In either of these cases, these CV methods do not come with rigorous theoretical justification. In this section we propose a general CV framework, which is a variant of $K$ fold CV. This framework provides a unified, theoretically principled and computationally efficient way to design CV versions for a variety of nonparametric estimation methods. 

We then use this framework to propose and study a cross validated version of Dyadic CART (DC), a classical regression tree method originally proposed in~\cite{donoho1997cart}. In a sense, DC can be thought of as an $\ell_0$ penalized version of Trend Filtering which is an (generalized) $\ell_1$ penalized least squares estimator. In~\cite{chatterjee2019adaptive}, DC has been shown to be a computationally faster and statistically competitive alternative to Trend Filtering and its multivariate versions such as the Total Variation Denoising estimator (proposed by~\cite{rudin1992nonlinear} and used heavily for image processing). This makes it natural for us to consider Dyadic CART alongside Trend Filtering in this paper. In spite of Dyadic CART being a classical nonparametric regression method and having been applied in various settings over the years; all the available theoretical results depend on a theoretical choice of the tuning parameter $\lambda$ which depends on unknown problem parameters. We again show our cross validated version is able to attain nearly the same risk bound as is known for the optimally tuned one.

Trend Filtering and Dyadic CART are the two prime examples considered in this paper where we apply our general CV framework. However, our CV framework is quite general and is potentially applicable to any other method which uses tuning parameters. To illustrate the generality and flexibility of our CV framework, we further consider two fundamental estimation methods, Singular Value Thresholding for Matrix Estimation and Lasso for high dimensional regression. We propose and study new cross validated versions of these fundamental methods. In the case of matrix estimation, we consider Singular Value Thresholding which is a canonical matrix estimation method; see~\cite{cai2010singular},~\cite{donoho2014minimax},~\cite{chatterjee2015matrix}.  We use our cross validation framework to derive a cross validated version of Singular Value Thresholding and provide rigorous adaptivity guarantees for it. Finally, for the case of the Lasso, our cross validated version can be thought of as the penalized counterpart of the estimator proposed in~\cite{chatterjee2015prediction} which cross validates constrained Lasso. %We believe our estimator has certain advantages over the one proposed in~\cite{chatterjee2015prediction}. 
We show our cross validated Lasso estimator enjoys both types of standard rates of convergence known for optimally tuned Lasso. %These are the so called slow rate and the fast rate. This is a new addition to the literature for cross validated Lasso. 
%We finally consider matrix estimation by Singular Value Thresholding which is a canonical matrix estimation method; see~\cite{cai2010singular},~\cite{donoho2014minimax},~\cite{chatterjee2015matrix}.  We use our cross validation framework to derive a cross validated version of Singular Value Thresholding and provide rigorous adaptivity guarantees for it. 

To summarize, this paper gives a general framework for cross validation and presents one general risk bound (Theorem~\ref{thm:main}) for any CV version of an estimation method which is built within our framework. Then we consider four different estimation methods, namely a) Trend Filtering, b) Dyadic CART, c) Singular Value Thresholding and d) Lasso. For each of these methods, we show how to construct a CV version within our framework. Next, we show how to apply Theorem~\ref{thm:main} to our CV versions and establish rate optimality and adaptivity which is only known for the optimally tuned analogues of these methods. Essentially, our results for these estimators look like the one below (stated informally),
\begin{theorem}[Informal]
	Let $\hat{\theta}$ be an optimally tuned estimation method (any one of the four stated above). Let $\hat{\theta}_{CV}$ be our CV version. Then, with high probability, 
	\begin{equation*}
	 MSE(\hat{\theta}_{CV},\theta^*) \leq  MSE(\hat{\theta},\theta^*) \:p(\log n)
	\end{equation*}
	where $\theta^*$ denotes the true signal, $MSE$ denotes the usual mean squared error and $p(\log n)$ is a (low degree) polynomial factor of $\log n$ where $n$ is the sample size.
	\end{theorem}

%{\color{red} change order?}

\textbf{Outline:} \textbf{Outline:} This paper is organized as follows. In Section~\ref{sec:cv} we describe and explain our cross validation framework in detail. We also give a general risk bound (see Theorem~\ref{thm:main}) for any CV estimator which falls under the scope of our framework in this section. We also provide a sketch of proof of Theorem~\ref{thm:main} in this section. In Section~\ref{sec:dc} we propose a CV version of Dyadic CART and establish an oracle risk bound for it which is only known for an optimally tuned Dyadic CART. One of the attractive aspects of Dyadic CART is fast computation and in this section we  similarly establish fast computation for our CV version by providing an algorithm in Section ~\ref{sec:compu}. In Section~\ref{sec:tf} we propose a CV version of Trend Filtering and establish both the so-called slow and fast rates known for optimally tuned Trend Filtering. 
In Sections~\ref{sec:svt} and ~\ref{sec:lasso} 
we propose CV versions of Singular Value Thresholding (SVT) and Lasso, and establish that they enjoy similar theoretical guarantees as are known for the optimally tuned versions. Section~\ref{sec:discuss} discusses some matters naturally related to the research in this article. Section~\ref{sec:simu} contains simulations done for the CV versions of Dyadic CART and Trend Filtering proposed here. 
Section ~\ref{sec:proofmain} 
contains the proof of our general risk bound (which is  Theorem~\ref{thm:main}). Sections~\ref{sec:dcproofs},~\ref{sec:tfproofs},~\ref{sec:svtproofs} and~\ref{sec:appendixlasso} 
contain the proofs of the risk bounds shown for Dyadic CART, Trend Filtering, SVT and Lasso respectively. %The proofs for Dyadic CART and Trend Filtering are somewhat lengthy. 
For the convenience of the reader, we have provided proof sketches at the beginning of Sections~\ref{sec:dcproofs} and~\ref{sec:tfproofs}.

\textbf{Notation:} Throughout the paper we use the usual $O(\cdot)$
notation to compare sequences. We write $a_n = O(b_n)$ if
there exists a constant $C > 0$ such that $a_n \leq C b_n$ for
all sufficiently large $n.$ We also use $a_n = \tilde{O}(b_n)$ to denote
$a_n = O(b_n(\log n)^C)$ for some $C > 0.$ The $O_{d}$ notation is the same as the $O$ notation except it signifies that the constant factor while comparing two sequences in $n$ may depend on the underlying dimension $d.$ For an event $A$, we will denote $\mathrm{1}(A)$ to denote the indicator random variable of the event $A.$

We use $C$ to denote a universal constant throughout the paper. This will be a positive constant independent of the problem parameters unless otherwise stated. The precise value of the constant $C$ may change from line to line. We use $[m]$ to denote the set of positive integers from $1$ to $m.$ For any vector $v \in \R^m$, we denote its $\ell_2$ norm to be $\|v\| = \sqrt{\sum_{i = 1}^{m} v_i^2}.$ Similarly, we use $\|v\|_0,\|v\|_{1}$ and $|v|_{\infty}$ to denote its $\ell_0,\ell_1$ and $\ell_{\infty}$ norms respectively. Also, for any subset $S \subset [m]$, we use $v_{S}$ to denote the vector in $\R^{|S|}$ obtained by restricting $v$ to the coordinates in $S.$ For any two vectors $v,v' \in \R^m$ we denote $\|v - v'\|^2$ by $SSE(v,v')$ where $SSE$ stands for sum of squared errors. We denote the set of all positive real numbers by $\R_{+}.$

%To give a flavor of our theoBelow we informally state

%In the next section, we describe our cross validation framework. %Such a cross validated singular value threshholding estimator appears to not exist in the current literature to the best of our knowledge. In the next section, we describe our cross validation framework. 
%{\color{red} introduce and explain the terminology completion version.}

%{\color{red} add a simulation section.}

%\input{general2}

%\documentclass[11pt]{article}

%\usepackage{amsthm,amsmath,amssymb}
%\usepackage[numbers]{natbib}

%\usepackage[colorlinks]{hyperref}
%\usepackage{hyperref}
%\hypersetup{
 %   colorlinks,%
  %  citecolor=black,%
  %  filecolor=black,%
   % linkcolor=black,%
   % urlcolor=black
%}
	\section{Cross Validation Framework}\label{sec:cv}

	%Our main focus is on developing theoretically tractable CV versions for modern fixed design nonparametric regression/signal denoising methods such as Trend Filtering, Image/Matrix denoising methods, et cetera. Now, generally speaking, there are several possible variations of CV (like $K$ fold, leave one out etc) that can be done for a given estimation method. For some of these methods (like Trend Filtering) there appears to be a particular CV version that is used heavily in practice (implemented in the Rpackage Genlasso~\cite{arnold2020package}). For other methods (like image denoising methods) there appears to be no one particular version that dominates but instead several different adhoc versions of CV are used in practice. In either of these cases, these CV methods do not come with rigorous theoretical justification. In this section we propose a general CV framework, which is a variant of $K$ fold CV. This framework provides a unified, theoretically principled and computationally efficient way to design CV versions for a variety of nonparametric estimation methods. 

	The precise setting we consider is that of signal denoising or fixed design regression where we observe $y = \theta^* + \epsilon$, where all these are $n \times 1$ vectors or vectorized matrices/tensors. $\theta^*$ is the true signal and $\epsilon$ is the noise vector consisting of i.i.d mean $0$ subgaussian noise with subgaussian norm $\sigma$. This model is sometimes called the \textit{subgaussian sequence model} and we will use the notation $y \sim Subg(\theta^*, \sigma^2)$ to mean that $y$ arises from this probabilistic model. The precise distribution of the \textit{errors} $\epsilon$ could be anything as long as subgaussianity is satisfied. The problem is to denoise or estimate the signal $\theta^*$ after observing $y.$ Many well known and popular methods to estimate $\theta^*$ in this model involve the use of tuning parameters. For such methods, we now lay out our general $K$ fold cross validation framework.

	%For any index set (ordered) $I \subseteq [n]$, and any vector $v \in \R^{n}$, $v_{I} \in \R^{|I|}$ represents the subvector of $v$ such that its $i^{th}$ element $$(v_{I})_i = v_{I_i}\ \text{for all}\ i = 1, \dots, |I|.$$
	%For any penalty function $\text{pen}(\theta)$ and tuning parameter $\lambda \geq 0$, the general estimation problem is
	%$$\argmin_{\theta \in \R^n} ||y - \theta||_2^2 + \lambda\ \text{pen}(\theta).$$
%\subsection{Notations}
%Let $[m]$ denote the set of positive integers from $1$ to $m.$ For any vector $v \in \R^m$, denote its $\ell_2$ norm to be $\|v\| = \sqrt{\sum_{i = 1}^{m} v_i^2}.$ We also denote its $\ell_1$ and $\ell_{\infty}$ norms by $\|v\|_1,\|v\|_{\infty}$ respectively. Also, for any subset $S \subset [m]$, let $v_{S}$ denote the vector in $\R^{|S|}$ obtained by restricting $v$ to the coordinates in $S.$ For any two vectors $v, v' \in \R^m$ we denote $\|v - v'\|^2$ by $SSE(v, v')$, where $SSE$ stands for sum of squared errors. Moreover, we denote $\frac{SSE(v, v')}{n}$ by $MSE(v, v')$ where $MSE$ stands for mean squared errors. We write $a_n = O(b_n)$   if there exists   constants $C>0$  and  $n_0 >0$ such that  $n \geq  n_0$ implies   that  $a_n \,\leq \,  C b_n$.   We also use the notation 
%$a_n = \tilde{O}(b_n)$ to indicate   that      $a_n    \,\leq\,  C  b_n g(\log n) $ for  $n\geq n_0$ where  $g(\cdot)$ is a polynomial function.
	
\subsection{\textbf{A General Framework of Cross Validation}}\label{sec:scheme}
%Let us set up some notations to be used throughout. Let $[m]$ denote the set of positive integers from $1$ to $m.$ For any vector $v \in \R^m$, denote its $\ell_2$ norm to be $\|v\| = \sqrt{\sum_{i = 1}^{m} v_i^2}.$ Also, for any subset $S \subset [m]$, let $v_{S}$ denote the vector in $\R^{|S|}$ obtained by restricting $v$ to the coordinates in $S.$ For any two vectors $v,v' \in \R^m$ we denote $\|v - v'\|^2$ by $SSE(v,v')$ where $SSE$ stands for sum of squared error.%Also we will abbreviate cross validation to CV throughout. 

The following $6$ general steps constitute our CV framework. This is a variant of $K \geq 2$ fold CV and is different from the traditional $K$ fold CV in some respects. Let $\hat{\theta}^{(\lambda)}$ be a given family of estimators (with tuning parameter $\lambda$) for which a CV version is desired. %For example, $\hat{\theta}^{(\lambda)}$ could be the Trend Filtering estimator or the Dyadic CART estimator with tuning $\lambda.$

\vspace{0.1cm}
\noindent\fbox{%
    \parbox{\textwidth}{%
    
\begin{enumerate}
	\item Choose the number of folds $K.$
	\item  Partition $[n]$ into $K$ disjoint index sets or folds $I_1, I_2, \dots, I_K$. We allow  this division to be done in a deterministic way or by using additional randomization. For any $j \in [K]$, denote $I^c_{j}$ to be the index set which excludes the indices in $I_j$, that is $I_{j}^c = [n] \setminus I_j.$
	%The elements of these two sets are denoted as
	%$$I = \{I_1, I_2, \dots, I_{|I|}\}, \quad I^c = \{I^c_1, I^c_2, \dots, I^c_{|I^c|}\}.$$
	\item For each $j \in [K]$ and any choice of a tuning parameter $\lambda_j > 0$, construct a version of $\hat{\theta}^{(\lambda_j)}$ which only depends on the data $y$ through the coordinates in $I_j^c$ or in other words is a function of $y_{I_j^{c}}.$ We denote this estimator by $\hat{\theta}^{(\lambda_j, I_j^c)} \in \R^n$. 
	
	%Formally,
	%	for any fixed $\lambda > 0$, $\hat{\theta}^{(\lambda, I_j^c)} = f_{\lambda}(y_{I_{j}^c})$
	%$$\hat{\theta}^{(\lambda, I_j^c)} = f_{\lambda}(y_{I_{j}^c}),\quad j \in [K],$$
	%for some measurable function $f_{\lambda}: \R^{|I_j^c|} \rightarrow \R^n$. The family of estimators $\hat{\theta}^{(\lambda, I_j^c)}$ can be suitably chosen depending on the method of estimation under consideration. 

	%For any $\lambda \geq 0$ and $j \in [K]$, define the estimate $\hat{\theta}^{(\lambda, I_j^c)} \in \R^{n}$ to be some measurable function of $y_{I_{j}^c}$ using the tuning parameter $\lambda$, which is used as a predictor of $\theta^*_{I_j}$. Formally,
	%$$\hat{\theta}^{(\lambda, I_j^c)} = f(y_{I_{j}^c}, \lambda),\quad j \in [K],$$
	%for some measurable function $f$. This predictor function $f$ can be suitably chosen depending on the problem under consideration. 
	%Similarly define $\hat{\theta}^{(\lambda, 1)} \in \R^{|I^c|}$ as
	%$$\hat{\theta}^{(\lambda, 1)} = f(y_I, \lambda).$$
	%For notational convenience, we denote the index of their elements same as the indices in the actual response vector, i.e.,
	%$$\hat{\theta}^{(\lambda, 1)} = \left(\left(\hat{\theta}^{(\lambda, 1)}_{i}\right)\right)_{i \in I^c}, \quad \hat{\theta}^{(\lambda, 2)} = \left(\left(\hat{\theta}^{(\lambda, 2)}_{i}\right)\right)_{i \in I}.$$
	\item For each $j \in [K]$, choose a finite set of possible candidate values of the tuning parameter $\lambda_j$, namely $\Lambda_j$. The set $\Lambda_j$ can be chosen deterministically or even in a data driven way as a function of $y_{I_j^{c}}$. %The simplest and generic way is to set $\Lambda_j$ to be a deterministic exponentially spaced grid $\{1,2,4,\dots,2^{K*}\}$ where $2^{K*} = \lceil \log_2(n) \rceil$; see Remark~\ref{rem:tuning}. An alternative way is to consider data driven choices of $\Lambda_j$ which can be implemented for convex optimization methods such as Trend Filtering (described in Section~\ref{sec:tf}).
	
	%A simple deterministic way is to set $\Lambda_j$ to be a exponentially spaced grid $\{1,2,4,\dots,2^{K*}$ where $2^{K*} = \lceil \log_2(n) \rceil$; see Remark~\ref{rem:tuning}. For convex optimization methods (like Trend Filtering) one can construct  data driven $\Lambda_j$ as explained later.
	
	\item For any $j \in [K]$, denote the total squared prediction error (as a function of $\lambda$) on the $j$th fold by
	$$CVERR_{j}(\lambda) = \left|\left|y_{I_j} - \hat{\theta}^{(\lambda, I_j^c)}_{I_j}\right|\right|^2.$$

	Define $\hat{\lambda}_j$ to be the candidate in $\Lambda_j$ for which the prediction error on the $j$th fold is the minimum, that is, $$\hat{\lambda}_j := \argmin_{\lambda \in \Lambda_j} CVERR_{j}(\lambda).$$
	%\begin{equation}\label{eq:ls}
	%\hat{\lambda}_j := \argmin_{\lambda \in \Lambda_j} \left|\left|y_{I_j} - \hat{\theta}^{(\lambda, I_j^c)}_{I_j}\right|\right|^2,\quad j \in [K].
	%\end{equation}
	Now define an intermediate estimator $\tilde{\theta} \in \R^n$ such that
	\begin{equation}\label{eq:defninter}
	\tilde{\theta}_{I_j} = \hat{\theta}^{(\hat{\lambda}_j, I_j^c)}_{I_j},\quad j \in [K].
	\end{equation}

	\item Define $\hat{\lambda} = \argmin_{\lambda \in \Lambda} \|\hat{\theta}^{(\lambda)} - \tilde{\theta}\|^2$ where $\Lambda$ is a deterministic set of candidate tuning parameter values to be chosen by the user. 
	Now define the final estimator $\hat{\theta}_{CV} \in \R^n$ to be %the original estimator $\hat{\theta}^{(\lambda)}$ with $\lambda = \hat{\lambda}$ as defined below: 
	\begin{equation}\label{eq:defnmain}
	\hat{\theta}_{CV} = \hat{\theta}^{(\hat{\lambda})}.
	\end{equation}
\end{enumerate}
}%
}
%{\color{red} Should we keep the notation  $\hat{\theta}_{CV}$? I think we can mention here that: later, in order to emphasize the specific problems of Lasso, Trend Filtering, Dyadic CART and Singular Value Thresholding, we denote $\hat{\theta}_{CV}$ by $\hat{\theta}_{LS}$, $\hat{\theta}_{TF}$, $\hat{\theta}_{DC}$ and $\hat{\theta}_{SV}$ respectively.}
%{\color{red} can we make a box around this?}
\vspace{0.1cm}

We now discuss about various aspects of our $K$ fold CV scheme and how it differs from the typical $K$ fold CV scheme. 
\begin{itemize}
	
	\item Dividing the dataset into folds is typically done randomly which is natural when the design is random. Since our focus is on fixed design methods or signal denoising methods our framework is a bit more general and allows deterministic folds as well. In our application to Trend Filtering we prefer using a simple deterministic strategy to create the folds. This avoids the use of extra randomization and makes our estimator conceptually simpler. In other applications such as Low Rank Matrix Estimation and Lasso it is not clear if there is a sensible way to create folds deterministically and thus we propose to create the folds randomly.

	\item In any cross validation procedure, one needs to predict on a part of the data based on observations from the rest of the data. In our framework, the way one does this is by constructing estimators which are a function of a strict subset of the data. In particular, for each fold $j \in [K]$, the user needs to define estimators $\hat{\theta}^{(\lambda, I_j^{c})} \in \R^n$ of $\theta^*$ which are functions only of $y_{I_j^{c}}$. Estimating the true signal $\theta^*$ based on only a subset of the data $y$ can be thought of as a \textit{completion} problem. Thus, we refer to the estimators $\hat{\theta}^{(\lambda, I_j^{c})}$ as \textit{completion} estimators. One can use these completion estimators to define the predictions on the $j$th fold as $\hat{\theta}^{(\lambda, I_j^{c})}_{I_j}.$ As an example, if the estimation method under consideration is Trend Filtering, then the user needs to design completion versions of Trend Filtering which are based only on a strict subset of the data. How exactly can one define these completion versions is problem specific and is described later. 
	%d) For each $j \in [K],$, the user needs to set a grid of possible values $\Lambda_j$, possibly depending on $y_{I_j^c}.$ We recommend a deterministic geometrically growing grid (described later).
	
	\item In our framework, for each fold $j \in [K]$, the user needs to build a finite set of candidate tuning values $\Lambda_j$ which is allowed to depend on $y_{I_j^c}.$ This gives quite a bit of flexibility to the user. For example, for Trend Filtering we can use a particular data driven choice of $\Lambda_j$ (see the discussion in Section~\ref{sec:alttf}). However, it should be said here that for all of our estimators, we find that setting $\Lambda_j = \Lambda = \{1,2,2^2,\dots,2^{N^*}\}$ with $N^* = O(\log n)$, a simple deterministic exponentially growing grid, is sufficient for our purposes.
	%{\color{red} Completion?}
	
	\item In traditional/typical $K$ fold CV, a single optimized tuning parameter is commonly chosen by taking $\Lambda_j = \Lambda$ and by minimizing the sum of prediction errors over all folds, that is, 
	\begin{equation}\label{eq:cvtypical}
	\hat{\lambda} = \argmin_{\lambda \in \Lambda} \sum_{j = 1}^{K} CVERR_j(\lambda).
	\end{equation}
	%$$\hat{\lambda} = \argmin_{\lambda \in \Lambda} \sum_{j = 1}^{K} CVERR_j(\lambda).$$
	Subsequently, we refer to this as the \textit{typical} or \textit{traditional} CV framework. In Section~\ref{sec:simu}, we give simulations comparing this typical CV framework to our framework where in contrast, we first construct $K$ optimized tuning parameters $\hat{\lambda}_j$ (one for each fold) which minimize the prediction error on each fold. We then construct an intermediate interleaved estimator $\tilde{\theta}$ by gluing together the optimized fits on each of the folds as in~\eqref{eq:defninter}. Finally, we then come up with a single optimized tuning parameter by minimizing the squared distance of $\hat{\theta}^{(\lambda)}$ (over a set $\Lambda$) to the intermediate fit $\tilde{\theta}$ as in~\eqref{eq:defnmain}. This seemingly roundabout way of choosing $\hat{\lambda}$ makes our cross validation scheme theoretically tractable; see our explanation in Section~\ref{sec:tract}.

%	\item In most of our applications, the set $\Lambda$ in Step $5$ is set to be an exponentially spaced finite set $\{1,2,4,\dots,2^{K^*}\}$ where $K^* = O(\log n).$ For certan convex optimization methods like Fused Lasso, $\Lambda$ can be taken to be the whole set of non negative numbers $\R^{+}$. This is because it would be possible to efficiently compute $\hat{\lambda}$ by computing the entire path of solutions $\hat{\theta}^{(\lambda)}$ over $\lambda \in \R^{+}$, {\color{red} see Section~\ref{sec:discuss} for more on this}.

	\item The main advantage of our variant of $K$ fold CV versus the traditional or typical version of $K$ fold CV is mathematical tractability (see Theorem~\ref{thm:main} below). We believe that ingredients of the theoretical analysis of our variant could be a stepping stone towards a theoretical analysis of other CV versions used in practice. Furthermore, our simulations suggest that our CV versions not only enjoy rigorous theoretical guarantees but are also practically useful, providing good finite sample performance. For example, we found that in our simulations for Trend Filtering (see Section~\ref{sec:simu}), the practical performance of our CV variant is very similar with the state of the art CV version implemented by the R package~\cite{arnold2020package}. %Our simulations suggest that our CV versions not only enjoy rigorous theoretical guarantees but are also practically useful, providing good finite sample performance.
	%The PI hopes that analyzing our CV framework would prove to be a stepping stone towards analyzing CV versions which are actually used in practice. 
\end{itemize}

\begin{comment}
\begin{itemize}
\item The user needs to set $K.$ It can be set to any small constant not less than $2.$

\item The user needs to partition the data set into $K$ folds or groups $I_1,I_2,\dots,I_K$. How to choose a good partition is going to be problem specific.
In some of our applications such as for Trend Filtering and Dyadic CART we propose to use an appropriate deterministic assignment strategy and in other applications such as Matrix Estimation and Lasso we will use a simple randomized assignment.

\item  The user needs to design a family of estimation methods using a tuning parameter $\lambda$ which is only a function of a strict subset of the data. In particular, for each $j \in [K]$, the user needs to define the estimators $\hat{\theta}^{(\lambda, I_j^{c})}$ which are functions only of $y_{I_j^{c}}$. This can be thought of as a completion problem. As an example, if the estimation method under consideration is Trend Filtering, then the user needs to design completion versions of Trend Filtering which are based only on a strict subset of the data.

\item For each $j \in [K],$, the user needs to set a grid of possible values $\Lambda_j$, possibly depending on $y_{I_j^c}.$ We can use a couple of general strategies (described later) here.
\end{itemize}
%\end{remark}
\end{comment}

\subsection{\textbf{A General Result}}

We now describe the main theoretical result underlying our cross validation framework. This result is our main tool and is used throughout the paper. This result bounds the squared error loss of the cross validated estimator $\hat{\theta}_{CV}$ defined in~\eqref{eq:defnmain}.

\begin{theorem}\label{thm:main}
%Given any deterministic or a valid randomized assignment of the index sets $I_1, I_2, \dots, I_K$, the cross validated estimator $\hat{\theta}$ satisfies for all $x \geq 0$, with probability not less than $1 - 2K \exp(-x^2/2 \sigma^2)$, the following inequality:
Let $\hat{\theta}^{(\lambda)}$ be a given family of estimators in the subgaussian sequence model for a tuning parameter $\lambda$ ranging in the set $\mathbb{R}_{+}$. Then the $K$ fold cross validated estimator $\hat{\theta}_{CV}$ defined in~\eqref{eq:defnmain} satisfies for all $x \geq 0$, with probability not less than $1 - 2K \exp(-x^2/2)$, the following inequality:
\begin{align*}
\|\hat{\theta}_{CV} - \theta^*\| \leq \min_{\lambda \in \Lambda} \|\hat{\theta}^{(\lambda)} - \theta^*\| + 4 \sum_{j \in [K]} \min_{\lambda_j \in \Lambda_j} \|\hat{\theta}^{(\lambda_j, I_j^c)}_{I_j} - \theta^*_{I_j}\| \:\:+\:\:8 \sqrt{2}  \sigma &\sum_{j \in [K]} \sqrt{\log |\Lambda_j|}\\ 
&+ 8 \sigma Kx.
\end{align*}
%Here, valid randomized assignment means that the randomization used to generate the assignment $I_1, I_2, \dots, I_K$ is independent of the data vector $y.$

\end{theorem}
%\begin{proof}
%Fix $j \in [K]$. For any arbitrary $\lambda_j \in \Lambda_j$, we have
%\begin{equation*}
%\left\|y_{I_j} - \hat{\theta}^{(\hat{\lambda}_j, I_j^c)}_{I_j}\right\|^2 \leq \left\|y_{I_j} - \hat{\theta}^{({\lambda}_j, I_j^c)}_{I_j}\right\|^2.
%\end{equation*}
%\end{theorem}
%{\color{red} A high prob version can easily be written down. Do we need that?}

%{\color{red} Add a remark about the single tuning parameter version. Maybe this can go at the discussions section at the end of the paper.}
We now explain the above theorem in more detail.
%{\color{red} Explain each of the terms above.}
\begin{itemize}

	\item The above theorem holds for all subgaussian error distributions. In the above theorem, the stated high probability event holds under the joint distribution of the errors $\epsilon$ and the (possibly) independently randomized assignment $I_1,\dots,I_k$.

\item Theorem~\ref{thm:main} bounds the root sum of squared error (RSSE) of the cross-validated estimator $\hat{\theta}_{CV}$ as a sum of four terms. The first term is $\min_{\lambda \in \Lambda} \sqrt{SSE(\hat{\theta}^{(\lambda)},\theta^*)}.$ This is basically the RSSE of the optimally tuned version of $\hat{\theta}^{(\lambda)}$ as long as $\Lambda$ is chosen to contain the theoretically optimal value of $\lambda.$ Clearly, this term is necessary as the CV version has to incur RSSE atleast as much as what is incurred by the optimally tuned version. For instance, considering the example of Trend Filtering, state of the art bounds for the RSSE are known under appropriate choices of the tuning parameter  (see~\cite{wang2016trend},~\cite{van2019prediction},~\cite{guntuboyina2020adaptive}). As long as $\Lambda$ is chosen containing these ideal choices of the tuning parameter; this term will scale exactly like the known bounds for Trend Filtering. The third term says that the dependence on the cardinality of $\Lambda_j$ in the bound in Theorem~\ref{thm:main} is logarithmic so as long as the cardinalities of $\Lambda_j$ are bounded above by a polynomial in $n$ our bound would only incur an additional $\log n$ term. It is not hard to ensure that the cardinality of $\Lambda_j$ is at most a polynomial in $n$ as will be shown in our applications. The fourth term gives a parametric $O(1/\sqrt{n})$ rate which is always going to be a lower order term.

	\item The second term appearing in the bound in Theorem~\ref{thm:main} is really the key term which arises due to cross validation. %This is simply the prediction errors $\min_{\lambda_j \in \Lambda_j} CVERR_j(\lambda_j).$
	Bounding this term becomes the central task in our applications. The second term behooves us, for each $j \in [K]$, to bound $SSE(\hat{\theta}^{(\lambda_j, I_j^c)}_{I_j}, \theta^*_{I_j})$ for \textit{\textbf{some good choice}} of the tuning parameter $\lambda_j \in \Lambda_j$. %Suppose $K = 2.$ Now, when $j = 1$, we can think of $I_2$ as the training data and the estimator $\hat{\theta}^{(\lambda_j, I_j^c)}$ is \textit{trained} on $I_2$. Moreover, $I_1$ here is the test data and the term $SSE(\hat{\theta}^{(\lambda_j, I_j^c)}_{I_j}, \theta^*_{I_j})$ can be thought of as the test error $T_1$ of $\hat{\theta}^{(\lambda_j, I_j^c)}.$ Similarly, when $j = 2$, the roles of $I_1,I_2$ as test/training data are swapped and the test error $SSE(\hat{\theta}^{(\lambda_j, I_j^c)}_{I_j}, \theta^*_{I_j})$ incurred can be denoted by $T_2.$ 
	%Theorem~\ref{thm:main} says that one needs to be able to give a good bound on the test errors $T_1$ and $T_2$ for any good choices of the tuning parameters $\lambda_1,\lambda_2$ in $\Lambda_1,\Lambda_2$ respectively. 
	%The user also needs to ensure that these good choices of $\lambda_1,\lambda_2$ lie within the finite sets $\Lambda_1,\Lambda_2$ without making their cardinalities too large. 

\item As per the earlier point, the main mathematical problem then facing us is to design completion estimators $\hat{\theta}^{(\lambda_j, I_j^c)}$ and bound the prediction errors $\min_{\lambda_j \in \Lambda_j} SSE(\hat{\theta}^{(\lambda_j, I_j^c)}_{I_j}, \theta^*_{I_j}).$ For instance, a first trivial step could be to write
$$\min_{\lambda_j \in \Lambda_j} SSE(\hat{\theta}^{(\lambda_j, I_j^c)}_{I_j}, \theta^*_{I_j}) \leq \min_{\lambda_j \in \Lambda_j} SSE(\hat{\theta}^{(\lambda_j, I_j^c)}, \theta^*)$$
where in the inequality we have just dropped the subscript $I_j.$ Now, the problem of bounding the R.H.S $\min_{\lambda_j \in \Lambda_j} SSE(\hat{\theta}^{(\lambda_j, I_j^c)}, \theta^*)$ looks similar to the problem of bounding the SSE of the original estimator $\hat{\theta}^{(\lambda)}$ with one major difference. The estimator $\hat{\theta}^{(\lambda_j, I_j^c)}$ is a completion estimator, meaning that it is a function only of $y_{I_j}^c$ in contrast with the original estimator (optimally tuned) $\hat{\theta}^{(\lambda)}$ which is based on the full data. Nevertheless, we will show that for several estimation methods, there exists a way to divide the data into folds $I_1,I_2,\dots,I_K$, design completion estimators $\hat{\theta}^{(\lambda_j, I_j)}$ for $j \in [K]$ so that the estimation errors $\min_{\lambda_j \in \Lambda_j} SSE(\hat{\theta}^{(\lambda_j, I_j^c)}, \theta^*)$ scale like the usual SSE (possibly with an extra multiplicative log factor) for the original estimation method with optimal tuning. In Sections~\ref{sec:tf},~\ref{sec:dc} we will propose some specific ways to do this for Trend Filtering and Dyadic CART respectively. A high level intuition why we can expect $\min_{\lambda_j \in \Lambda_j} SSE(\hat{\theta}^{(\lambda_j, I_j^c)}, \theta^*)$ to have same rates of convergence as $\min_{\lambda \in \Lambda} SSE(\hat{\theta}^{(\lambda)}, \theta^*)$ is the following. Observe that $\hat{\theta}^{(\lambda_j, I_j^c)}$ is based on $y_{I_j}^c$ which has on the order of $\left(n - \frac{n}{K}\right)$ data points if $I_1,I_2,\dots,I_K$ is chosen to have roughly equal size. On the other hand, $\hat{\theta}^{(\lambda)}$ is based on the full dataset with $n$ points. A good estimator based on $\left(n - \frac{n}{K}\right)$ representative samples should have the same rate of convergence as a good estimator based on all the $n$ samples with at most worse constants (since $K$ is a small constant). 
\end{itemize}

%{\color{red} change the wordings above?}

\begin{comment}
\begin{remark}
	{\color{red}Coming back to the example of Trend Filtering as an illustration, state of the art  bounds are known under appropriate choices of the tuning parameter  (see~\cite{wang2016trend},~\cite{van2019prediction},\cite{guntuboyina2020adaptive}). However, \textit{\textbf{if one can design versions of Trend Filtering which are based only on a subset of the data containing a constant fraction of the full data %such as $I_1$ or $I_2$%
			 and still maintain these known guarantees for similar appropriate choices of the tuning parameter}} contained within $\Lambda$ then we can recover these guarantees for our cross validated estimator up to an additional factor of $\sqrt{\frac{\log |\Lambda|}{n}}.$ Typically, this term will be of a lower order as $\Lambda$ can be chosen to have cardinality growing like $o(poly(n))$ as explained later. This is the typical operational implication of Theorem~\ref{thm:main}....is this needed?} 
	\end{remark}
\end{comment}

%{\color{red} I think the above remark can be shortened or compactified. In the discussion section say that this thm is useful only when $K$ is constant. So it is not useful for leave one out cross validation for instance. A different theory would be needed for that and we leave it as a future research.} 

\begin{remark}
	Theorem~\ref{thm:main} is useful only when $K$ is constant and not growing with $n$. So it is not useful for leave one out cross validation for instance where $K = n.$ A different theory would be needed for that and we leave it as a topic for future research.
	\end{remark}

\subsection{\textbf{Why is our CV Estimator theoretically tractable?}}\label{sec:tract}
In this section we explain what makes our CV Estimator theoretically tractable. In particular, we give a proof sketch of Theorem~\ref{thm:main}.

\bigskip
\textbf{Proof Sketch:}
\begin{itemize}
	\item \textbf{Step $1$:}
	A simple argument (see Lemma~\ref{lem:extra}) %~\ref{lem:extra}) 
	shows 
		\begin{equation*}
	\|\hat{\theta}_{CV} - \theta^*\| \leq 2 \|\tilde{\theta} - \theta^*\| + \min_{\lambda \in \Lambda} \|\hat{\theta}^{(\lambda)} - \theta^*\|.
	\end{equation*}
	
	Therefore, this step reduces our problem to bounding the squared error of the intermediate estimator $\tilde{\theta}.$
	
		\item \textbf{Step $2$:}
		
		We note that 
		$$\|\tilde{\theta} - \theta^*\|^2 = \sum_{j = 1}^{K}  SSE(\hat{\theta}^{(\hat{\lambda}_j, I_j^c)}_{I_j}, \theta^*_{I_j}).$$
		
		Therefore, this step reduces our problem to bounding $SSE(\hat{\theta}^{(\hat{\lambda}_j, I_j^c)}_{I_j}, \theta^*_{I_j})$ for $j = 1$ say, as the same argument can be used for all $j \in [K]$. 
	
	\item \textbf{Step $3$:}
	
	At this point, we make the crucial observation that conditionally on the (possibly random) assignment of folds $I_1,\dots,I_K$ and the noise variables $\epsilon_{I_j^{c}}$ (on all folds except the $j$th fold), the estimator $\hat{\theta}^{(\hat{\lambda}_j, I_j^c)}_{I_j}$ can be seen as a least squares estimator (see Step $5$ in Section~\ref{sec:scheme}) over the finite set $\{\hat{\theta}^{(\lambda_j, I_j^c)}_{I_j}: \lambda_j \in \Lambda_j\}.$ Note that this finite set becomes non random, once we condition on $I_1,\dots,I_K$ and $\epsilon_{I_j^{c}}.$ This is because in our framework both the estimator $\hat{\theta}^{(\lambda_j, I_j^c)}$ and the set $\Lambda_j$ can only be functions of $y_{I_j^c}$. This allows us to use an oracle risk bound for a least squares estimator over a finite set (see Lemma~\ref{lem:fls}) %~\ref{lem:fls}) 
	to conclude a conditional high probability statement for all $x \geq 0$,
	\begin{align*}
	P\Bigg(SSE(\hat{\theta}^{(\hat{\lambda}_j, I_j^c)}_{I_j}, \theta^*_{I_j}) \leq 2 \min_{\lambda_j \in \Lambda_j} \|\hat{\theta}^{(\lambda_j, I_j^c)}_{I_j} - \theta^*_{I_j}\| + 4 \sqrt{2} \sigma &\sqrt{\log |\Lambda_j|} + 4x\Bigg| I_1,,\dots,I_K,\epsilon_{I_j^{c}}\Bigg) \\&\geq  1 -  2 \exp(-x^2/2 \sigma^2). 
	\end{align*}

	\item \textbf{Step $4$:}
	Note that the probability on the R.H.S in the above conditional high probability statement does not depend on the conditioned variables. Hence, we realize that we can actually drop the conditioning in the above statement which, along with the display in Step $1$, then furnishes the statement in Theorem~\ref{thm:main}.
\end{itemize}

Let's compare with the typical/traditional $K$ fold CV framework where $\Lambda_j = \Lambda$ and the final CV estimator is $\hat{\theta}^{(\hat{\lambda})}$ where $\hat{\lambda}$ is chosen according to \eqref{eq:cvtypical}.
%\begin{equation}\label{eq:cvusual}
%\hat{\lambda} = \argmin_{\lambda \in \Lambda} \sum_{j = 1}^{K} CVERR_j(\lambda)
%\end{equation}
To the best our knowledge, there is yet no general analysis available for this typical CV framework. There does not seem to be a way to invoke a conditional least squares estimator intepretation here as we are minimizing over the sum of prediction errors over all folds at once. In general, the above method seems harder to analyze and we leave this for future research.

%The fact that we first construct $K$ optimized tuning parameters $\hat{\lambda}_j$ (one for each fold) minimizing the prediction error on each fold is critical to have this conditional least squares estimator interpretation which allows us to use least squares theory and makes the theoretical analysis tractable.

%We then construct an intermediate interleaved estimator $\tilde{\theta}$ by gluing together the optimized fits on each of the folds as in~\eqref{eq:defninter}. 
\newpage

%We then construct an intermediate interleaved estimator $\tilde{\theta}$ by gluing together the optimized fits on each of the folds as in~\eqref{eq:defninter}. 

%\input{dyadiccart}

\section{Dyadic CART}\label{sec:dc}

%\subsection{A Brief Review of Dyadic CART}
%So far we focussed on a univariate signal denoising method, namely trend filtering. It is possible to use our CV scheme for multivariate signal denoising methods as well. In this paper, we consider Dyadic CART as an illustrative multivariate signal denoising method.  
\subsection{\textbf{Background and Related Literature}}
%One class of decision trees for which an optimal tree can be computed very fast, in low to moderate dimensions, is the class of \textit{dyadic decision trees}. These trees are constructed from recursive dyadic partitioning. 
The Dyadic CART estimator is a computationally feasible decision tree method proposed first in~\cite{donoho1997cart} in the context of regression on a two-dimensional grid design. This estimator optimizes a penalized least squares criterion over the class of \textit{dyadic decision trees}. Subsequently, several papers have used ideas related to dyadic partitioning for regression, classification and density estimation, e.g., see~\cite{nowak2004estimating},~\cite{scott2006minimax},~\cite{blanchard2007optimal},~\cite{willett2007multiscale}. %Recently, the paper~\cite{chatterjee2019adaptive} generalized the Dyadic CART estimator to general dimensions and to higher orders. In this paper, we will focus on the original version of Dyadic CART which fits a rectangular piecewise constant function

The two main facts about Dyadic CART are 

\begin{itemize}
	\item The Dyadic CART estimator attains an oracle risk bound; e.g see Theorem $2.1$ in~\cite{chatterjee2019adaptive}. This oracle risk bound can then be used to show that the Dyadic CART estimator is minimax rate optimal (up to small log factors) for several function classes of interest.
	
	\item The Dyadic CART estimator can be computed very fast by a bottom up dynamic program with computational complexity linear in the sample size, see Lemma $1.1$ in~\cite{chatterjee2019adaptive}.
\end{itemize}

These two properties of the Dyadic CART make it a very attractive signal denoising method. 
However, this oracle risk bound is satisfied only when a tuning parameter is chosen to be larger than a threshhold which depends on the unknown noise variance of the error distribution. In practice, an user is naturally led to cross validate this tuning parameter. To the best of our knowledge, there has been no rigorous study done so far on Dyadic CART when the tuning parameter is chosen by cross validation. Our goal here is to propose a cross validated version of Dyadic CART in general dimensions which retain the two properties stated above. We now set up notations, define the Dyadic CART estimator more precisely and state the existing oracle risk bound. 

%The author showed that it is possible to compute this estimator by a fast bottom up dynamic program which has linear time computational complexity $O(n \times n)$ for a $n \times n$ grid. Moreover, the author showed that Dyadic CART satisfies an oracle risk bound which in turn was used to show that it is adaptively minimax rate optimal over classes of anisotropically smooth bivariate functions. Ideas in this paper were later used in~\cite{nowak2004estimating} in the context of adaptively estimating piecewise Holder smooth functions. In~\cite{chatterjee2019adaptive}, the Dyadic CART estimator has been generalized to general dimensions and of higher orders. We now precisely define the $r$th order Dyadic CART estimator in dimension $d$ which fits a polynomial of degree $r$ on each rectangle of an optimal dyadic decision tree. 

\subsection{\textbf{Notations and Definitions}}

%{\red{$L_{d,n}$ is sometimes written as $L_{n,d}$ Correct this everywhere!}}
%In this section, we formally define the problem we are addressing in this paper. 
Let us denote the $d$ dimensional lattice with $N$ points by $L_{d,n} \coloneqq \{1,\dots,n\}^d$ where 
$N = n^{d}.$ %Throughout this paper we will consider the standard fixed design setting 
%where we treat the $N$ design points as fixed and located on the $d$ dimensional grid/lattice 
%$L_{d,n}.$ One may think of the design points embedded in $[0,1]^d$ and of the form 
%$\frac{1}{n}(i_1,\dots,i_d)$ where $(i_1,\dots,i_d) \in L_{d,n}$. 
The lattice design is quite 
commonly used for theoretical studies in multidimensional nonparametric function estimation 
(see, e.g.~\cite{nemirovski2000topics}) and is also the natural setting for 
certain applications such as image denoising, matrix/tensor estimation. 
%The first step towards estimating the entire function $f$ is to be able to estimate it at the design points only. This task then reduces to estimating $d$ dimensional arrays as we now explain. %It is worth mentioning here that in all cases that we are aware of, characterization of problem difficulty such as optimal rates of convergence remain similar for both fixed design and uniformly random design. 
Letting $\theta^*$ denote the true signal, our observation 
model becomes $$y = \theta^* + \epsilon,$$
 where 
$y,\theta^*,\epsilon$ are real valued functions on $L_{d,n}$ and 
hence are $d$ dimensional arrays. Furthermore, $\epsilon$ is a noise array consisting of i.i.d subgaussian errors with an unknown subgaussian norm $\sigma > 0.$
%For an estimator $\hat{\theta}$, we will evaluate 
%its performance by the usual fixed design expected mean 
%squared error $$\MSE(\hat{\theta},\theta^*) \coloneqq 
%\frac{1}{N}\:\E_{\theta^*} \|\hat{\theta} - \theta^*\|^2.$$ 
%Here $\|.\|$ refers to the usual Euclidean norm of an array 
%where we treat an array as a vector in $\R^N.$ 

For any $a < b \in \Z_{+}$, let us define the interval of positive integers $[a,b] := \{i \in \Z_{+}: a \leq i \leq b\}$ where $\Z_{+}$ denotes 
the set of all positive integers. For a positive integer $n$ we 
also denote the set $[1,n]$ by just $[n].$ A subset $R 
\subset L_{d,n}$ is called an \textit{axis aligned 
	rectangle} if $R$ is a product of 
intervals, i.e. $R = \prod_{i = 1}^{d} [a_i,b_i].$ 
Henceforth, we will just use the word rectangle to denote an 
axis aligned rectangle. Let us define a \textit{rectangular 
	partition} of $L_{d,n}$ to be a set of rectangles 
$\mathcal{R}$ such that (a) the rectangles in $\mathcal{R}$ 
are pairwise disjoint and (b) $\cup_{R \in \mathcal{R}} R = 
L_{d,n}.$

%Recall that a multivariate polynomial of degree at most $r \geq 0$ is a finite linear combination of the monomials $\Pi_{i = 1}^{d} (x_i)^{r_i}$ satisfying $\sum_{i = 1}^{d} r_i 
%\leq r.$ It is thus clear that they form a linear space of dimension 
%$K_{r,d} \coloneqq {r + d - 1 \choose d - 1}.$ Let us now define the set of discrete multivariate polynomial arrays as 
%\begin{align*}\label{eq:polydef}
%\mathcal{F}^{(r)}_{d,n} = \big\{\theta \in \R^{L_{d,n}}: \theta(i_1,\dots,i_d) = &f(i_1/n,\dots,i_d/n)\:\:\forall (i_1,\dots,i_d) \in [n]^d \\&
%\text{for some polynomial $f$ of degree at most $r$}\big\}.
%\end{align*}

For a given rectangle $R \subset L_{d,n}$ and any $\theta \in \R^{L_{d,n}}$ let us denote the array obtained by restricting $\theta$ to $R$ by $\theta_{R}.$ %We say that $\theta$ is a degree $r$ polynomial on the rectangle $R$ if $\theta_{R} = \alpha_{R}$ for some $\alpha \in \mathcal{F}^{(r)}_{d,n}.$ 
For a given array $\theta \in \R^{L_{d,n}}$, let \textit{$k(\theta)$ denote the smallest positive integer $k$ such that a set of $k$ rectangles $R_1,\dots,R_k$ form a rectangular partition of $L_{d,n}$ and the restricted array $\theta_{R_i}$ is a constant array.} %Here, $[k]$ denotes the discrete interval of integers from $1$ to $k.$ 
In other words, $k(\theta)$ is the cardinality of the minimal rectangular partition of $L_{d,n}$ such that $\theta$ is piecewise constant on the partition.

\subsubsection{\textbf{Description of Dyadic CART}}
Let us consider a generic discrete interval $[a,b].$ We define a \textit{dyadic split} of the interval to be a split of the interval $[a,b]$ into two equal intervals. To be concrete, the interval $[a,b]$ is split into the intervals $[a,a - 1 + \ceil{(b - a + 1)/2}]$ and $[a + 
\ceil{(b - a + 1)/2}, b].$ Now consider a generic rectangle $R = \prod_{i = 1}^{d} [a_i,b_i].$ A \textit{dyadic split} of the rectangle $R$ involves the choice of a coordinate $1 \leq j \leq d$ to be split and then the $j$th interval in the product defining the rectangle $R$ undergoes a dyadic split. Thus, a dyadic split of $R$ produces two sub rectangles $R_1$ and $R_2$ where $R_2 = R \cap R_1^{c}$ and $R_1$ is of the following form for some $j \in [d]$,
\begin{equation*}
R_1 = \prod_{i = 1}^{j - 1} [a_i,b_i] \times [a_j ,a_j - 1 + \ceil{(b_j - a_j + 1) / 2}] \times \prod_{i = j + 1}^{d} [a_i,b_i].
\end{equation*}

Starting from the trivial partition which is just $L_{d,n}$ itself, %we can now recursively create more and more refined partitions of $L_{n,d}$ by dyadically splitting rectangles. For instance in the first step, 
we can create a refined partition by dyadically splitting $L_{d,n}.$  This will result in a partition of $L_{d,n}$ into two rectangles. We can now keep on dividing recursively, generating new partitions. In general, if at some stage we have the partition $\Pi = (R_1,\dots,R_k)$, we can choose any of the rectangles $R_i$ and dyadically split it to get a refinement of $\Pi$ with $k + 1$ nonempty rectangles. \textit{A recursive dyadic partition} (RDP) is any partition reachable by such successive dyadic splitting. Let us denote the set of all recursive dyadic partitions of $L_{d,n}$ as $\mathcal{P}_{\rdp,d,n}.$ Indeed, a natural way of encoding any RDP of $L_{d,n}$ is by a binary tree where each nonleaf node is labeled by an integer in $[d].$ This labeling corresponds to the choice of the coordinate that was used for the split. %Let us denote the set of all RDP of $L_{d,n}$ by $\mathcal{P}_{\rdp,n,d}.$

%This can be seen as follows. Let the root represent the full set $L_{n,d}.$ If the first step of partitioning is done by dividing in half along coordinate $i$ then label the root vertex by $i \in [d].$ The two children of the root now represent the subsets of $L_{n,d}$ given by $[n]^{i - 1} \times [n/2] \times [n]^{d - i}$ and its complement. Depending upon the coordinate of the next split, these vertices can now also be labelled. Let us denote the set of all RDP of $L_{d,n}$ by $\mathcal{P}_{\rdp,n,d}.$

%We can now consider a constrained version of $\hat{\theta}^{(r)}_{\all,\lambda}$ which only optimizes over $\mathcal{P}_{\rdp,d,n}$ instead of optimizing over $\mathcal{P}_{\all,d,n}.$ 

For a given array $\theta \in \R^{L_{d,n}}$, let \textit{$k_{\rdp}(\theta)$ denote the smallest positive integer $k$ such that a set of $k$ rectangles $R_1,\dots,R_k$ form a recursive dyadic partition of $L_{d,n}$ and the restricted array $\theta_{R_i}$ is a constant array for all $1 \leq i \leq k.$} %Here, $[k]$ denotes the discrete interval of integers from $1$ to $k.$ 
In other words, $k_{\rdp}(\theta)$ is the cardinality of the minimal recursive dyadic partition of $L_{d,n}$ such that $\theta$ is constant on every rectangular partition.

By definition, we have for any $\theta \in \R^{L_{d,n}}$,
$$k(\theta) \leq k_{\rdp}(\theta).$$
%{\color{red} Insert pics from oscar?}

We can now define the Dyadic CART estimator for a tuning parameter $\lambda > 0$,
\begin{equation}\label{eq:dcdefn}
\hat{\theta}^{(\lambda)} \coloneqq \argmin_{\theta \in \R^{L_{d,n}}} \big(\|y - \theta\|^2 + \lambda k_{\rdp}(\theta)\big).
\end{equation}

Equivalently, we can also define $\hat{\theta}^{(\lambda)} = P_{S_{\hat{\pi}^{(\lambda)}}} y$, where $\hat{\pi}^{(\lambda)}$ is a data dependent partition defined as 
\begin{equation*}
\hat{\pi}^{(\lambda)} := \argmin_{\pi \in \mathcal{P}_{\rdp,d,n}} \big(\|y - P_{S_\pi}y\|^2 + \lambda  |\pi|\big).
\end{equation*}
In the above, for any $\pi \in \mathcal{P}_{\rdp,d,n}$, $|\pi|$ denotes the number of rectangles constituting $\pi$, $S_\pi$ denotes the subspace of $\R^{L_{d,n}}$ which consists of all arrays which are constant on every rectangle of $\pi$ and $P_{S_\pi}$ denotes the orthogonal projection matrix on that subspace. The discrete optimization problem in the last display can be solved by a dynamic programming algorithm in $O_{d}(N)$ time which makes fast computation of Dyadic CART possible.

\subsubsection{\textbf{Existing Oracle Risk Bound and its Implications}}

We now state the oracle risk bound satisfied by the Dyadic CART estimator.

\begin{theorem}\label{thm:adapt}[Theorem $2.1$ in~\cite{chatterjee2019adaptive}]
	
	 Suppose the error vector $\epsilon$ is gaussian with mean $0$ and covariance matrix $\sigma^2 I$. Then there exists an absolute constant $C > 0$ such that if we set $\lambda \geq C \sigma^2 \log N$, then we have the following risk bound
	%\begin{equation*}
	%\E \|\hat{\theta}_{\lambda} - \theta^*\|^2 \leq \inf_{\theta \in \Theta} \big(2\:\|\theta - \theta^*\|^2 + \lambda k(\theta)\big) + 576 \sigma^2.and by setting $\lambda \geq C d\:\sigma^2\:\log n$ we have the following risk bound for $\hat{\theta}_{\rdp,\lambda}$:
	\begin{equation*}
	\E \|\hat{\theta}^{(\lambda)} - \theta^*\|^2 \leq \inf_{\theta \in \R^N} \left[3\:\|\theta - \theta^*\|^2 + 2\lambda\:k_{\rdp}(\theta)\right] + C \sigma^2.
	\end{equation*}
	%where $a \in \{\rdp,\hier,\all\}.$%thus giving risk bounds for three different estimators.
	%For the same choice of $\lambda$ and $C$, we also have
	%\begin{equation*}
	%\E \|\hat{\theta}_{\hier,\lambda} - \theta^*\|^2 \leq \inf_{\theta \in \R^{L_{n,d}}} \big(2\:\|\theta - \theta^*\|^2 + \lambda k_{\hier}(\theta)\big) + C\:\sigma^2.
	%\end{equation*}
\end{theorem}

%{\color{red} Set $\delta$ equal to something? Setting $\delta = 1/2$ gives us certain expressions which are similar to the case of $\hat{\theta}^{(\lambda, I)}$, e.g., $R(\theta^*, \lambda)$.}

In the case when $d = 2$ the above oracle risk bound had already appeared in the original paper~\cite{donoho1997cart}, albeit up to an extra log factor. This oracle risk bound actually can be used to show that the Dyadic CART estimator adaptively attains near optimal rates of convergence for several function classes of interest. For instance, the above oracle risk bound was used in~\cite{donoho1997cart} to show that Dyadic CART is minimax rate optimal over several bivariate anistropic smoothness classes of functions. In~\cite{chatterjee2019adaptive}, the above oracle risk bound was used to show %two facts. 
%Firstly, the univariate Dyadic CART estimator of order $r$ attains both the slow rate and the fast rate attained by Trend Filtering, up to an extra log factor. Secondly, 
that the Dyadic CART estimator is minimax rate optimal over the class of bounded variation signals in general dimensions and thus matches the known rates of convergence attained by the Total Variation Denoising estimator; see~\cite{hutter2016optimal},~\cite{sadhanala2016total}. It has been explained in detail in~\cite{chatterjee2019adaptive} how the above oracle risk bound (along with its fast computation) puts forward Dyadic CART as a computationally faster alternative to Trend Filtering and its multidimensional versions while essentially retaining (and even improving in some aspects) its statistical benefits.

The main question we consider here is the following.

\textit{\textbf{Can a cross validated version of Dyadic CART still attain the oracle risk bound in Theorem~\ref{thm:adapt}?}}

To the best of our knowledge, the above question is unanswered as of now. We answer this question in the affirmative in this paper. Since our cross validated version of Dyadic CART will also satisfy a result very similar to the oracle risk bound as in Theorem~\ref{thm:adapt} it will essentially inherit all the known results for the usual Dyadic CART mentioned above. 

\subsection{\textbf{Description of the CVDCART estimator}}
We will follow our general scheme of defining cross validated estimators as laid out in Section~\ref{sec:scheme}. In what follows and in the descriptions of all our CV estimators in the later sections as well, we will only describe the first three steps to avoid repetition. The first three steps describe the number of folds, the construction of the folds and the construction of the completion estimators respectively. The last three steps are common to all of them and are the same as the last three steps laid out in the general framework in Section~\ref{sec:scheme}. In all of our CV estimators, we take our grid of candidate tuning values $$\Lambda = \Lambda_j = \{1, 2, 2^2, 2^3,\dots, 2^{N^*}\}$$ for an appropriate $N^*$.

Let $\hat{\theta}^{(\lambda)}$ be the family of Dyadic CART estimators with tuning parameter $\lambda \geq 0$ as defined in~\eqref{eq:dcdefn}. We denote the final resulting cross validated Dyadic CART estimator by $\hat{\theta}_{CVDC}$.

\begin{enumerate}
	\item Set $K = 2$. 
	\vspace{0.2in}
	\item  We divide $L_{d,n}$ randomly into two folds/subsets $I_1,I_2$ as follows:
	Let $W \in \R^{L_{d,n}}$ be a random array consisting of i.i.d $Bernoulli(1/2)$ entries. Now define
	$$ I_1 = \{(i_1,\dots,i_d) \in L_{d,n}: W(i_1,\dots,i_d) = 1\}.$$ The set $I_2$ is just the complement of $I_1$ in $L_{d,n}.$
		\vspace{0.2in}
	\item Next, we define the estimators for $j \in \{1,2\}$,
	\begin{equation}\label{eq:defncross}
    \hat{\theta}^{(\lambda, I_j)} = \argmin_{\theta \in \R^{L_{d,n}}} ||y_{I_j} - \theta_{I_j}||^2 + \lambda k_{\rdp}(\theta).
    \end{equation}
    Note that $\hat{\theta}^{(\lambda, I_j)}$ is a completion version of Dyadic CART because it only depends on $y_{I_j}.$
    \begin{comment}
	\vspace{0.2in}
	\item Consider an exponentially spaced finite grid of possible values of the tuning parameter $\lambda_j$, namely $\Lambda_j = \{1, 2, 2^2, 2^3,\dots, 2^{N^*}\}$, for $j \in \{1,2\}$, where the choice of $N^*$ will be determined later. Now define $\hat{\lambda}_j$ to be the candidate in $\Lambda_j$ for which the prediction error $CVERR_{j}$ is the minimum, that is,
	$$\hat{\lambda}_j := \argmin_{\lambda \in \Lambda_j} \left|\left|y_{I_j} - \hat{\theta}^{(\lambda, I_j^c)}_{I_j}\right|\right|^2,\quad j \in \{1,2\}.$$ 
	Now define an intermediate estimator $\tilde{\theta} \in \R^{L_{d,n}}$ such that
	\begin{equation}\label{eq:defninterdc}
	\tilde{\theta}_{I_j} = \hat{\theta}^{(\hat{\lambda}_j, I_j^c)}_{I_j},\quad j \in \{1,2\}.
	\end{equation}

	\item Define $$\hat{\lambda} = \argmin_{\lambda \in \Lambda} \|\hat{\theta}^{(\lambda)} - \tilde{\theta}\|^2$$ where $\Lambda = \{1, 2, 2^2, 2^3,\dots, 2^{N^*}\}.$
	 %the original estimator $\hat{\theta}^{(\lambda)}$ with $\lambda = \hat{\lambda}$ as defined below: 
	 Finally, our estimator (CVDCART) $\hat{\theta}_{CVDC}$ is defined such that
	\begin{equation}\label{eq:defn}
	\hat{\theta}_{CVDC} = \hat{\theta}^{(\hat{\lambda})}.
	\end{equation}
	\vspace{0.2in}
	%\item Finally, our estimator (CVDCART) $\hat{\theta}$ is defined such that
	%\begin{equation}\label{eq:cvtfdefn}
	%\hat{\theta}_{I_j} = \hat{\theta}^{(\hat{\lambda}_j, I_j^c)}_{I_j},\quad j \in \{1,2\}.
	%\end{equation}
	\end{comment}
	
\end{enumerate}

%{\color{red} Say what is $N^*$? Make the notation $\hat{\theta}_{DC}?$ I have written something regarding $N^*$ and changed $\hat{\theta}$ to $\hat{\theta}_{DC}$.}

\subsection{\textbf{Computation of the CVDCART estimator}}
The major step in computing $\hat{\theta}_{CVDC}$ is to compute $\hat{\theta}^{(\lambda,I_j)}$ for $j = \{1,2\}.$ %It turns out that $\hat{\theta}^{(\lambda,I_j)}$ can be computed by a similar botton up dynamic program which is used to compute the original version of Dyadic CART; see Lemma $1.1$ in~\cite{chatterjee2019adaptive}. 
We present a lemma below stating the computational complexity of $\hat{\theta}_{CVDC}.$

\begin{lemma}\label{lem:comp_cvdc}
Let $I$ denote the set $I_1$ or $I_2.$ The computational complexity of the completion estimators $\hat{\theta}^{(\lambda,I)}$, i.e the number of elementary operations involved in computing $\hat{\theta}^{(\lambda,I)}$ is bounded by $C 2^d d N$ for some absolute constant $C > 0$. Therefore, the overall computational complexity of $\hat{\theta}_{CVDC}$ is bounded by $C N^* 2^d d N.$ Since $N^*$ can be taken to be $O(\log N)$ (as explained later), the overall computational complexity becomes $O_{d}(N \log N)$  in this case.

%\begin{remark}
%	The $O_{d}$ notation is the same as the $O$ notation except it signifies that the constant factor while comparing two sequences in $n$ may depend on the underlying dimension $d.$
%	\end{remark}
%{\color{red} Make clear $O$, $\tilde{O}$ and $O_{d}$ notations at the beginning.}

\begin{remark}
	The above lemma ensures that our CVDCART estimator can also be computed in near linear time in the sample size.
\end{remark}

%The computational complexity of $\hat{\theta}^{(\lambda,I)}$, i.e the number of elementary operations involved in computing $\hat{\theta}^{(\lambda,I)}$ is bounded by $C 2^d d N$ for some absolute constant $C > 0$.
%\begin{comment}
%\begin{equation*}
%\begin{cases}
%C 2^d d N\:\:\:\text{for}\:\:\:\:r = 0\\
%C 2^d d^{3r} N\:\:\:\text{for}\:\:\: r = 1.
%\end{cases}
%\end{equation*}
%\end{comment}
	\end{lemma}

In Section~\ref{sec:compu} %~\ref{sec:compu}, 
we describe a bottom up dynamic programming based algorithm to compute $\hat{\theta}^{(\lambda,I)}$. The underlying idea behind this algorithm is similar to the original algorithm given in~\cite{donoho1997cart} to compute the original version of Dyadic CART based on the full data. The description of the algorithm in Section~\ref{sec:compu} %~\ref{sec:compu} 
also clarifies what is the computational complexity of the algorithm, thereby proving Lemma~\ref{lem:comp_cvdc}.

\subsection{\textbf{Specification of $\hat{\theta}^{(\lambda,I)}$}}
There may be multiple solutions to the optimization problem defined in~\eqref{eq:defncross}. For our main result (which is Theorem~\ref{thm:main_cv_dc}) to hold, we need to add one more specification which will complete the definition of $\hat{\theta}^{(\lambda,I)}$ for any given subset $I \subset L_{d,n}.$ Below, we use the notation $\overline{y}_{I}$ to denote the mean of all entries of $y$ in $I$.

In Section~\ref{sec:compu} %~\ref{sec:compu}, 
it is shown that the optimization problem in~\eqref{eq:defncross} can be solved by first solving the following discrete optimization problem over the space of all recursive dyadic partitions

\begin{equation}\label{eq:repre}
\hat{\pi} = \argmin_{\pi \in \mathcal{P}_{\rdp,d,n}} \big[\sum_{R \in \pi} \mathrm{1}(R \cap I \neq \emptyset) \sum_{u \in R \cap I} (y_{u} - \overline{y}_{R \cap I})^2 + \lambda |\pi|\big]
\end{equation}
where the notation $\sum_{R \in \pi}$ means summing over all the constituent rectangles $R$ of $\pi$ and $|\pi|$ denotes the number of constituent rectangles of the partition $\pi.$

It is then shown that a solution $\hat{\theta}^{(\lambda,I)}$ to the optimization problem~\eqref{eq:defncross} is a piecewise constant array over the optimal partition $\hat{\pi}$. For all $u \in R$, for each constituent rectangle $R$ of $\hat{\pi}$ such that the set $R \cap I$ is non empty, 
\begin{equation}\label{eq:defnonempty}
\hat{\theta}^{(\lambda,I)}_{u} = \overline{y}_{R \cap I}.
\end{equation}

It is possible for a constituent rectangle $R$ of the optimal partition $\hat{\pi}$ to not contain any data point from $I$, i.e, the set $R \cap I$ is empty. In that case, $\hat{\theta}^{(\lambda,I)}$ can take any constant value within $R$ and still be an optimal solution to the optimization problem~\eqref{eq:defncross}. In such a case, for all $u \in R$, we set 
\begin{equation}\label{eq:defempty}
\hat{\theta}^{(\lambda,I)}_{u} = \overline{y}_{I}
\end{equation}

This fully specifies the estimator $\hat{\theta}^{(\lambda,I)}$ which is a valid completion version of the Dyadic CART estimator being a function of $y_{I}$ only.

\begin{remark}
	Note that the above specification still does not mean that $\hat{\theta}^{(\lambda,I)}$ is uniquely defined. This is because $\hat{\pi}$ in~\eqref{eq:repre} is not necessarily uniquely defined. However, as long as we take a solution $\hat{\pi}$ to the optimization problem in~\eqref{eq:repre} and then construct $\hat{\theta}^{(\lambda,I)}$ satisfying both~\eqref{eq:defnonempty} and~\eqref{eq:defempty}, Theorem~\ref{thm:main_cv_dc} holds. 
	\end{remark}

%The proof of this lemma is similar to that of Lemma $1.1$ in~\cite{chatterjee2019adaptive} with some differences in the details. For the sake of completeness, we give a proof in Section~\ref{}.

%{\color{red} Remark about non uniqueness. Remark about empty rectangles and estimating by the mean on the empty rectangles. }

%{\color{red} To do for A: Mention $\hat{\pi}$ and then $\hat{\theta}.$ Give the computational complexity lemma and prove it with algorithm in the appendix. Fill in this lemma}

%{\color{red} In this paper, we have used single $|$ for $l_{\infty}$ and $l_1$ norms. Also, there is no subscript for the case of $l_1$ norm, and there is a possibility that it could be confused with the mod/absolute value function.}

\subsection{\textbf{Main Result for the CVDCART estimator}}
Now we state an oracle risk bound for our proposed CVDCART estimator in general dimensions. Before that, let us define the following quantities
\begin{align*}
R(\theta^*, \lambda) &:= \inf_{\theta \in \R^{L_{d,n}}} \big(3 \|\theta - \theta^*\|^2 + 2 \lambda k_{\rdp}(\theta)\big), \quad \text{and}\\
V(\theta^*) &:= \max_{u,v \in L_{d,n}} |\theta^*_{u} - \theta^*_{v}|.
\end{align*}

%The term $R(\theta^*,\lambda)$ is 

\begin{theorem}\label{thm:main_cv_dc}
Fix any $\alpha \geq 1$ and any $\delta > 0$. There exists an absolute constant $C > 0$ such that if we set our grid $\Lambda_j = \Lambda = \{1, 2, 2^2, \dots, 2^{N^*}\}$ for $j = \{1, 2\}$, satisfying
$$2^{N^*} > C\sigma^2\log N,$$
then we have the following bound
\begin{align*}
 \frac{1}{N} \|\hat{\theta}_{CVDC} - \theta^*\|^2 \leq \frac{C_3}{N}\Bigg\{\big(R(\theta^*, C \sigma^2 \log n)
  + \alpha \sigma^2 \log N\big) &\left(\frac{V(\theta^*)}{\sigma} + \sqrt{\log N}\right)^2\\
  &\qquad \qquad  + \sigma^2\log(4N^*/\delta)\Bigg\}
\end{align*}
with probability at least $1 - \delta - C_1 \log (V(\theta^*) \sqrt{N})N^{-C_2 \alpha}$, where $C,C_1, C_2,C_3$ are absolute positive constants which may only depend on the underlying dimension $d.$

\end{theorem}

It is now worthwhile discussing some aspects of Theorem~\ref{thm:main_cv_dc}.

\begin{enumerate}
	\item 
	The above theorem basically ensures that the mean squared error (MSE) of our cross validated estimator $\hat{\theta}_{CVDC}$ can also essentially be bounded (up to additive and multiplicative log factors) by the desired factor $R(\theta^*, \sigma^2 \log N)$. The same bound holds for the optimally tuned version of Dyadic CART as is stated in Theorem~\ref{thm:adapt}. The only essential difference is that our bound contains an extra multiplicative factor $(\frac{V(\theta^*)}{\sigma} + \sqrt{\log N})^2.$ The term $V(\theta^*)$ captures the range of the underlying signal. For realistic signals, the range should stay bounded. In these cases, this extra multiplicative factor would then be a logarithmic factor. %Clearly, one has to incur extra MSE for the cross validated version compared to the optimally tuned version. Theorem~\ref{thm:main_cv_dc} ensures that 
	
	\medskip
	
	\item Theorem~\ref{thm:main_cv_dc} ensures that our cross validated estimator $\hat{\theta}_{CVDC}$ (up to log factors) enjoys a similar oracle risk bound as the optimally tuned version of Dyadic CART. Therefore, the CVDCART estimator essentially inherits all the known statistical risk bounds for Dyadic CART. In particular, CVDCART estimator would be minimax rate optimal (up to log factors) for several function/signal classes of interest such as anisotropically smooth functions (see~\cite{donoho1997cart}), piecewise constant signals on arbitrary rectangular partitions when $d \leq 2$ and signals with finite bounded variation (see~\cite{chatterjee2019adaptive}). 
	
	\medskip
	
	\item The existing risk bound Theorem~\ref{thm:adapt} says that the optimal tuning parameter choice is $C \sigma^2 \log N$ for some absolute constant $C > 0$. Theorem~\ref{thm:main_cv_dc} holds as long as such a choice of $\lambda$ is included in $\Lambda.$ This suggests that we  would like to have the grid of choices $\Lambda$ to be dense enough so that we do not miss the optimal tuning value range. On the other hand, we would like to have a sparse grid $\Lambda$ because the computational complexity scales like the cardinality of $\Lambda$ times the complexity of computing one Dyadic CART estimator for a given $\lambda.$
	
	\medskip
	Observe that our risk bound (in particular $R(\theta^*,\lambda)$) scales proportionally with $\lambda$. This implies that missing the optimal tuning value by a factor of $2$ means that we pay at most $2$ times the MSE of the ideally tuned version. This fact allows us to take a geometrically growing grid $\Lambda = \{1,2,2^2,\dots,2^{N^*}\}.$ Selecting such a sparse grid $\Lambda$ then has obvious computational benefits. The only disadvantage here is that $N^*$ then becomes like a tuning parameter to be set by the user. However, in practice and in theory, this seems to be a minor issue. Theorem~\ref{thm:main_cv_dc} holds if $2^{N^*}$ is larger than the theoretically recommended choice of $\lambda = C \sigma^2 \log N$. Hence, plugging in even a gross overestimate $\sigma_{over}$ of $\sigma$ and choosing $N^* = C' \left(\log \log N + \log \sigma_{over}\right)$ for a large enough constant $C'$ would suffice for any realistic value of $\sigma.$ In our simulations we simply take $N^* = \log_{2} N.$

\end{enumerate}

%{\color{red} Rewrite this theorem and clean up the constants some more. Mention a computational comp lemma. Remark why this theorem is important. How should one choose $N^*$? How to read the above theorem. }

%{\color{red} maybe introduce some terminology such as CV errors and so in the general framework and then integrate this in the proof.}

%{\color{red} Rewrite this theorem and clean up the constants some more. Mention a computational comp lemma. Remark why this theorem is important. How should one choose $N^*$? How to read the above theorem. }

%{\color{red} maybe introduce some terminology such as CV errors and so in the general framework and then integrate this in the proof.}

%\input{trendfiltering}
\section{Trend Filtering}\label{sec:tf}
\subsection{\textbf{Background and Related Work}}

Trend Filtering, proposed by~\cite{kim2009ell_1}, is a univariate nonparametric regression method that has become popular recently; see~\cite{tibshirani2020divided} for a comprehensive overview. For a given integer $r \geq 1$ and any tuning parameter $\lambda \geq 0$, the $r^{th}$ order trend filtering estimator $\hat{\theta}^{(r)}_{\lambda}$ is defined as the minimizer of the sum of squared errors when we penalize the sum of the absolute $r^{th}$ order discrete derivatives of the signal. Formally, given a data vector $y$,
\begin{equation}\label{eq:tfdefn}
\hat{\theta}^{(r)}_{\lambda} := \left(\argmin_{\theta \in \R^{n}} \frac{1}{2}||y - \theta||^2 + \lambda n^{r - 1}||D^{(r)}\theta||_1\right)
\end{equation}
where $D^{(1)}\theta := \left(\theta_2 - \theta_1, \theta_3 - \theta_2, \dots, \theta_n - \theta_{n-1}\right)$ and $D^{(r)}\theta$, for $r \geq 2$, is recursively defined as $D^{(r)}\theta := D^{(1)}D^{(r-1)}\theta$. For any positive integer $r \geq 1$, let us now define the $r^{th}$  order total variation of a vector $\theta$ as follows:
\begin{equation}
\mathrm{TV}^{(r)}(\theta) = n^{r - 1} \|D^{(r)}(\theta)\|_{1}
\end{equation}
where $\|\cdot\|_1$ denotes the usual $\ell_1$ norm of a vector. 
\begin{remark}
	The $n^{r - 1}$ term in the above definition is a normalizing factor and is written following the convention adopted in the trend filtering literature; see for instance the terminology of canonical scaling introduced in~\cite{sadhanala2016total}. If we think of $\theta$ as evaluations of a $r$ times differentiable function $f:[0,1] \rightarrow  R$ on the grid $(1/n,2/n\dots,n/n)$ then the Riemann approximation to the integral $\int_{[0,1]} \vert f^{(r)}(t)\vert  dt$ is precisely equal to $\mathrm{TV}^{(r)}(\theta).$ Here $f^{(r)}$ denotes the $r$th derivative of $f.$ Thus, for natural instances of $\theta$, the reader can imagine that $\mathrm{TV}^{(r)}(\theta) = O(1).$ 
\end{remark}

A continuous version of these trend filtering estimators, where discrete derivatives are replaced by continuous derivatives, was proposed much earlier in the statistics literature 
by~\cite{mammen1997locally} under the name {\em locally adaptive regression splines}. 
By now, there exists a body of literature studying the risk properties of trend filtering under squared error loss. There exists two strands of risk bounds for trend filtering in the literature focussing on two different aspects.

Firstly, it is known that for any $r \geq 1,$ a well tuned trend filtering estimator $\hat{\theta}^{(r)}_{\lambda}$ attains a MSE bound $O\left(\frac{(TV^{(r)}(\theta^*))^{1/r}}{n}\right)^{2r/2r + 1}.$ This bound is minimax rate optimal over the space $\{\theta \in  \R^n: \TV^{(r)}(\theta) \leq V\}$ for a given $V > 0$ and has been shown in~\cite{tibshirani2014adaptive} and~\cite{wang2014falling} building on earlier results by~\cite{mammen1997locally}. A standard terminology in this field terms this $O(n^{-2r/2r + 1})$ rate as the \textit{slow rate}.

Secondly, it is also known that an ideally tuned Trend Filtering (of order $r$) estimator can adapt to $\|D^{r}(\theta)\|_0$, the number of non zero elements in the $r^{th}$  order differences, \textit{under some assumptions on $\theta^*$}. Such a result has been shown in~\cite{guntuboyina2020adaptive} (for the constrained version of Trend Filtering of all orders) and~\cite{ortelli2019prediction} (for the penalized version of Trend Filtering with $r \leq 4$). In this case, the Trend Filtering estimator of order $r$ attains the near parametric $\tilde{O}(\|D^{(r)}(\theta)\|_0/n)$ rate which can be much faster than the $O(n^{-2r/2r + 1})$ rate. Standard terminology in this field terms this as the \textit{fast rate}.

The big problem is that the results described above are shown to hold only under theoretical choices of the tuning parameter. These choices depend on unknown problem parameters and hence cannot be directly implemented in practice. Moreover different tuning is needed to achieve slow or fast rates.  A square root version of Trend Filtering was proposed by~\cite{ortelli2021oracle} to mitigate this issue. It has been shown that the ideal choice of the tuning parameter for the square root version does not depend on the noise variance. However, the tuning parameter still needs to be set to a particular unspecified constant (differently depending on whether slow or fast rates are desired) and thus does not really solve this problem.

Therefore, a version of Trend Filtering which chooses the tuning parameter in a data driven way and attains both slow and fast rates is highly desirable. This naturally leads us to consider cross validation. Practical usage of Trend Filtering almost always involves cross validation to choose $\lambda$; e.g, see~\cite{politsch2020trend}. However, no theoretical properties are known for a cross validated version of Trend Filtering. We attempt to fill this gap in the literature by proposing a cross validated version of Trend Filtering based on our general framework. %Our hope is to show both a \textit{slow} rate and a \textit{fast} rate for our cross validated Trend Filtering (CVTF) estimator 
Our goal here is to show that our cross validated version nearly (atmost up to log factors) attains the risk (both the slow rate and the fast rate) of the ideally tuned versions.
%{\color{red} a comment on square root tf?}
\subsection{\textbf{Description of the CVTF Estimator}}\label{sec:cvtf}
We will again follow our general scheme of defining cross validated estimators as laid out in Section~\ref{sec:scheme}.  Fix any $r \geq 1$ and let $\hat{\theta}^{(r)}_{\lambda}$ be the family of $r^{th}$ order Trend Filtering estimators with tuning parameter $\lambda \geq 0$ as defined in~\eqref{eq:tfdefn}. We denote the final resulting cross validated trend filtering estimator of order $r \geq 1$ by $\hat{\theta}^{(r)}_{CVTF}.$%We use the notation $[m]$ for a positive integer $m$ to denote the set of positive integers $\{1,2,\dots,m\}.$
\begin{enumerate}
	\item Set $K = r + 1$. 
	
	\item  Divide $[n]$ deterministically into $K$ disjoint index sets (ordered) $I_1, I_2, \dots, I_K$ as follows. Let $n_0 = \floor{\frac{n}{K}}$. Then for any $j \in [K]$, define
	$$I_j = \{Kt + j \leq n : t = 0, 1, \dots, n_0\}.$$
	%The elements of these two sets are denoted as
	%$$I = \{I_1, I_2, \dots, I_{|I|}\}, \quad I^c = \{I^c_1, I^c_2, \dots, I^c_{|I^c|}\}.$$
	In words, data  points $K$ positions apart are placed into the same fold. 
	\item For all $j \in [K]$, define  $\tilde{y}(I_j^c) \in \R^n$ by interpolating $y_{I_j^c}$ as follows:
	$$\tilde{y}(I_j^c)_i =
	\begin{cases}
	y_i\ \quad &\text{if}\ i \in I_j^c\\
	\sum_{l = 1}^r (-1)^{l+1}{r \choose l} y_{i + l}        &\text{if}\ i \in I_j\ \text{and}\ i+r \leq n\\ %&\text{if}\ i = Kt_* + j\ \text{for some}\ t_* \leq n_0-1\\
	\sum_{l = 1}^r (-1)^{l+1}{r \choose l} y_{i - l}        &\text{if}\ i \in I_j\ \text{and}\ i+r > n %&\text{if}\ i = Kn_0 + j \leq n
	\end{cases}.$$
	In words, $\tilde{y}(I_j^c)$ is defined in such a way that within the index set $I_j^c$ it is same as $y$ but for any index in $I_j$ it is linearly interpolated from the neighbouring indices of $y$ in $I_j^c$ by a $r^{th}$ order polynomial interpolation scheme. %Moreover, the method of extrapolation is: for any index in $i \in I_j$, we always extrapolate from its right neighborhood as long as it contains $r$ elements, that is, $i + r \leq n$. Otherwise, we extrapolate from its left neighborhood. The main reason for considering such an extrapolation is explained in Remark {\color{red} blah}.%Moreover, the method of extrapolation is: for any index in $i \in I_j$, we always extrapolate from its right neighborhood as long as it contains $r$ elements, that is, $i + r \leq n$. Otherwise, we extrapolate from its left neighborhood.	

	Next, we define the completion estimators
	$$\hat{\theta}^{(\lambda, I_j^c,r)} := \argmin_{\theta \in \R^n} \frac{1}{2}\left\|\tilde{y}(I_j^c) - \theta\right\|^2 + \lambda n^{r-1}\left\|D^{(r)}\theta\right\|_1, \quad j \in [K].$$
	Note that $\hat{\theta}^{(\lambda, I_j^c,r)}$ is a valid completion version as it is a function of $y_{I_j^c}$ only.

	\begin{comment}
	\item Consider an exponentially spaced finite grid of possible values of the tuning parameter $\lambda$, namely $\Lambda_j = \{1,2,2^2,2^3,\dots,2^{N^*}\}.$ For any $j \in [K]$, define $\hat{\lambda}_j$ to be the candidate in $\Lambda_j$ for which the prediction error is the minimum, that is,
	$$\hat{\lambda}_j := \argmin_{\lambda \in \Lambda_j} \left|\left|y_{I_j} - \hat{\theta}^{(\lambda, I_j^c,r)}_{I_j}\right|\right|^2,\quad j \in [K].$$ %A different and completely data driven choice of $\Lambda_j$ and $\hat{\lambda}_j$ is explained in Section~\ref{sec:alttf}}.

	Now define an intermediate estimator $\tilde{\theta} \in \R^{n}$ such that
	\begin{equation}\label{eq:defnintertf}
	\tilde{\theta}_{I_j} = \hat{\theta}^{(\hat{\lambda}_j, I_j^c,r)}_{I_j},\quad j \in [K].
	\end{equation}

	\item Define $$\hat{\lambda} = \argmin_{\lambda \in \Lambda} \|\hat{\theta}^{(r)}_{\lambda} - \tilde{\theta}\|^2$$ where $\Lambda = \{1,2,2^2,2^3,\dots,2^{N^*}\}.$
	%the original estimator $\hat{\theta}^{(\lambda)}$ with $\lambda = \hat{\lambda}$ as defined below: 
	Finally, our estimator (CVTF) $\hat{\theta}_{CVTF}^{(r)}$ is defined such that
	\begin{equation}\label{eq:defn}
	\hat{\theta}_{CVTF}^{(r)} = \hat{\theta}^{(r)}_{\hat{\lambda}}.
	\end{equation}
	\end{comment}
\end{enumerate}

\begin{remark}
	A different yet valid choice of $\Lambda_j$ and $\hat{\lambda}_j$, $j \in [K]$, $\Lambda$ and $\hat{\lambda}$ is described in Section~\ref{sec:alttf}.
	\end{remark}

%\subsection{Results of~\cite{ortelli2019prediction}}
%For the state of the art risk bounds on constrained/penalized Trend Filtering, we refer to~\cite{guntuboyina2020adaptive} and~\cite{ortelli2019prediction} respectively. The CVTF estimator (of order $r \geq 1$) is a $r + 1$ fold crossvalidated version of the penalized version of Trend Filtering (of order $r$).  Here, we attempt to prove a risk bound for the CVTF estimator which is similar to the one presented in~\cite{ortelli2019prediction}. Therefore, before stating our risk bound for the CVTF estimator, we first describe the results of~\cite{ortelli2019prediction} precisely. 

%Our risk bound and its proof is based on the style of analysis of~\cite{ortelli2019prediction}. To describe our risk bound 

%The state of the art results for the usual penalized version of Trend Filtering have been shown in~\cite{ortelli2019prediction} where both the \textit{slow rate} and the ~\textit{fast rate} have been shown to hold under ideal choices of the tuning parameter. 
\subsection{\textbf{Main Results for the CVTF Estimator}}
Below we state both the slow rate and the fast rate results for the proposed CVTF estimator. 

%For $r \geq 1$, let $\cD = [n-r]$ and $\cS$ be an index set such that $\cS = \{t_1, \dots, t_s\} \subseteq \cD$, where $1 \leq t_1 < \dots < t_s \leq n-r$. Also, let $t_0 := 0$ and $t_{s+1} := n-r+1$. Next, we define $n_i := t_i - t_{i-1}$, $i \in [s+1]$ and $n_{\rm \max} := \max_{i \in [s+1]} n_i$.
%Furthermore, let us define the sign vector $q^* \in \{-1, +1\}^s$ containing the signs of the elements in $(D^{(r)}\theta^*)_{\cS}$, that is,
%$$q^*_{i} := sign(D^{(r)}\theta^*)_{t_i}, \quad i \in [s],$$
%and the index set $$\cS^{\pm} := \{2 \leq i \leq s : q^*_{i}q^*_{i-1} = -1\} \cup \{1, s+1\}.$$
%{\color{red} new notation with $r$ th order. somewhere we have to say we drop the index r.}
\begin{theorem}\label{thm:slowratetv}[Slow Rate]

    Fix any $r \geq 1$ and any $\delta > 0.$ There exists a constant $C_r$ only depending on $r$ such that if we take our grid $\Lambda = \{1,2,2^2,2^3,\dots,2^{N^*}\}$ satisfying
	\begin{equation*}
	2^{N^*} \geq C_r \sigma \big(n \log n\big)^{1/(2r + 1)} 
	\end{equation*}
	then we have the following bound with probability at least $1 - \delta$,
	%\begin{align*}
	%\sqrt{SSE(\hat{\theta}_{cvlasso},\theta^*)} \leq 4 \sqrt{32} |\beta^*|_1 (M n \log p)^{1/4} + 4 (2 \log 1/\delta)^{1/4} + 8\sqrt{2} \sigma \sqrt{\log 2^{K^*}} + 8 \sqrt{2 \sigma^2 \log 4/\delta}
	%\end{align*}
	\begin{align*}
	 \frac{1}{n} \|\hat{\theta}^{(r)}_{CVTF} - \theta^*\|^2 \leq \frac{2C_rV^*}{n^r} \left|D^{(r)}\theta^*\right|_{\infty} + C_r \sigma^2\left(n^{-\frac{2r}{2r+1}} (V^*\log (n/\delta))^{\frac{1}{2r+1}} + \frac{N^* + \log(1/\delta)}{n}\right),
	\end{align*}
	where $V^* = n^{r-1} \left\|D^{(r)}\theta^*\right\|_1 = \TV^{(r)}(\theta^*)$.
\end{theorem}

%{\color{red} The exponent of $V$ does not look right.}
%\end{theorem}
%\begin{comment}
%\begin{theorem}
%For any $u, v > 0$, if $\lambda$ is chosen such that 
%\begin{equation}\label{eq:lamtv}
%\lambda = O\left(n^{\frac{1}{2r+1}}(\log n + u)^{\frac{1}{2r+1}}\right),
%\end{equation}
%and $n^{r-1} \left\|D^{(r)}\theta^*\right\|_1 = V^*$,
%then there exists positive constants $C_1$, $C_2$ such that with probability at least $1 - e^{-C_1u} - e^{-C_2v}$, for any $j \in [K]$,
%\begin{align*}
%SSE\left(\hat{\theta}^{(\lambda, I_j^c)}_{I_j}, \theta^*_{I_j}\right) = 2\left\|D^{(r)}\theta^*\right\|_{\infty}\left\|D^{(r)}\theta^*\right\|_{1} + O\left(n^{\frac{1}{2r+1}} (V^*(\log n + u))^{\frac{1}{2r+1}}\right) + 8\sigma^2 v
%\end{align*}
%\end{theorem}
%\end{comment}

\begin{theorem}\label{thm:fastratetv}[Fast Rate]

Fix any $1 \leq r \leq 4$ and any $\delta > 0$. Let $s =\|D^{(r) } \theta^*\|_0 $ and $\cS=\{  j \,:\,  (D^{(r) } \theta^*)_j  \neq 0   \}$. Then $\cS$ can be represented as $\cS = \{t_1, \dots, t_s\} \subseteq [n-r]$, where $1 \leq t_1 < \dots < t_s \leq n-r$. Also, let $t_0 := 0$ and $t_{s+1} := n-r+1$. Next, we define $n_i := t_i - t_{i-1}$, $i \in [s+1]$ and $n_{\rm \max} := \max_{i \in [s+1]} n_i$. 

Define the sign vector $q^* \in \{-1, +1\}^s$ containing the signs of the elements in $(D^{(r)}\theta^*)_{\cS}$, that is, for every $i \in [s]$,
$q^*_{i} := sign(D^{(r)}\theta^*)_{t_i},$
and the index set $$\cS^{\pm} := \{2 \leq i \leq s : q^*_{i}q^*_{i-1} = -1\} \cup \{1, s+1\}.$$
Suppose $\theta^*$ satisfies the following minimum length assumption, for a constant $c > 1$,
$$n_{\rm \max} \leq c n_i \ \text{and}\ n_i \geq r(r+2)\ \text{for all}\ i \in \cS^{\pm}.$$
Then there exists a constant $C_r$ only depending on $r$ such that if we take our grid $\Lambda = \{1, 2^1, 2^2, 2^3, \dots, 2^{N^*}\}$ satisfying 
\begin{equation*}
	2^{N^*} \geq C_r \sigma s^{-(2r-1)/2} \sqrt{n \log n}
	\end{equation*}
	then we have the following bound with probability at least $1 - \delta$,
	\begin{align*}
	 \frac{1}{n} \|\hat{\theta}^{(r)}_{CVTF} - \theta^*\|^2 &\leq \frac{2C_rs}{n} \left|D^{(r)}\theta^*\right|_{\infty}^2 + C_r \sigma^2\left(\frac{s}{n} \log n \log(n/\delta) + \frac{N^* + \log (1/ \delta)}{n}\right).
	\end{align*}
%	where $\cS^c = [n-r] \setminus \cS$.
\end{theorem}
%{\color{red} Can we simplify the statement of this proof please? Too much notations?}

We now make some remarks to explain certain aspects of the above theorems.

\begin{enumerate}
	\item We have presented both our slow rate and the fast rate theorem following the notations and presentation style adopted by~\cite{ortelli2019prediction} in Theorem $1.1$ in their paper. We have done this mainly because our proofs rely on the results developed by~\cite{ortelli2019prediction} and also to remain consistent with the existing literature. The two theorems above ensure that the $\hat{\theta}_{CVTF}$ estimator essentially attains the slow rate and the fast rate (both implied by Theorem $1.1$ in~\cite{ortelli2019prediction}) known for an ideally tuned penalized Trend Filtering estimator. The main difference in both our bounds are the extra additive terms involving $|D^{(r)}\theta^*|_{\infty}.$ However, as we explain below, this is typically a lower order term.

	\medskip

	\item Both the bounds above involve the term $|D^{(r)}\theta^*|_{\infty}.$ Note that under the canonical scaling where $V^* = O(1)$, we have $|D^{(r)}\theta^*|_{\infty} \leq \|D^{(r)}\theta^*\|_{1} = O(n^{1 - r}).$ This means that the terms involving $|D^{(r)}\theta^*|_{\infty}$ in our bounds can again be considered to be of a lower order for all $r \geq 1$ under realistic regimes of $V^*.$ 
	
	\medskip

  \item In light of the above two remarks, under the canonical scaling, the bound in Theorem \ref{thm:slowratetv} can be read as scaling like the near minimax rate $\tilde{O}(n^{-\frac{2r}{2r + 1}})$ and the bound in Theorem \ref{thm:fastratetv} scales like the near parametric rate $\tilde{O}(|D^{(r)}\theta^*|_{0}\:\: n^{-1})$ up to additive lower order terms. Thus, our bounds show that the CVTF estimator attains the slow rate and the fast rate, up to log factors, and hence does not suffer too much in comparison to ideally tuned trend filtering estimators, atleast in the context of rates of convergence.  
  
  \medskip
  
  \item We only state Theorem~\ref{thm:fastratetv} for $r \in \{1,2,3,4\}$ and the assumptions on $\theta^*$ in Theorem~\ref{thm:fastratetv} are identical to the assumptions made in Theorem $1.1$ of~\cite{ortelli2019prediction}. This is because our proof is based on the proof technique employed by~\cite{ortelli2019prediction}, as explained in Section~\ref{sec:tvsketch}. %~\ref{sec:tvsketch}.  
  The fast rate result in~\cite{ortelli2019prediction} also is shown to hold for $r \in \{1,2,3,4\}$. To the best of our knowledge, a complete proof of the fast rate for penalized trend filtering of order $r > 4$ is not yet available in the literature. In contrast, fast rates have been established for an ideally tuned constrained trend filtering of all orders; see~\cite{guntuboyina2020adaptive}. It is possible to develop a cross validated version of the constrained trend filtering using our framework and show that it will then enjoy fast rates for all orders $r \geq 1.$ However, in this paper we prefer considering the penalized version due to its popularity and computational benefits. %Moreover, arguably, in most practical applications, the case when $r \in \{1,2,3,4\}$ is of interest. 

  \medskip
  
  \item The assumption $n_{max} \leq c n_i$ for a constant $c > 1$ means that the length of each of the blocks in $\cS^{\pm} $ are within a constant factor of each other. This kind of minimum length assumption is standard and is also known to be necessary for fast rates to hold; see Remark $2.4$ in~\cite{guntuboyina2020adaptive}. Note that such a minimum length assumption is needed only for the blocks in $\cS^{\pm}$ and not for all blocks.
  For example, when $r = 1$, the blocks in $\cS^{\pm}$ are either the first and last constant pieces of $\theta^*$ or those constant pieces of $\theta^*$ which constitute a local maxima stretch or a local minima stretch. %This fact was already noted in Remark $2.4$ in~\cite{guntuboyina2020adaptive}.  %Our fast rate is $O(\frac{s}{n} (\log n)^2).$ It is possible that the exponent of $(\log n)$ can be improved to be $1$ instead of $2$; this is known in the case when $r =1$; see Theorem in~\cite{guntuboyina2020adaptive}. This extra log factor is an artifact of the proof technique of~\cite{ortelli2019prediction} which we adopt here. Since saving a log factor is not a major focus for us here we leave this possible improvement 

\medskip
	\item 	The ideal choice of the tuning parameter (as shown in Theorem $1.1$ in~\cite{ortelli2019prediction}) depends on whether we desire the slow rate or the fast rate. However, both these choices scale (with $n$) like $n^{\alpha}$ for some $0 < \alpha < 1.$ Therefore, as long as $2^{N^*}$ is chosen to be larger than these idealized choices, both the theorems presented above will hold. By construction,  $\Lambda$ contains an exponentially growing grid which means that $N^*$ can be chosen so that it grows logarithmically in $n,\sigma$. Therefore, in the regime where $\sigma$ stays bounded away from $\infty$, the term $\sigma^2  \frac{N^* + \log (1/ \delta)}{n} = O(\log n/n)$ appearing in both of the above theorems is a lower order term. In practice, one can choose $2^{N^*} = n$, for instance, which will satisfy the required condition for realistic sample sizes $n$ and $\sigma.$

\end{enumerate}

\section{Singular Value Thresholding for Matrix Estimation}\label{sec:svt}
Singular Value Thresholding (SVT) is a fundamental matrix estimation and completion method; see~\cite{cai2010singular},~\cite{donoho2014minimax},~\cite{chatterjee2015matrix}. It is known that Singular Value thresholding is an \textit{all purpose} matrix estimation method and performs well in a wide variety of structured matrix estimation problems; see~\cite{chatterjee2015matrix}. However, the existing guarantees for this estimator depend on a thresholding parameter being chosen to be larger than a cutoff value which depends on the noise variance. In practice, the choice of the threshold matters in regards to the finite sample performance of the SVT estimator; see Section $5$ (simulations) in~\cite{chatterjee2019estimation} where the authors were investigating the SVT estimator in the context of estimating Nonparametric Bradley Terry Matrices. Thus, it is of both theoretical and practical interest in using a cross validated version of the SVT estimator. To the best of our knowledge, a theoretical analysis of a CV version of SVT is not available in the literature. Our goal here is to demonstrate that our CV framework is well suited to develop a cross validated version of this fundamental estimator.

\subsection{\textbf{Background and Related Literature}}
The literature on SVT is vast. For our purposes here, we will just consider one particular result known for an optimally tuned SVT. We will then develop a CV version of SVT and  show that this particular result continues to hold for our CV version of SVT as well. We consider the basic denoising setting where we observe 
$$y = \theta^* + \epsilon$$
where $\theta^*$ is an underlying $n \times n$ signal matrix and $\epsilon$ is a $n \times n$ noise matrix consisting of i.i.d subgaussian errors with unknown subgaussian norm $\sigma.$ Consider the data matrix $y$ and consider its singular value decomposition
$$y = \sum_{i = 1}^{n} s_i u_i v_i^{t}.$$

Let $S_{\lambda} = \{i: |s_i| > \lambda\}$ be the set of thresholded singular values of $y$ with threshold level $\lambda > 0.$ Define the estimator 
\begin{equation}\label{eq:svtdef}
\hat{\theta}^{(\lambda)} = \sum_{i \in S_{\lambda}} s_i u_i v_i^{t}.
\end{equation}
This is how a standard version of the SVT estimator is defined. Under this setting, the following lemma can be traced back at least to Lemma $3$ in~\cite{shah2016stochastically}. It is probable that this result is even older. We use the notation $\|M\|_{op}$ to denote the operator norm of a $n \times n$ matrix $M.$

\begin{lemma}\label{lem:deter}

	For any fixed $\eta > 0$, if the threshold $\lambda$ is chosen such that $\lambda = (1 + \eta) \tau$ with $\tau > \|\epsilon\|_{op}$ then we have the following inequality:
	\begin{equation*}
	\|\hat{\theta}^{(\lambda)} - \theta^*\|^2 \leq 8 (1 + \eta)^2 \left[\sum_{j = 1}^{n} \min\{\tau^2, \sigma_j^2(\theta^*)\}\right],
	\end{equation*}
	where $\sigma_j(\theta^*)$ is the $j$th largest singular value (in absolute value) of $\theta^*$.
	
	%$\sigma_1^2(\theta^*) \geq \sigma_2(\theta^*) \geq \dots \geq \sigma_n^2(\theta^*)$ are the $n$ singular values of $\theta^*$, ranked in order of their absolute values.
\end{lemma}

\begin{remark}
	The above lemma is purely a deterministic inequality; there is no notion of randomness here.
	\end{remark}

The following is a standard bound on the maximum singular value of a random subgaussian matrix quoted from~\cite{vershynin2018high}.

\begin{theorem}[Theorem $4.4.5$ from~\cite{vershynin2018high}]\label{thm:op}
  Let $\epsilon$ be an $n \times n$ matrix whose entries are independent mean $0$ subgaussian random variables with subgaussian norm at most $\sigma.$ Then there exists an absolute constant $C > 0$ such that for any $t > 0$, we have 
  \begin{equation*}
  P\left(\|\epsilon\|_{op} \leq C \sigma (\sqrt{n} + t)\right) \geq 1 - 2 \exp(-t^2).
  \end{equation*}
 \end{theorem} 
  
Combining Lemma~\ref{lem:deter} and Theorem~\ref{thm:op} (after plugging in $t = \sqrt{n})$ immediately yields the following theorem.

\begin{theorem}\label{thm:svt}
		There exists an absolute constant $C > 0$ such that if the threshold $\lambda$ is chosen satisfying $\lambda = C\: \sigma \sqrt{n}$, then the following inequality holds with probability at least $1 - 2 \exp(-n)$,
		\begin{equation*}
		\frac{\|\hat{\theta}^{(\lambda)} - \theta^*\|^2}{n^2} \leq \frac{C}{n^2} \left[\sum_{j = 1}^{n} \min\{n \sigma^2, \sigma_j^2(\theta^*)\}\right]
		\end{equation*}
\end{theorem}

This theorem reveals the adaptive nature of the SVT estimator. This is because the right hand side is a deterministic quantity which only depends on the true signal $\theta^*.$ Intuitively, this term can be thought of as describing the spectral complexity of $\theta^*.$ The above risk bound can be used to derive the rates of convergence of the SVT estimator for several different types of classes of matrices of interest. We mention two standard classes below. For more interesting matrix classes where SVT can be applied; see~\cite{chatterjee2015matrix}.

\begin{enumerate}
	\item \textit{Low Rank Matrices}: If $\theta^*$ has rank $k$, then $\sigma_{j}(\theta^*) = 0$ for $j > k$ and hence we obtain a bound on the MSE which is $C \frac{k \sigma^2}{n}.$ It is well known that this is the minimax rate of estimation for the class of $n \times n$ matrices of rank $k.$ 
	
	\medskip
	\item \textit{Nonparametric Bradley Terry Matrices}: For a general structured class of matrices, one can typically show by an approximation theoretic argument that the singular values decay at a certain rate, even if they do not become exactly $0$ as in the exact low rank case. For example,~\cite{shah2016stochastically} showed that the right hand side in the above theorem scales like $\frac{C}{\sqrt{n}}$ for the class of $n \times n$ Nonparametric Bradley Terry Matrices. These matrices are monotone in both row and column, up to an unknown permutation and arise in modeling of pairwise comparison data; see~\cite{shah2016stochastically},~\cite{chatterjee2019estimation}.
\end{enumerate}

The important point to note here is that $\lambda$ needs to be set proportional to $\sigma\sqrt{n}$ for the above theorem to hold. Since $\sigma$ is typically unknown and the constant $C$ is unspecified it is natural to cross validate over $\lambda.$ Therefore, it would be highly desirable for a cross validated version (where the tuning parameter is chosen in a data driven way) of the SVT estimator to also satisfy a risk bound of the form given in Theorem~\ref{thm:svt}. We propose such an estimator in the next section.

\subsection{\textbf{Description of the CVSVT estimator}}
We will follow our general scheme of defining cross validated estimators as laid out in Section~\ref{sec:scheme}. Let $\hat{\theta}^{(\lambda)}$ be the family of Singular Value thresholding estimators with threshold parameter $\lambda \geq 0$ as defined in~\eqref{eq:svtdef}. We denote the final resulting cross validated SVT estimator as $\hat{\theta}_{CVSVT}.$

\begin{enumerate}
	\item Set $K = 2$. 
	\item  Divide $[n] \times [n]$ into $I_1,I_2$ randomly as follows. Each entry $(i,j) \in [n] \times [n]$ belongs to $I_1$ or $I_2$ with probability $1/2$ independently of other entries. 
	\item Let us denote $I_1$ by $I$ and $I_2$ by $I^c.$ Define the $n \times n$ binary matrix $W$ which takes the value $1$ on the entries in $I$ and $0$ elsewhere. Now define  $\tilde{y}(I) \in \R^{n \times n}$ as follows:
	$$\tilde{y}(I) = 2 y \circ W$$
	where $\circ$ denotes the operation of entrywise multiplication of two matrices of the same size. Similarly, define 
	$$\tilde{y}(I^c) = 2 y \circ (1 - W).$$
	Thus, $\tilde{y}(I),\tilde{y}(I^c)$ are matrices obtained by zeroing out entries of $y$ (corresponding to entries in $I$ or $I^c$) and then doubling it.
	
	Next, we define the completion estimator $\hat{\theta}^{(\lambda, I)}$ to be the 
	SVT estimator applied to the matrix $\tilde{y}(I)$ with threshold $\lambda.$ Define $\hat{\theta}^{(\lambda, I^c)}$ similarly using the matrix $\tilde{y}(I^c).$ Note that by definition, $\hat{\theta}^{(\lambda, I)}$ is a function only of $y_{I}$ and hence is a valid completion estimator.

    \begin{comment}
	\item Consider an exponentially spaced finite grid of possible values of the tuning parameter $\lambda_j$, namely $\Lambda_j = \{1, 2, 2^2, 2^3,\dots, 2^{N^*}\}$, for $j \in \{1,2\}$. Now define $\hat{\lambda}_j$ to be the candidate in $\Lambda_j$ for which the test error is the minimum, that is,
	$$\hat{\lambda}_j := \argmin_{\lambda \in \Lambda_j} \left|\left|y_{I_j} - \hat{\theta}^{(\lambda, I_j^c)}_{I_j}\right|\right|^2,\quad j \in \{1,2\}.$$ 
	Now define an intermediate estimator $\tilde{\theta} \in \R^{L_{d,n}}$ such that
	\begin{equation}\label{eq:defninterdc}
	\tilde{\theta}_{I_j} = \hat{\theta}^{(\hat{\lambda}_j, I_j^c)}_{I_j},\quad j \in \{1,2\}.
	\end{equation}

	\item Define $$\hat{\lambda} = \argmin_{\lambda \in \Lambda} \|\hat{\theta}^{(\lambda)} - \tilde{\theta}\|^2,$$ where $\Lambda = \{1, 2, 2^2, 2^3,\dots, 2^{N^*}\}.$
	%the original estimator $\hat{\theta}^{(\lambda)}$ with $\lambda = \hat{\lambda}$ as defined below: 
	Finally, our estimator (CVSVT) $\hat{\theta}_{CVSVT}$ is defined such that
	\begin{equation}\label{eq:defn}
	\hat{\theta}_{CVSVT} = \hat{\theta}^{(\hat{\lambda})}.
	\end{equation}
	\vspace{0.2in}
	%\item Finally, our estimator (CVDCART) $\hat{\theta}$ is defined such that
	%\begin{equation}\label{eq:cvtfdefn}
	%\hat{\theta}_{I_j} = \hat{\theta}^{(\hat{\lambda}_j, I_j^c)}_{I_j},\quad j \in \{1,2\}.
	%\end{equation}
	\end{comment}
\end{enumerate}
\begin{remark}
	The main difference with  Dyadic CART and Trend Filtering here is in the way we construct the completion estimator. We essentially construct $\tilde{y}(I)$ as an unbiased estimator of $\theta^*$ by randomly doubling or zeroing out each entry of $y$ and then perform SVT on $\tilde{y}(I)$ to create our completion version of the SVT estimator. This idea of randomly zeroing out and inflating other entries to preserve unbiasedness is not new and appears in several matrix completion papers. We call this particular method of creating completion estimators as the zero doubling method. This method is quite generic and can be used for several other signal denoising methods, see Section~\ref{sec:methods} for more on this. %In particular, this could have been used for Dyadic CART and Trend Filtering as well, see Section~\ref{sec:methods} for more on this. 
	\end{remark}

%{\color{red} maybe no need to mention the final steps as they are the same.}

%{\color{red} Say what is $N^*$? Make the notation $\hat{\theta}_{DC}?$ I have written something regarding $N^*$ and changed $\hat{\theta}$ to $\hat{\theta}_{DC}$.}

\subsection{\textbf{Main Result}}
\begin{theorem}\label{thm:cvsvt}
Fix any $\delta > 0$. There exists an absolute constant $C > 0$ such that if we set our grid $\Lambda_j = \Lambda = \{1, 2, 2^2, \dots, 2^{N^*}\}$ for $j = \{1, 2\}$, satisfying
$$2^{N^*} > C \left(|\theta^*|_{\infty} + \sigma\right) \sqrt{n},$$ then the following inequality holds with probability atleast $1 - 2 \exp(-n) - \delta$,
\begin{equation*}
\frac{\|\hat{\theta}_{CVSVT} - \theta^*\|^2}{n^2} \leq \frac{C}{n^2} \left[\sum_{j = 1}^{n} \min\{n \left(|\theta^*|_{\infty} + \sigma\right)^2, \sigma_j^2(\theta^*)\}\right] + \frac{C}{n^2} \sigma^2\log(N^*/\delta).
\end{equation*}
\end{theorem}

We now make some remarks about this theorem.

	\begin{remark}
    The above theorem ensures that our CVSVT estimator also enjoys the adaptive risk bound given in Theorem~\ref{thm:svt} with the only difference being that $\sigma$ is replaced by the term $|\theta^*|_{\infty} + \sigma.$ In realistic scenarios, the term $|\theta^*|_{\infty}$ should remain bounded even if $n$ grows. Therefore, our CVSVT estimator essentially inherits all the implications of Theorem~\ref{thm:svt} for various structured subclasses of matrices.
\end{remark}
	
	\begin{remark}
As mentioned before, setting $N^*$ so that $2^{N^*} > C \left(|\theta^*|_{\infty} + \sigma\right) \sqrt{n}$ is a mild requirement as it is not usually difficult to set upper bounds on the values of $\sigma$ and $|\theta^*|_{\infty}$. Note that, $N^*$ would scale logarithmically in $n$, $\sigma$ and $|\theta^*|_{\infty}$. 
\end{remark}

\section{Discussion}\label{sec:discuss}
In this section we discuss some naturally related matters.

\subsection{Other Signal Denoising Methods}
In this paper, we applied our CV framework to produce CV versions of Trend Filtering, Dyadic CART, Singular Value Thresholding, and Lasso (see Section~\ref{sec:lasso}) %~\ref{sec:lasso} in the supplementary file) 
with the main focus being on the first two estimators whereas the latter two are considered to further illustrate the generality of the framework. Our CV framework is based on a general principle and should be looked upon as providing a general recipe to develop \textit{theoretically tractable} CV versions of potentially any other estimator that uses a tuning parameter, for example, the Total Variation Denoising estimator proposed by~\cite{rudin1992nonlinear} (also see~\cite{hutter2016optimal},~\cite{sadhanala2016total},~\cite{chatterjee2019new}), the Hardy Krauss estimator (see~\cite{fang2021multivariate},~\cite{ortelli2020adaptive}), the Optimal Regression Tree estimator proposed in~\cite{chatterjee2019adaptive}, a higher dimensional version of Trend Filtering of order $2$ proposed in~\cite{ki2021mars} and many more. As a starting point, considering the Zero Doubling method for constructing the completion estimators along with a geometrically doubling grid of candidate tuning values should be useable in these problems.

\subsection{Three Different Methods for Creating Completion Estimators}\label{sec:methods}

One of the main ingredients of the proposed CV framework is the construction of the completion estimators, where the user needs to build a version of the estimator of interest depending only on a subset of the data, namely $y_{I}.$ In this paper, we have considered three different strategies for constructing the completion estimators. In Dyadic CART and Lasso, we restrict the squared error term in the optimization objective to only be summed over the subset $I$, and let us call this method \textit{Restricted Optimization} (RO). In Trend Filtering, we first construct an interpolated data vector $\tilde{y}(I) \in \R^n$ by interpolating on the subset of indices $I^c$ based on $y_{I}$, and then feed $\tilde{y}(I)$ into the Trend Filtering optimization objective; let us call this method \textit{Interpolate then Optimize} (IO). In Singular Value Thresholding, we zero out entries in $I^c$ and double the entries in $I$ to create a new matrix $\tilde{y}(I)$, and then use the singular value thresholding operator on $\tilde{y}(I)$; let us call this method \textit{Zero Doubling} (ZD). To summarize, RO, IO and ZD are three different ways to construct completion estimators, among which the user can try any one or even come up with a different method in some given problem. For example, the ZD method is extremely generic and could have been applied in Trend Filtering or Dyadic CART as well, however, we did not use it because in our numerical experiments we observed that RO (resp. IO) was performing better (in MSE) than ZD for Dyadic CART (resp. Trend Filtering) by a factor of $2$ or $3$. %Similarly, IO performed better than ZD in the case of Trend Filtering. %It might be interesting to undertake a rigorous comparison of these different methods for a particular estimator although it is not clear to us how to do so. 

\subsection{Comparison with the R Package~\cite{arnold2020package} Cross Validation Version of Trend Filtering}\label{sec:compare}
It is instructive to compare the proposed CVTF estimator with the CV version of Trend Filtering implemented in the R package ~\cite{arnold2020package} (see the R command \textit{cv.trendfilter}). In particular, it is worth noting the similarities and the differences in the two CV algorithms. Simulations comparing the finite sample performance of both these algorithms is given in Section~\ref{sec:simu}.

\begin{enumerate}
	\item We construct the folds in the same way as in the R package: after deciding the number of folds $K$, the folds are created by placing every $K$th index into the same fold. However, unlike the R package, for $r$th order Trend Filtering we specifically set $K  = r + 1$ to enable our interpolation scheme to create $\tilde{y}(I_j^c)$, $j \in [K]$ in such a way that further allows us to obtain the two inequalities stated in Step $3$ of the proof sketch in Section~\ref{sec:tvsketch}. %~\ref{sec:tvsketch}.

	%\medskip
	
	\item In the R package, the predictions for a given fold $I_j$ is made by first performing trend filtering on the shortened data vector $y_{I_j^c} \in \R^{|I_j^c|}$ and then at any point in $I_j$, the prediction is given by averaging the fits at its two neighbors (guaranteed to be in a different fold). However, we follow the reverse order by first applying a specific polynomial interpolation scheme on $y_{I_j^c}$, and then performing trend filtering to this interpolated data vector $\tilde{y}(I_j^{c}) \in \R^n$, and finally obtain the completion estimators $\hat{\theta}^{(\lambda, I_j^{c})}$, based on which our prediction at any point in $I_j$ is given by the fit $\hat{\theta}^{(\lambda, I_j^{c})}$ at that point.

	\item The main point of difference of the two methods is how they choose the final data driven value of the tuning parameter $\hat{\lambda}.$ The R package implementation uses~\eqref{eq:cvtypical} to choose $\hat{\lambda}$ whereas our method is different as has been explained before.
	
	\item The grid of candidate tuning values $\Lambda_j$ and $\Lambda$ are chosen in a fully data driven way in the R package implementation as explained in Section~\ref{sec:alttf} below. However, we prefer to simply set $\Lambda_j = \Lambda = \{1,2,2^2,\dots,2^{N^*}\}$ for a large enough $N^*$ such that $2^{N^*} = O(n).$ We could also mimick the R package implementation in choosing $\Lambda_j$ and $\Lambda$ and this will be perfectly in accordance with our CV framework as also explained in Section~\ref{sec:alttf}.
		
\end{enumerate}

\subsubsection{An Alternative Way to Construct $\hat{\theta}_{CVTF}^{(r)}$}\label{sec:alttf}
Recall that the last few steps to construct the CVTF estimator is identical to the corresponding steps to construct the CVDCART estimator. However, unlike Dyadic CART, Trend Filtering is based on convex optimization which brings with it its inherent advantages. For the discussion below let us fix $r = 1.$ %It is known (see~\cite{}) that the entire path of solutions (for all $\lambda > 0$) can be computed for Trend Filtering of any order $r \geq 1.$ 
Infact, (see~\cite{hoefling2010path}) for Trend Filtering of order $1$, also known as Fused Lasso, the solution (as a function of $\lambda$) is piecewise constant with a finite number of pieces. Moreover, the entire path of solutions (for all $\lambda > 0$) can be computed in $O(n \log n)$ time and the number of distinct solutions is always bounded by $O(n).$ 
%It is known (see~\cite{}) that the entire path of solutions (for all $\lambda > 0$) can be computed for Trend Filtering of any order $r \geq 1.$ Moreover, the solution (as a function of $\lambda$) is piecewise linear and convex with a finite number of knots. Infact, the Rpackage implementing CV for Trend Filtering makes use of the above fact. The grid of candidate tuning parameter values in this package is simply taken to be the (finite) set of knots in the entire path of solutions. 
In the Rpackage, CV for Fused Lasso is implemented by using these facts, and the grid of candidate tuning parameters is simply chosen to be a (random and finite) set of tuning values $\lambda$'s that correspond to the set of all possible solutions.

We can mimick the Rpackage implementation in the last few steps and still stay within our CV framework which gives us a different CV version of Fused Lasso. This is because we can define $\Lambda_j$, $j \in [K] = [r+1] = [2]$, to be a finite set of tuning values, one for each of the distinct solutions of the following optimization problem 
$$\min_{\theta \in \R^n} \frac{1}{2}\left\|\tilde{y}(I_j^c) - \theta\right\|^2 + \lambda n^{r-1}\left\|D^{(r)}\theta\right\|_1.$$
This is allowed in our framework because this set only depends on $\tilde{y}(I_j^c)$ by definition. Similarly, we can define $\Lambda$ to be a set of tuning values, one for each of the distinct solutions of the full optimization problem 
$$\min_{\theta \in \R^n} \frac{1}{2}\left\|y - \theta\right\|^2 + \lambda n^{r-1}\left\|D^{(r)}\theta\right\|_1.$$
Under these choices of $\Lambda_j,\Lambda$, our slow rate and fast rate theorems are still valid because of the following. 
In view of Theorem~\ref{thm:main}, we need to bound $\min_{\lambda \in \Lambda} SSE(\hat{\theta}^{(r)}_{\lambda}, \theta^*)$ and $\min_{\lambda \in \Lambda_j} SSE(\hat{\theta}^{(\lambda,I,r)}_{I^c}, \theta^*_{I^c})$, where $I = I_j$ for $j \in [K]$. Note that, $\min_{\lambda \in \Lambda} SSE(\hat{\theta}^{(r)}_{\lambda}, \theta^*) = \min_{\lambda \in \R} SSE(\hat{\theta}^{(r)}_{\lambda}, \theta^*)$ and thus, we can use known bounds for the ideally tuned versions. Bounding $\min_{\lambda \in \Lambda_j} SSE(\hat{\theta}^{(\lambda, I, r)}_{I^c}, \theta^*_{I^c})$ can be similarly accomplished for this data driven choice of $\lambda_j$ by again noting that, $\min_{\lambda \in \Lambda_j} SSE(\hat{\theta}^{(\lambda,I,r)}_{I^c}, \theta^*_{I^c}) = \min_{\lambda \in \R} SSE(\hat{\theta}^{(\lambda,I,r)}_{I^c}, \theta^*_{I^c})$, and then simply following our existing proof. The only point to further consider would be the term involving $\log |\Lambda_j|$ because now it is random unlike previous. Here, we can invoke the result of~\cite{hoefling2010path} and use a deterministic bound (scaling like $O(n))$ on $|\Lambda_j|.$ Thus, the term involving $\log |\Lambda_j|$ can be bounded by $O({\log n}/{n})$, a lower order term.

The advantage of this particular variant is that the entire procedure is fully data driven and one does not even need to set the value of $N^*$ as before. This furnishes a truly completely data driven cross validated Fused Lasso estimator which attains both the slow rate and the fast rate. To the best of our knowledge, such a version of Fused Lasso did not exist in the literature before our work here.

One can also use this approach and choose $\Lambda_j$ and $\Lambda$ similarly, for Trend Filtering of any general order $r \geq 1$. This is because it is known (see Section $6.2$ in~\cite{tibshirani2011solution}) that the entire path of solutions (for all $\lambda > 0$) can again be computed for Trend Filtering of any order $r \geq 1.$ Moreover, the solution (as a function of $\lambda$) is piecewise linear and convex with a finite number of knots. Therefore, one can simply take $\Lambda_j, j \in [K]$ and $\Lambda$ to be the finite set of knots of the appropriate optimization problems. This would produce fully data driven CV (within our framework) versions of Trend Filtering of general order $r \geq 1$, and to obtain a theoretical guarantee one can then use Theorem~\ref{thm:main}. The only missing part is that a deterministic bound on $\log |\Lambda_j|$ is not known for $r > 1$, however, from numerical experiments we are led to conjecture that the number of distinct Trend Filtering solutions (of any order, w.r.t $\lambda$) grows at most polynomially with $n$. If this conjecture were true, we can then again conclude that the term involving $\log |\Lambda_j|$ is of a lower order term. We prefer to write our theorem for the current version because a) we have a complete proof of a risk bound for all orders $r \geq 1$, b) practically, the parameter $N^*$ is not difficult to set and in our simulations both these versions perform similarly when we set $N^* = \log_2 n.$

%\subsection{Square Root Trend Filtering}

\subsection{Heavy Tailed Errors}
The proof of our main result in Theorem~\ref{thm:adapt} relies heavily on the errors being subgaussian. It would be interesting to explore the robustness of our CV framework to heavy tailed errors. In particular, can one develop CV algorithms for corresponding quantile versions of Dyadic CART and Trend Filtering (see~\cite{hernan2021risk} and~\cite{padilla2021quantile}) using our framework? We leave this question for future research.

%\subsection{Other CV Methods}

\section{Simulations}\label{sec:simu}

\subsection{Dyadic CART}

We conduct a simulation study to observe the performance of the proposed CVDCART estimator in three different scenarios each corresponding to a different true signal $\theta^*$ with dimension $d = 2$. In every scenario, we vary the sample sizes $n = 128, 256, 512$, generate the errors from $N(0, 1)$, and estimate the MSE by $100$ Monte Carlo replications. Moreover, these results are compared with the traditional CV version of Dyadic CART, where we consider the same exact folds and the same completion estimator as in CVDCART, and choose $\hat{\lambda}$ according to~\eqref{eq:cvtypical}. The results are presented in Table~\ref{tab:dc}, where the Monte Carlo standard errors in the estimation of the MSEs are reported in the parentheses next to the corresponding estimates, and we observe that the performance of CVDCART is almost similar to the traditional version with a minute difference of order $10^{-3}$.  The plots corresponding to Scenarios $1$, $2$ and $3$ are provided in Figures~\ref{fig:dc_box}, \ref{fig:dc_circle} and \ref{fig:dc_sinu} respectively. %where the first diagram refers to the true signal, the second one to the noisy signal and the third one to the estimated signal by CVDCART. 
\begin{enumerate}
\item Scenario 1 [Rectangular Signal]: The true signal $\theta^*$ is such that for every $(i_1, i_2) \in L_{2, n}$, %we have
$$\theta^*_{(i_1, i_2)} = 
\begin{cases}
1 \quad &\text{if} \;\; n/3 \leq i_1, i_2 \leq 2n/3\\
0 \quad &\text{otherwise}
\end{cases}.$$
%The corresponding plots are shown in Figure~\ref{fig:dc_box} when $n = 256$.
\begin{figure}[H]
\centering
\includegraphics[scale = 0.25]{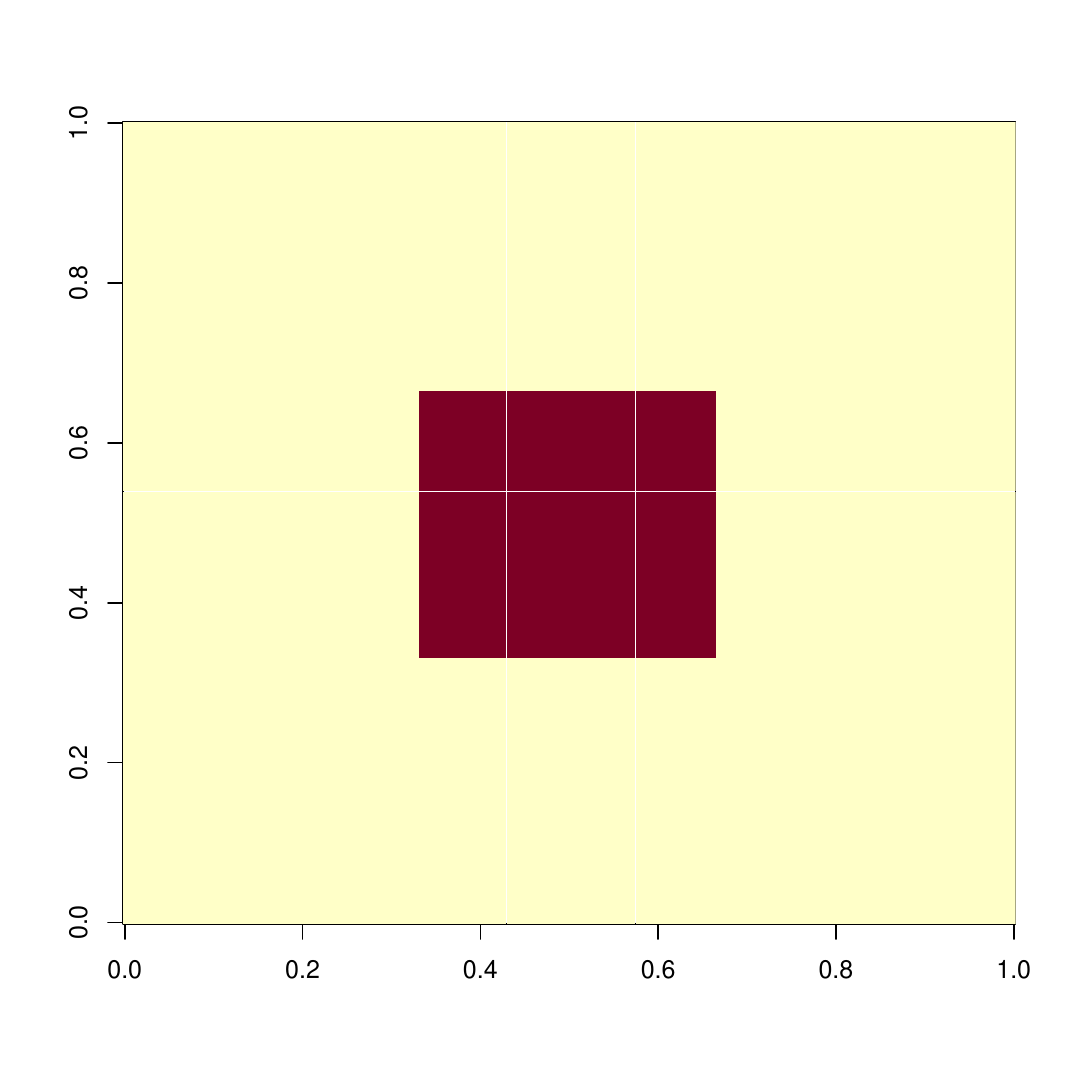} \includegraphics[scale = 0.25]{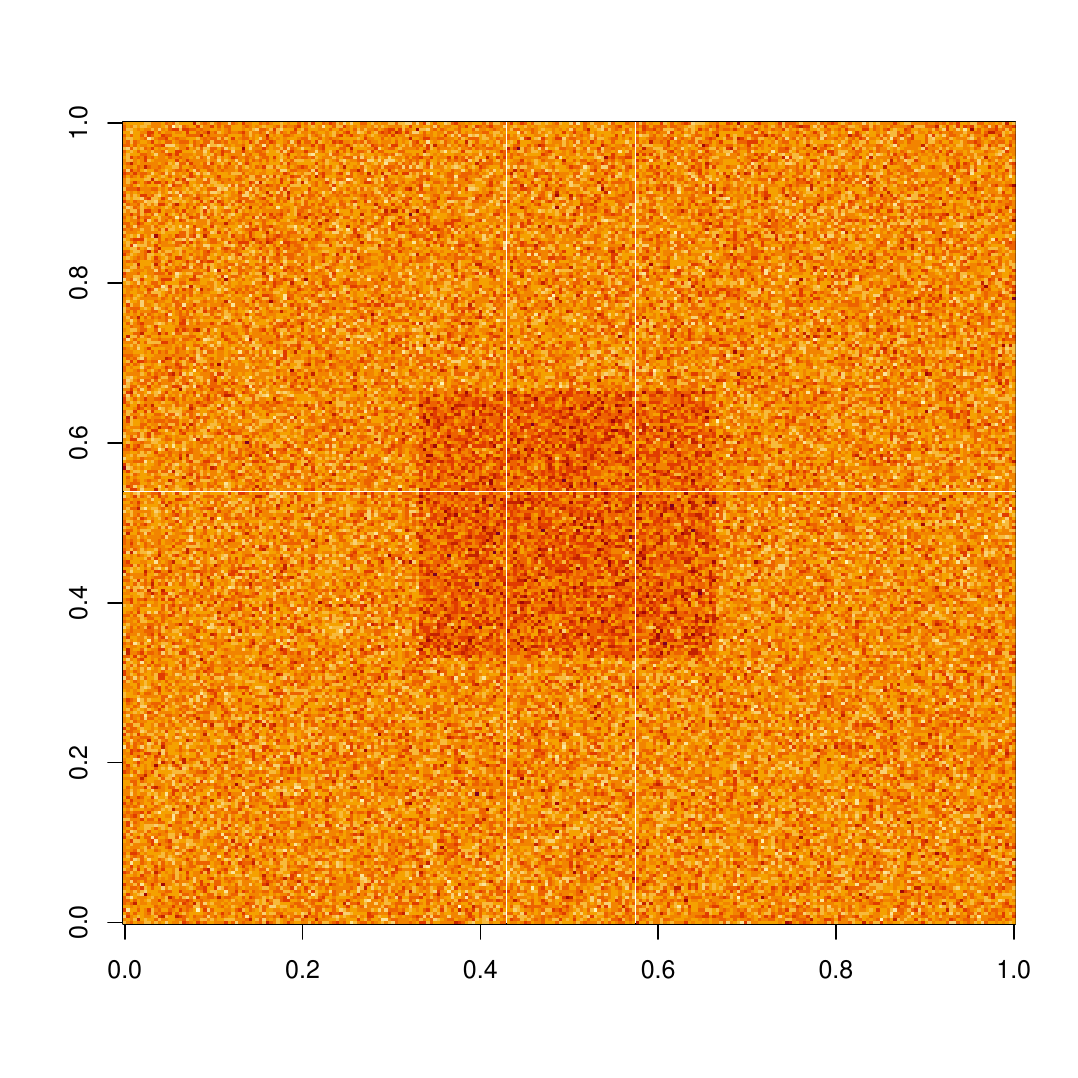} \includegraphics[scale = 0.25]{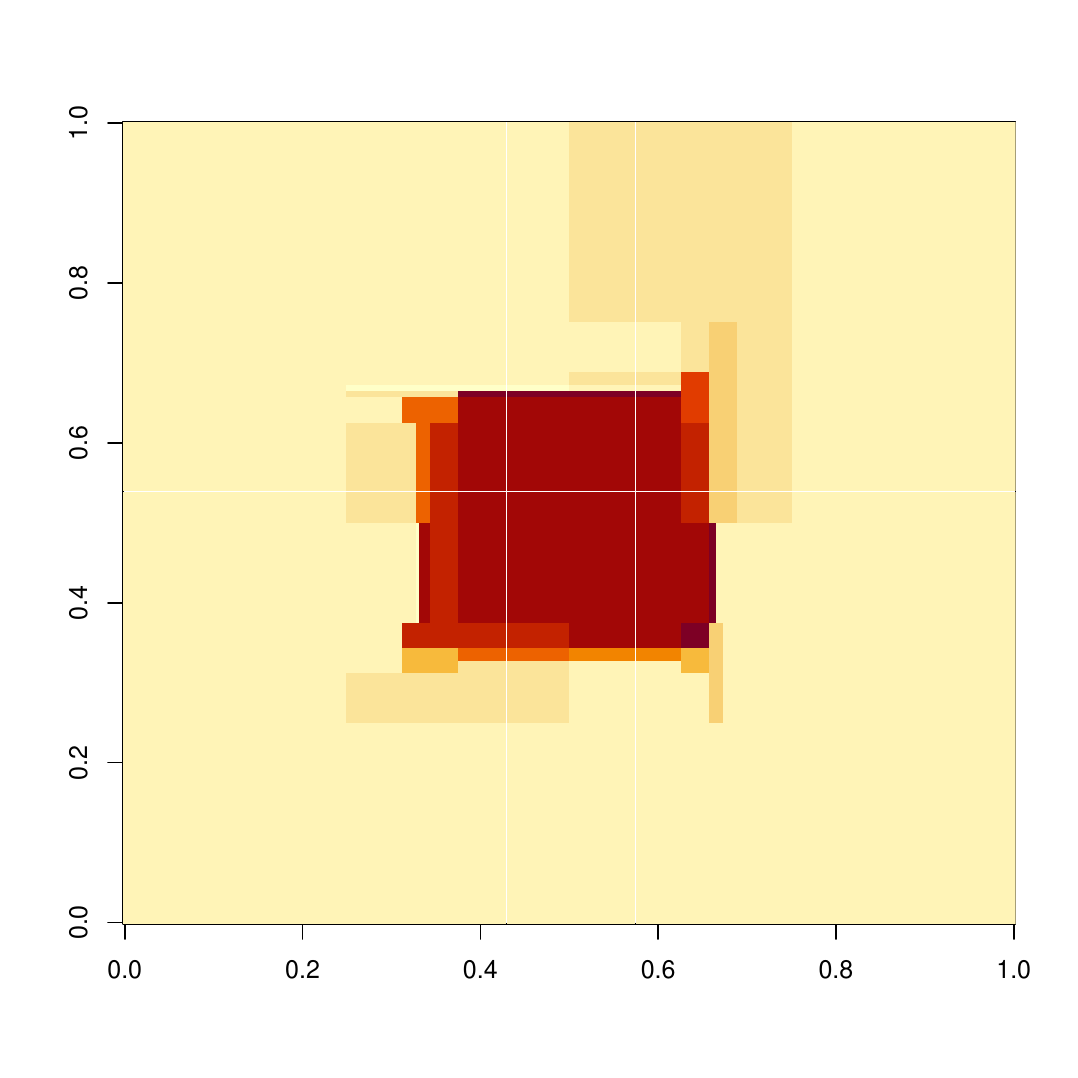}
%\caption{The first diagram refers to the true signal, the second one to the noisy signal and the third one to the estimated signal by the CV Dyadic CART estimator.}
%\caption{Estimation in Scenario 1.}
\caption{The true signal, the noisy signal and the estimated signal by CVDCART.}
\label{fig:dc_box}
\end{figure}

\item Scenario 2 [Circular Signal]: The true signal $\theta^*$ is such that for every $(i_1, i_2) \in L_{2, n}$, %we have
$$\theta^*_{(i_1, i_2)} = 
\begin{cases}
1 \quad &\text{if} \;\; \sqrt{(i_1 - n/2)^2 + (i_2 - n/2)^2} \leq n/4\\
0 \quad &\text{otherwise}
\end{cases}.$$
%The corresponding plots are shown in Figure~\ref{fig:dc_circle} when $n = 256$.
\begin{figure}[H]
\centering
\includegraphics[scale = 0.25]{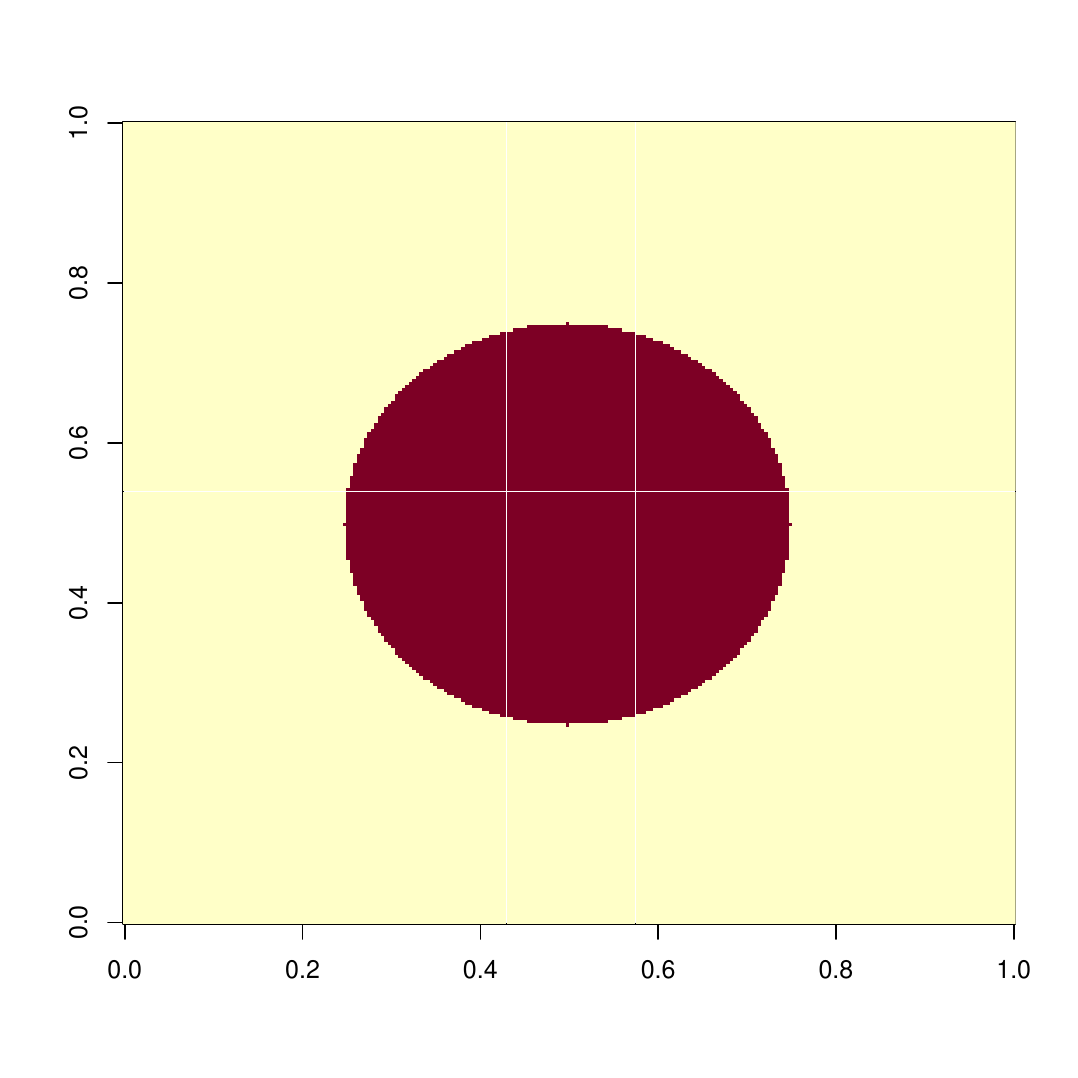} \includegraphics[scale = 0.25]{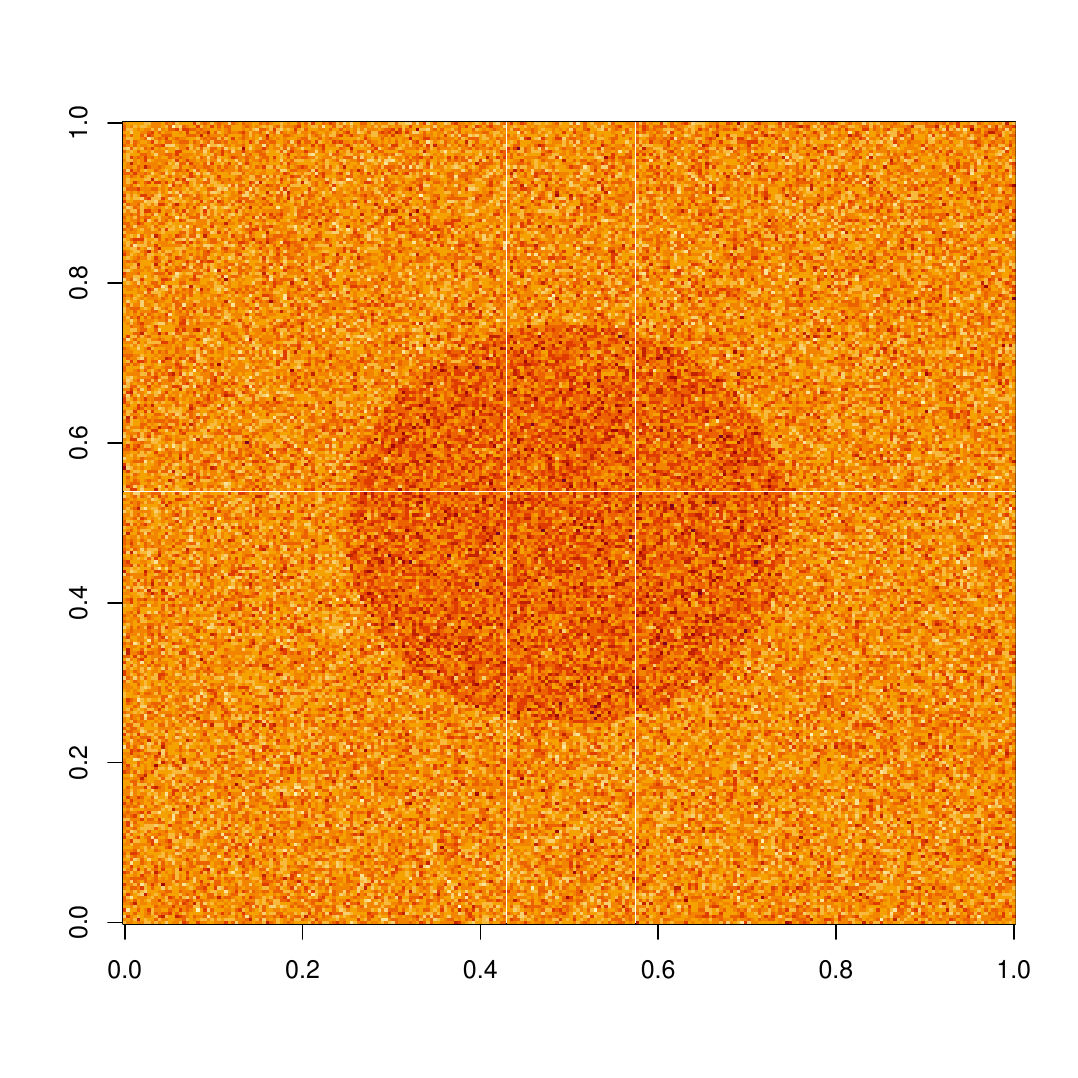} \includegraphics[scale = 0.25]{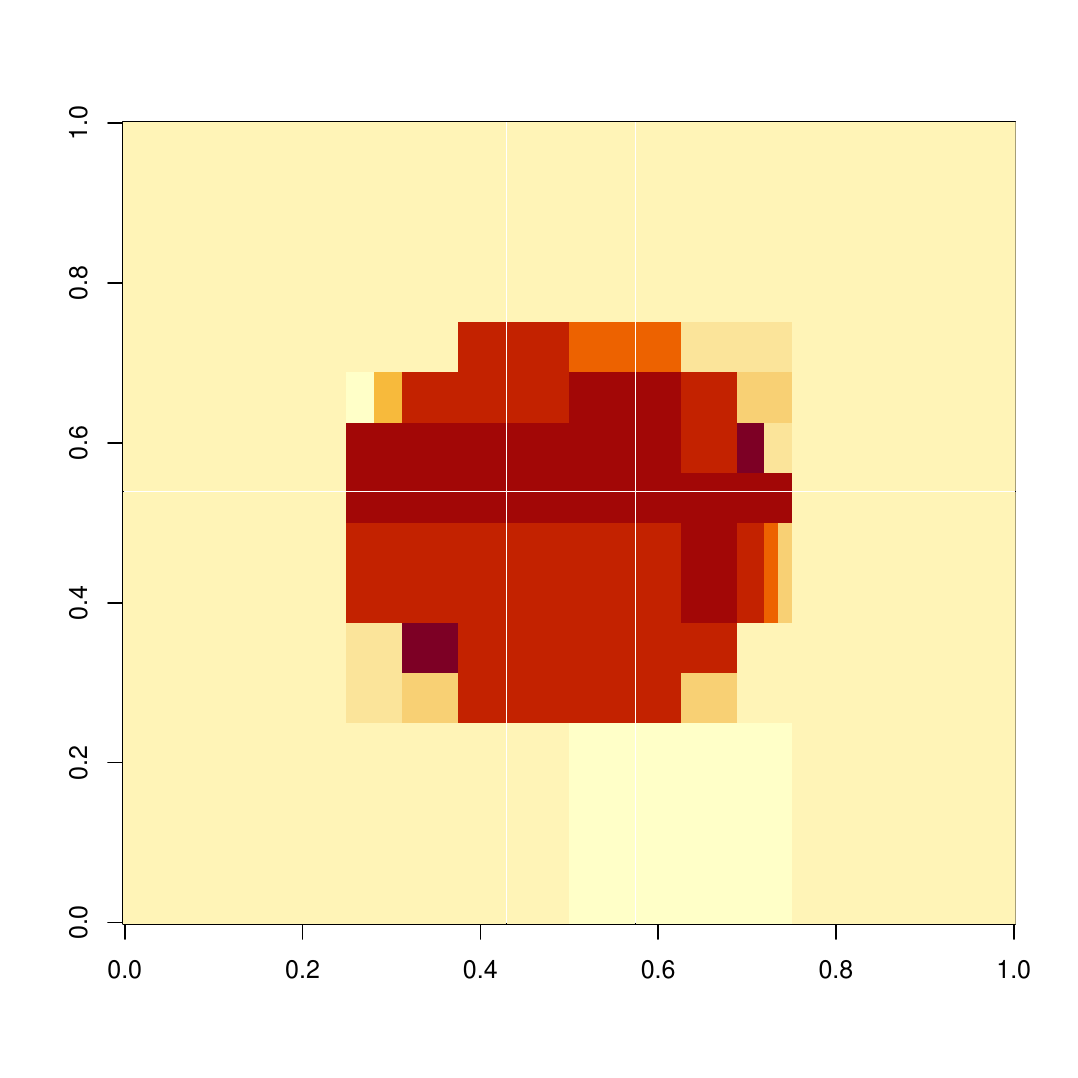}
%\caption{Estimation in Scenario 2.}
\caption{The true signal, the noisy signal and the estimated signal by CVDCART.}
\label{fig:dc_circle}
\end{figure}

\item Scenario 3 [Smooth Signal]: The true signal $\theta^*$ is such that for every $(i_1, i_2) \in L_{2, n}$, we have $\theta^*_{(i_1, i_2)} = f\left(i_1/n, i_2/n\right)$, where
$$f(x, y) = 20\exp\left(-5\{(x - 1/2)^2 + (y - 1/2)^2 - 0.9(x - 1/2)(y - 1/2)\}\right), \;\; 0 \leq x, y \leq 1.$$
%(\mathrm{1}(x \in [0, 1/3])) + -36(x-1/2-1/\sqrt{12})(x-1/2+\sqrt{12})\mathrm{1}(x \in [1/3, 2/3]) + 18(x - 1)^2(\mathrm{1}(x \in [2/3, 1])),$$
%The corresponding plots are shown in Figure~\ref{fig:dc_sinu} when $n = 256$.
\begin{figure}[H]
\centering
\includegraphics[scale = 0.25]{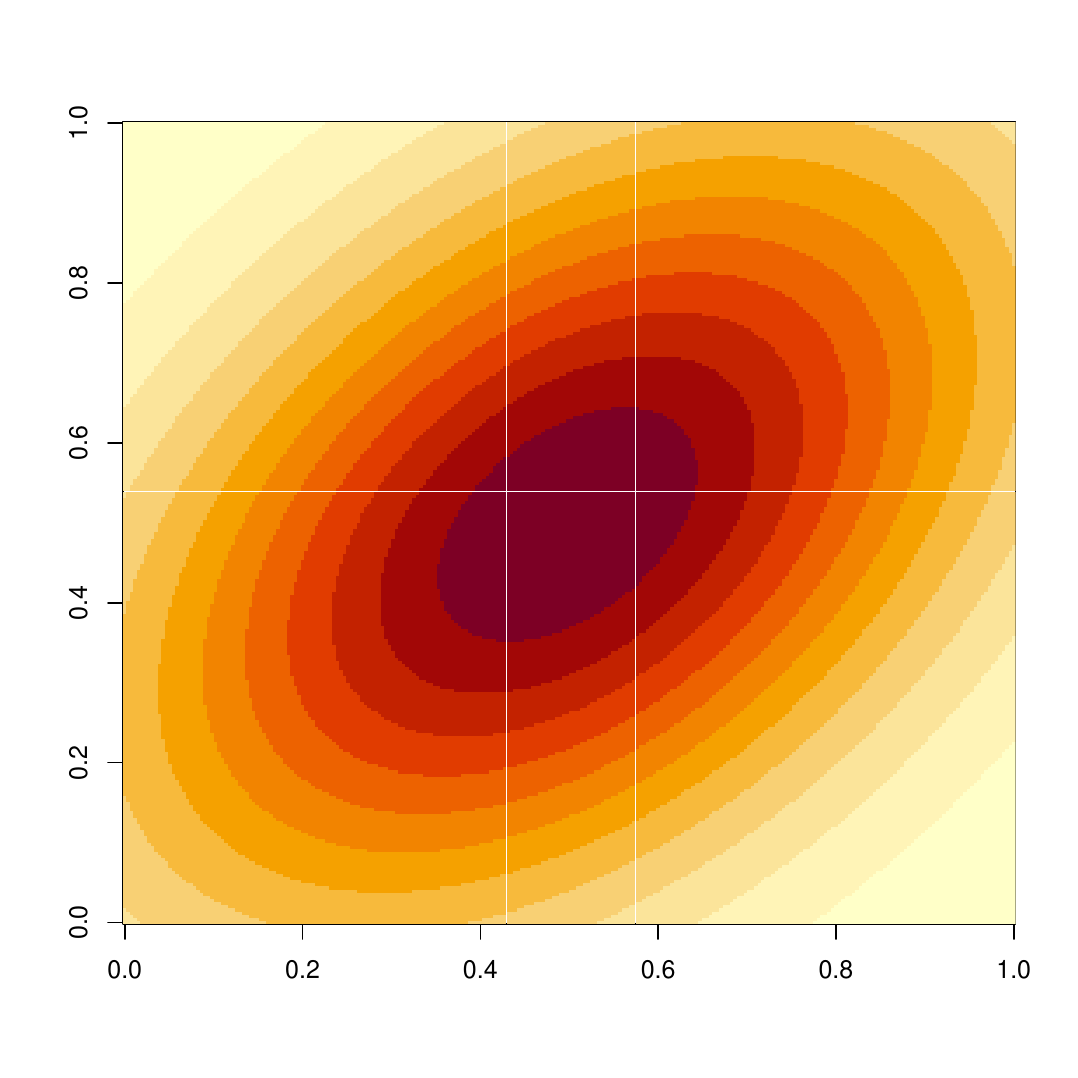} \includegraphics[scale = 0.25]{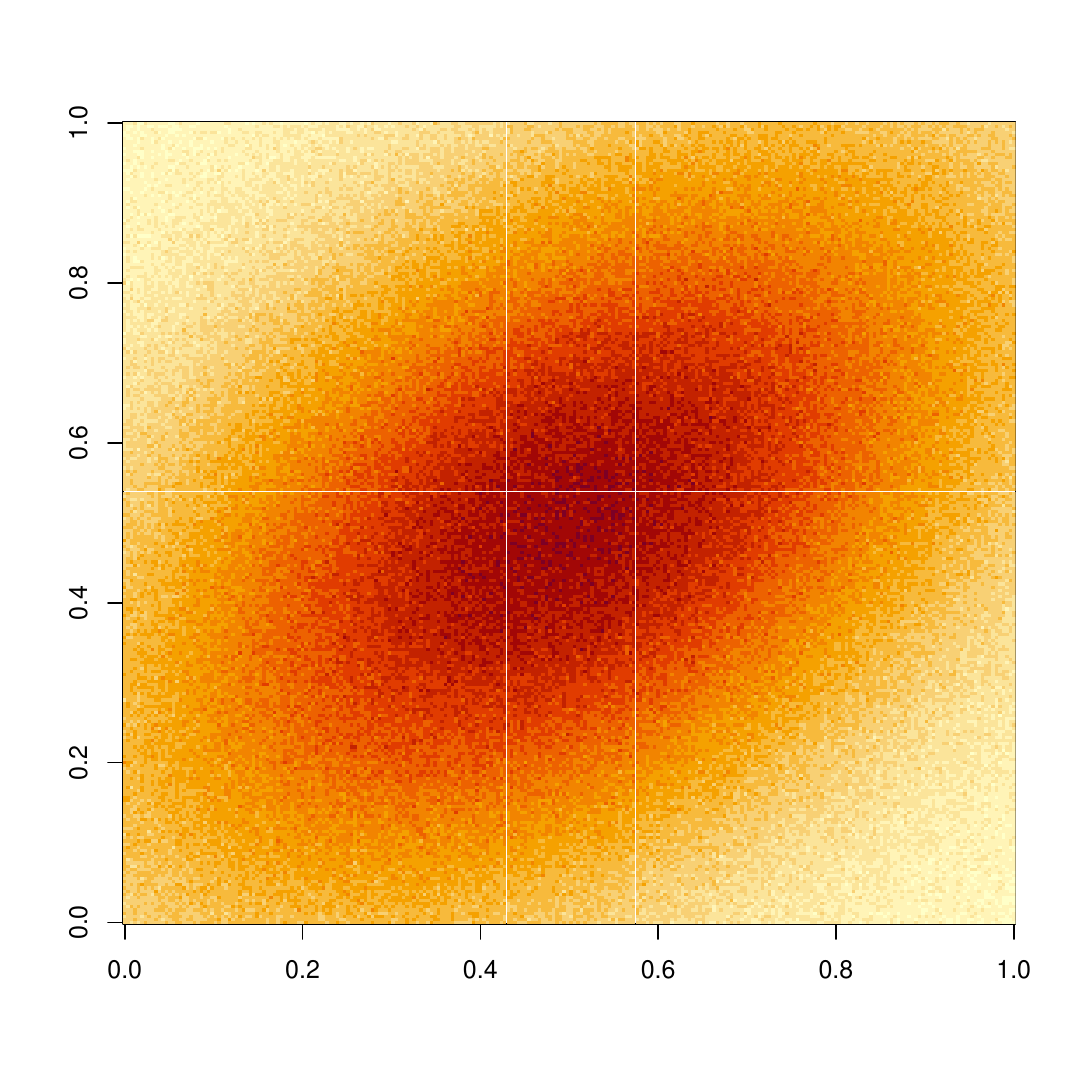} \includegraphics[scale = 0.25]{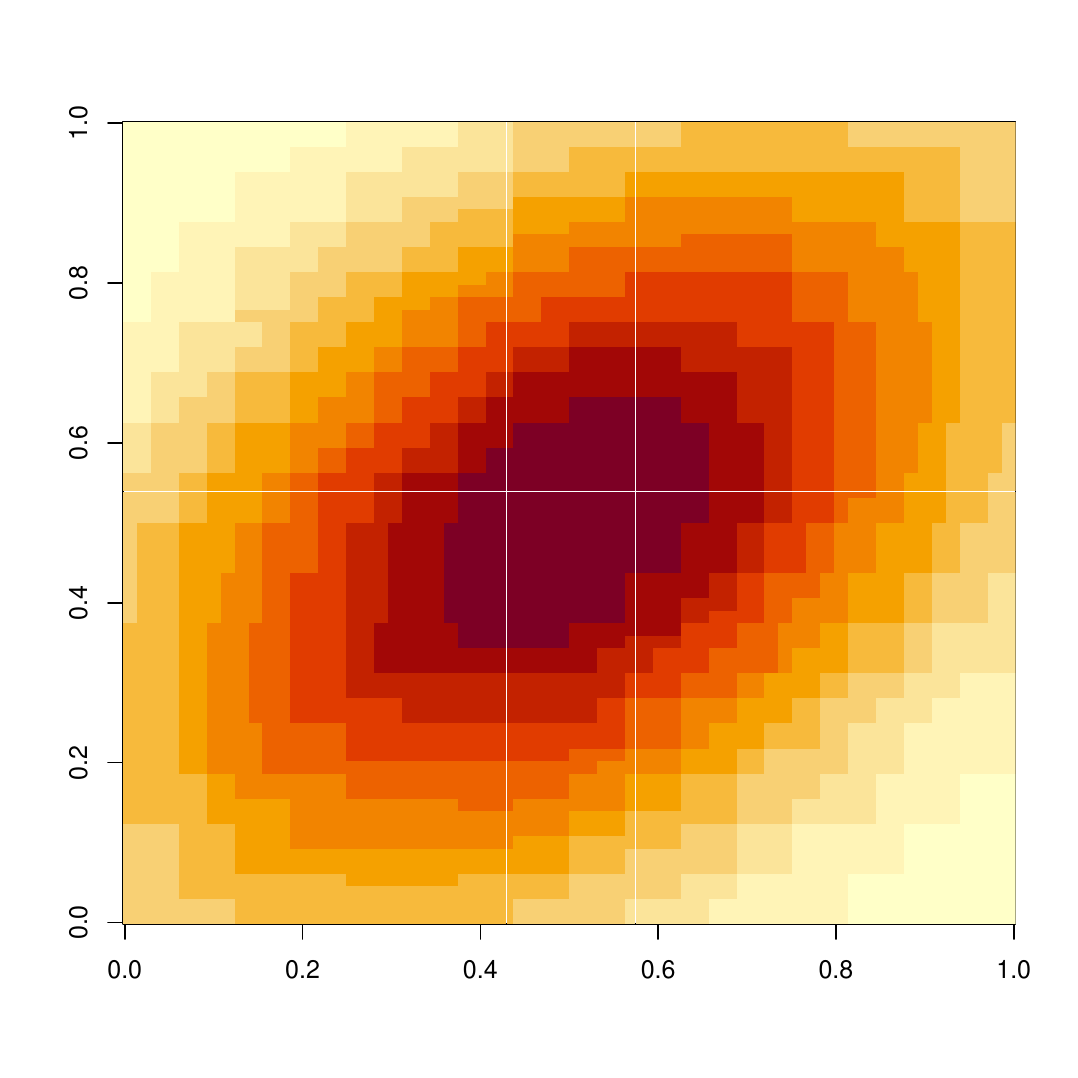}
%\caption{The first diagram refers to the true signal, the second one to the noisy signal and the third one to the estimated signal by the CV Dyadic CART estimator.}
\caption{The true signal, the noisy signal and the estimated signal by CVDCART.}
\label{fig:dc_sinu}
\end{figure}
\end{enumerate}
\begin{table}[H]
\caption{Comparison of MSEs between the traditional CV and the proposed CVDCART estimator in different scenarios}
\label{tab:dc}
\begin{tabular}{c|cc|cc|cc}
\hline
&\multicolumn{2}{c|}{Scenario 1}  &\multicolumn{2}{c|}{Scenario 2} & \multicolumn{2}{c}{Scenario 3}\\
\hline
$n$ &Traditional CV                    & CVDCART  &Traditional CV                     & CVDCART  &Traditional CV                     &CVDCART\\
\hline
128 & 0.011 (0.0013) & 0.019 (0.0032) & 0.018 (0.0014) & 0.021 (0.0018) & 0.0007 ($10^{-4}$) & 0.0005 ($10^{-4}$) \\
256 & 0.003 (0.0004) & 0.005 (0.0009) & 0.009 (0.0011) & 0.014 (0.0008) & 0.0004 ($10^{-5}$) & 0.0004 ($10^{-5}$) \\
512 & 0.001 (0.0001) & 0.001 (0.0003) & 0.005 (0.0001) & 0.006 (0.0003) & 0.0003 ($10^{-5}$) & 0.0003 ($10^{-5}$) \\
%1024 & 0.023 & 0.022 & 0.007 \\
\hline
\end{tabular}
\end{table}

\subsection{Trend Filtering}

We conduct a simulation study to observe the performance of the proposed CVTF estimator and compare it with three different model selection methods: the traditional CV method (with the same number of folds, the same completion estimator as in the proposed CVTF, but $\hat{\lambda}$ is chosen according to~\eqref{eq:cvtypical}), the CV Trend Filtering method implemented in the R package \textit{genlasso}~\cite{arnold2020package}, and the model selection by Stein's Unbiased Risk Estimation (SURE) for Trend Filtering, see~\cite{dfLasso}. In particular, to implement SURE for Trend Filtering we use \textit{trendfilter} command in the R package \textit{genlasso}~\cite{arnold2020package}.  

%{\color{red} description of SURE.}

\begin{figure}[H]
\centering
\includegraphics[scale = 0.26]{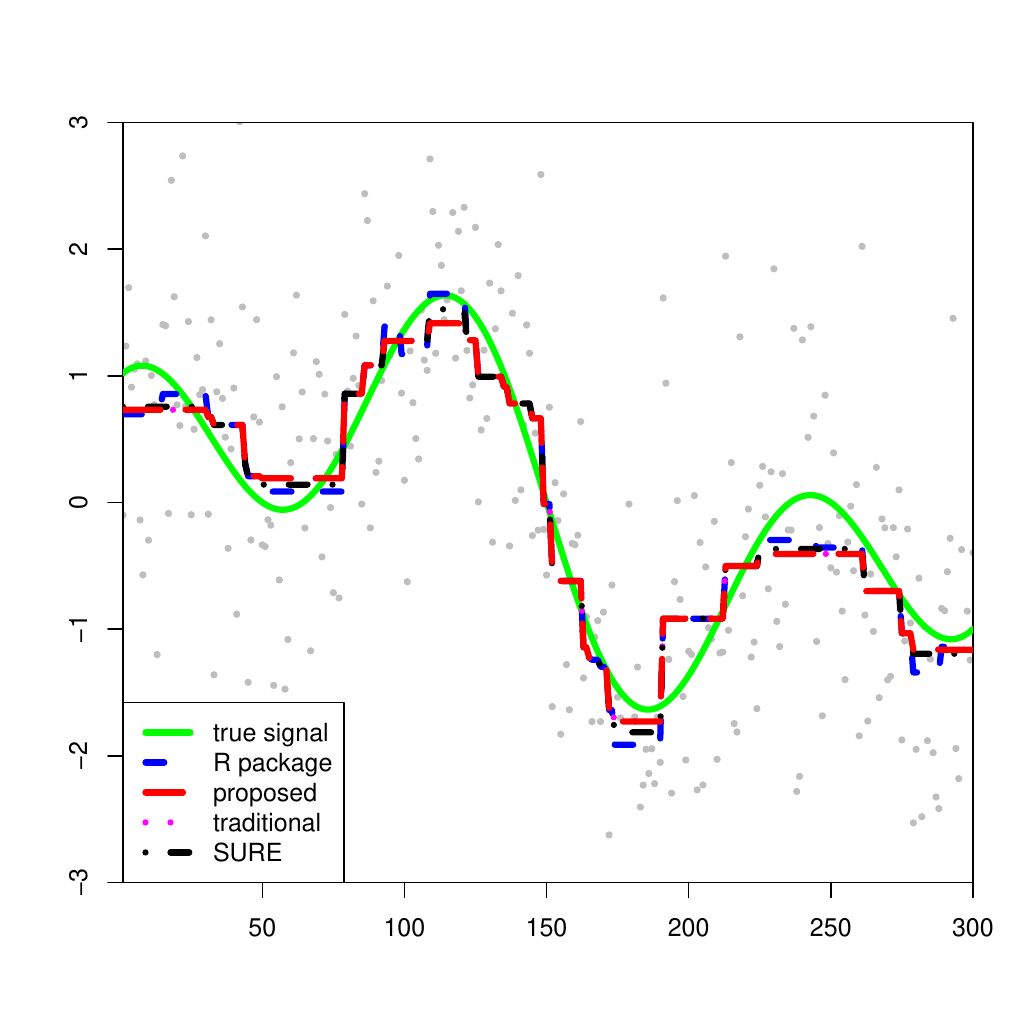} \includegraphics[scale = 0.26]{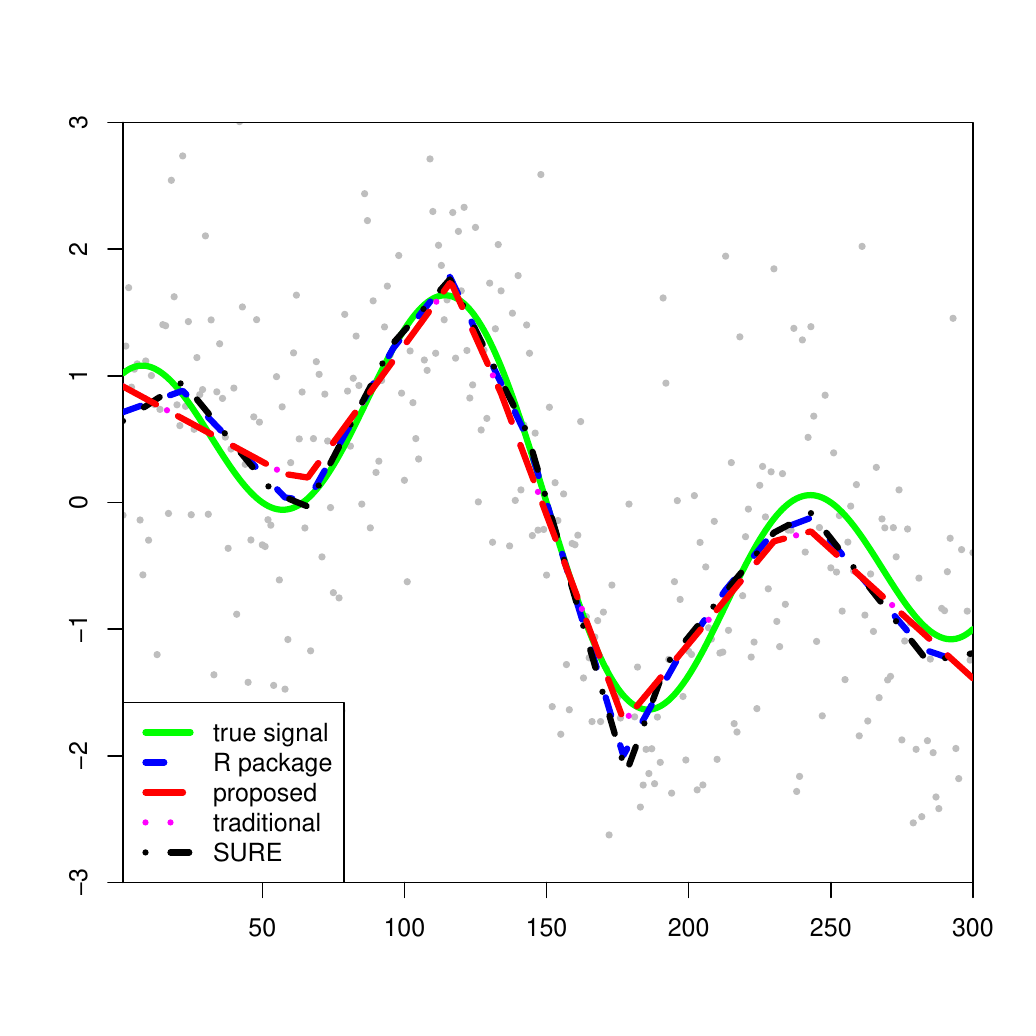} \includegraphics[scale = 0.26]{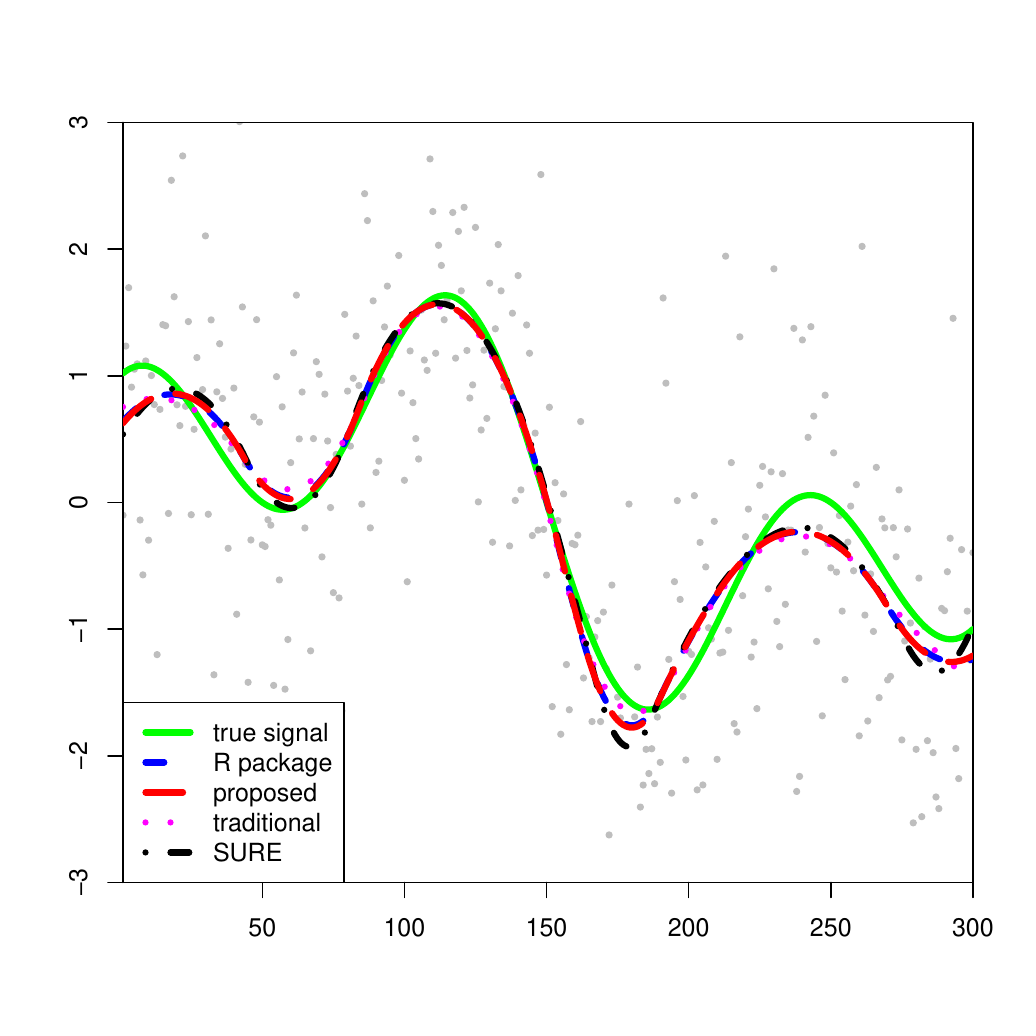}
\caption{The fits of order 1, the fits of order 2 and the fits of order 3 in Scenario 4.}
\label{fig:tv2}
\end{figure}

The studies are carried out in four different scenarios each corresponding to a different true signal $\theta^* \in \R^n$, where for any $i \in [n]$, we have $\theta^*_i = f(i/n)$ for some function $f : [0,1] \to \R$, specified below. In every scenario, we vary the sample sizes $n = 300, 600, 1200, 2400$, generate the errors from $N(0, 1)$, and estimate the MSE by 100 Monte Carlo replications. Furthermore, the comparison of MSE between the aforementioned model selection methods in Scenario $1$, $2$ and $3$ are presented in Table~\ref{tab:tv1}, \ref{tab:tv2} and \ref{tab:tv3} respectively, where the Monte Carlo standard errors in the estimation of the MSEs are reported in the parentheses next to the corresponding estimates. It is important to note that, implementation of SURE requires gaussianity of the errors as well as the knowledge of their variance, however, these are not required for the other CV methods. 
\begin{enumerate}
\item Scenario 1 [Piecewise Constant Signal]: We consider the piecewise constant function 
$$f(x) = 2(\mathrm{1}(x \in [1/5, 2/5])) + \mathrm{1}(x \in [2/5, 3/5]) - \mathrm{1}(x \in [3/5, 4/5]) + 2(\mathrm{1}(x \in [4/5, 1])),$$
and apply Trend Filtering of order $r = 1$. %The corresponding plot is shown in the first diagram of Figure~\ref{fig:tv1} when $n = 300$.

\begin{table}[H]
\caption{Comparison of MSE in Scenario 1.}
\label{tab:tv1}
\begin{tabular}{c|cccc}
\hline
% &\multicolumn{2}{c|}{Scenario 1}  &\multicolumn{2}{c|}{Scenario 2} & \multicolumn{2}{c}{Scenario 3}\\
%\hline
$n$ &RPackage  & CVTF  &SURE     & Traditional CV        \\
\hline
300 & 0.081 (0.030) & 0.073 (0.020) & 0.076 (0.030) & 0.074 (0.020)\\
600 & 0.042 (0.010) & 0.040 (0.010) & 0.040 (0.010) & 0.039 (0.010)\\
1200 & 0.022 (0.008) & 0.022 (0.009) & 0.023 (0.007) & 0.022 (0.008)\\
2400 & 0.013 (0.003) & 0.012 (0.003) & 0.014 (0.004) & 0.012 (0.004)\\
\hline
\end{tabular}
\end{table}
\begin{comment}
\begin{table}[H]
\caption{Comparison of MSE in Scenario 1.}
\label{tab:tv1}
\begin{tabular}{c|cccc}
\hline
% &\multicolumn{2}{c|}{Scenario 1}  &\multicolumn{2}{c|}{Scenario 2} & \multicolumn{2}{c}{Scenario 3}\\
%\hline
$n$ &RPackage  & CVTF  &SURE     & Traditional CV        \\
\hline
300 & 0.081 & 0.073 & 0.076 & 0.074\\
600 & 0.042 & 0.040 & 0.040 & 0.039\\
1200 & 0.022 & 0.022 & 0.023 & 0.022\\
2400 & 0.013 & 0.012 & 0.014 & 0.012\\
\hline
\end{tabular}
\end{table}
\end{comment}
\item Scenario 2 [Piecewise Linear Signal]: We consider the piecewise linear function 
$$f(x) = 6x(\mathrm{1}(x \in [0, 1/3])) + (-12x + 6)\mathrm{1}(x \in [1/3, 2/3]) + (x - 8/3)(\mathrm{1}(x \in [2/3, 1])),$$
and apply Trend Filtering of order $r = 2$. %The corresponding plot is shown in the second diagram of Figure~\ref{fig:tv1} when $n = 300$.
\begin{table}[H]
\caption{Comparison of MSE in Scenario 2.}
\label{tab:tv2}
\begin{tabular}{c|cccc}
\hline
% &\multicolumn{2}{c|}{Scenario 1}  &\multicolumn{2}{c|}{Scenario 2} & \multicolumn{2}{c}{Scenario 3}\\
%\hline
$n$ &RPackage  & CVTF  &SURE     & Traditional CV        \\
\hline
300 & 0.029 (0.021) & 0.029 (0.017) & 0.032 (0.023) & 0.028 (0.016)\\
600 & 0.015 (0.010) & 0.013 (0.006) & 0.016 (0.010) & 0.013 (0.006)\\
1200 & 0.007 (0.004) & 0.006 (0.003) & 0.008 (0.004) & 0.007 (0.004)\\
2400 & 0.004 (0.003) & 0.003 (0.002) & 0.004 (0.004) & 0.003 (0.002)\\
\hline
\end{tabular}
\end{table}
\begin{comment}
\begin{table}[H]
\caption{Comparison of MSE in Scenario 2.}
\label{tab:tv2}
\begin{tabular}{c|cccc}
\hline
% &\multicolumn{2}{c|}{Scenario 1}  &\multicolumn{2}{c|}{Scenario 2} & \multicolumn{2}{c}{Scenario 3}\\
%\hline
$n$ &RPackage  & CVTF  &SURE     & Traditional CV        \\
\hline
300 & 0.029 & 0.029 & 0.032 & 0.028\\
600 & 0.015 & 0.013 & 0.016 & 0.013\\
1200 & 0.007 & 0.006 & 0.008 & 0.007\\
2400 & 0.004 & 0.003 & 0.004 & 0.003\\
\hline
\end{tabular}
\end{table}
\end{comment}
\item Scenario 3 [Piecewise Quadratic Signal]: We consider the piecewise quadratic function 
$$f(x) = 
\begin{cases}
18x^2 \quad &\text{if}\;\; x \in [0, 1/3]\\
-36(x-1/2-1/\sqrt{12})(x-1/2+\sqrt{12}) \quad &\text{if}\;\; x \in [1/3, 2/3]\\
18(x-1)^2 \quad &\text{if}\;\; x \in [2/3, 1]
\end{cases},
$$
%(\mathrm{1}(x \in [0, 1/3])) + -36(x-1/2-1/\sqrt{12})(x-1/2+\sqrt{12})\mathrm{1}(x \in [1/3, 2/3]) + 18(x - 1)^2(\mathrm{1}(x \in [2/3, 1])),$$
and apply Trend Filtering of order $r = 3$. %The corresponding plot is shown in the third diagram of Figure~\ref{fig:tv1} when $n = 300$.
\begin{comment}
\begin{table}[H]
\caption{Comparison of MSE in Scenario 3.}
\label{tab:tv3}
\begin{tabular}{c|cccc}
\hline
% &\multicolumn{2}{c|}{Scenario 1}  &\multicolumn{2}{c|}{Scenario 2} & \multicolumn{2}{c}{Scenario 3}\\
%\hline
$n$ &RPackage & CVTF  &SURE     & Traditional CV        \\
\hline
300 & 0.025 & 0.034 & 0.032 & 0.025\\
600 & 0.014 & 0.019 & 0.016 & 0.014\\
1200 & 0.007 & 0.014 & 0.007 & 0.010\\
2400 & 0.003 & 0.009 & 0.003 & 0.008\\
\hline
\end{tabular}
\end{table}
\end{comment}
\begin{table}[H]
\caption{Comparison of MSE in Scenario 3.}
\label{tab:tv3}
\begin{tabular}{c|cccc}
\hline
% &\multicolumn{2}{c|}{Scenario 1}  &\multicolumn{2}{c|}{Scenario 2} & \multicolumn{2}{c}{Scenario 3}\\
%\hline
$n$ &RPackage & CVTF  &SURE     & Traditional CV        \\
\hline
300 & 0.025 (0.015) & 0.034 (0.017) & 0.032 (0.027) & 0.025 (0.013)\\
600 & 0.014 (0.012) & 0.019 (0.008) & 0.016 (0.014) & 0.014 (0.007)\\
1200 & 0.007 (0.004) & 0.014 (0.005) & 0.007 (0.004) & 0.010 (0.004)\\
2400 & 0.003 (0.002) & 0.009 (0.002) & 0.003 (0.002) & 0.008 (0.002)\\
\hline
\end{tabular}
\end{table}
The fits corresponding to Scenario $1$, $2$ and $3$ are shown in %the 1st, 2nd and 3rd diagram of 
Figure~\ref{fig:tv1} when $n = 300$.
\begin{figure}[H]
\centering
\includegraphics[scale = 0.26]{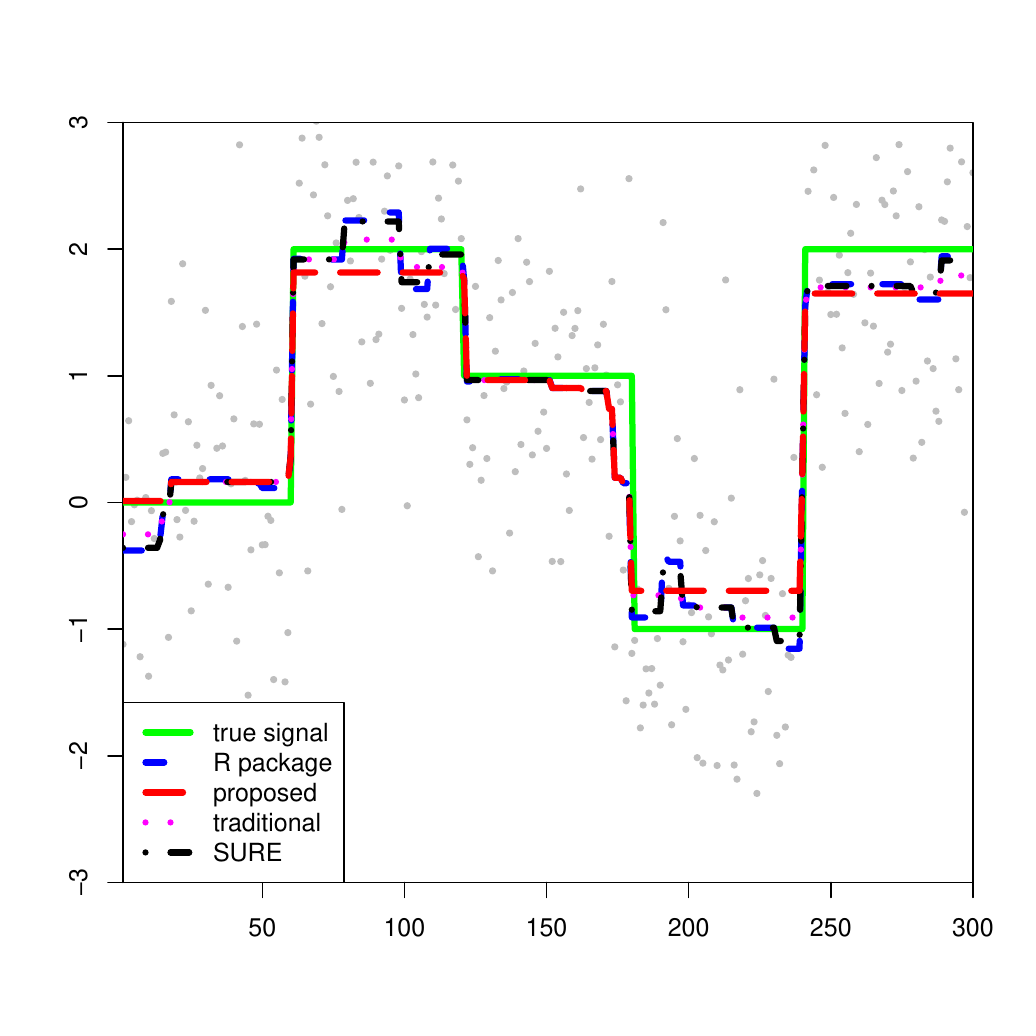} \includegraphics[scale = 0.26]{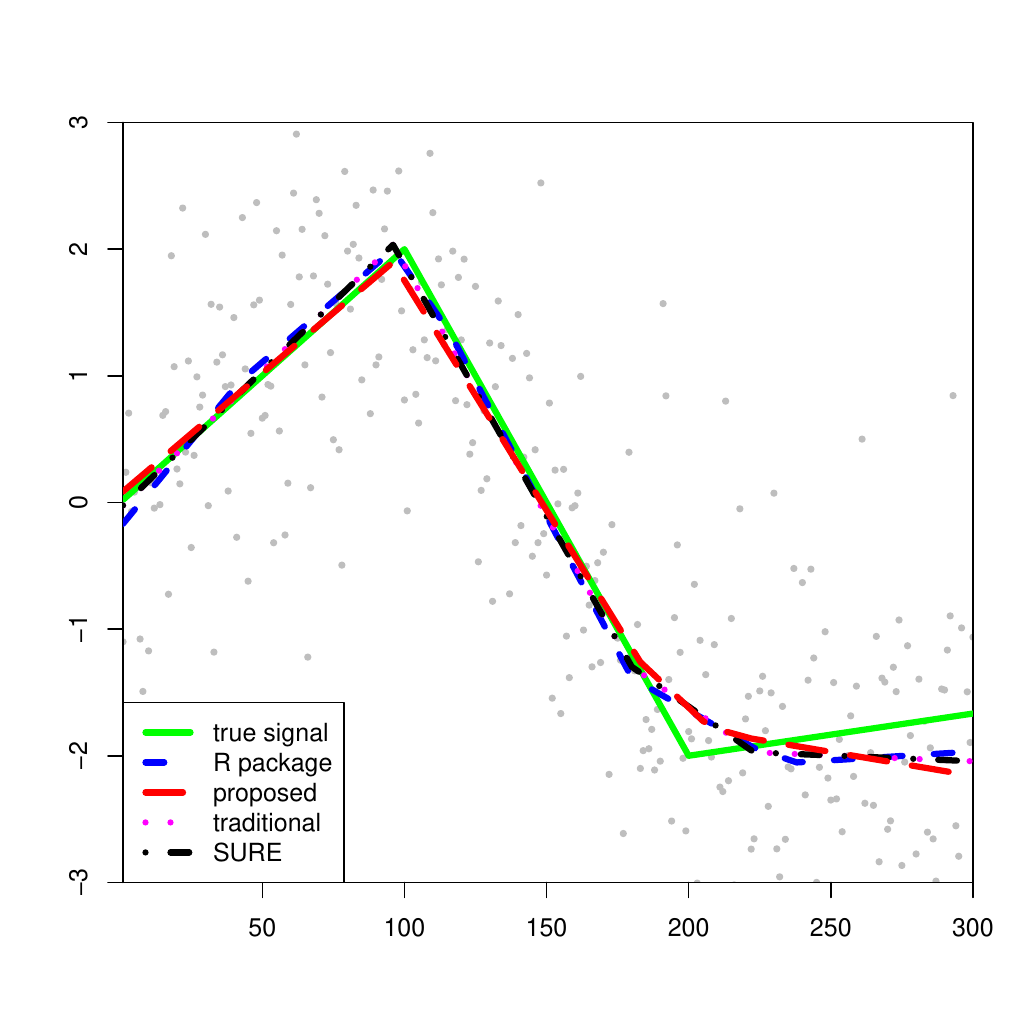} \includegraphics[scale = 0.26]{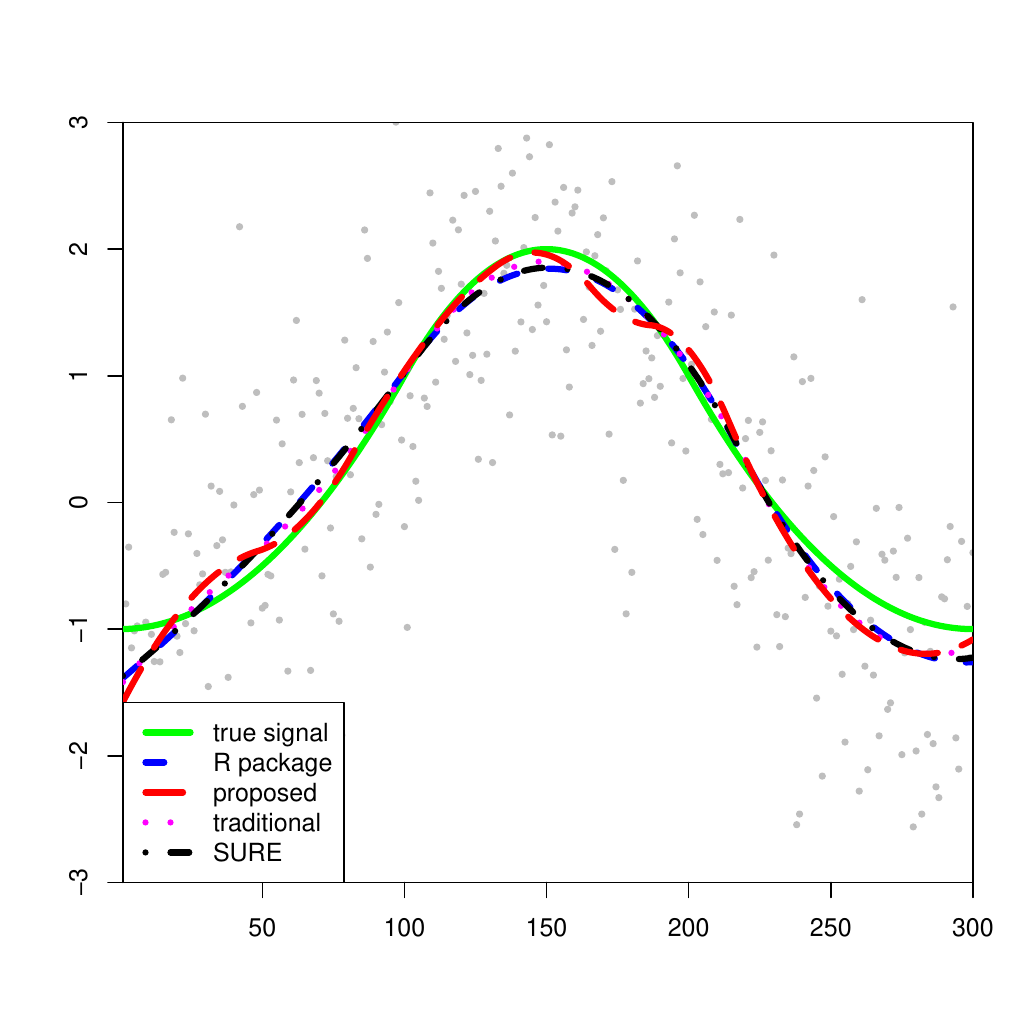}
\caption{The fits of order 1, the fits of order 2 and the fits of order 3 in Scenario 1, 2 and 3 respectively.}
\label{fig:tv1}
\end{figure}

\item Scenario 4 [Smooth Sinusoidal Signal]: We consider the smooth sinusoidal function 
$$f(x) = \sin 2 \pi x + \cos 5 \pi x,$$
and apply Trend Filtering of order $r = 1, 2, 3$. The corresponding plots are shown in Figure~\ref{fig:tv2} when $n = 300$.
\end{enumerate}

\begin{comment}
\begin{table}[H]
\caption{Comparison of MSE between R package and CVTF in different scenarios.}
\label{tab:tv}
\begin{tabular}{c|cc|cc|cc}
\hline
 &\multicolumn{2}{c|}{Scenario 1}  &\multicolumn{2}{c|}{Scenario 2} & \multicolumn{2}{c}{Scenario 3}\\
\hline
$n$ &RPackage Fit                     & CVTF Fit &RPackage Fit                     & CVTF Fit &RPackage Fit                     & CVTF Fit\\
\hline
300 & 0.077 & 0.071 & 0.029 & 0.029 &0.025 & 0.034\\
600 & 0.042 & 0.038 & 0.015 & 0.013 & 0.014 & 0.019\\
1200 & 0.023 & 0.022 & 0.007 & 0.006 & 0.007 & 0.006\\
2400 & 0.013 & 0.012 & 0.004 & 0.003 & 0.003 & 0.003\\
\hline
\end{tabular}
\end{table}
\end{comment}

To summarize, we observe from Tables~\ref{tab:tv1}, \ref{tab:tv2} and \ref{tab:tv3} that the proposed CVTF estimator is adequately competitive with the other CV versions as well as the model selection method by SURE. In fact, in almost every case, surprisingly the CVTF estimator performs slightly better than the SURE method even though SURE requires the knowledge of the error variance. Moreover, its performance is also favourably competitive when compared with the method in the existing R package. On the other hand, in most of the cases, the traditional CV method incurs marginally less (of order $10^{-3}$) MSE than the CVTF estimator. The results here are fully reproducible and our code is available upon request. 

\section{Proof of General Risk Bound for CV Estimator}\label{sec:proofmain}

In this section we will prove Theorem~\ref{thm:main}.
Recall the definition of the intermediate estimator $\tilde{\theta}$ in~\eqref{eq:defninter}. We start by stating the following lemma. 
\begin{lemma}\label{lem:extra}
	The following inequality is true:
	\begin{equation*}
	\|\hat{\theta}_{CV} - \theta^*\| \leq 2 \|\tilde{\theta} - \theta^*\| + \min_{\lambda \in \Lambda} \|\hat{\theta}^{(\lambda)} - \theta^*\|.
	\end{equation*}
\end{lemma}

\begin{proof}
	Recall the definition $\hat{\lambda} = \argmin_{\lambda \in \Lambda} \|\hat{\theta}^{(\lambda)} - \tilde{\theta}\|^2$ and the notation $\hat{\theta}^{(\hat{\lambda})} = \hat{\theta}_{CV}.$ We have by triangle inequality, 
	\begin{equation*}
	\|\hat{\theta}^{(\hat{\lambda})} - \theta^*\| \leq \|\hat{\theta}^{(\hat{\lambda})} - \tilde{\theta}\| + \|\tilde{\theta} - \theta^*\|.
	\end{equation*}
	Further, for any fixed $\lambda \in \Lambda$ we have
	\begin{equation*}
	\|\hat{\theta}^{(\hat{\lambda})} - \tilde{\theta}\| \leq \|\hat{\theta}^{(\lambda)} - \tilde{\theta}\| \leq \|\hat{\theta}^{(\lambda)} - \theta^*\| + \|\tilde{\theta} - \theta^*\|
	\end{equation*}
	where the first inequality follows from the definition of $\hat{\lambda}$ and the second inequality again follows from the triangle inequality. 
\end{proof}

We next state a proposition bounding the root squared error of $\tilde{\theta}.$

\begin{proposition}\label{prop:inter}
	Let $\hat{\theta}^{(\lambda)}$ be a given family of estimators in the subgaussian sequence model for a tuning parameter $\lambda$ ranging in the set $\mathbb{R}_{+}$. Then the intermediate estimator $\tilde{\theta}$ defined in~\eqref{eq:defninter} satisfies for all $x \geq 0$, with probability not less than $1 - 2K \exp(-x^2/2 \sigma^2)$, the following inequality:
	\begin{align*}
	\|\tilde{\theta} - \theta^*\| \leq 2 \sum_{j \in [K]} \min_{\lambda_j \in \Lambda_j} \|\hat{\theta}^{(\lambda_j, I_j^c)}_{I_j} - \theta^*_{I_j}\| \:\:+\:\:4 \sqrt{2}  \sigma \sum_{j \in [K]} \sqrt{\log |\Lambda_j|} + 4Kx.
	\end{align*}
\end{proposition}

Note that combining Proposition~\ref{prop:inter} with Lemma~\ref{lem:extra} finishes the proof of Theorem~\ref{thm:main}. Therefore the remaining task is to prove Proposition~\ref{prop:inter}. Towards this end, the first thing to do is to state an oracle inequality for a least squares estimator over a finite parameter space in the subgaussian sequence model. The topic of oracle inequalities for least squares estimators is classical and its origins date back atleast to~\cite{LiBarron} and~\cite{BarronBirgeMassart}. The fact that such oracle inequalities hold for least squares is well known; also see Chapter $2$ of~\cite{rigollet2015high}. For us, the following version of oracle inequality is the main tool in the proof of Proposition~\ref{prop:inter}.

\begin{lemma}[Oracle Inequality for Finite Least Squares]\label{lem:fls}
	Fix any positive integer $m.$ Let $y \sim Subg(\theta^*, \sigma^2)$ where $y \in \R^m.$ Let $\Theta \subset \R^m$ be a finite parameter space not necessarily containing $\theta^*.$ Let $\hat{\theta}$ be the least squares estimator over the finite parameter space $\Theta$, that is, 
	\begin{equation*}
	\hat{\theta} = \argmin_{\theta \in \Theta} \|y - \theta\|^2.
	\end{equation*}
	Then we have for all $x \geq 0$,
	\begin{equation*}
	P\big(\|\hat{\theta} - \theta^*\| \leq 2 \min_{\theta \in \Theta} \|\theta - \theta^*\| + 4 \sqrt{2} \sigma (\sqrt{\log |\Theta|}) + 4x\big) \geq 1 -  2 \exp(-x^2/2 \sigma^2). 
	\end{equation*}
\end{lemma}

\begin{proof}
	For any arbitrary $\theta \in \Theta$, we have
	\begin{equation*}
	\|y - \hat{\theta}\|^2 \leq \|y - \theta\|^2.
	\end{equation*}
	From here, writing $y = \theta^* + \epsilon$ and expanding the squares one obtains the following pointwise inequality
	%\begin{equation}\label{eq:basic}
	%\|\hat{\theta} - \theta^*\|^2 \leq \|\theta - \theta^*\|^2 + 2 \langle \epsilon, \hat{\theta} - \theta \rangle =  \underbrace{\|\theta - \theta^*\|^2}_{T_1} + \underbrace{2 \langle \epsilon, \hat{\theta} - \theta^* \rangle}_{T_{21}} + \underbrace{2 \langle \epsilon, \theta^* - \theta \rangle}_{T_{22}} 
	%\end{equation}
	
	\begin{align*}\label{eq:basic}
	\|\hat{\theta} - \theta^*\|^2 \leq \|\theta - \theta^*\|^2 + 2 \langle \epsilon, \hat{\theta} - \theta \rangle =  \underbrace{\|\theta - \theta^*\|^2}_{T_1} + \underbrace{2 \langle \epsilon, \hat{\theta} - \theta^* \rangle}_{T_{21}} + \underbrace{2 \langle \epsilon, \theta^* - \theta \rangle}_{T_{22}} 
	\end{align*}
	
	%Note that the first term $T_1$ is  $SSE(\hat{\theta}^{(\lambda_j, I_j^c)}_{I_j}, \theta^*_{I_j})$. %the prediction error with the tuning parameter $\lambda_j$ based on the data $y_{I_j^c}$.
	
	%We can rewrite $T_2$ as 
	%$$T_2 = 2 \left< \epsilon_{I_j},\:\hat{\theta}^{(\hat{\lambda}_j,I_j^c)}_{I_j} - \hat{\theta}^{(\lambda_j, I_j^c)}_{I_j} \right> = \underbrace{2 \left< \epsilon_{I_j},\:\hat{\theta}^{(\hat{\lambda}_j,I_j^c)}_{I_j} - \theta^*_{I_j} \right>}_{T_{21}} + \underbrace{2 \left< \epsilon_{I_j},\:\theta^*_{I_j} - \hat{\theta}^{(\lambda_j, I_j^c)}_{I_j} \right>}_{T_{22}}.$$
	%Therefore, we have a pointwise inequality that,
	%$$SSE(I_j) \leq T_1 + T_{21} + T_{22},$$
	This leads to the fact that for any $t \geq 0$,
	\begin{equation}\label{eq:1}
	\|\hat{\theta} - \theta^*\| \leq t + \mathbbm{1}\left\{\|\hat{\theta} - \theta^*\| > t\right\}\frac{T_1 + T_{21} + T_{22}}{\|\hat{\theta} - \theta^*\|}\\
	= t + A_1 + A_2 + A_3,
	\end{equation}
	where
	$$A_1 := \mathbbm{1}\left\{\|\hat{\theta} - \theta^*\| > t\right\}\frac{T_1}{\|\hat{\theta} - \theta^*\|},$$
	$$A_2 := \mathbbm{1}\left\{\|\hat{\theta} - \theta^*\| > t\right\}\frac{T_{21}}{\|\hat{\theta} - \theta^*\|},$$
	$$A_3 := \mathbbm{1}\left\{\|\hat{\theta} - \theta^*\| > t\right\}\frac{T_{22}}{\|\hat{\theta} - \theta^*\|}.$$
	%Now we upper bound the expectation of each of these three terms in the RHS. First, we will work conditioning on $\epsilon_{I_j^c}$ and then take the final expectation with respect to $\epsilon_{I_j^c}$.
	Now we bound each of these three terms in the RHS. Note that, since $T_1$ is non-negative we have
	$$A_1 \leq \frac{T_1}{t}.$$
	%and since $T_1$ is a function of $\epsilon_{I_j^c}$, conditional expectation keeps it unchanged, that is,
	%$$\E \left[A_1| \epsilon_{I_j^c}\right] \leq \frac{1}{t}\E\left[T_1| \epsilon_{I_j^c}\right] = \frac{T_1}{t}.$$
	
	The second term
	\begin{align*}
	A_2 \leq \frac{|T_{21}|}{\|\hat{\theta} - \theta^*\|}
	= 2\frac{ \left|\left< \epsilon,\:\hat{\theta} - \theta^* \right>\right|}{\|\hat{\theta} - \theta^*\|} \leq 2\max_{\theta \in \Theta} \frac{ \left|\left< \epsilon,\:\theta - \theta^* \right>\right|}{\|\theta - \theta^*\|}.
	\end{align*}
	Now note that the RHS in the above expression is the maximum of a finite number of sub-gaussian random variables (with absolute value) each of subgaussian norm atmost $\sigma.$ Therefore, by %a standard tail inequality for a maximum of finitely many subgaussian random variables; see Lemma $5.2$ in~\cite{van2014probability}, 
	Lemma \ref{lem:max_subG}
	we have
	\begin{equation*}
	P\left(\max_{\theta \in \Theta} \frac{ \left|\left< \epsilon,\:\theta - \theta^* \right>\right|}{\|\theta - \theta^*\|} > \sqrt{2 \sigma^2 \log |\Theta|} + x\right) \leq 2 \exp(-x^2/2 \sigma^2) \:\:\: \forall x \geq 0.
	\end{equation*}
	
	%{\red Include it in the appendix.}
	
	%Also, using the fact that $\epsilon_{I_j}$ is independent of $\epsilon_{I_j^c}$ we obtain
	%\begin{align*}
	%\E\left[A_2| \epsilon_{I_j^c}\right] \leq 2\sigma\sqrt{\log|\Lambda_j|}. 
	%\end{align*}
	%In case, the assignment of $I_1,I_2,\dots,I_k$ is done using randomization independent of $\epsilon$, the fact that $\epsilon_{I_j}$ is independent of $\epsilon_{I_j^c}$ remains true conditionally on the assigment. Therefore, this entire argument holds verbatim under an extra conditioning on the assignment as well. 
	The third term
	\begin{align*}
	A_3 \leq \frac{|T_{22}|}{t} \leq \frac{2}{t} \frac{ |\langle \epsilon, \theta^* - \theta \rangle|}{\|\theta^* - \theta\|} \|\theta^* - \theta\| \leq \frac{2}{t} \sqrt{T_1} \max_{\theta \in \Theta} \frac{ |\langle \epsilon, \theta^* - \theta \rangle|}{\|\theta^* - \theta\|}.
	\end{align*}
	%We can now bound $\max_{\theta \in \Theta} \frac{ |\langle \epsilon, \theta^* - \theta \rangle|}{\|\theta^* - \theta\|}$ again as before. 
	
	%Since $\frac{ \langle \epsilon, \theta^* - \theta \rangle}{\|\theta^* - \theta\|}$ is subgaussian with subgaussian norm at most $\sigma$ we again have 
	%\begin{equation*}
	%P(|\frac{ \langle \epsilon, \theta^* - \theta \rangle}{\|\theta^* - \theta\|}| > \sqrt{2 \sigma^2} + x) \leq \exp(-x^2/2 \sigma^2) \:\:\: \forall x \geq 0.
	%\end{equation*}
	%Therefore, we have
	%\begin{align*}
	%\E\left[A_3| \epsilon_{I_j^c}\right] \leq \frac{1}{t}\E\left[|T_{22}|| \epsilon_{I_j^c}\right] \leq \frac{1}{t}C\sigma||\hat{\theta}^{(\lambda_j, I_j^c)}_{I_j} - \theta^*_{I_j}|| = C\sigma \frac{\sqrt{T_1}}{t},
	%\end{align*}
	%for some constant $C$.
	Plugging these bounds for $A_1,A_2,A_3$ in~\eqref{eq:1} we have the following inequality:
	\begin{align*}
	\|\hat{\theta} - \theta^*\| \leq \min_{\theta \in \Theta} \inf_{t \geq 0} \big[t + \frac{T_1}{t} + 2(1 + \frac{\sqrt{T_1}}{t}) \max_{\theta \in \Theta} \frac{ |\langle \epsilon, \theta^* - \theta \rangle|}{\|\theta^* - \theta\|}\big].
	\end{align*}
	Now setting $t = \sqrt{T_1}$ in the R.H.S gives us
	\begin{align*}
	\|\hat{\theta} - \theta^*\| \leq \min_{\theta \in \Theta} \big[2 \|\theta - \theta^*\| + 4\max_{\theta \in \Theta} \frac{ |\langle \epsilon, \theta^* - \theta \rangle|}{\|\theta^* - \theta\|}\big].
	\end{align*}
	
	Now plugging in the tail bound for $\max_{\theta \in \Theta} \frac{ |\langle \epsilon, \theta^* - \theta \rangle|}{\|\theta^* - \theta\|}$ we obtain for any $x \geq 0$,

	\begin{align*}
	P\big(\|\hat{\theta} - \theta^*\| \leq 2 \min_{\theta \in \Theta} \|\theta - \theta^*\| + 4 \sqrt{2} \sigma \sqrt{\log |\Theta|} + 4x) \geq 1 - 2 \exp(-x^2/2 \sigma^2). 
	\end{align*}
	%Note that the above argument holds for any $\theta \in \Theta$ and the high probability set does not depend on $\theta.$ Thus we can further minimize the bound over $\theta \in \Theta.$ Now setting $t = \sqrt{T_1}$ we finish the proof. 
\end{proof}

We are now ready to prove Proposition~\ref{prop:inter}.

\begin{proof}[Proof of Proposition~\ref{prop:inter}]

	By definition of the intermediate estimator $\tilde{\theta}$ as in~\eqref{eq:defninter}, we have 
	\begin{equation*}
	\|\tilde{\theta} - \theta^*\| \leq \sum_{j = 1}^{K} \|\tilde{\theta}_{I_j} - \theta^*_{I_j}\| = \sum_{j = 1}^{K} \|\hat{\theta}^{(\hat{\lambda}_j,I_j^c)}_{I_j} - \theta^*_{I_j}\|.
	\end{equation*}
	Therefore, it is enough to bound $\|\tilde{\theta}_{I_j} - \theta^*_{I_j}\|$ for any fixed $j \in [K].$

	Fix a $j \in [K].$ We will now argue conditionally on the assignment of folds $\{I_1,I_2,\dots,I_K\}$ and $\epsilon_{I_j^{c}}.$ Since the (finite) set of candidate tuning values $\Lambda_j$ itself can only depend on $\epsilon_{I_j^{c}}$ and each of the estimators in the set $\Theta_j = \{\hat{\theta}^{(\lambda, I_j^c)}_{I_j}: \lambda \in \Lambda_j\}$ are measurable functions of $\epsilon_{I_j^{c}}$ we can treat this set $\Theta_j$ as a fixed finite set of vectors after conditioning on $\{I_1,I_2,\dots,I_K\}$ and $\epsilon_{I_j^{c}}.$ Now recall that $$\hat{\lambda}_j := \argmin_{\lambda \in \Lambda_j} \left|\left|y_{I_j} - \hat{\theta}^{(\lambda, I_j^c)}_{I_j}\right|\right|^2.$$ Therefore, we can view $\tilde{\theta}_{I_j} = \hat{\theta}^{(\hat{\lambda}_j,I_j^c)}_{I_j}$ as a least squares estimator for $\theta^*_{I_j}$ over this finite set $\Theta_j.$ Since the coordinates of $\epsilon$ are i.i.d and the assignment $\{I_1,I_2,\dots,I_K\}$ is either deterministic or done using randomization independent of $\epsilon$ therefore the conditional distribution of $\epsilon_{I_j}$ given the assignment $\{I_1,I_2,\dots,I_K\}$ and $\epsilon_{I_j^{c}}$ is same as the unconditional distribution of the random vector $(\epsilon_1,\dots,\epsilon_{|I_j|})$. Therefore, we are exactly in the setting of Lemma~\ref{lem:fls} and we can apply the oracle inequality for finite least squares to obtain a conditional probability statement for all $x \geq 0$,

	\begin{align*}
	&P\big(\|\tilde{\theta}_{I_j} - \theta^*_{I_j}\| \leq 2 \min_{\lambda_j \in \Lambda_j} \|\hat{\theta}^{(\lambda_j, I_j^c)}_{I_j} - \theta^*_{I_j}\| + 4 \sqrt{2} \sigma \sqrt{\log |\Lambda_j|} + 4x\Bigg| I_1,,\dots,I_K,\epsilon_{I_j^{c}}\big) \geq  \\& 1 -  2 \exp(-x^2/2 \sigma^2). 
	\end{align*}

	Since the upper bound on the probability in the R.H.S above does not depend on $\{I_1,I_2,\dots,I_K\},\epsilon_{I_j^{c}}$ we can drop the conditioning to deduce the unconditional probability statement for all $x \geq 0$,
	\begin{align*}
	P\big(\|\tilde{\theta}_{I_j} - \theta^*_{I_j}\| \leq 2 \min_{\lambda_j \in \Lambda_j} \|\hat{\theta}^{(\lambda_j, I_j^c)}_{I_j} - \theta^*_{I_j}\| + 4 \sqrt{2} \sigma \sqrt{\log |\Lambda_j|} + 4x\big) \geq 1 -  2 \exp(-x^2/2 \sigma^2). 
	\end{align*}
	
	We can now consider a union bound over $j \in [K]$ to finish the proof. 
\end{proof}

\section{Computation of CVDCART}\label{sec:compu}
In this section, we describe an algorithm to compute a completion version of the Dyadic CART estimator defined as
\begin{equation*}
\hat{\theta} \coloneqq \argmin_{\theta \in \R^{L_{d,n}}} ||y_{I} - \theta_{I}||^2 + \lambda k_{\rdp}(\theta)
\end{equation*}
where $I$ is any given arbitrary subset of $L_{d,n}$ and $\lambda > 0.$

The first observation to make is the following fact.
\begin{equation}\label{eq:repres}
\min_{\theta \in \R^{L_{d,n}}} \big(||y_{I} - \theta_{I}||^2 + \lambda k_{\rdp}(\theta)\big)
= \min_{\pi \in \mathcal{P}_{\rdp,d,n}} \big(\min_{\theta \in S_\pi} ||y_{I} - \theta_{I}||^2 + \lambda |\pi|\big)
\end{equation}

Take $\hat{\theta}$, a minimizer of the L.H.S in~\eqref{eq:repres}. Let $\hat{\pi}$ be the minimal recursive dyadic partition correspoding to $\hat{\theta}.$ Then $\hat{\theta} \in S_{\hat{\pi}}.$ Clearly, L.H.S is at most the R.H.S in~\eqref{eq:repres} because the R.H.S takes minimum over all recursive dyadic partitions $\pi$ and and all $\theta \in S_{\pi}.$

Now, consider $\hat{\pi}$, a minimizer of the R.H.S in~\eqref{eq:repres}. Also, let $\hat{\theta} = \argmin_{\theta \in S_{\hat{\pi}}} ||y_{I} - \theta_{I}||^2.$ It is clear that $\hat{\theta}$ is a piecewise constant array such that, within every rectangle $R$ constituting $\hat{\pi}$, it takes the mean value of the entries of $y_I$ within $R$, if $R \cap I \neq \emptyset$; otherwise, it can take any arbitrary value. In particular, 

\begin{align*}
\min_{\theta \in S_\pi} ||y_{I} - \theta_{I}||^2 + \lambda |\pi| = \sum_{R \in \pi} \mathrm{1}(R \cap I \neq \emptyset) \sum_{u \in R \cap I} (y_{u} - \overline{y}_{R \cap I})^2 + \lambda |\pi| = Cost(\pi).
\end{align*}

Now, note that $k_{\rdp}(\hat{\theta}) = |\hat{\pi}|.$ This is because if there exists a recursive dyadic partition $\pi'$ such that $\hat{\theta} \in S_{\pi'}$ and $|\pi'| < |\hat{\pi}|$ then this violates the fact that $\hat{\pi}$ minimizes the R.H.S in~\eqref{eq:repres}. Therefore,  the R.H.S in~\eqref{eq:repres} equals $ ||y_I - \hat{\theta}_I||^2 + \lambda k_\rdp(\hat{\theta})$ and hence is not smaller than the L.H.S in~\eqref{eq:repres}. This shows the correctness of~\eqref{eq:repres}.

Next, we describe the algorithm to compute $\hat{\theta}$. It is clear that, we need to compute the minimum:
$$OPT(L_{d, n}) := \min_{\pi \in \mathcal{P}_{\rdp, d, n}} \min_{\theta \in S_\pi} ||y_{I} - \theta_{I}||^2 + \lambda |\pi| = \min_{\pi \in \mathcal{P}_{\rdp, d, n}} Cost(\pi),$$
and find the optimal partition $\hat{\pi}$.

Now, for any given rectangle $R \subset L_{d, n}$ we can define the minimum of the corresponding subproblem restricted to $R$,
\begin{equation}\label{eq:OPT_R}
OPT(R) := \min_{\pi \in \mathcal{P}_{\rdp, R}} \min_{\theta \in S_\pi(R)} ||y_{R \cap I} - \theta_{R \cap I}||^2 + \lambda |\pi|,
\end{equation}
where we are now optimizing only over the class of recursive dyadic partitions of the rectangle $R$ denoted by $\mathcal{P}_{\rdp, R}$ and for any such partition $\pi \in \mathcal{P}_{\rdp, R}$, let $S_\pi(R)$ denote the subspace of $\R^{R}$ which consists of all arrays which are constant on every rectangle of $\pi$.
Furthermore, we can again write
%note that, the minimizer $\argmin_{\theta \in S_\pi(R)} ||y_{R \cap I} - \theta_{R \cap I}||^2$ is a piecewise constant array such that, within every rectangle $R'$ constituting $\pi$, it takes the mean value of the entries of $y_{R \cap I}$ within $R'$, if $R' \cap I \neq \emptyset$; otherwise, it  takes any arbitrary value. This implies,
$$\min_{\theta \in S_\pi(R)} ||y_{R \cap I} - \theta_{R \cap I}||^2 + \lambda |\pi| = 
\sum_{R' \in \pi} \mathrm{1}(R' \cap I \neq \emptyset) \sum_{u \in R' \cap I} (y_{u} - \overline{y}_{R' \cap I})^2 + \lambda |\pi|,$$
and let the quantity in the right hand side be denoted as $Cost(\pi; R)$.
Therefore, from \eqref{eq:OPT_R}, we can write
$$OPT(R) = \min_{\pi \in \mathcal{P}_{\rdp, R}} Cost(\pi; R).$$

A key point to note here is that the objective function $Cost(\pi; R)$ enjoys an additively separable property. By this, we mean that if we know the first split of the optimal partition $\hat{\pi}$ which minimizes the above problem; then we can separately solve the problem in the resulting two disjoint sub rectangles of $R.$

%enjoys the "additive nature" over any partition of $R$ into two disjoint rectangles. Indeed, if $\pi$ is a recursive dyadic partition of $R$ such that the first split of $\pi$ gives rise to two disjoint rectangles $R_1,R_2$ with $R = R_1 \cup R_2$, then $\pi$ can be written as $\pi = \pi_1 \cup \pi_2$, where $\pi_1$ and $\pi_2$ are partitions of $R_1$ and $R_2$ respectively, and we have
% \begin{align*}
% Cost(\pi; R) &= \sum_{R' \in \pi} \mathrm{1}(R' \cap I \neq \emptyset) \sum_{u \in R' \cap I} (y_{u} - \overline{y}_{R' \cap I})^2 + \lambda |\pi|\\
%  &= \sum_{R' \in \pi_1} \mathrm{1}(R' \cap I \neq \emptyset) \sum_{u \in R' \cap I} (y_{u} - %\overline{y}_{R' \cap I})^2 + \lambda |\pi_1|\\ 
%  &\qquad + \sum_{R' \in \pi_2} \mathrm{1}(R' \cap I \neq \emptyset) \sum_{u \in R' \cap I} (y_{u} - \overline{y}_{R' \cap I})^2 + \lambda |\pi_2|\\
%  &= Cost(\pi_1; R_1) + Cost(\pi_2; R_2).
%  \end{align*}
This property implies the following \textit{dynamic programming principle} for computing $OPT(R)$:
\begin{equation*}
OPT(R) = \min_{(R_1, R_2)} \left\{OPT(R_1) + OPT(R_2), \quad \mathrm{1}(R \cap I \neq \emptyset) \sum_{u \in R \cap I}(y_u -\overline{y}_{R \cap I})^2 + \lambda\right\}.
\end{equation*}
In the above, $(R_1, R_2)$ ranges over all possible \textit{nontrivial} dyadic splits of $R$ into two disjoint rectangles. %Consequently, in order to compute $OPT(R)$ and obtain the optimal recursive dyadic partition, the first step is to obtain the corresponding first split of $R$, which let be denoted by $SPLIT(R)$.
The minimizer of the above problem lets us obtain the optimal first split of $R$, which we denote by $SPLIT(R)$.

Now it is important to make some observations. For any rectangle $R$, the number of non trivial dyadic splits possible is at most $d$, one for each dimension. Any split of $R$ creates two disjoint sub rectangles $R_1$ and $R_2$. Suppose we know $OPT(R_1)$ and $OPT(R_2)$ for $R_1$, $R_2$ arising out of each possible split. Then, to compute $SPLIT(R)$ we have to compute the minimum of the sum of $OPT(R_1)$ and $OPT(R_2)$ for each possible split as well as the number $\mathrm{1}(R \cap I \neq \emptyset) \sum_{u \in R \cap I}(y_u -\overline{y}_{R \cap I})^2$, which corresponds to not splitting $R$ at all. Thus, we need to compute the minimum of at most $d+1$ numbers.

From Lemma \ref{lem:dyadic}, the number of distinct dyadic rectangles of $L_{d, n}$ is at most $2^dN$. Any rectangle $R$ has dimension $n_1 \times n_2 \times \cdots \times n_d$. Let us denote the number $n_1 + \cdots + n_d$ by $Size(R)$. Now we are ready to describe the main scheme of computing the optimal partition.

For each rectangle $R$, the goal is to store $SPLIT(R)$ and $OPT(R)$. We do this inductively on $Size(R)$. We will make a single pass/visit through all distinct rectangles $R \subset L_{d, n}$, in increasing order of $Size(R)$. Thus, we will first start with all $1 \times 1 \times \cdots \times 1$ rectangles of size equals $d$. Then we visit rectangles of size $d+1$, $d+2$ all the way to $nd$. Fixing the size, we can choose some arbitrary order in which we visit the rectangles.

For $1 \times 1 \times \cdots \times 1$ rectangles, computing $SPLIT(R)$ and $OPT(R)$ is trivial. Consider a generic step where we are visiting some rectangle $R$. Note that we have already computed $OPT(R')$ for all rectangles $R'$ with $Size(R') < Size(R)$. Since a possible split of $R$ generates two rectangles $R_1$, $R_2$ of strictly smaller size, to compute $OPT(R_1) + OPT(R_2)$ we just need to sum two previously computed numbers and store it. We do this for each possible split to get a list of at most $d$ numbers. Moreover, we also compute $\mathrm{1}(R \cap I \neq \emptyset) \sum_{u \in R \cap I}(y_u -\overline{y}_{R \cap I})^2$ (described later) and add this number to the list. Finally, we take the minimum of these $d+1$ numbers. In this way, we obtain $OPT(R)$ and $SPLIT(R)$.

The number of basic operations needed per rectangle here is $O(d)$. Since there are at most $2^dN$ many rectangles in all, the total computational complexity of the overall inductive scheme scales like $O(2^ddN)$.

To compute $\mathrm{1}(R \cap I \neq \emptyset) \sum_{u \in R \cap I}(y_u -\overline{y}_{R \cap I})^2$ for every rectangle $R$, we can again induct on size in increasing order. Define $SUM(R)$ to be the sum of entries of $y_{R \cap I}$ and $SUMSQ(R)$ to be the sum of squares of entries of $y_{R \cap I}$. One can easily keep storing $SUM(R)$ and $SUMSQ(R)$ in a similar bottom up fashion, visiting rectangles in increasing order of $Size(R)$. This requires constant number of basic operations per rectangle $R$. Once we have computed $SUM(R)$ and $SUMSQ(R)$, we can then calculate $\mathrm{1}(R \cap I \neq \emptyset) \sum_{u \in R \cap I}(y_u -\overline{y}_{R \cap I})^2$. Thus, this inductive scheme requires lower order computation.

Once we finish the above inductive scheme, we have stored $SPLIT(R)$ for every rectangle $R$. We can now start going topdown, starting from the biggest rectangle which is $L_{d, n}$ itself. We can recreate the full optimal partition by using $SPLIT(R)$ to split the rectangles at every step. Once the full optimal partition $\hat{\pi}$ is obtained, computing $\hat{\theta}$ just involves computing the mean value of the entries of $y_I$ within the rectangles. If a rectangle does not contain any entry of $y_I$, we simply compute the mean value of the entries of $y_I$, i.e., $\overline{y}_I$. It can be checked that this step requires lower order computation as well.

\section{Proofs for Dyadic CART}\label{sec:dcproofs}

\medskip
The main goal of this section is to prove Theorem~\ref{thm:main_cv_dc}.
\subsection{\textbf{Sketch of Proof of Theorem~\ref{thm:main_cv_dc}}}
For the convenience of the reader, we first present a sketch of proof of Theorem~\ref{thm:main_cv_dc}. This sketch is divided into several steps and is meant to convey the essential aspects of our proof strategy.

Let $\hat{\theta}^{(\lambda)}$ denote the usual Dyadic CART estimator based on the full data array $y$ as defined in~\eqref{eq:dcdefn}. 
In view of Theorem \ref{thm:main}, it is enough to bound $\min_{\lambda \in \Lambda} SSE(\hat{\theta}^{(\lambda)}, \theta^*)$ and $\min_{\lambda \in \Lambda} SSE(\hat{\theta}^{(\lambda,I)}_{I^c}, \theta^*_{I^c})$, where $I = I_1$ or $I_2$.

The existing bound in Theorem~\ref{thm:adapt} gives the desired upper bound on $SSE(\hat{\theta}^{(\lambda)}, \theta^*)$ as long as $\lambda$ is chosen to be not smaller than $C \sigma^2 \log N$ where $C$ is some absolute constant. Since $2^{N^*}$ is assumed to be strictly larger than $C \sigma^2 \log N$, therefore $\Lambda$ contains this good choice $C \sigma^2 \log N$ of the tuning parameter. In particular, by construction of $\Lambda$, there exists a $\lambda^* \in \Lambda$ satisfying $C \sigma^2 \log n < \lambda^* \leq 2 C \sigma^2 \log n.$ For our purpose, this $\lambda^*$ can be thought of as a sufficiently good choice of the tuning parameter.

Clearly, $\min_{\lambda \in \Lambda} SSE(\hat{\theta}^{(\lambda)}, \theta^*) \leq SSE(\hat{\theta}^{(\lambda^*)}, \theta^*).$ We can now use the existing bound in Theorem~\ref{thm:adapt} to bound $SSE(\hat{\theta}^{(\lambda^*)}, \theta^*)$ for the usual Dyadic CART estimator by generalizing the proof of Theorem \ref{thm:adapt} given in~\cite{chatterjee2019adaptive} to general subgaussian errors. Therefore, the desired bound for $\min_{\lambda \in \Lambda} SSE(\hat{\theta}^{(\lambda)}, \theta^*)$ follow more or less directly from existing results.

The main new task for us here is bounding $\min_{\lambda \in \Lambda} SSE(\hat{\theta}^{(\lambda,I)}_{I^c}, \theta^*_{I^c})$ where $I = I_1$ or $I = I_2.$ Again, since $\min_{\lambda \in \Lambda} SSE(\hat{\theta}^{(\lambda,I)}_{I^c}, \theta^*_{I^c}) \leq SSE(\hat{\theta}^{(\lambda^*,I)}_{I^c}, \theta^*_{I^c})$ it is sufficient for us to bound $SSE(\hat{\theta}^{(\lambda^*,I)}_{I^c}, \theta^*_{I^c})$ which in turn is trivially upper bounded by $SSE(\hat{\theta}^{(\lambda^*,I)}, \theta^*)$.

We now outline the main steps in our proof which bounds $SSE(\hat{\theta}^{(\lambda^*,I)}, \theta^*)$. 
Before that we say a few words about notation at this point. In this proof sketch, we abuse notation and write $A \lesssim B$, which basically means $A \leq B$ after ignoring multiplicative and additive logarithmic factors. This is done to increase readability and interpretability of our bounds within this proof sketch. Also, we write $A (\leq^{+}) B$ to mean that $A \leq \max\{B,0\}.$

\begin{enumerate}
	\item \textbf{Preliminary Localizations}

	We first establish a few preliminary localization properties of the completion estimator $\hat{\theta}^{(\lambda^*,I)}.$ Let us introduce a few notations at this point. Let us denote $L^* = V(\theta^*) + \sqrt{\log N}$ where we recall that $V(\theta^*) = \max_{u,v \in L_{d,n}} |\theta^*_u - \theta^*_v|.$ For realistic signals we would expect $V(\theta^*) = O(1)$ and not grow with $n$ to $\infty.$ Therefore, the reader can think of $L^*$ as a $\tilde{O}(1)$ term, growing at most logarithmically. %Let us also define the function $R: \R^{L_{d,n}} \rightarrow \R$ as follows:
	%\begin{equation*}
	%R(\theta) = \inf_{\alpha \in \R^{L_{d,n}}} \big(\big)
	%\end{equation*} 

	Now define the events  
	\begin{equation}\label{eq:events}
	\begin{split}
	A_1 &= \left\{|\hat{\theta}^{(\lambda^*,I)} - \theta^*|_{\infty} \lesssim L^*\right\},\\
	A_2 & = \left\{||\theta^*_I - \hat{\theta}^{(\lambda^*, I)}_I||^2 + \lambda^* k_\rdp(\hat{\theta}^{(\lambda^*, I)}) \lesssim R(\theta^*,\lambda^*) \right\}.
	%A_2 &= \left\{\|\hat{\theta}^{(\lambda,I)}_{I} - \theta^*_{I}\|^2 \leq R(\theta^*, \lambda)\right\},\quad \text{and}\\
	%A_3 &= \left\{k_{\rdp}(\hat{\theta}^{(\lambda,I)}) \leq \frac{2}{\lambda} R(\theta^*, \lambda)\right\}.
	\end{split}
	\end{equation}
	We first show that both the events $A_1$ and $A_2$ hold with high probability. This is formalized in Proposition~\ref{prop:dcevents}. In the rest of the steps, we will assume that the events $A_1$ and $A_2$ hold.
	
	\item \textbf{Total Error = In Sample Error + Out of Sample Error}

	We need to bound $SSE(\hat{\theta}^{(\lambda^*,I)}, \theta^*)$; let's term this as the total error. %From the previous step, we have a bound for $SSE(\hat{\theta}^{(\lambda^*,I)}_{I}, \theta^*_{I})$ (under the event $A_2$) which can be thought of as an \textit{in sample error}. 
	We now decompose the total error as twice the in sample error plus the difference of out of sample and in sample errors where by \textit{out of sample error} we mean $SSE(\hat{\theta}^{(\lambda^*,I)}_{I^c}, \theta^*_{I^c})$ and by \textit{in sample error} we mean $SSE(\hat{\theta}^{(\lambda^*,I)}_{I}, \theta^*_{I})$.

	%{\color{red} change every $\lambda$ to $\lambda^*$}
	\begin{align*}
	&SSE(\hat{\theta}^{(\lambda^*,I)}, \theta^*) = SSE(\hat{\theta}^{(\lambda^*,I)}_{I}, \theta^*_{I}) + SSE(\hat{\theta}^{(\lambda^*,I)}_{I^c}, \theta^*_{I^c}) = \\& 2 SSE(\hat{\theta}^{(\lambda^*,I)}_{I}, \theta^*_{I}) +  SSE(\hat{\theta}^{(\lambda^*,I)}_{I^c}, \theta^*_{I^c}) -  SSE(\hat{\theta}^{(\lambda^*,I)}_{I}, \theta^*_{I}).
	\end{align*}
	It would be convenient if we further rewrite the earlier display after division by the square root of $SSE(\hat{\theta}^{(\lambda^*,I)}, \theta^*)$ to obtain
	
	\begin{align}\label{eq:loc}
	\|\hat{\theta}^{(\lambda^*,I)} - \theta^*\| &= 2 \frac{\|\hat{\theta}^{(\lambda^*,I)}_{I} - \theta^*_{I}\|^2}{\|\hat{\theta}^{(\lambda^*,I)} - \theta^*\|} + \frac{\|\hat{\theta}^{(\lambda^*,I^c)}_{I^c} - \theta^*_{I^c}\|^2 - \|\hat{\theta}^{(\lambda^*,I)}_{I} - \theta^*_{I}\|^2}{\|\hat{\theta}^{(\lambda^*,I)} - \theta^*\|} \nonumber \\&\leq \underbrace{2 \|\hat{\theta}^{(\lambda^*,I)}_{I} - \theta^*_{I}\|}_{T_1} + \underbrace{\frac{\|\hat{\theta}^{(\lambda^*,I^c)}_{I^c} - \theta^*_{I^c}\|^2 - \|\hat{\theta}^{(\lambda^*,I)}_{I} - \theta^*_{I}\|^2}{\|\hat{\theta}^{(\lambda^*,I)} - \theta^*\|}}_{T_2}.
	\end{align}
	
	Since the in sample error term $T_1$ has the desired oracle risk bound (under the event $A_2$) from Step $1$ itself, our focus henceforth is on bounding the normalized difference of out of sample and in sample errors which is the term $T_2.$ 
	
	\item \textbf{Reduction to Bounding Rademacher Averages}

	To bound $T_2$, we express the numerator in $T_2$ as a Rademacher average. \begin{align*}
	\|\hat{\theta}^{(\lambda^*,I)}_{I^c} - \theta^*_{I^c}\|^2 - \|\hat{\theta}^{(\lambda^*,I)}_{I} - \theta^*_{I}\|^2 = \sum_{u \in L_{d, n}} \eta_{u} (\hat{\theta}^{(\lambda^*,I)}_{u} - \theta^*_u)^2,
	\end{align*}
	where $\eta_u = \mathrm{1}(u \in I^c) - \mathrm{1}(u \in I)$ are i.i.d $\pm1$ rademacher random variables (same as the entries of $W$ or $-W$ depending on whether $I = I_1$ or $I = I_2.$). This is the step where we really use the fact that our folds are random and chosen in a uniformly i.i.d manner. The above display gives 
	\begin{equation*}
	T_2 = \sum_{u \in L_{n,d}} \frac{\eta_{u} (\hat{\theta}^{(\lambda^*,I)}_{u} - \theta^*_u)^2}{\|\hat{\theta}^{(\lambda^*,I)} - \theta^*\|}.
	\end{equation*}

	\item \textbf{Peeling Step}

	At this point we are still not quite ready to bound the R.H.S in the last display in Step $3$, primarily because of the presence of the square in the exponent of  $(\hat{\theta}^{(\lambda^*,I)}_{u} - \theta^*_u).$ If there were no square, $T_2$ would have been in "standard" form and could have been directly bounded by a Rademacher complexity term which in turn could have been bounded by standard techniques as laid out in the proof of Theorem $2.1$  in~\cite{chatterjee2019adaptive}. Henceforth, our effort will be to "drop the square". This will require us to use standard techniques from the theory of bounding maxima of stochastic processes, carefully and in appropriate order.

	We first implement the so called peeling step where we peel on the value of $\|\hat{\theta}^{(\lambda^*,I)}- \theta^*\|.$ 
	
	\begin{align*}
	&T_2 \leq \sum_{l = 1}^{L} \sum_{u \in L_{n,d}} \frac{\eta_{u} (\hat{\theta}^{(\lambda^*,I)}_{u} - \theta^*_u)^2}{\|\hat{\theta}^{(\lambda^*,I)} - \theta^*\|} \mathrm{1}\left(2^l < \|\hat{\theta}^{(\lambda^*,I)}- \theta^*\| \leq 2^{l + 1}\right) (\leq^{+}) \\& \sum_{l = 1}^{L} \underbrace{\frac{1}{2^l} \sum_{u \in L_{d,n}} \eta_{u} (\hat{\theta}^{(\lambda^*,I)}_{u} - \theta^*_u)^2 \mathrm{1}\left(2^l < \|\hat{\theta}^{(\lambda^*,I)}- \theta^*\| \leq 2^{l + 1}\right)}_{T_{2,l}}.
	\end{align*}

	Since the event $A_1$ holds, $|\hat{\theta}^{(\lambda^*,I)} - \theta^*|_{\infty} \lesssim L^*$ and hence $\|\hat{\theta}^{(\lambda^*,I)} - \theta^*\|^2 \lesssim N (L^*)^2.$ Therefore, in the above sum, $l$ only needs to go up to $L \lesssim \log \left(N (L^*)^2\right) = \tilde{O}(1).$ 
	
	%{\color{red} Mention about the possible values of $l.$ Also the last ineq is up to positive part. Explain that.}
	
	\item \textbf{Reduction to Bounding Suprema}

	In view of the last display in the last step, we need to bound $T_{2,l}$ for each $l$ from $1$ to $L = \tilde{O}(1)$. At this point, we "sup out" the random variable $\hat{\theta}^{(\lambda^*,I)}$ using the high probability localization properties established in Step $1.$ Specifically, for any fixed $l \in [L]$ we can write
	
	\begin{align*}
	&T_{2,l}\\ &= \frac{1}{2^l} \sum_{u \in L_{d, n}} \eta_{u} (\hat{\theta}^{(\lambda^*,I)}_{u} - \theta^*_u)^2 \mathrm{1}\Bigg(2^l < \|\hat{\theta}^{(\lambda^*,I)}- \theta^*\| \leq 2^{l + 1}, k_{\rdp}(\hat{\theta}^{(\lambda^*,I)}) \lesssim \frac{R(\theta^*,\lambda^*)}{\lambda^*}, \\
	&\qquad \qquad \qquad \qquad \qquad \qquad \qquad \qquad \qquad \qquad \qquad \qquad |\hat{\theta}^{(\lambda^*,I)} - \theta^*|_{\infty} \lesssim L^*\Bigg) \\
	&\leq 2L^* \sup_{\substack{\pi \in \mathcal{P}_{\rdp,d,n}, \\ |\pi| \lesssim \frac{R(\theta^*,\lambda^*)}{\lambda^*}}} \: \frac{1}{2^l}\:\underbrace{\sup_{\substack{v \in S_{\pi}: 2^l < \|v - \theta^*\| \leq 2^{l + 1},\\ |v - \theta^*|_{\infty} \lesssim L^*}} \frac{1}{2L^*} \sum_{u \in L_{d, n}} \eta_{u} (v_{u} - \theta^*_u)^2}_{X_{\pi,l}}.
	\end{align*}
	where the first equality follows because the events $A_1,A_2$ hold. Here $S_{\pi}$ refers to the subspace of arrays in $\R^{L_{d,n}}$ which are piecewise constant on each of the rectangles of $\pi.$

	%{\color{red} Change $M,L$ to $R(\theta^*)$ }
	
	\item \textbf{Contraction Principle for  Rademacher Averages}
	
	%In the lat display, the random variables $X_{\pi,l}$ were introduced for any given recursive dyadic partition $\pi \in \mathcal{P}_{\rdp,d,n}$ and $l \in L.$ 
	At this point we can indeed "drop the square" by using the so-called contraction principle for Rademacher averages, see Theorem $2.2$ in~\cite{koltchinskii2011oracle}. This is allowed because the square function is Lipschitz in $[-L^*,L^*]$ with lipschitz constant $2L^*.$ Note that we can localize to $[-L^*,L^*]$ precisely because of the $\ell_{\infty}$ localization property we established in Step $1.$ In particular, an in probability version of the contraction principle allows us to conclude that the random variable $X_{\pi,l}$ (defined in the last display in Step $5$) is stochastically upper bounded by the random variable $Y_{\pi,l}$ which is defined as 
	$$ Y_{\pi,l} := \sup_{\substack{v \in S_{\pi}: \|v- \theta^*\| \leq 2^{l + 1}, \\ |v- \theta^*|_{\infty} \leq L^*}} \sum_{u \in {L_{d,n}}} \eta_u (v_{u} - \theta^*_{u}).$$
	Note that there is no square in the exponent of $v - \theta^*$ in $Y_{\pi,l}$ anymore.

	This in turn essentially implies that we can write 
	$$T_{2,l} \leq 2L^* \sup_{\substack{\pi \in \mathcal{P}_{\rdp,d,n}, \\ |\pi| \leq \frac{R(\theta^*,\lambda^*)}{\lambda^*}}} \frac{1}{2^l} Y_{\pi,l}.$$
	
	Note that there is no square in the exponent of $v - \theta^*$ in $Y_{\pi,l}$ anymore. So we have succesfully dropped the square while going from $X_{\pi,l}$ to $Y_{\pi,l}.$ This entire step is carried out within the proof of Proposition~\ref{prop:sqrad}.

	%{\color{red} Check these steps again!}	

	\item \textbf{Bounding Suprema}

	The final step is to observe that the right hand side of the last display is a Rademacher complexity term in the "standard form". Moreover, the $\frac{1}{2^l}$ factor effectively cancels out with the $2^{l + 1}$ radius of the $\ell_2$ ball over which we are taking the supremum in the term $Y_{\pi,l}.$ Hence, we can give a bound (for any $l \in [L]$) by essentially similar proof techniques as laid out in the proof of Theorem $2.1$  in~\cite{chatterjee2019adaptive}. This is carried out by using the so called basic inequality and applying standard bounds on suprema of maxima of subgaussian random variables along with basic estimates of the cardinality of $\mathcal{P}_{\rdp,d,n}.$ This yields a bound 
	\begin{equation*}
	T_{2,l} \lesssim L^* \sqrt{\frac{R(\theta^*,\lambda^*)}{\lambda^*}}
	\end{equation*}
	which holds with high probability. Since there are only $\tilde{O}(1)$ many $l$'s which need to accounted for, a simple union bound applied to the above display then furnishes with high probability
	\begin{equation*}
	T_{2} \lesssim L^* \sqrt{\frac{R(\theta^*,\lambda^*)}{\lambda^*}}.
	\end{equation*}
	This display gives us the desired bound for $T_2$ and thus finishes the proof.

\end{enumerate}

\subsection{\textbf{Detailed Proof of Theorem~\ref{thm:main_cv_dc}}}
We now give the detailed proof.
Throughout the proofs in this section, we will use $C,C_1,C_2$ to denote constants which may only depend on the dimension $d$ but not the sample size $N$, the true signal $\theta^*$ or the distribution of the errors. Also, the precise values of these constants may change from line to line. Generically, we will denote any partition in $\mathcal{P}_{\rdp,d,n}$ by $\pi.$ We will then denote by $S_\pi$ the subspace of $\R^{L_{d,n}}$ consisting of all arrays which are constant on every subrectangle of $\pi.$

%{\color{red} Some basic notations? Constant line to line change.}

The proof of Theorem~\ref{thm:main_cv_dc} is organized in the following way. First, we state two propositions which are crucial ingredients in this proof. After stating these two propositions, we finish the proof of Theorem~\ref{thm:main_cv_dc}.

\begin{proposition}\label{prop:dcevents}

Let $I$ denote any fixed subset of $L_{d,n}.$ %the random set of indices $I_1$ or $I_2.$  
Recall that the estimator $\hat{\theta}^{(\lambda,I)}$ is defined as
\begin{equation*}
\hat{\theta}^{(\lambda, I)} \coloneqq \argmin_{\theta \in \R^{L_{d,n}}} ||y_{I} - \theta_{I}||^2 + \lambda k_{\rdp}(\theta).
\end{equation*}

There exists an absolute constant $C$ such that if $\lambda \geq C \sigma^2 \log N$ then for any $\alpha > 1$, $$P(A_1 \cap A_2) \geq 1 - C N^{-C \alpha},$$ where the events $A_1,A_2$ are defined as follows: 
%{\color{red} fix this probability bound and introduce an $\alpha.$}
	
\begin{equation}\label{eq:events}
\begin{split}
A_1 &= \left\{|\hat{\theta}^{(\lambda,I)} - \theta^*|_{\infty} \leq V(\theta^*) + (2 + \sqrt{\alpha}) \sigma \sqrt{\log N}\right\},\\
A_2 & = \left\{||\hat{\theta}^{(\lambda, I)}_I - \theta^*_I||^2 + \lambda k_\rdp(\hat{\theta}^{(\lambda, I)}) \leq R(\theta^*,\lambda) + C\alpha \sigma^2 \log N \right\}
%A_2 &= \left\{\|\hat{\theta}^{(\lambda,I)}_{I} - \theta^*_{I}\|^2 \leq R(\theta^*, \lambda)\right\},\quad \text{and}\\
%A_3 &= \left\{k_{\rdp}(\hat{\theta}^{(\lambda,I)}) \leq \frac{2}{\lambda} R(\theta^*, \lambda)\right\}.
\end{split}
\end{equation}

Since, the probability lower bound $1 - C N^{-C \alpha}$ does not depend on the subset $I$, the same conclusion is true even if $I$ is now a random subset (chosen independently of $\epsilon$), such as the random folds $I_1$ or $I_2.$ 
\end{proposition}

%{\color{red} write this for a general subset $I.$}

\begin{proposition}\label{prop:sqrad}
	Suppose $\eta \in \R^{L_{d,n}}$ is a random array consisting of i.i.d Rademacher random variables. For any integer $l \geq 0$ and $L, M > 0$, define the random variable 
	\begin{align*}
	T_{l} := \sum_{u \in L_{d, n}} \frac{\eta_{u} (\hat{\theta}^{(\lambda,I)}_{u} - \theta^*_u)^2}{\|\hat{\theta}^{(\lambda,I)} - \theta^*\|} \mathrm{1}\Bigg(2^l < \|\hat{\theta}^{(\lambda,I)}- \theta^*\| \leq 2^{l + 1}, &k_{\rdp}(\hat{\theta}^{(\lambda,I)}) \leq M,\\
	&|\hat{\theta}^{(\lambda,I)} - \theta^*|_{\infty} \leq L\Bigg).
	\end{align*}
	Then there exists absolute positive constants $C_{1}$ and $C_{2}$ such that for any $\alpha \geq 1$, %and $l \geq 0$,
	\begin{equation*}
	P\left(T_{l} \geq (C_{1} \sqrt{\alpha} +C_{2}) L \sqrt{M \log N}\right) \leq N^{-\alpha}.
	\end{equation*}
	\end{proposition}

Now we will finish the proof of Theorem~\ref{thm:main_cv_dc} assuming these propositions hold.

\begin{proof}[Proof of Theorem~\ref{thm:main_cv_dc}]
	%{\color{red} First say that it is enough to bound the MSE of $\tilde{\theta}$ as defined in~\eqref{eq:defninterdc}. For this it is enough to bound MSE of $\hat{\theta}^{(\lambda,I)}$ where $I = I_1$ or $I_2.$}
	Let $\hat{\theta}^{(\lambda)}$ denote the usual Dyadic CART estimator based on the full data array $y$ as defined in~\eqref{eq:dcdefn}. 
	In view of Theorem \ref{thm:main}, it is enough to bound $SSE(\hat{\theta}^{(\lambda)}, \theta^*)$ and $\min_{\lambda \in \Lambda} SSE(\hat{\theta}^{(\lambda,I)}, \theta^*)$, where $I = I_1$ or $I_2$.%and $\lambda \geq C\sigma^2\log N$.

    We can use the existing bound in Theorem~\ref{thm:adapt} to bound $SSE(\hat{\theta}^{(\lambda)}, \theta^*)$ for the usual Dyadic CART estimator (under a Gaussian assumption on the errors). In fact, Proposition~\ref{prop:dcevents} generalizes Theorem~\ref{thm:adapt} by giving bounds for a completion version of Dyadic CART (and for general subgaussian errors). By setting $I$ equal to the entire set $L_{d,n}$, it is seen that Theorem~\ref{thm:adapt} is a special case of Proposition~\ref{prop:dcevents}. Thus, Proposition~\ref{prop:dcevents} itself also implies that if $\lambda \geq C \sigma^2 \log N$ then for any $\alpha > 1$, there exists an absolute constant $C$ such that with probability at least $1 - C N^{-C \alpha}$, we have
    $$SSE(\hat{\theta}^{(\lambda)}, \theta^*) \leq R(\theta^*,\lambda) + C \alpha \sigma^2 \log N.$$
   
%Our goal is to bound $\|\hat{\theta}^{(\lambda,I)} - \theta^*\|.$ 
Our main goal now is to bound $SSE(\hat{\theta}^{(\lambda,I)}, \theta^*)$, for $I = I_1$ or $I_2$, under the choice of $\lambda \geq C \sigma^2 \log N.$ Consider the events $A_1$ and $A_2$ as defined in \eqref{eq:events}. Furthermore, define the event
$A_3 := \{\|\hat{\theta}^{(\lambda,I)} - \theta^*\| > 1\}.$ In the calculation below we will assume that the event $A_1 \cap A_2 \cap A_3$ holds because $A_1 \cap A_2$ hold with high probability and on $A_3^c$ we anyway have $\|\hat{\theta}^{(\lambda,I)} - \theta^*\| \leq 1$.

We first write
\begin{align*}
\|\hat{\theta}^{(\lambda,I)} - \theta^*\|^2 &= 2 \|\hat{\theta}^{(\lambda,I)}_{I} - \theta^*_{I}\|^2 + \|\hat{\theta}^{(\lambda,I)}_{I^c} - \theta^*_{I^c}\|^2 - \|\hat{\theta}^{(\lambda,I)}_{I} - \theta^*_{I}\|^2\\ &= 2 \|\hat{\theta}^{(\lambda,I)}_{I} - \theta^*_{I}\|^2 + \sum_{u \in L_{d, n}} \eta_{u} (\hat{\theta}^{(\lambda,I)}_{u} - \theta^*_u)^2
\end{align*}
where the array $\eta \coloneqq \eta_u = \mathrm{1}(u \in I^c) - \mathrm{1}(u \in I)$ consists of i.i.d $\pm1$ rademacher random variables (same as $W$ or $-W$ depending on whether $I = I_1$ or $I = I_2.$).
	%{\red I think in the beginning we should define that $1(A)$ defines the indicator variable of the event $A$.}
%{\color{red} To do for Anamitra: Are we using this above notation for indicator? If not make it consistent throughout the paper.}

Dividing both sides by $\|\hat{\theta}^{(\lambda,I)} - \theta^*\|$ (nonzero on $A_3$) we can write
\begin{align}\label{eq:loc}
\|\hat{\theta}^{(\lambda,I)} - \theta^*\| &= 2 \frac{\|\hat{\theta}^{(\lambda,I)}_{I} - \theta^*_{I}\|^2}{\|\hat{\theta}^{(\lambda,I)} - \theta^*\|} + \sum_{u \in L_{n,d}} \frac{\eta_{u} (\hat{\theta}^{(\lambda,I)}_{u} - \theta^*_u)^2}{\|\hat{\theta}^{(\lambda,I)} - \theta^*\|} \nonumber \\&\leq \underbrace{2 \|\hat{\theta}^{(\lambda,I)}_{I} - \theta^*_{I}\|}_{T_1} + \underbrace{\sum_{u \in L_{n,d}} \frac{\eta_{u} (\hat{\theta}^{(\lambda,I)}_{u} - \theta^*_u)^2}{\|\hat{\theta}^{(\lambda,I)} - \theta^*\|}}_{T_2}.
\end{align}

%In the above display we are implicitly assuming that $\|\hat{\theta}^{(\lambda,I)} - \theta^*\| > 0$ as if $\|\hat{\theta}^{(\lambda,I)} - \theta^*\| = 0$ then there is nothing to bound. 

%We have 
%$$\|\hat{\theta}_{I^{c}} - \theta^*_{I^c}\| \leq \|\hat{\theta} - \theta^*\| \leq \underbrace{2 \frac{\|\hat{\theta}_{I} - \theta^*_{I}\|^2}{\|\hat{\theta}- \theta^*\|}}_{T_1} + \underbrace{\frac{1}{\|\hat{\theta}- \theta^*\|} \sum_{u \in \R^{L_{d,n}}} \eta_u (\hat{\theta}_{u} - \theta^*_{u})^2}_{T_2}.$$

Now on the event $A_2$ we have 
\begin{equation}\label{eq:part1}
T_1 \leq 2 \sqrt{R(\theta^*, \lambda) + C \alpha \sigma^2 \log N}.
\end{equation}
Next, we proceed to bound $T_2.$

%Henceforth, within this proof, let us denote the quantity $3V(\theta^*) + C \sigma \sqrt{\log N}$ by $L.$ 
The event $A_1$ in particular implies that $ \|\hat{\theta}^{(\lambda,I)} - \theta^*\|^2 \leq N (L^*)^2$ where we denote $$L^* = V(\theta^*) + (2 + \sqrt{\alpha}) \sigma \sqrt{\log N}.$$%where $$L^* := 3V(\theta^*) + (2 + \sqrt{\alpha}) \sigma \sqrt{\log N}.$$
 Therefore, on the event $A_1 \cap A_2 \cap A_3$, we can write $T_{2} = \sum_{l = 0}^{\log_2 (L^*\sqrt{N}) \:-\:1} T_{2,l}$, where 
\begin{align*}
T_{2,l} := T_2 \;\mathrm{1}\Bigg(2^l < \|\hat{\theta}^{(\lambda,I)}- \theta^*\| \leq 2^{l + 1}, k_{\rdp}(\hat{\theta}^{(\lambda,I)}) &\leq \frac{1}{\lambda} (R(\theta^*, \lambda) + C\alpha \sigma^2 \log N),\\
&\qquad \qquad |\hat{\theta}^{(\lambda,I)} - \theta^*|_{\infty} \leq L^*\Bigg).
\end{align*}
To bound $T_2$, we now need to bound $T_{2,l}$ for each $l.$ We can now use Proposition~\ref{prop:sqrad} to deduce that there exists absolute constants $C_{1}$ and $C_{2}$ such that for any $\alpha \geq 1$,
\begin{equation}\label{eq:smallprob1}
P(A_{4,l}^c) \leq \frac{1}{N^{\alpha}}
\end{equation}
where the event $A_{4,l}$ is defined as 
\begin{equation*}
A_{4,l} := \left\{T_{2,l} \leq (C_{1} \sqrt{\alpha} +C_{2}) L^* \sqrt{\frac{1}{\lambda}\left(R(\theta^*, \lambda) + 2\alpha \sigma^2 \log N\right)\log N}\right\}.
\end{equation*}

Define the event 
\begin{equation*}
A_{4} := \left\{\cap_{l = 0}^{\log_2 (L^*\sqrt{N}) \:-\:1} A_{4,l}\right\}
\end{equation*}

Since $T_{2,l}$ cannot be positive for two distinct values of $l$, we can conclude that under the event $A_1 \cap A_2 \cap A_3 \cap A_4$, 
\begin{equation*}
%P(T_2 = \sum_{l = 0}^{\log_2 (NL^2) \:-\:1} T_{2,l} \geq (C_{1,d} \sqrt{\alpha} +C_{2,d}) \sqrt{k^* \log N}) \leq \frac{\log NL^2}{N^{C''}}.
T_2 \leq (C_{1} \sqrt{\alpha} +C_{2}) L^* \sqrt{\frac{1}{\lambda}\left(R(\theta^*, \lambda) + 2\alpha \sigma^2 \log N\right)\log N}.
\end{equation*}

By combining~\eqref{eq:loc},~\eqref{eq:part1} and the last display we can conclude that under the event $A_1 \cap A_2 \cap A_3 \cap A_4$ there exists constants $C_{1}$ and $C_{2}$ such that for any $\alpha \geq 1$, 
\begin{equation*}%\label{eq:part2}
%\|\hat{\theta}^{(\lambda,I)} - \theta^*\| \leq (C_{1} \sqrt{\alpha} +C_{2})L \sqrt{k_{\rdp}(\theta^*) \log N}.
\|\hat{\theta}^{(\lambda,I)} - \theta^*\| \leq \sqrt{R(\theta^*, \lambda) + C \alpha \sigma^2 \log N} \left\{2 + (C_{1} \sqrt{\alpha} +C_{2})L^* \sqrt{\frac{\log N}{\lambda}}\right\}.
\end{equation*}
This further implies that on the event $A_1 \cap A_2 \cap A_4$ there exists constants $C_{1}$ and $C_{2}$ such that for any $\alpha \geq 1$, 
\begin{equation}\label{eq:ineqdeterm}
%\|\hat{\theta}^{(\lambda,I)} - \theta^*\| \leq 1 + (C_{1} \sqrt{\alpha} +C_{2}) L \sqrt{k_{\rdp}(\theta^*) \log N}.
\|\hat{\theta}^{(\lambda,I)} - \theta^*\| \leq 1 + \sqrt{R(\theta^*, \lambda) + C \alpha \sigma^2 \log N} \left\{2 + (C_{1} \sqrt{\alpha} +C_{2})L^* \sqrt{\frac{\log N}{\lambda}}\right\}.
\end{equation}

Note that by~\eqref{eq:smallprob1} and a union bound argument we have
$$P(A_4) \geq 1 - \frac{\log_2 (L^*\sqrt{N})}{N^{\alpha}}.$$

Therefore, by Proposition \ref{prop:dcevents} and another union bound we have
$$P(A_1 \cap A_2 \cap A_4) \geq 1 - \frac{\log_2 (L^*\sqrt{N})}{N^{\alpha}} - C N^{-C \alpha},$$
for an absolute constant $C$.
The last display alongwith~\eqref{eq:ineqdeterm} shows that if $\lambda \geq C \sigma^2 \log N$ then with probability atleast $1 - C_1 \log (L^*\sqrt{N})N^{-C_2 \alpha}$ we have

\begin{equation*}
SSE(\hat{\theta}^{(\lambda, I)}, \theta^*) \leq C_3 \Bigg\{\big(R(\theta^*, \lambda) + \alpha \sigma^2 \log N\big)\alpha (V(\theta^*) + \sigma\sqrt{\log N})^2 \frac{\log N}{\lambda}\Bigg\}
\end{equation*}
where $C_3$ is an absolute constant.

Now, by construction of $\Lambda$ and the fact that $C \sigma^2 \log n < 2^{N^{*}}$ there exists a $\lambda^* \in \Lambda$ satisfying $C \sigma^2 \log n < \lambda^* \leq 2 C \sigma^2 \log n.$ Therefore, the previous display further implies that with probability atleast $1 - C_1 \log (L^*\sqrt{N})N^{-C_2 \alpha}$ we have

\begin{align*}
&\min_{\lambda \in \Lambda} SSE(\hat{\theta}^{(\lambda, I)}, \theta^*) \leq SSE(\hat{\theta}^{(\lambda^*, I)}, \theta^*)\\ &\leq C_3 \Bigg\{\big(R(\theta^*, \lambda^*) + \alpha \sigma^2 \log N\big)\alpha (V(\theta^*) + \sigma\sqrt{\log N})^2 \frac{\log N}{\lambda^*}\Bigg\}\\ &\leq C_3 \Bigg\{\big(R(\theta^*, 2 C \sigma^2 \log n) + \alpha \sigma^2 \log N\big) (\frac{V(\theta^*)}{\sigma} + \sqrt{\log N})^2\Bigg\},
\end{align*}
where in the last inequality we used the fact that $C \sigma^2 \log n < \lambda^* \leq 2 C \sigma^2 \log n.$
The last display finishes the proof. 
%{\color{red} Check whether you like this. If so, then change the statement of the theorem. Check the constants. Its a mess now.}

%The above display means that under the event $A_1 \cap A_2 \cap A_3 \cap A_5$, we have 
%This finishes the proof. 
\end{proof}

It now remains to prove Proposition~\ref{prop:dcevents} and Proposition~\ref{prop:sqrad}.

\subsection{Proof of Proposition~\ref{prop:dcevents}}

\begin{proof}[Proof of Proposition~\ref{prop:dcevents}]
	Fix $I \subset L_{d,n}$.

	\textbf{Part $1$:} In this part, we will prove the first assertion that the event $A_1$ holds with high probability where we recall
	$$A_1 = \left\{|\hat{\theta}^{(\lambda,I)} - \theta^*|_{\infty} \leq 3V(\theta^*) + (2 + \sqrt{\alpha}) \sigma \sqrt{\log N}\right\}.$$
	
	For any $u \in L_{d,n}$ we can represent 
	\begin{equation*}
	\hat{\theta}^{(\lambda,I)}_{u} = \overline{y}_{R \cap I} \mathrm{1}(R \cap I \neq \emptyset) + \overline{y}_{I} \mathrm{1}(R \cap I = \emptyset). 
	\end{equation*}
	where $R$ is the constituent rectangle of the optimal partition $\hat{\pi}^{(\lambda, I)}$ defining $\hat{\theta}^{(\lambda, I)}$ containing $u$, see~\eqref{eq:defnonempty} and~\eqref{eq:defempty}.
	%{\color{red} what is this lemma? Mention somewhere that the optimal partition may have empty cells.}
	
	Let's consider the case when $R \cap I$ is non empty as the other case can be done similarly. %If the reader is wondering, it is indeed possible that $R \cap I$ is empty. {\color{red} See discussion blah for more on this.}

	%It is possible that this rectangle $R$ is empty; in the sense that $
	
	We have
	\begin{align*}
	\hat{\theta}^{(\lambda,I)}_{u} - \theta^*_u =  \overline{y}_{R \cap I} - \theta^*_u = \overline{\theta^*}_{R \cap I} - \theta^*_u + \overline{\epsilon}_{R \cap I} 
	\end{align*}
	and thus, 
	\begin{align*}
	|\hat{\theta}^{(\lambda,I)}_{u} - \theta^*_u| \leq |\overline{\theta^*}_{R \cap I} - \theta^*_u| + \max_{\substack{R \subset L_{d,n}\\ R \:\:\text{is a dyadic rectangle}}} |\overline{\epsilon}_{R \cap I}| \leq &\\ V(\theta^*) + \max_{\substack{R \subset L_{d,n}\\ R \:\:\text{is a dyadic rectangle}}} |\overline{\epsilon}_{R \cap I}|.
	\end{align*}

	Now, the fact that for any dyadic rectangle $R$, $\overline{\epsilon}_{R \cap I}$ is sub-Gaussian with sub-Gaussian norm bounded by $\sigma^2$ and also the fact that the number of dyadic rectangles is bounded by $O_d(N)$ (see Lemma~\ref{lem:dyadic}) allow us to apply the standard result~Lemma \ref{lem:max_subG} about finite maxima of subgaussians. Thus, we obtain
	$$\max_{\substack{R \subset L_{d,n}\\ R \:\:\text{is a dyadic rectangle}}} |\overline{\epsilon}_{R \cap I}| \leq 2\sigma\sqrt{\log N} + \sigma \sqrt{\alpha \log N} = (2 + \sqrt{\alpha})\sigma\sqrt{\log N}$$
	with probability at least $1 - 2N^{-\alpha/2}$. This finishes the first part of this proof.

	%{\color{red} Write two lemmas. One is about representation and the other is about max of rectangles.}

	\textbf{Part $2$:} 
 	%In this part we show that the event $A_2$ holds with high probability. In particular, our goal is to show that there is an ev
 %	\begin{equation*}
 %	\frac{1}{2}||\theta^*_I - \hat{\theta}^{(\lambda, I)}_I||^2 + \frac{\lambda}{2} k_\rdp(\hat{\theta}^{(\lambda, I)}) \leq \frac{3}{2}||\theta_I - \theta^*_I||^2 + \lambda k_\rdp(\theta) + L(\epsilon_I, \lambda)
 %	\end{equation*}	
	
	From the definition of $\hat{\theta}^{(\lambda,I)}$ in~\eqref{eq:defncross}, for any fixed but arbitrary $\theta \in \R^{L_{d,n}}$, we have
	$$||y_I - \hat{\theta}^{(\lambda, I)}_I||^2 + \lambda k_{\rdp}(\hat{\theta}^{(\lambda, I)}) \leq ||y_I - \theta_I||^2 + \lambda k_{\rdp}(\theta).$$
	%Letting $\theta = \theta^*$, we have
	%$$||y_I - \hat{\theta}^{(\lambda, I)}_I||^2 + \lambda k(\hat{\theta}^{(\lambda, I)}) \leq ||y_I - \theta^*_I||^2 + \lambda k(\theta^*).$$
	Since $y = \theta^* + \epsilon$, this simplifies to
	\begin{align*}
	%\label{cont.}
	&||\theta^*_I - \hat{\theta}^{(\lambda, I)}_I||^2 + \lambda k_\rdp(\hat{\theta}^{(\lambda, I)}) \leq ||\theta_I^* - \theta_I||^2 + 2\left<\hat{\theta}^{(\lambda, I)}_I - \theta_I, \epsilon_I\right> + \lambda k_\rdp(\theta)\\
	&\leq ||\theta^*_I - \theta_I||^2 + 4 \left(\left<\frac{\hat{\theta}^{(\lambda, I)}_I - \theta_I}{||\hat{\theta}^{(\lambda, I)}_I - \theta_I||}, \epsilon_I\right>\right)^2 + \frac{1}{4}||\hat{\theta}^{(\lambda, I)}_I - \theta_I||^2 + \lambda k_\rdp(\theta)\\ \notag
	&\leq ||\theta^*_I - \theta_I||^2 + 4 \left(\left<\frac{\hat{\theta}^{(\lambda, I)}_I - \theta_I}{||\hat{\theta}^{(\lambda, I)}_I - \theta_I||}, \epsilon_I\right>\right)^2 + \frac{1}{2} ||\hat{\theta}^{(\lambda, I)}_I - \theta^*_I||^2 + \frac{1}{2} ||\theta_I - \theta^*_I||^2 + \lambda k_\rdp(\theta),
	\end{align*}
	where the second inequality follows from the elementary inequality, $2ab \leq 4a^2 + \frac{1}{4}b^2$ and the last inequality follows from the fact that $(a + b)^2 \leq 2 a^2 + 2 b^2$ for any two arbitrary real numbers $a$ and $b.$

	%Note that we are implicitly assuming that $\hat{\theta}^{(\lambda, I) \neq \theta.$ On the set $\hat{\theta}^{(\lambda, I) \neq \theta$, we anyway have 

	Note that we implicitly assumed that $\hat{\theta}^{(\lambda, I)}_{I} \neq \theta_{I}$ in the previous display. Simplifying further, we write
	\begin{align}\label{eq:propbasic1}
	&\frac{1}{2}||\theta^*_I - \hat{\theta}^{(\lambda, I)}_I||^2 + \frac{\lambda}{2} k_\rdp(\hat{\theta}^{(\lambda, I)}) \leq \frac{3}{2}||\theta_I - \theta^*_I||^2 + \lambda k_\rdp(\theta) + L(\epsilon_I, \lambda) \mathrm{1}(\hat{\theta}^{(\lambda, I)}_I \neq \theta_{I})
	\end{align}
	
	where 
	\begin{equation}\label{eq:L_eps}
	L(\epsilon_I, \lambda) := 4 \left(\left<\frac{\hat{\theta}^{(\lambda, I)}_I - \theta_I}{||\hat{\theta}^{(\lambda, I)}_I - \theta_I||}, \epsilon_I\right>\right)^2 - \frac{\lambda}{2} k_\rdp(\hat{\theta}^{(\lambda, I)}).
	\end{equation}
	Now, for any $\alpha > 1$, we show in Lemma \ref{lem:L_eps} that there exists an absolute constant $C > 0$ such that if $\lambda \geq C \sigma^2 \log n$ then 
	$$P(\sup_{\theta \in \R^{L_{d,n}}} L(\epsilon_I, \lambda) \mathrm{1}(\hat{\theta}^{(\lambda, I)}_I \neq \theta_{I}) > \alpha \sigma^2 \log N) \leq 2 N^{-C\alpha}.$$
	Combining the previous display with~\eqref{eq:propbasic1} finishes the proof. 
\end{proof}

In order to finish the proof of Proposition~\ref{prop:dcevents}, we now need to provide the proofs of the lemmas used within Proposition~\ref{prop:dcevents}.

\begin{lemma}\label{lem:L_eps}
	Recall the definition of $L(\epsilon_I, \lambda)$ from \eqref{eq:L_eps}. There exists a positive absolute constant $C$ such that if $\lambda \geq C\sigma^2\log N$, then for any $\alpha > 1$, 
	$$P(\sup_{\theta \in \R^{L_{d,n}}} L(\epsilon_I, \lambda) \mathrm{1}(\hat{\theta}^{(\lambda, I)}_I \neq \theta_{I}) > C \alpha \sigma^2 \log N) \leq C N^{-C \alpha}.$$ 
\end{lemma}
\begin{proof}
	Recall the definition of $L(\epsilon_I, \lambda)$ from \eqref{eq:L_eps}.
	By taking supremum over the possible values of $\hat{\theta}^{(\lambda, I)}$, we obtain
	\begin{align*}
	\sup_{\theta \in \R^{L_{d,n}}} L(\epsilon_I, \lambda) \mathrm{1}(\hat{\theta}^{(\lambda, I)}_I \neq \theta_{I}) \leq \max_{k \in [N]}\left[4 \sup_{\substack{S \in \cS_\rdp,\\ Dim(S) = k}} \sup_{\substack{v \in S \\ v_I \neq \theta_I}}\left(\left<\frac{v_I - \theta_I}{||v_I - \theta_I||}, \epsilon_I\right>\right)^2 - \frac{\lambda}{2} k\right].
	\end{align*}
	
	Here $\cS_\rdp$ refers to the collection of subspaces $\{S_{\pi}: \pi \in \mathcal{P}_{\rdp,d,n}\}.$

	%{\color{red} Where is $ \cS_\rdp$ defined? Explain this above line a bit?}
	
	%We will now argue conditionally on $I.$ 
	For any subspace $S \in \cS_\rdp$, we can define a corresponding subspace of $\R^{|I|}$,
	$$S_{I} := \{v_{I}: v \in S\}.$$ Note that, we must have $Dim(S_{I}) \leq Dim(S).$ Armed with this observation, we can further write
	\begin{align*}
	\sup_{\theta \in \R^{L_{d,n}}} L(\epsilon_I, \lambda) \mathrm{1}(\hat{\theta}^{(\lambda, I)}_I \neq \theta_{I}) \leq \max_{k \in [N]} \left[4 \underbrace{\sup_{\substack{S \in \cS_\rdp,\\ Dim(S_I) = k}} \sup_{\substack{v \in S_I \\ v \neq \theta_I}}\left(\left<\frac{v - \theta_I}{||v - \theta_I||}, \epsilon_I\right>\right)^2}_{W_k} - \frac{\lambda}{2} k\right].
	\end{align*}
	
	We will now bound the random variables $W_k$ for each $k \in [N].$ For any $k \in [N]$, fix a subspace $S \in \cS_\rdp$ such that $Dim(S_I) = k$. 
	
	By Lemma~\ref{lem:HSsubG} (stated and proved after this proof) we can assert that for any $u > 0$, there exist absolute positive constants $c, C$ such that with probability at least $1 - 2e^{-Cu}$,
	$$\left(\sup_{\substack{v \in S_I\\ v \neq \theta_I}}\left<\frac{v - \theta_I}{||v - \theta_I||}, \epsilon_I\right>\right)^2 \leq c \sigma^2 (k + 1) + 4 \sigma^2 u.$$

We can now use a union bound argument along with the fact that $|S \in \cS_\rdp: Dim(S) = k| \leq n^{c'k}.$ This cardinality bound follows because this cardinality is clearly at most the number of distinct rectangles in $L_{d,n}$ raised to the power $k$ and the number of distinct rectangles of $L_{d,n}$ is at most $N^2.$ Hence, for some absolute constant $c' > 0$, we obtain for all $u > 0$,
	$$ P(W_k > c \sigma^2 (k + 1) + 4 \sigma^2 u) \leq 2 \exp(c' k \log n - c u).$$ Equivalently, by reparametrizing $u =  2\frac{c'}{c} k \log n  + t$ we can write for all $t > 0$, 
	$$P(W_k > c \sigma^2 (k + 1) + 8 \sigma^2 \frac{c'}{c} k \log n + 4 \sigma^2 t) \leq 2 \exp(-c' k \log n - c t).$$

	The above now implies the existence of an absolute constant $C_1$ such that if $\lambda \geq C_1 \sigma^2 \log n$, then for all $t > 0$,
	$$P(4 W_k - \frac{\lambda}{2} k > 4 \sigma^2 t) \leq 2 \exp(-c' k \log n - c t).$$
	
	A further union bound implies that for an absolute constant $C_2 > 0$,
	$$P(L(\epsilon_I, \lambda) > 4 \sigma^2 t) \leq \sum_{k = 1}^{N} 2 \exp(-c' k \log n - c t) \leq C_2 \exp(-c t).$$

	Setting $t = \alpha \log N$ finishes the proof.%and noticing that since the above bound does not depend on $I$ one can now integrate over the random choice of $I$ finishes the proof.
\end{proof}

The following lemma is essentially the same as Lemma $9.1$ in~\cite{chatterjee2019adaptive}, the only difference being we would need this lemma to hold for a general subgaussian random variable instead of a Gaussian random variable as was done in~\cite{chatterjee2019adaptive}.
\begin{lemma}\label{lem:HSsubG}
	Let $Z \in \R^n$ be a random vector with independent sub-Gaussian entries with mean $0$ and sub-Gaussian norm $\sigma^2$, $S$ be a subspace in $\R^n$ and $\theta \in \R^n$ be any fixed vector.  Then there exist absolute constants $c$ and $C$ such that for any $u > 0$, with probability at least $1 - 2e^{-C u}$ we have
	%$$\left(\sup_{v \in S, v \neq \theta}\left<\frac{v - \theta}{||v - \theta||}, Z\right>\right)^2 \leq 2 Dim(S) + 4(1 + u^2).$$
	$$\left(\sup_{\substack{v \in S\\ v \neq \theta}}\left<\frac{v - \theta}{||v - \theta||}, Z\right>\right)^2 \leq c\sigma^2 \big(Dim(S) + 1\big) + 4 \sigma^2u.$$
\end{lemma}
\begin{proof}
	Let $P_S$ be the orthogonal projection matrix on to the subspace $S$ and $P_{S'}$ be the orthogonal projection matrix on to the subspace $S'$, where $S'$ is the one dimensional subspace spanned by the vector $(I - P_S)\theta$.
	Then we can write,
	\begin{align}\notag
	&\left|\sup_{v \in S, v \neq \theta}\left<\frac{v - \theta}{||v - \theta||}, Z\right>\right| = \left|\sup_{v \in S, v \neq \theta} \left<\frac{v - P_S\theta -(I - P_S)\theta}{\sqrt{||v -P_S\theta||^2 + ||(I - P_S)\theta||^2}}, Z\right>\right|\\ \notag
	&\leq \left|\sup_{v \in S, v \neq \theta}\left<\frac{v - P_S\theta}{\sqrt{||v -P_S\theta||^2 + ||(I - P_S)\theta||^2}}, Z\right>\right|\\ \label{bound_inner_pdt}
	&\qquad \qquad \qquad + \left|\sup_{v \in S, v \neq \theta}\left<\frac{(I - P_S)\theta}{\sqrt{||v -P_S\theta||^2 + ||(I - P_S)\theta||^2}}, Z\right>\right|\\ \notag
	&\leq \left|\sup_{v \in S, ||v|| \leq 1}\left<Z, v\right>\right| + \left|\sup_{v \in S', ||v||\leq 1} \left<Z, v\right>\right|,
	\end{align}
	%Now we will give high probability bound on both the terms above separately.
	Note that, for any subspace $S$,
	$$\sup_{v \in S, ||v|| \leq 1}\left<Z, v\right> = ||P_SZ||.$$
	Therefore, from \eqref{bound_inner_pdt},
	\begin{equation}\label{bound_inner_pdt_sq}
	\left(\sup_{\substack{v \in S\\ v \neq \theta}}\left<\frac{v - \theta}{||v - \theta||}, Z\right>\right)^2 \leq 2||P_SZ||^2 + 2||P_{S'}Z||^2 = 2Z^TP_SZ + 2Z^TP_{S'}Z.
	\end{equation}
	Now we will give a high probability bound on both of the above terms separately.
	From the Hanson-Wright concentration inequality, see Theorem 6.2.1 in \cite{vershynin2018high}, we obtain that, for any matrix $A_{m \times n}$ and $u > 0$, %large enough,
	\begin{equation}\label{HSineq}
	P\left(|Z^TAZ - \E(Z^TAZ)| > u\right) \leq 2\exp\left( - C \min\left\{\frac{u}{\sigma^2||A||_{OP}},\frac{u^2}{\sigma^4 \|A\|_{F}^2}\right\}\right),
	\end{equation}
	for some $C > 0$. Here, the notation $\|A\|_{OP}$ and $\|A\|_{F}$ refers to the operator norm and the frobenius norm of the matrix $A$ respectively. In our case, when $A = P_S$, we have $||A||_{OP} = 1$ and $\|A\|_{F}^2 = Dim(S).$ It can be checked that if $u > \sigma^2 Dim (S)$ the first term inside the brackets in the R.H.S above dominates.

	Now letting $A = P_S$, we have $\E(Z^TAZ) = trace(P_S\Sigma)$, where $\Sigma$ is the diagonal covariance matrix of $Z$. %Since $\{Z_i : i \in [n]\}$ are independent with sub-Gaussian norm $\sigma^2$, we have $\Sigma_{ij} = 0$ for any $i \neq j$, and  $i, j \in [n]$ and for every $i \in [n]$, $\Sigma_{ii} \leq c\sigma^2$, for some $c > 0$. 
	Therefore, we have
	\begin{align*}
	\E(Z^TAZ) = trace(P_S\Sigma) &\leq trace(P_S)\;\max_{i \in [n]} \Sigma_{ii} \leq c\sigma^2 Rank(P_S) = c\sigma^2Dim(S)
	\end{align*}
	for some absolute constant $c > 0$. Similarly by letting $A = P_{S'}$, we have $trace(P_{S'}\Sigma) \leq c\sigma^2$ since $Dim(S') = 1$. %Moreover, for both $A = P_S, P_{S'}$, we have $||A||_{OP} = 1$.
	Therefore, by \eqref{HSineq} and the arguments outlined thereafter, we can actually assert that for any $u > 0$, both the following events happen
	\begin{align*}
	Z^TP_SZ &\leq c\sigma^2Dim(S) + u\\
	Z^TP_{S'}Z &\leq c\sigma^2 + u
	\end{align*}
	with probability at least $1 - 2\exp\left(-C\frac{u}{\sigma^2}\right)$. Finally, by using \eqref{bound_inner_pdt_sq}, the proof is complete.
\end{proof}

\subsection{Proof of Proposition~\ref{prop:sqrad}}

\begin{proof}[Proof of Proposition~\ref{prop:sqrad}]
Recall that for any partition $\pi \in \mathcal{P}_{\rdp,d,n}$ we denote by $S_\pi$ the subspace of $\R^{L_{d,n}}$ which consists of all arrays which are constant on every subrectangle of $\pi.$ Now observe that we can write
\begin{align*}
&T_{l} = \sum_{u \in L_{d, n}} \frac{\eta_{u} (\hat{\theta}^{(\lambda,I)}_{u} - \theta^*_u)^2}{\|\hat{\theta}^{(\lambda,I)} - \theta^*\|} \mathrm{1}\Big(2^l < \|\hat{\theta}^{(\lambda,I)}- \theta^*\| \leq 2^{l + 1}, k_{\rdp}(\hat{\theta}^{(\lambda,I)}) \leq M,\\
&\qquad \qquad \qquad \qquad \qquad \qquad \qquad \qquad |\hat{\theta}^{(\lambda,I)} - \theta^*|_{\infty} \leq L\Big) \\
&\leq \max\left\{0,\sup_{\substack{\pi \in \mathcal{P}_{\rdp,d,n}, \\ |\pi| \leq M}} \: \sup_{\substack{v \in S_{\pi}: 2^l < \|v - \theta^*\| \leq 2^{l + 1},\\ |v - \theta^*|_{\infty} \leq L}} \sum_{u \in L_{d, n}} \frac{\eta_{u} (v_{u} - \theta^*_u)^2}{\|v - \theta^*\|}\right\} \\
&\leq 
2L \; \max\left\{0,\sup_{\substack{\pi \in \mathcal{P}_{\rdp,d,n},\\ |\pi| \leq M}}  \frac{1}{2^l} \sup_{\substack{v \in S_{\pi}: 2^l < \|v - \theta^*\| \leq 2^{l + 1}, \\ |v - \theta^*|_{\infty} \leq L}} \sum_{u \in L_{d, n}} \frac{1}{2L} \eta_{u} (v_{u} - \theta^*_u)^2\right\}.
\end{align*}

%We have
%\begin{equation*}
%T_{2,l} \leq 2L \max\{0,\sup_{\pi: |\pi| \lesssim k^*} \frac{1}{2^l} \sup_{v \in \pi: \|v- \theta^*\| \leq 2^{l + 1}, \|v- \theta^*\|_{\infty} \leq L} \sum_{u \in \R^{L_{d,n}}} \frac{1}{2L}\eta_u (\hat{\theta}_{u} - \theta^*_{u})^2\}
%\end{equation*}
For any fixed subspace $S_{\pi}$ indexed by $\pi \in \mathcal{P}_{\rdp,d,n}$ and any integer $l \geq 0$, let us define the random variable $$X_{\pi,l} := \sup_{\substack{v \in S_{\pi}: \|v- \theta^*\| \leq 2^{l + 1},\\ |v- \theta^*|_{\infty} \leq L}} \sum_{u \in {L_{d,n}}} \frac{1}{2L}\eta_u (v_{u} - \theta^*_{u})^2.$$

Concurrently, let us also define the random variable
$$ Y_{\pi,l} := \sup_{\substack{v \in S_{\pi}: \|v- \theta^*\| \leq 2^{l + 1}, \\ |v- \theta^*|_{\infty} \leq L}} \sum_{u \in {L_{d,n}}} \eta_u (v_{u} - \theta^*_{u}).$$

It is easier to control $Y_{\pi,l}$ than $X_{\pi,l}$ since $Y_{\pi,l}$ is a rademacher complexity term. So we will bound $Y_{\pi,l}$ and then conclude a similar bound for $X_{\pi,l}$ by arguing that $X_{\pi,l}$ is dominated by $Y_{\pi,l}$ in a certain sense. This is accomplished by using the contraction principle for rademacher complexity (Theorem $2.2$ in~\cite{koltchinskii2011oracle}) and a concentration inequality for convex lipxhitz functions of independent rademachers (see Theorem $6.10$ in~\cite{boucheron2013concentration}.) 
%The bound for $X_{\pi,l}$ is stated in the following lemma.

Altogether, in order to obtain a high probability bound for $X_{\pi, l}$ the above stated strategy is employed in Lemma~\ref{lem:xbound} and Lemma~\ref{lem:old}, which together imply that there is an absolute constant $C$ such that for all $u \geq 1$,
\begin{equation*}
P\left(\frac{1}{2^{l + 1}} X_{\pi,l} \geq  C \sqrt{M} + \frac{u}{2^l}\right) \leq \exp\left(-\frac{u^2}{32 \cdot 2^{2l}}\right).
\end{equation*}
Setting $\frac{u}{2^l} = (C_{1} \sqrt{\alpha} +C_{2}) \sqrt{M  \log N}$ in the above display for large enough constants $C_{1}$ and $C_{2}$ and using a union bound argument now gives us
\begin{equation}\label{eq:xbound}
P\left(\sup_{\substack{\pi \in \mathcal{P}_{\rdp, d, n}, \\ |\pi| \leq M}} \frac{1}{2^l} X_{\pi,l} \geq (C_{1} \sqrt{\alpha} +C_{2}) \sqrt{M \log N}\right) \leq  \exp \left(O_{d}\left(M  \log N\right) - \frac{u^2}{2^{2l}}\right) \leq \frac{1}{N^{\alpha}}.
\end{equation}
In the above display, we used the fact that the number of dyadic rectangles of $L_{d,n}$ is at most $2^d N$, see Lemma~\ref{lem:dyadic}, and hence the number of recursive dyadic partitions $\pi$ with $|\pi| \leq M$ is at most $(2^d N)^M = \exp\left(O_{d}\left(M \log N\right)\right)$ where $O_d$ hides the constants depending on $d.$ 

Recall that for any $l$, 
\begin{equation*}
T_{l} \leq 2L \max\left\{0,\sup_{\substack{\pi \in \mathcal{P}_{\rdp, d, n},\\  |\pi| \leq M}} \frac{1}{2^l} X_{\pi,l}\right\}.
\end{equation*}
Therefore, by~\eqref{eq:xbound} we have for each $l$,
\begin{equation*}
P\left(T_{l} \geq (C_{1} \sqrt{\alpha} +C_{2}) L \sqrt{M \log N}\right) \leq \frac{1}{N^{\alpha}},
\end{equation*}
and this completes the proof.
\end{proof}

In order to finish the proof of Proposition~\ref{prop:sqrad}, we now need to provide proofs of the lemmas used within the above proof.

The following result provides a concentration inequality for $Y_{\pi,l}.$
\begin{lemma}\label{lem:conlipconc}
	For all $u \geq 1$, the random variable $Y_{\pi,l}$ satisfies the following concentration inequality:
	\begin{equation}\label{eq:conlipconc}
	P\left(Y_{\pi,l} \geq \mathbb{E} Y_{\pi,l} + u\right) \leq \exp\left(-\frac{u^2}{8 \cdot 2^{2l}}\right).
	\end{equation}
\end{lemma}

\begin{proof}%[Proof of Lemma~\ref{lem:conlipconc}]
	The random variable $Y_{\pi,l}$ is a convex lipschitz function of $\{\eta_u: u \in L_{d, n}\}$ with lipschitz constant at most $2^{l + 1}.$ Convexity follows because $Y_{\pi,l}$ is a supremum of linear functions. The lipschitz property follows by similar arguments as in Lemma $A.2$ in~\cite{chatterjeel2019adaptive}. Since $Y_{\pi,l}$ is a convex lipschitz function of $\{\eta_u: u \in L_{d, n}\}$ the concentration inequality in~\eqref{eq:conlipconc} follows from Theorem $6.10$ in~\cite{boucheron2013concentration}. 
\end{proof}
The next lemma states that $X_{\pi,l}$ is stochastically dominated by $Y_{\pi,l}$ in the following sense:
\begin{lemma}\label{lem:stochdom}
	We have 
	$\mathbb{E} F\left(\frac{1}{2} X_{\pi,l}\right) \leq \mathbb{E} F(Y_{\pi, l})$ for all non decreasing convex functions $F: \mathbb{R} \rightarrow \mathbb{R}.$
\end{lemma}

\begin{proof}
	We observe that the function $x \rightarrow x^2$ is $2L$ lipschitz when $|x| \leq L.$ The result now follows directly from the Contraction principle for Rademacher Complexity; see Theorem $2.2$ in~\cite{koltchinskii2011oracle}.
\end{proof}

The following lemma is a modified version of Panchenko's lemma; see Lemma $7.6$ in~\cite{van2014probability}, which is useful to provide a probability tail bound for $X_{\pi, l}$.
\begin{lemma}[Panchenko's Lemma]\label{lem:panchenko}
	Suppose $X,Y$ are real valued random variables such that $\mathbb{E} F(X) \leq \mathbb{E} F(Y)$ for every non decreasing convex function $F: \mathbb{R} \rightarrow \mathbb{R}.$ Suppose also the following inequality is true for some $a > 0, c > 0$ and every $u > 0$,
	\begin{equation*}
	P(Y \geq a + u) \leq \exp(-c u^2).
	\end{equation*}
	Then the following inequality is also true for all $u > 1$,
	\begin{equation*}
	P(X \geq a + u) \leq \exp(-c u^2/4).
	\end{equation*}
\end{lemma}

\begin{proof}
	Take the non decreasing convex function $F(x) = (x - t)_{+}$ for any $t \in \R.$ Note that $\E F(X) = \E (X - t)_{+} = \int_{t}^{\infty} P(X \geq s)\:ds.$ We now obtain for any $t \geq a$,
	\begin{align*}
	&\int_{t}^{\infty} P(X \geq s)\:ds \leq \int_{t}^{\infty} P(Y \geq s)\:ds \leq \int_{t}^{\infty} \exp(-c (s - a)^2) ds = \\& \int_{t - a}^{\infty} \exp(-c u^2) du = P(N(0,\sigma^2 = \frac{1}{2c}) \geq t - a) \leq \exp(-c (t - a)^2).
	\end{align*}
	where the first two inequalities follow due to our assumptions and the last inequality follows by the standard Mills ratio upper bound to Gaussian tails.
	
	Theefore, this implies that
	\begin{align*}
	&P(X \geq a + u) \leq \frac{1}{u - u'} \int_{a + u'}^{a + u} P(X \geq s) ds \leq \frac{1}{u - u'} \int_{a + u'}^{\infty} P(X \geq s) ds \leq \frac{1}{u - u'} \exp(-c u'^2).
	\end{align*}
	Setting $u' = u/2$ we get 
	\begin{align*}
	&P(X \geq a + u) \leq \frac{2}{u} \exp(-c \frac{u^2}{4}) \leq 2 \exp(-c \frac{u^2}{4})
	\end{align*}
	where the last inequality is true because $u > 1.$
\end{proof}

Finally, we have the following result that provides a probability tail bound for $X_{\pi, l}$.
\begin{lemma}\label{lem:xbound}
	For all $u \geq 1$, we have the following probability tail bound
	\begin{equation*}%\label{eq:xconc}
	P\left(\frac{1}{2} X_{\pi,l} \geq \E Y_{\pi,l} + u\right) \leq \exp\left(-\frac{u^2}{32  \cdot 2^{2l}}\right).
	\end{equation*}
\end{lemma}
\begin{proof}
	First, by Lemma \ref{lem:conlipconc}, we obtain the following concentration inequality for $Y_{\pi, l}$, that for every $u > 0$,
	\begin{equation*}
	P\left(Y_{\pi,l} \geq \mathbb{E} Y_{\pi,l} + u\right) \leq \exp\left(-\frac{u^2}{8 \cdot 2^{2l}}\right).
	\end{equation*}
	Next, by Lemma \ref{lem:stochdom}, we obtain the fact that $X_{\pi, l}$ is stochastically dominated by $Y_{\pi, l}$. 
	Therefore, we can use Panchenko's lemma; see Lemma~\ref{lem:panchenko}, %in Section~\ref{sec:auxiliary}, 
	to conclude that for all $u \geq 1$,
	\begin{equation}\label{eq:xconc}
	P\left(\frac{1}{2} X_{\pi,l} \geq \E Y_{\pi,l} + u\right) \leq \exp\left(-\frac{u^2}{32 \cdot 2^{2l}}\right).
	\end{equation}
\end{proof}

Now, in view of Lemma \ref{lem:xbound}, we need to provide an upper bound of $\E Y_{\pi,l}$, which is done in the following lemma. 

\begin{lemma}\label{lem:old}
	For any $\pi \in \mathcal{P}_{\rdp,d,n}$ with $Dim(S_{\pi}) \leq M$ we have
	\begin{equation*}
	\E Y_{\pi,l} \leq C 2^l \sqrt{M},
	\end{equation*}
	for some $C > 0$.
\end{lemma}
\begin{proof}
	Recall that $$ Y_{\pi,l} = \sup_{\substack{v \in S_{\pi}: \|v- \theta^*\| \leq 2^{l + 1},\\ |v- \theta^*|_{\infty} \leq L}} \sum_{u \in {L_{d,n}}} \eta_u (v_{u} - \theta^*_{u}).$$
	%Also recall that $\pi \in \mathcal{P}_{\rdp,d,n}$ such that $Dim(S_{\pi}) \leq \sqrt{\frac{2}{\lambda} R(\theta^*)}.$
	We can write
	\begin{align*}
	\E Y_{\pi,l} \leq  \E \sup_{v \in S_{\pi}: \|v- \theta^*\| \leq 2^{l + 1}} \sum_{u \in {L_{d,n}}} \eta_u (v_{u} - \theta^*_{u}) \leq \sqrt{\pi/2}\:\: \E \sup_{v \in S_{\pi}: \|v- \theta^*\| \leq 2^{l + 1}} \sum_{u \in {L_{d,n}}} Z_u (v_{u} - \theta^*_{u})
	\end{align*}
	where $Z$ is an array consisting of i.i.d $N(0,1)$ entries. In the last display we bounded the rademacher complexity by gaussian complexity which is a well known result; see Page $132$ in~\cite{wainwright2019high}. Now note that we can write
	\begin{align*}
	&\quad\; \E \sup_{v \in S_{\pi}: \|v- \theta^*\| \leq 2^{l + 1}} \sum_{u \in {L_{d,n}}} Z_u (v_{u} - \theta^*_{u})\\ &= \E \sup_{v \in S_{\pi}: \|v- \theta^*\| \leq 2^{l + 1}} \sum_{u \in {L_{d,n}}} Z_u (v_{u} - (O_{\pi}\theta^*)_{u}) + \E  \sum_{u \in {L_{d,n}}} Z_u (\theta^*_u - (O_{\pi}\theta^*)_{u})\\ &= \E \sup_{v \in S_{\pi}: \|v\| \leq 2^{l + 1}} \sum_{u \in {L_{d,n}}} Z_u v_{u} = 2^{l + 1} \E \sup_{v \in S_{\pi}: \|v\| \leq 1} \sum_{u \in {L_{d,n}}} Z_u v_{u}\\ &= 2^{l + 1} \sqrt{Dim(S_{\pi})} \leq 2^{l + 1} \sqrt{M}.
	\end{align*}
	In the display above, $O_\pi$ refers to the orthogonal projection matrix for the subspace $S_\pi$ and the last equality follows from standard facts about projection of a standard gaussian vector onto a subspace. 
\end{proof}

\subsection{Auxiliary Lemmas}\label{sec:auxiliary}

\begin{lemma}\label{lem:dyadic}
	The number of dyadic rectangles of $L_{d,n}$ is at most $2^d N$.
\end{lemma}
\begin{proof}
	Any dyadic rectangle must have each side length dyadic or equivalently, it is a product of dyadic rectangles. In any dimension, the number of dyadic intervals is at most $2n$ which then furnishes the lemma as $(2n)^d = 2^d N.$

	In one dimension, we can count the number of dyadic intervals as follows. Start bottom up and there are at most $n$ singletons which are dyadic intervals of length $1.$ Then there are at most $n/2$ dyadic intervals of length $2$ and so on. Hence the total number of dyadic intervals is at most $n + n/2 + n/4 + \dots \leq 2n.$ 
\end{proof}

The following result is a standard tail inequality for a maximum of finitely many subgaussian random variables; also see Lemma $5.2$ in~\cite{van2014probability}. 
\begin{lemma}\label{lem:max_subG}
	Let $Z_i$ be a sub-Gaussian random variable with sub-Gaussian norm $\sigma_i^2$, where $i \in [n]$. Then for all $x \geq 0$,
	$$P\left(\max_{1 \leq i \leq n} |Z_i| > \sqrt{2\sigma^2 \log n} + x\right) \leq 2\exp(-x^2/2\sigma^2),$$
	where $\sigma^2 := \max_{1 \leq i \leq n} \sigma_i^2$.
\end{lemma}
\begin{proof}
	We have, for any $t > 0$,
	\begin{align*}
	P\left(\max_{1 \leq i \leq n} |Z_i| > t\right) &= P\left(\bigcup_{i = 1}^n \{|Z_i| > t\}\right)\\ 
	&\leq \sum_{i = 1}^n P(|Z_i| > t) \leq \sum_{i = 1}^n 2 \exp(-t^2/2\sigma_i^2)\\
	&\leq \sum_{i = 1}^n 2 \exp(-t^2/2\sigma^2) = 2n\exp(-t^2/2\sigma^2),
	\end{align*}
	and the proof follows.
\end{proof}

\section{Proofs for Trend Filtering}\label{sec:tfproofs}

\subsection{Sketch of Trend Filtering Proofs}\label{sec:tvsketch}

%\subsection{\textbf{Sketch of Proof}}\label{sec:tvsketch}
For the convenience of the reader, we first present a sketch of proof of Theorems~\ref{thm:slowratetv} and~\ref{thm:fastratetv}. %This sketch is divided into several steps and is meant to convey the essential aspects of our proof strategy.

Fix $r \geq 1.$ Let $\hat{\theta}^{(r)}_{\lambda}$ be the usual $r$th order Trend Filtering estimators based on the full data as defined in~\eqref{eq:tfdefn}.  Also, to reduce notational clutter let us denote the completion estimators $\hat{\theta}^{(\lambda,I,r)}$ by simply $\hat{\theta}^{(\lambda,I)}$ where $I = I_j$ for $j \in [K] = [r + 1].$

In view of Theorem \ref{thm:main}, it suffices to bound two main quantities: $\min_{\lambda \in \Lambda} SSE(\hat{\theta}^{(r)}_{\lambda}, \theta^*)$ and $\min_{\lambda \in \Lambda} SSE(\hat{\theta}^{(\lambda,I)}_{I^c}, \theta^*_{I^c})$.
State of the art existing bounds for penalized Trend Filtering are available in~\cite{ortelli2019prediction}. In particular, Theorem $1.1$ in~\cite{ortelli2019prediction} gives the desired upper bounds (both slow and fast rates) on $SSE(\hat{\theta}^{(r)}_{\lambda}, \theta^*)$ as long as $\lambda$ is chosen appropriately depending on the unknown $\sigma$ and some properties of $\theta^*.$ We choose $\Lambda_j = \Lambda = \{1,2,2^2,\dots,2^{N^*}\}$ for all $j \in [K]$ and assume that $2^{N^*}$ exceeds these theoretically optimal choices of $\lambda.$ Therefore, by construction of $\Lambda$, there exists a $\lambda^* \in \Lambda$ which (up to a factor of $2$) scales like the theoretically recommended choice in Theorem $1.1$.

Clearly, $\min_{\lambda \in \Lambda} SSE(\hat{\theta}^{(\lambda)}, \theta^*) \leq SSE(\hat{\theta}^{(\lambda^*)}, \theta^*).$ We can now use the existing bound given by Theorem $1.1$ in~\cite{ortelli2019prediction} to bound $SSE(\hat{\theta}^{(\lambda^*)}, \theta^*)$ for the usual Trend Filtering estimator by generalizing its proof to hold for general subgaussian errors. Therefore, the desired bound for $\min_{\lambda \in \Lambda} SSE(\hat{\theta}^{(\lambda)}, \theta^*)$ again follows pretty much directly from the existing result Theorem $1.1$ in~\cite{ortelli2019prediction}.

As in the Dyadic CART proof, the main new task for us here is to bound the quantity:  $\min_{\lambda \in \Lambda} SSE(\hat{\theta}^{(\lambda,I)}_{I^c}, \theta^*_{I^c})$, where $I = I_j$ for $j \in [K].$ For each $j \in [r + 1]$ our eventual bound will be the same and hence our final bound on the MSE of the CVTF estimator will merely be $(r + 1)$ times the bound for $\min_{\lambda \in \Lambda_j} SSE\left(\hat{\theta}^{(\lambda, I^c)}_{I}, \theta^*_{I}\right)$ when $I= I_1.$ Therefore, we can just consider $I = I_1$ below.

Again, since $\min_{\lambda \in \Lambda} SSE(\hat{\theta}^{(\lambda,I)}_{I^c}, \theta^*_{I^c}) \leq SSE(\hat{\theta}^{(\lambda^*,I)}_{I^c}, \theta^*_{I^c})$ it is sufficient for us to bound $SSE(\hat{\theta}^{(\lambda^*,I)}_{I^c}, \theta^*_{I^c})$ which in turn is trivially upper bounded by $SSE(\hat{\theta}^{(\lambda^*,I)}, \theta^*)$.

We now outline the main steps in our proof which bounds $SSE(\hat{\theta}^{(\lambda^*,I)}, \theta^*)$. %For the convenience of the reader, we now divide the proof sketch into various steps. 

\begin{enumerate}
	%	\item Step 1: By Theorem~\ref{thm:main}, it suffices to bound for any $j \in [r + 1]$, the terms $$\min_{\lambda \in \Lambda_j} SSE\left(\hat{\theta}^{(\lambda, I_j^c)}_{I_j}, \theta^*_{I_j}\right).$$
	
	%	 For each $j \in [r + 1]$ our eventual bound will be the same and hence the leading order term of our final bound on the MSE of the CVTF estimator will merely be $r + 1$ times the bound for $\min_{\lambda \in \Lambda_j} SSE\left(\hat{\theta}^{(\lambda, I_j^c)}_{I_j}, \theta^*_{I_j}\right)$ when $j = 1.$ Firstly, we write the trivial inequality
	%	\begin{equation*}
	%	SSE\left(\hat{\theta}^{(\lambda, I_j^c)}_{I_j}, \theta^*_{I_j}\right) \leq SSE\left(\hat{\theta}^{(\lambda, I_j^c)}, \theta^*\right).
	%	\end{equation*} 
	%	Therefore, it suffices to bound $SSE\left(\hat{\theta}^{(\lambda, I_j^c)}, \theta^*\right)$ for any good choice of $\lambda$ in $\Lambda_j.$ 

	\item \textbf{KKT Condition:} %The estimator $\hat{\theta}^{(\lambda, I_j^c)}$ is exactly the penalized trend filtering estimator when the data vector is $\tilde{y}(I_j^c).$ Since the Trend Filtering optimization objective function is convex, it is natural to write down the KKT condition. 
	
	Recall the definition of $\tilde{y}(I_j^c)$ in the description of the CVTF estimator (see  Step $3$ in Section~\ref{sec:cvtf}).
	Define $$\tilde{\theta}(I_j^c) := E\left[\tilde{y}(I_j^c)\right],\ \text{and}\ \ \tilde{\epsilon}(I_j^c) := \tilde{y}(I_j^c) - \tilde{\theta}(I_j^c).$$ 
	%Therefore, $\hat{\theta}^{(\lambda, I_j^c)}$ should be a good estimator of $\tilde{\theta}(I_j^c).$

	Since $\tilde{y}(I_j^c) = \tilde{\theta}(I_j^c) + \tilde{\epsilon}(I_j^c)$, we can think of $\tilde{y}(I_j^c)$ as a noisy version of  $\tilde{\theta}(I_j^c).$ Therefore, $\hat{\theta}^{(\lambda, I_j^c)}$, being the trend filtering estimator applied to $\tilde{y}(I_j^c)$ should estimate the mean of $\tilde{y}(I_j^c)$, which is $\tilde{\theta}(I_j^c)$, very well. This suggests decomposing $SSE\left(\hat{\theta}^{(\lambda, I_j^c)}, \theta^*\right)$ into a sum of two errors; one is $SSE\left(\hat{\theta}^{(\lambda, I_j^c)}, \tilde{\theta}(I_j^c)\right)$ and the other is $SSE\left(\tilde{\theta}(I_j^c),\theta^*\right)$. The first term is simply the usual mean squared error of the usual Trend Filtering estimator and the second is a non-stochastic approximation error term. %capturing the squared distance of the true signal $\theta^*$ and the the mean of $\tilde{y}(I_j^c)$, which is $\tilde{\theta}(I_j^c)$. 
	This decomposition is enabled for us once we write down the KKT conditions for the Trend Filtering convex opptimization objective (see Lemma~\ref{lem:basic}). This yields

	%However, we want to estimate $\theta^*$. This suggests that we should be able to write 

	%The KKT conditions for the trend filtering convex optimization objective yields (see Lemma\ref{lem:basic})

	%We can write 
	\begin{align*}
	&SSE\left(\hat{\theta}^{(\lambda, I_j^c)}, \theta^*\right)\\
	&\leq \underbrace{2\left<\tilde{\epsilon}(I_j^c), \hat{\theta}^{(\lambda, I_j^c)} - \theta^*\right> + 2\lambda n^{r-1}\left(\|D^{(r)}\theta^*\|_1 - \|D^{(r)}\hat{\theta}^{(\lambda, I_j^c)}\|_1\right)}_{T_1} + \underbrace{\left\|\theta^* -  \tilde{\theta}(I_j^c)\right\|^2}_{T_2}.
	\end{align*}
	Here the term $T_1$ can be thought of as a further upper bound to $SSE\left(\hat{\theta}^{(\lambda, I_j^c)}, \tilde{\theta}(I_j^c)\right).$
	We now bound the two terms $T_1$ and $T_2$ separately.
	
	\item \textbf{Bounding $T_1$ by suitably modifying arguments in~\cite{ortelli2019prediction}}:

	Note that $\tilde{y}(I_j^c) \in \R^n$ is linearly extrapolated from the vector $y_{I_j^c}.$ Therefore, there exist matrices $H^{j} \in \R^{n \times n-|I_j|}$ such that $$\tilde{y}(I_j^c) = H^j y_{I_j^c} \:\:\:\:j \in [K].$$ The exact form of these matrices $H^j $ is explicitly described later in the detailed proof. 
	%Now, the term $T_1$ is simply the MSE of a trend filtering estimator with tuning $\lambda$, where the data vector is $\tilde{y}(I_j^c) = \tilde{\theta}(I_j^c) + \tilde{\epsilon}(I_j^c)$. 
	Now, a similar term to $T_1$ appears in the analysis of SSE of the trend filtering estimator
	%This suggests that to bound $T_1$ we should be able to use existing risk bounds for penalized trend filtering such as those 
	in~\cite{ortelli2019prediction}. So we could apply the entire proof machinery developed 
	in~\cite{ortelli2019prediction} to bound $T_1.$ The only potential problem is that while ~\cite{ortelli2019prediction} considered the case where the error vector is i.i.d gaussian with variance $\sigma^2$, the error vector here is $\tilde{\epsilon}(I_j^c) = H^j \epsilon_{I_j^c}.$ Thus, in our setting $\tilde{\epsilon}(I_j^c)$ is a linear transformation of a vector of i.i.d subgaussian $\sigma$ random variables for a general subgaussian distribution. In particular, in our setting the error variables are no longer even independent.

	By using appropriate concentration inequalities which hold for general subgaussian errors and the fact that the operator norm of $H^j$ stays bounded away from $\infty$; (see Lemma~\ref{ABopnorm}), we show that the entire proof machinery of~\cite{ortelli2019prediction} to bound $T_1$ can be adapted to our setting as well. We state and prove Propositions~\ref{prop:slowratetv} and~\ref{prop:fastratetv} to bound $T_1$ which can be thought of as extensions (to general subgaussian error distribution and dependent errors) of the two bounds in Theorem $2.2$ in~\cite{ortelli2019prediction} proved for the standard i.i.d Gaussian errors setting. This means that up to constant factors, our slow rate and fast rate bounds on $T_1$ match the corresponding bounds given in~\cite{ortelli2019prediction} under their recommended choice of the tuning parameter, which $\lambda^*$ is taken to be.
	
	\item \textbf{Bounding the Approximation Error Term $T_2$}:

	To bound the deterministic quantity $T_2$, we show that 
	\begin{enumerate}
		\item $\left\|\theta^* -  \tilde{\theta}(I_j^c)\right\|^2 \leq \left|D^{(r)}\theta^*\right|_{\infty}\left\|D^{(r)}\theta^*\right\|_{1}.$
		\item $\left\|\theta^* -  \tilde{\theta}(I_j^c)\right\|^2 \leq \left|D^{(r)}\theta^*\right|_{\infty}^2\left\|D^{(r)}\theta^*\right\|_{0}.$
	\end{enumerate}
	The first bound is used for the slow rate and the second bound is used for the fast rate. These bounds work well for our purposes because they turn out to be of lower order than the bounds for $T_1$ under realistic choices of the true signal $\theta^*.$ Therefore, obtaining this particular form of the upper bounds is important for us. This is ensured by our $r^{th}$ order polynomial interpolation scheme (step $3$ of our description of the CVTF estimator). Other simpler interpolation schemes (like simple two neighbor averaging for all orders $r \geq 1$ as used in the Rpackage~\cite{arnold2020package}; see Section~\ref{sec:compare} for more on this) may potentially be used but with our particular interpolation scheme, it becomes possible to prove inequalities (see Lemma~\ref{lem:firstterm}) like in the above display in a fairly simple fashion. 
\end{enumerate}
%{\color{red} Change $r$th to $r^{th}$.}

\subsection{Detailed Proofs for Trend Filtering}\label{sec:tfproofsdet}
First, we will state two propositions which go a long way towards proving Theorem~\ref{thm:slowratetv} and Theorem~\ref{thm:fastratetv}. After stating these two propositions, we finish the proofs of the two theorems assuming these propositions. In these proofs, we will again drop the superscript $r$ and denote the completion estimators $\hat{\theta}^{(\lambda,I,r)}$ by simply $\hat{\theta}^{(\lambda,I)}$ where $I = I_j$ for $j \in [K] = [r + 1].$

\subsection{Proposition~\ref{prop:slowratetv} and Proposition~\ref{prop:fastratetv}}

 For $r \geq 1$, let $\cD = [n-r]$. The following two propositions hold for a fixed but arbitrary index set $\cS \subseteq \cD$. We fix $\cS = \{t_1, \dots, t_s\} \subseteq \cD$, where $1 \leq t_1 < \dots < t_s \leq n-r$. Hence, $s$ denotes the cardinality of the subset $\cS.$ Also, let $t_0 := 0$ and $t_{s+1} := n-r+1$. Next, we define $n_i := t_i - t_{i-1}$, $i \in [s+1]$ and $n_{\rm \max} := \max_{i \in [s+1]} n_i$. One could think of $\cS$ as forming a partition of $[n - r]$ into $s + 1$ blocks where the $i$th block has right end point $t_{i}$ and left end point $t_{i - 1}.$ Then $n_i$ is the length of the $i$ th block and $n_{\rm \max}$ is the length of the longest block. Also as in the previous section, since $r \geq 1$ is fixed, we reduce notational clutter by denoting the completion estimators $\hat{\theta}^{(\lambda,I,r)}$ as simply $\hat{\theta}^{(\lambda,I)}$ where $I = I_j$ for $j \in [K] = [r + 1].$

% {\color{red} Figure out what to do with what notation to use to indicate the order $r.$}

\begin{proposition}\label{prop:slowratetv}
	For any $u, v > 0$, if $\lambda$ is chosen such that 
	\begin{equation}\label{eq:lamtv}
	\lambda \geq C_r \frac{\sigma}{n^{r-1}} \left(\frac{n_{\rm max}}{2}\right)^{\frac{2r-1}{2}}\sqrt{2\log(2(n-s-r)) + 2u},
	\end{equation}
	then there exist positive constants $C_1$, $C_2$ such that with probability at least $1 - e^{-C_1u} - e^{-C_2v}$, for any $j \in [K]$,
	\begin{align*}
	SSE\left(\hat{\theta}^{(\lambda, I_j^c)}_{I_j}, \theta^*_{I_j}\right) &\leq 2\left|D^{(r)}\theta^*\right|_{\infty}\left\|D^{(r)}\theta^*\right\|_{1} + 8\lambda n^{r-1} \left\|D^{(r)}\theta^*\right\|_1\\
	&\qquad +4\left(\sqrt{C_r(s+r)} + \sigma\sqrt{v}\right)^2
	\end{align*}
	for some positive constant $C_r$ which depends only on $r$.
\end{proposition}

\begin{proposition}\label{prop:fastratetv}
	Fix $r \in \{1,2,3,4\}.$ Define the sign vector $q^* \in \{-1, 0, +1\}^s$ containing the signs of the elements in $(D^{(r)}\theta^*)_{\cS}$, that is, for every $i \in [s]$,
	$q^*_{i} := sign(D^{(r)}\theta^*)_{t_i},$
	and the index set $$\cS^{\pm} := \{2 \leq i \leq s : q^*_{i}q^*_{i-1} = -1\} \cup \{1, s+1\}.$$
	Suppose $n_i \geq r(r+2)$ for all $i \in \cS^{\pm}$.
	For any $u, v > 0$, if $\lambda$ is chosen such that 
	$$\lambda \geq C_r \frac{\sigma}{n^{r-1}} \left(\frac{n_{\rm max}}{2}\right)^{\frac{2r-1}{2}}\sqrt{2\log(2(n-s-r)) + 2u},$$ then there exist positive constants $C_1, C_2$ such that with probability at least $1 - e^{-C_1u} - e^{-C_2v}$, for any $j \in [K]$,
	\begin{align*}
	SSE\left(\hat{\theta}^{(\lambda, I_j^c)}_{I_j}, \theta^*_{I_j}\right) &\leq 2\left|D^{(r)}\theta^*\right|_{\infty}^2\left\|D^{(r)}\theta^*\right\|_{0} + 8\lambda n^{r-1} \left\|(D^{(r)}\theta^*)_{\cD \setminus \cS}\right\|_1\\
	&\qquad +4\left(\sqrt{C_r(s+r)} + \sigma \sqrt{v} + \lambda n^{r-1}\Gamma_{\cS}\right)^2,
	\end{align*}
	where
	$$\Gamma_{\cS}^2 := \tilde{C}_r\left\{\sum_{i \in \cS^{\pm}}\frac{1 + \log n_i}{n_i^{2r-1}} + \sum_{i \notin \cS^{\pm}}\frac{1 + \log n_i}{n_{\rm max}^{2r-1}}\right\},$$
	for some positive constants $C_r$ and $\tilde{C}_r$ which depend only on $r$.
\end{proposition}

We now explain the role of these propositions in our proof. Both these propositions hold for any fixed subset $\cS \subset [n -  r].$ Proposition~\ref{prop:slowratetv} will be used to prove the slow rate theorem and Proposition~\ref{prop:fastratetv} will be used to prove the fast rate theorem. These two propositions can be thought of as extensions (to general subgaussian error distribution and dependent errors of the kind we are interested in) of the two bounds in Theorem $2.2$ in~\cite{ortelli2019prediction} proved for the standard i.i.d Gaussian errors setting.

Note that the bounds in these two propositions involve cardinalities of the blocks of $\cS$ and $s = |\cS|$ itself. When we apply these propositions we just need to choose the set $\cS$ appropriately. For the fast rate theorem, the choice of $\cS$ is clear, it is just the set of indices where $D^{(r)}\theta^*$ is non zero. For the slow rate theorem, we choose  $\cS$ to consist of blocks of equal size where the number of blocks is obtained by optimizing a trade off of terms. We will now finish the proofs of Theorem~\ref{thm:slowratetv} and Theorem~\ref{thm:fastratetv} assuming these two propositions.

%Basically, we will choose 

\subsection{Proof of Theorem \ref{thm:slowratetv}}
In view of Theorem \ref{thm:main}, it is enough to bound the quantities: $\min_{\lambda \in \Lambda} SSE(\hat{\theta}^{(\lambda)}, \theta^*)$ and $\min_{\lambda \in \Lambda} SSE(\hat{\theta}^{(\lambda,I)}, \theta^*)$, where $I = I_j$ for $j \in [K].$ The desired bounds (both slow rate and fast rate) for $\min_{\lambda \in \Lambda} SSE(\hat{\theta}^{(\lambda)}, \theta^*)$ follow from Theorem $1.1$ in~\cite{ortelli2019prediction}. The final bounds are going to be exactly similar to the bounds we obtain here for $SSE(\hat{\theta}^{(\lambda)}, \theta^*)$. Infact, Theorems~\ref{thm:slowratetv} and~\ref{thm:fastratetv} generalize Theorem $1.1$ in~\cite{ortelli2019prediction} by giving bounds for a completion version of Trend Filtering (and for general subgaussian errors). By setting $I_j^{c}$ equal to the entire set $[n]$ and following our proof, it can be checked that Theorem $1.1$ in~\cite{ortelli2019prediction} is a special case of Theorems~\ref{thm:slowratetv} and~\ref{thm:fastratetv} combined. %Thus, Theorems~\ref{thm:slowratetv} and~\ref{thm:fastratetv}~\ref{prop:dcevents} itself also implies that if $\lambda \geq C \sigma^2 \log N$ then for any $\alpha > 1$, there exists an absolute constant $C$ such that with probability at least $1 - C N^{-C \alpha}$, we have
%$$SSE(\hat{\theta}^{(\lambda)}, \theta^*) \leq R(\theta^*,\lambda) + C \alpha \sigma^2 \log N.$$

%Our goal is to bound $\|\hat{\theta}^{(\lambda,I)} - \theta^*\|.$ 
Our main goal therefore is to bound $SSE(\hat{\theta}^{(\lambda,I)}, \theta^*)$, where $I = I_j$ for $j \in [K]$, under the optimal choices of $\lambda$, which are different for slow rate and fast rate. However, since we assume that $2^{N^*}$ contains both these optimal choices, this will essentially give us bounds on $\min_{\lambda \in \Lambda} SSE(\hat{\theta}^{(\lambda,I)}, \theta^*)$.

In particular, for the slow rate, it suffices to prove that, 
for any $u, v > 0$, if $\lambda$ is chosen such that 
\begin{equation}\label{eq:lamtvthm1}
\lambda \geq C_r \sigma \left(n^{1/(2r+1)}(\log n + u)^{1/(2r+1)}\right),
\end{equation}
and $n^{r-1} \left\|D^{(r)}\theta^*\right\|_1 = V^*$,
then there exists positive constants $C_1$, $C_2$ such that with probability at least $1 - e^{-C_1u} - e^{-C_2v}$, for any $j \in [K]$,
\begin{align*}
MSE\left(\hat{\theta}^{(\lambda, I_j^c)}_{I_j}, \theta^*_{I_j}\right) \leq \frac{2C_r}{n}\left|D^{(r)}\theta^*\right|_{\infty}\left\|D^{(r)}\theta^*\right\|_{1} + C_r \sigma^2 \left(n^{-\frac{2r}{2r+1}} (V^*(\log n + u))^{\frac{1}{2r+1}} + \frac{v}{n}\right),
\end{align*}
where $C_r$ is a constant only depending on $r$.
In order to prove the above result, we will use Proposition \ref{prop:slowratetv}, which we are going to prove later. Consider $\cS$ in Proposition \ref{prop:slowratetv} such that each block is of equal size, which implies $n_{\rm max} \approx \frac{n}{s}$. Then, by choosing
$$\lambda \geq C_r \sigma \left(\frac{1}{n^{r-1}} \left(\frac{n}{s}\right)^{\frac{2r-1}{2}} \sqrt{\log n + u}\right),$$
there exists positive constants $C_1$, $C_2$ such that with probability at least $1 - e^{-C_1u} - e^{-C_2v}$, for any $j \in [K]$, 
\begin{align*}
SSE\left(\hat{\theta}^{(\lambda, I_j^c)}_{I_j}, \theta^*_{I_j}\right) &\leq 2\left|D^{(r)}\theta^*\right|_{\infty}\left\|D^{(r)}\theta^*\right\|_{1} + 8\lambda V^*
 + 4\left(\sqrt{C_r(s+r)} + \sigma\sqrt{v}\right)^2\\
&\leq 2\left|D^{(r)}\theta^*\right|_{\infty}\left\|D^{(r)}\theta^*\right\|_{1} + 8\lambda V^* + 8C_r(s + r) + 8\sigma^2 v,
\end{align*}
where the right hand side is a function of $s$. In order to optimize this with respect to $s$, we equate
$$\frac{1}{n^{r-1}} \left(\frac{n}{s}\right)^{\frac{2r-1}{2}} (\sqrt{\log n + u}) V^* = s,$$
which leads to $s = n^{\frac{1}{2r+1}} (V^*(\log n + u))^{\frac{1}{2r+1}}$. This completes the proof.

%In order to present the fast rate result, we introduce some additional notations as follows.

\subsection{Proof of Theorem \ref{thm:fastratetv}}
For the fast rate, it suffices to prove that, 
for any $u, v > 0$, if $\lambda$ is chosen such that 
\begin{equation}\label{eq:lamtvthm}
\lambda \geq C_r \sigma s^{-(2r-1)/2} \sqrt{n(\log n + u)},
\end{equation}
then there exists positive constants $C_1$, $C_2$ such that with probability at least $1 - e^{-C_1u} - e^{-C_2v}$, for any $j \in [K]$,
\begin{align*}
MSE\left(\hat{\theta}^{(\lambda, I_j^c)}_{I_j}, \theta^*_{I_j}\right) &\leq \frac{2C_r}{n}\left|D^{(r)}\theta^*\right|_{\infty}^2\left\|D^{(r)}\theta^*\right\|_{0} + 8\lambda n^{r-2} \left\|(D^{(r)}\theta^*)_{\cD \setminus \cS}\right\|_1\\ 
&\quad \quad \quad + C_r \sigma^2 \left(\frac{s}{n} \log n (\log n+u) + \frac{v}{n}\right),
\end{align*}
where $C_r$ is a constant only depending on $r$.
In order to prove the above result, we will use Proposition \ref{prop:fastratetv}, which we are going to prove later. Apply Proposition \ref{prop:fastratetv} for the particular $\cS = \{  j \,:\,  (D^{(r) } \theta^*)_j  \neq 0   \}$ which satisfies the length assumption. Then, by choosing
$$\lambda \geq C_r \sigma \left(\frac{1}{n^{r-1}} \left(\frac{n}{s}\right)^{\frac{2r-1}{2}} \sqrt{\log n + u}\right) =  C_r \sigma \left(s^{-(2r-1)/2} \sqrt{n(\log n + u)}\right),$$
there exists positive constants $C_1$, $C_2$ such that with probability at least $1 - e^{-C_1u} - e^{-C_2v}$, for any $j \in [K]$, 
\begin{align*}
SSE\left(\hat{\theta}^{(\lambda, I_j^c)}_{I_j}, \theta^*_{I_j}\right) &\leq 2\left|D^{(r)}\theta^*\right|_{\infty}^2\left\|D^{(r)}\theta^*\right\|_{0} + 8\lambda n^{r-1} \left\|(D^{(r)}\theta^*)_{\cS^c}\right\|_1 + C_r \sigma^2 \left(v + \lambda^2n^{2r - 2}\Gamma_{\cS}^2\right).
\end{align*}
In the third term above, we use the fact that $\frac{n}{2s} \leq n_{\rm max} \leq c n_i$ for all $i \in \cS^{\pm}$ to bound $\Gamma_{\cS}^2$ as follows
\begin{align*}
\Gamma_{\cS}^2 &= \tilde{C}_r\left\{\sum_{i \in \cS^{\pm}}\frac{1 + \log n_i}{n_i^{2r-1}} + \sum_{i \notin \cS^{\pm}}\frac{1 + \log n_i}{n_{\rm max}^{2r-1}}\right\}\\
&\leq C_r\left\{\sum_{i \in \cS^{\pm}}\frac{\log n}{(n/s)^{2r-1}} + \sum_{i \notin \cS^{\pm}}\frac{\log n}{(n/s)^{2r-1}}\right\}
= C_r \frac{s^{2r}}{n^{2r-1}} \log n, 
\end{align*}
and this completes the proof.

It now remains to give the proofs of Proposition \ref{prop:slowratetv} and Proposition \ref{prop:fastratetv} which we provide in the next section.

\subsection{Proofs of Proposition \ref{prop:slowratetv} and Proposition \ref{prop:fastratetv}}

Fix any $r \geq 1$, $K = r + 1$, $j \in [K]$, and $\lambda > 0$. 
For notational simplicity we drop $I_j^c$ from the following notations: $\tilde{\theta} \equiv \tilde{\theta}(I_j^c)$, $\tilde{y} \equiv \tilde{y}(I_j^c)$ and $\tilde{\epsilon} \equiv \tilde{\epsilon}(I_j^c)$.

\subsubsection{Basic Inequality}
To begin with, we prove the following basic inequality.
\begin{lemma}\label{lem:basic}
	For all $\theta \in \R^{n}$,
	%\begin{align*}
	%\left\|\hat{\theta}^{(\lambda, I_j^c)} - \tilde{\theta}(I_j^c)\right\|^2 + \left\|\hat{\theta}^{(\lambda, I_j^c)} - \theta\right\|^2 &\leq \left\|\theta -  \tilde{\theta}(I_j^c)\right\|^2\\ 
	%&+ 2\left<\tilde{\epsilon}(I_j^c), \hat{\theta}^{(\lambda, I_j^c)} - \theta\right> + 2n^{r-1}\lambda\left(\|D^{(r)}\theta\|_1 - \|D^{(r)}\hat{\theta}^{(\lambda, I_j^c)}\|%_1\right).
	%\end{align*}
	\begin{align*}
	\left\|\hat{\theta}^{(\lambda, I_j^c)} - \tilde{\theta}\right\|^2 + \left\|\hat{\theta}^{(\lambda, I_j^c)} - \theta\right\|^2 &\leq \underbrace{\left\|\theta -  \tilde{\theta}\right\|^2}_{T_2}\\ 
	&+ \underbrace{2\left<\tilde{\epsilon}, \hat{\theta}^{(\lambda, I_j^c)} - \theta\right>}_{T_{11}} + \underbrace{2\lambda n^{r-1}\left(\|D^{(r)}\theta\|_1 - \|D^{(r)}\hat{\theta}^{(\lambda, I_j^c)}\|_1\right)}_{T_{12}}.
	\end{align*}
\end{lemma}
\begin{proof}
	The above result follows from the KKT conditions, see Lemma A.1 in~\cite{ortelli2019prediction}.
\end{proof}

%\begin{remark}
%In the sketch of the proof, we plugged $\theta = \theta^*$
%\end{remark}

\subsubsection{Upper Bound on $T_{11}$}
Next, we state and prove Lemma \ref{lem:boundonerr} which will help in providing an upper bound of the random variable $T_{11} = \left<\tilde{\epsilon}, \hat{\theta}^{(\lambda, I_j^c)} - \theta\right>$. %In fact, Lemma \ref{lem:boundonerr} will give a bound for $left<\tilde{\epsilon}, \hat{\theta}^{(\lambda, I_j^c)} - \theta\right>$ which holds for all $\theta$ and 

%Before that, we prove the next two lemmas which are going to be useful in proving Lemma \ref{lem:boundonerr}. 

%{\color{red} where are we saying that proof holds for $r \leq 4.$?}

At this point, we introduce some additional notations which we are going to use in the rest of the results. Let $\cN_{-\cS}$ denote the following subspace
$$\cN_{-\cS} := \{\theta \in \R^n : (D^{(r)}\theta)_{\cD \setminus \cS} = 0\}$$
and $\cN_{-\cS}^{\perp}$ is its orthogonal complement. It is not difficult to observe that $dim(\cN_{-\cS}) = r+s$. For any $\theta \in \R^n$, let $\theta_{\cN_{-\cS}}$ and $\theta_{\cN_{-\cS}^{\perp}}$ be its projection onto the subspaces $\cN_{-\cS}$ and $\cN_{-\cS}^{\perp}$ respectively. Furthermore, we denote $\Psi^{-\cS} \in \R^{n \times (n-r-s)}$ to be the matrix for which the following relation holds:
$$\theta_{\cN_{-\cS}^{\perp}} = \Psi^{-\cS}(D^{(r)}\theta)_{\cD \setminus \cS}.$$
For any $k \in [n-r-s]$, the $k^{th}$ column of $\Psi^{-\cS}$ is denoted by $\Psi_{k}^{-\cS} \in \R^n$. Now we are ready to prove Lemma \ref{lem:boundonerr}.

\begin{lemma}\label{lem:boundonerr}
For any $u, v > 0$, suppose $\lambda$ is chosen such that 
\begin{equation}\label{eq:lamtv1}
\lambda \geq \frac{\eta(u)}{n^{r-1}}\max_{k \in [n-r-s]} \left\|\Psi_{k}^{-\cS}\right\|,
\end{equation}
where $\eta(u) := \sigma\left\|{H^j}^T\right\|_{\rm OP}\sqrt{2\log(2(n-r-s)) + 2u}$ and $w \in \R^{n-r-s}$ is defined as 
\begin{equation}\label{eq:w1}
w_k := \frac{\eta(u)}{\lambda n^{r-1}}\left\|\Psi_k^{-\cS}\right\|,\quad k \in [n-r-s],
\end{equation}
then with probability at least $1 - e^{-C_1u} - e^{-C_2v}$,
\begin{align*}
\left<\tilde{\epsilon}, \theta\right> \leq \left(\sqrt{||H^j||_{\rm OP}(r+s)} + \sigma\sqrt{v}\right)||\theta|| + \lambda n^{r-1}  \sum_{k = 1}^{n-r-s} \left|w_k ((D^{(r)}\theta)_{\cD \setminus \cS})_k\right| \:\:\:\forall \theta \in \R^n%\left\|w_{-S} \times \left(D^{(r)}\theta\right)_{-S}\right\|_1,
\end{align*}
for some positive constants $C_1$, $C_2$.
\end{lemma}
\begin{proof}
We can write
$$\left<\tilde{\epsilon}, \theta\right> = T_1 + T_2,$$
where $T_1 = \left<\tilde{\epsilon}, \theta_{\cN_{-S}}\right>$ and $T_2 = \left<\tilde{\epsilon}, \theta_{\cN_{-S}^{\perp}}\right>$. Let 
$$\cV := \left\{\|\tilde{\epsilon}_{\cN_{-S}}\| \leq \sqrt{||H^j||_{\rm OP}(r+s)} + \sigma\sqrt{v}\right\}$$
Now, on $\cV$, we have
\begin{align*}
T_1 = \left<\tilde{\epsilon}, \theta_{\cN_{-S}}\right>
=  \left<\tilde{\epsilon}_{\cN_{-S}}, \theta_{\cN_{-S}}\right>
&\leq \|\tilde{\epsilon}_{\cN_{-S}}\|\|\theta_{\cN_{-S}}\|\\
&\leq \|\tilde{\epsilon}_{\cN_{-S}}\|\|\theta\|
\leq \left(\sqrt{||H^j||_{\rm OP}(r+s)} + \sigma\sqrt{v}\right)||\theta||. 
\end{align*}
Also, by Lemma~\ref{lem:lem1forerr} (a result about norms of projection of a matrix times a subgaussian random vector), $P(\cV) \geq 1 - e^{-C_2v}$ for some positive constant $C_2$. Next, let
$$\cU := \left\{\max_{k \in [n-r-s]}\left|\left<\tilde{\epsilon}, \frac{\Psi^{-S}_k}{\|\Psi^{-S}_k\|}\right>\right| \leq \eta(u)\right\}.$$ %= \left\|{H^j}^T\right\|_{\rm OP}\sqrt{2\log(2|\cD \setminus \cS|) + 2u}
Then, on $\cU$, we have
\begin{align*}
T_2 &= \left<\tilde{\epsilon}, \theta_{\cN_{-S}^{\perp}}\right>\\
&= \left<\tilde{\epsilon}, \Psi^{-\cS}\left(D^{(r)}\theta\right)_{\cD \setminus \cS}\right>\\
&= \left<\tilde{\epsilon}, \sum_{k = 1}^{n-r-s}\Psi^{-\cS}_k((D^{(r)}\theta)_{\cD \setminus \cS})_k\right>\\
&= \sum_{k =1}^{n-r-s}\left<\tilde{\epsilon}, \Psi^{-\cS}_k\right>((D^{(r)}\theta)_{\cD \setminus \cS})_k\\
&=\lambda n^{r-1} \sum_{k =1}^{n-r-s}\left<\tilde{\epsilon}, \frac{\Psi^{-\cS}_k}{\|\Psi^{-S}_k\|\eta(u)}\right>\frac{\eta(u)}{\lambda n^{r-1}}\|\Psi^{-S}_k\|((D^{(r)}\theta)_{\cD \setminus \cS})_k\\
&=\lambda n^{r-1} \sum_{k =1}^{n-r-s}\left<\tilde{\epsilon}, \frac{\Psi^{-\cS}_k}{\|\Psi^{-S}_k\|\eta(u)}\right>w_k((D^{(r)}\theta)_{\cD \setminus \cS})_k\\
&\leq \lambda n^{r-1} \sum_{k =1}^{n-r-s} \left|\left<\tilde{\epsilon}, \frac{\Psi^{-\cS}_k}{\|\Psi^{-S}_k\|\eta(u)}\right>\right|\left|w_k((D^{(r)}\theta)_{\cD \setminus \cS})_k\right|\\
&\leq \lambda n^{r-1} \sum_{k = 1}^{n-r-s} \left|w_k ((D^{(r)}\theta)_{\cD \setminus \cS})_k\right|.
\end{align*}
Also, by Lemma ~\ref{lem:lem2forerr} (a standard result about maxima of subgaussians), $P(\cU) \geq 1 - e^{-C_1u}$ for some positive constant $C_1$. Combining the above results, the proof is complete. Lemmas~\ref{lem:lem1forerr} and~\ref{lem:lem2forerr} which have been used in this proof are stated and proved in Section~\ref{sec:auxtf}.
\end{proof}

\begin{comment}
Since our extrapolation is linear we can express $\tilde{y}(I_j^c)$ as a linear transformation of $y_{I_j^c}$:
$$\tilde{y}(I_j^c) = H^jy_{I_j^c}$$ for some matrix $H^j \in \R^{n \times n-|I_j|}$ which can be defined according to the extrapolation stated before. Moreover, the method of extrapolation is: for any index in $i \in I_j$, we always extrapolate from its right neighborhood as long as it contains $r$ elements, that is, $i + r \leq n$. Otherwise, we extrapolate from its left neighborhood.	For example, let us consider $j = 1$ and we show how $H^1$ looks like in the following.
{\color{red} Reword the above?}
\end{comment}

\subsubsection{Upper Bound on $T_2$ when $\theta = \theta^*$}
Next, we prove a lemma that will provide an upper bound of the quantity $T_2 = \left\|\theta -  \tilde{\theta}\right\|^2$, when $\theta = \theta^*$.
\begin{lemma}\label{lem:firstterm}
We have
\begin{enumerate}
\item $\left\|\theta^* -  \tilde{\theta}\right\|^2 \leq \left|D^{(r)}\theta^*\right|_{\infty}\left\|D^{(r)}\theta^*\right\|_{1}.$
\item $\left\|\theta^* -  \tilde{\theta}\right\|^2 \leq \left|D^{(r)}\theta^*\right|_{\infty}^2\left\|D^{(r)}\theta^*\right\|_{0}.$
\end{enumerate}
\end{lemma}
\begin{proof}
Note that, by the construction of $\tilde{y}$, we have
\begin{align*}
&\left\|\theta^* -  \tilde{\theta}\right\|^2 = \sum_{i \in I_j}\left(\theta^*_i - \tilde{\theta}_i\right)^2\\
&=\sum_{i \in I_j : i+r \leq n}\left(\theta^*_i - \sum_{l = 1}^r(-1)^{l+1}{r \choose l}\theta^*_{i + l}\right)^2 + \sum_{i \in I_j : i+r > n}\left(\theta^*_i - \sum_{l = 1}^r(-1)^{l+1}{r \choose l}\theta^*_{i - l}\right)^2\\
&\leq \left\|D^{(r)}\theta^*\right\|^2.
\end{align*}
Now note that one can further upper bound $\left\|D^{(r)}\theta^*\right\|^2$ in either way:
$$\left\|D^{(r)}\theta^*\right\|^2 \leq \left|D^{(r)}\theta^*\right|_{\infty}\left\|D^{(r)}\theta^*\right\|_{1}, \quad \left\|D^{(r)}\theta^*\right\|^2 \leq \left|D^{(r)}\theta^*\right|_{\infty}^2\left\|D^{(r)}\theta^*\right\|_{0}.$$
\end{proof}

\begin{proof}[Proof of Proposition~\ref{prop:slowratetv}]

Using Lemmas \ref{lem:basic} and \ref{lem:boundonerr}, we have for any $\theta \in \R^n$, if $\lambda$ is chosen according to \eqref{eq:lamtv1} and $w$ is defined according to \eqref{eq:w1} then with probability at least $1 - e^{-C_1u} - e^{-C_2v}$,
\begin{align}\notag
&\left\|\hat{\theta}^{(\lambda, I_j^c)} - \tilde{\theta}\right\|^2 + \left\|\hat{\theta}^{(\lambda, I_j^c)} - \theta\right\|^2\\ \notag
&\leq \left\|\theta -  \tilde{\theta}\right\|^2
+ 2\left<\tilde{\epsilon}, \hat{\theta}^{(\lambda, I_j^c)} - \theta\right> + 2\lambda n^{r-1}\left(\|D^{(r)}\theta\|_1 - \|D^{(r)}\hat{\theta}^{(\lambda, I_j^c)}\|_1\right)\\ \notag
&\leq \left\|\theta -  \tilde{\theta}\right\|^2
+ 2 \left(\sqrt{||H^j||_{\rm OP}(r+s)} + \sigma\sqrt{v}\right)\left\|\hat{\theta}^{(\lambda, I_j^c)} - \theta\right\|\\ \label{eq:slowandfast}
&+ \underbrace{2\lambda n^{r-1} \left\{\sum_{k = 1}^{n-r-s} \left|w_k ((D^{(r)}(\hat{\theta}^{(\lambda, I_j^c)} - \theta))_{\cD \setminus \cS})_k\right| + \left(\|D^{(r)}\theta\|_1 - \|D^{(r)}\hat{\theta}^{(\lambda, I_j^c)}\|_1\right)\right\}}_{T_{3}}.
\end{align}
By definition, $|w|_{\infty} \leq 1$. Thus, replacing $w_k$ by $1$ for all $k \in [n-r-s]$, the third term of \eqref{eq:slowandfast}, namely $T_3$ can be further bounded by
\begin{align*}
&2\lambda n^{r-1}\left\{ \| (D^{(r)}(\hat{\theta}^{(\lambda, I_j^c)} - \theta))_{\cD \setminus \cS}\|_1 + \left(\|D^{(r)}\theta\|_1 - \|D^{(r)}\hat{\theta}^{(\lambda, I_j^c)}\|_1\right)\right\}\\
&\leq 2\lambda n^{r-1}\left\{ \| D^{(r)}(\hat{\theta}^{(\lambda, I_j^c)} - \theta)\|_1 + \left(\|D^{(r)}\theta\|_1 - \|D^{(r)}\hat{\theta}^{(\lambda, I_j^c)}\|_1\right)\right\}\\
&\leq 2\lambda n^{r-1}\left\{\|D^{(r)}\hat{\theta}^{(\lambda, I_j^c)}\|_1 + \|D^{(r)}\theta\|_1 + \|D^{(r)}\theta\|_1 - \|D^{(r)}\hat{\theta}^{(\lambda, I_j^c)}\|_1\right\}\\
&= 4\lambda n^{r-1}\|D^{(r)}\theta\|_1.
\end{align*}
Imposing the above bound in \eqref{eq:slowandfast}, and plugging in $\theta = \theta^*$ we have
\begin{align*}
&SSE\left(\hat{\theta}^{(\lambda, I_j^c)}_{I_j}, \theta^*_{I_j}\right) = \left\|\hat{\theta}^{(\lambda, I_j^c)}_{I_j} - \theta^*_{I_j}\right\|^2 \leq \left\|\hat{\theta}^{(\lambda, I_j^c)} - \theta^*\right\|^2\\
&\leq \left\|\hat{\theta}^{(\lambda, I_j^c)} - \tilde{\theta}\right\|^2 + \left\|\hat{\theta}^{(\lambda, I_j^c)} - \theta^*\right\|^2\\
&\leq \left\|\theta^* -  \tilde{\theta}\right\|^2
+ 2 \left(\sqrt{||H^j||_{\rm OP}(r+s)} + \sigma\sqrt{v}\right) \left\|\hat{\theta}^{(\lambda, I_j^c)} - \theta^*\right\| + 4\lambda n^{r-1}\|D^{(r)}\theta^*\|_1\\
&\leq \left|D^{(r)}\theta^*\right|_{\infty}\left\|D^{(r)}\theta^*\right\|_{1} + 4\lambda n^{r-1} \left\|D^{(r)}\theta^*\right\|_1\\
&\qquad +2\left(\sqrt{C_r(r + s)} + \sigma\sqrt{v}\right)^2 + \frac{1}{2}\left\|\hat{\theta}^{(\lambda, I_j^c)} - \theta^*\right\|^2.
\end{align*} 
The last inequality is obtained by bounding the first term by using the inequality (1) in Lemma~\ref{lem:firstterm}, and the second term by using the inequality $2ab \leq 2a^2 + \frac{b^2}{2}$ and Lemma \ref{ABopnorm} together for some positive constant $C_r$ that depends only on $r$.

All that remains is to check whether our choice of $\lambda$ is valid. We use Lemma 3.1 in ~\cite{ortelli2019prediction}, which gives us
$$\max_{k \in [n-r-s]} \left\|\Psi_{k}^{-\cS}\right\|^2 \leq \left(\frac{n_{\rm max}}{2}\right)^{2r-1}.$$
Plugging the above we obtain an upper bound on the right hand side of \eqref{eq:lamtv1} in Lemma \ref{lem:boundonerr} as follows.
\begin{align*}
\frac{\eta(u)}{n^{r-1}}\max_{k \in [n-r-s]} \left\|\Psi_{k}^{-\cS}\right\| &=  \frac{\sigma}{n^{r-1}}\left\|{H^j}^T\right\|_{\rm OP}\sqrt{2\log(2(n-r-s)) + 2u}\max_{k \in [n-r-s]} \left\|\Psi_{k}^{-\cS}\right\|\\
&\leq C_r \frac{\sigma}{n^{r-1}} \left(\frac{n_{\rm max}}{2}\right)^{\frac{2r-1}{2}} \sqrt{2\log(2(n-s-r)) + 2u},
\end{align*}
where the last inequality follows from Lemma \ref{ABopnorm} for some positive constant $C_r$ that depends only on $r$.
Thus, the requirement on the tuning parameter is met when $\lambda$ is chosen according to \eqref{eq:lamtv}.
Hence, the proof is complete.
\end{proof}

\subsubsection{Another Upper Bound on $T_3$}
In order to prove  Proposition~\ref{prop:fastratetv}, we again bound $T_3$, but this time in a different manner. This is to prove the fast rate theorem and this part uses the proof technique via interpolating vectors pioneered by~\cite{ortelli2019prediction}.

We quote the next three lemmas from~\cite{ortelli2019prediction}.
\begin{lemma}\label{lem:thirdtermfast1}
We have
\begin{align*}
&\sum_{k = 1}^{n-r-s} \left|w_k ((D^{(r)}(\hat{\theta}^{(\lambda, I_j^c)} - \theta))_{\cD \setminus \cS})_k\right| + \left(\|D^{(r)}\theta\|_1 - \|D^{(r)}\hat{\theta}^{(\lambda, I_j^c)}\|_1\right)\\
&\leq 2\left\|(D^{(r)}\theta)_{\cD \setminus \cS}\right\|_1  + \left\|(D^{(r)}\theta)_{\cS}\right\|_1 - \left\|(D^{(r)}\hat{\theta}^{(\lambda, I_j^c)})_{\cS}\right\|_1\\
&\qquad - \sum_{k = 1}^{n-r-s}\left|(1 - w_k) ((D^{(r)}(\hat{\theta}^{(\lambda, I_j^c)} - \theta))_{\cD \setminus \cS})_k\right|.
\end{align*}
\end{lemma}
\begin{proof}
See the proof of Theorem 2.2 in ~\cite{ortelli2019prediction}.
\end{proof}

Let $q \in \{-1, +1\}^s$ be the sign vector containing the signs of the elements in $(D^{(r)}\theta)_{\cS}$, that is,
$$q_k = sign((D^{(r)}\theta)_{\cS})_k, \quad k \in [s],$$
and $z \in \R^{n-r-s}$ be any vector such that $|z|_{\infty} \leq 1$.
Now, following~\cite{ortelli2019prediction} we define the quantity called noisy effective sparsity
$$\Gamma(q, z) := \max\left\{{q}^T(D^{(r)}\theta')_{\cS} - \sum_{k = 1}^{n-r-s}\left|(1 - z_k) ((D^{(r)}\theta')_{\cD \setminus \cS})_k\right| : \|\theta'\| = 1\right\}.$$
Next lemma will help us to improve the upper bound in Lemma~\ref{lem:thirdtermfast1} in terms of the noisy effective sparsity.
\begin{lemma}\label{lem:thirdtermfast2}
\begin{align*}
&\left\|(D^{(r)}\theta)_{\cS}\right\|_1 - \left\|(D^{(r)}\hat{\theta}^{(\lambda, I_j^c)})_{\cS}\right\|_1 - \sum_{k = 1}^{n-r-s}\left|(1 - w_k) ((D^{(r)}(\hat{\theta}^{(\lambda, I_j^c)} - \theta))_{\cD \setminus \cS})_k\right|\\ 
&\leq \Gamma(q, w)\left\|\hat{\theta}^{(\lambda, I_j^c)} - \theta\right\|
\end{align*}
\end{lemma}
\begin{proof}
See the proof of Lemma A.3 in ~\cite{ortelli2019prediction}.
\end{proof}

\begin{proof}[Proof of Proposition~\ref{prop:fastratetv}]
We start from \eqref{eq:slowandfast} and plug in $\theta = \theta^*$. Now, using Lemmas ~\ref{lem:thirdtermfast1} and ~\ref{lem:thirdtermfast2}, we have
\begin{align*}
&SSE\left(\hat{\theta}^{(\lambda, I_j^c)}_{I_j}, \theta^*_{I_j}\right) = \left\|\hat{\theta}^{(\lambda, I_j^c)}_{I_j} - \theta^*_{I_j}\right\|^2 \leq \left\|\hat{\theta}^{(\lambda, I_j^c)} - \theta^*\right\|^2\\
&\leq \left\|\hat{\theta}^{(\lambda, I_j^c)} - \tilde{\theta}\right\|^2 + \left\|\hat{\theta}^{(\lambda, I_j^c)} - \theta^*\right\|^2\\
&\leq \left\|\theta^* -  \tilde{\theta}\right\|^2
+ 2 \left(\sqrt{||H^j||_{\rm OP}(r+s)} + \sigma\sqrt{v}\right)\left\|\hat{\theta}^{(\lambda, I_j^c)} - \theta^*\right\|\\
&\qquad + 2\lambda n^{r-1} \left\{\sum_{k = 1}^{n-r-s} \left|w_k ((D^{(r)}(\hat{\theta}^{(\lambda, I_j^c)} - \theta^*))_{\cD \setminus \cS})_k\right| + \left(\|D^{(r)}\theta^*\|_1 - \|D^{(r)}\hat{\theta}^{(\lambda, I_j^c)}\|_1\right)\right\}\\
&\leq \left|D^{(r)}\theta^*\right|_{\infty}^2\left\|D^{(r)}\theta^*\right\|_{0} + 2 \left(\sqrt{2c'_r(r+s)} + \sigma\sqrt{v}\right)\left\|\hat{\theta}^{(\lambda, I_j^c)} - \theta^*\right\|\\
&\qquad +2\lambda n^{r-1}\left\{2\left\|(D^{(r)}\theta^*)_{\cD \setminus \cS}\right\|_1 + \Gamma(q^*, w)\left\|\hat{\theta}^{(\lambda, I_j^c)} - \theta^*\right\|\right\}\\
&= \left|D^{(r)}\theta^*\right|_{\infty}^2\left\|D^{(r)}\theta^*\right\|_{0} + 4\lambda n^{r-1}\left\|(D^{(r)}\theta^*)_{\cD \setminus \cS}\right\|_1\\
&\qquad + 2\left(\sqrt{2c'_r(r+s)} + \sigma\sqrt{v} + \lambda n^{r-1} \Gamma(q^*, w)\right)\left\|\hat{\theta}^{(\lambda, I_j^c)} - \theta^*\right\|,
\end{align*}
where in the third inequality, we bound the first term using the inequality (2) in Lemma \ref{lem:firstterm}.
Now, %{\color{red} since $n_i \geq r(r+2)$ for all $i \in \cS^{\pm}$}, 
we can bound the term $\Gamma(q^*, w)$ by using the proof technique via \textit{interpolating vector} (illustrated in~\cite{ortelli2019prediction}, see Definition $2.3$ and Lemma $2.4$ there). Specifically, we use the result which appears as the first display in page $17$ right before equation $10$ in ~\cite{ortelli2019prediction}, to obtain
$$\Gamma(q^*, w) \leq \Gamma_{\cS} = \tilde{C}_r\left\{\sum_{i \in \cS^{\pm}}\frac{1 + \log n_i}{n_i^{2r+1}} + \sum_{i \notin \cS^{\pm}}\frac{1 + \log n_i}{n_{\rm max}^{2r+1}}\right\},$$
for some constant $\tilde{C}_r$ as long as (see Eq. (10), (11) in ~\cite{ortelli2019prediction}) the requirement for tuning parameter $\lambda$ is strengthened to: for an appropriate constant $c''_r \geq 1$,
$$\lambda \geq c''_r \frac{\sigma}{n^{r-1}}\left(\frac{n_{\rm max}}{2}\right)^{\frac{2r-1}{2}} \left\|{H^j}^T\right\|_{\rm OP} \sqrt{2\log(2(n-r-s)) + 2u}.$$
Finally, let $C_r = 2c_r'c''_r$ to finish the proof.
\end{proof}

\subsection{Auxiliary Lemmas}\label{sec:auxtf}

\begin{lemma}\label{lem:lem1forerr}For any subspace $S \subseteq \R^n$, let $\tilde{\epsilon}_S$ be the projection of $\tilde{\epsilon}$ on $S$. For any $u > 0$, with probability at least $1 - 2e^{-Cu/2}$,
	$$\left\|\tilde{\epsilon}_S\right\| \leq \sqrt{||H^j||_{\rm OP}}\sqrt{Dim(S)} + \sigma\sqrt{u},$$
	where $||\cdot||_{\rm OP}$ denotes the operator norm of a matrix and $C$ is some positive constant.
\end{lemma}
\begin{proof}
	Let $P_S$ denote the projection matrix of subspace $S$. Then $$\tilde{\epsilon}_S = P_S\tilde{\epsilon} = P_SH^j\epsilon_{I_j^c}.$$
	From the Hanson-Wright concentration inequality, see Theorem 6.2.1 in \cite{vershynin2018high}, we obtain that, for any matrix $A_{m \times n}$ and $u > 0$, if $Z = (Z_1, \dots, Z_n) \overset{iid}{\sim} Subg(0, \sigma^2)$ then 
	$$P\left(\left|||AZ|| - ||A||_{\rm HS}\right| > u\right) \leq 2\exp\left( - C\frac{u^2}{2\sigma^2}\right),$$
	for some $C > 0$. Here, the notation $\|A\|_{\rm HS}$ refers to the Hilbert-Schmidt norm of the matrix $A$. 
	Letting $A = P_SH^j$, we have
	\begin{align*}
	\left\|A\right\|_{\rm HS} &= \left\|P_SH^j\right\|_{\rm HS}\\
	&=\sqrt{trace((P_SH^j)^T(P_SH^j))}\\
	&=\sqrt{trace({H^j}^TP_SH^j)}
	\end{align*}
	Let $r_0 = rank({H^j}^TP_SH^j)$ and the non-zero eigenvalues are $\{\lambda_i, i \in [r_0]\}$. Also, for any matrix $A$, let $\lambda_{\rm max}(A)$ denote its largest eigenvalue. Therefore,
	\begin{align*}
	trace({H^j}^TP_SH^j) &= \sum_{i = 1}^{r_0} \lambda_i\\
	&\leq \lambda_{\rm max}({H^j}^TP_SH^j) rank({H^j}^TP_SH^j)\\
	&\leq \sup_{x}\frac{x^T{H^j}^TP_SH^jx}{||x||} rank(P_S)\\
	&=\sup_{x}\frac{x^T{H^j}^TP_SH^jx}{||H^jx||}\frac{||H^jx||}{||x||} \times rank(P_S)\\
	&\leq \lambda_{\rm max}(P_S)||H^j||_{\rm OP}\times rank(P_S)\\
	&=||H^j||_{\rm OP}\times Dim(S).
	\end{align*}
	This completes the proof.
\end{proof}

\begin{lemma}\label{lem:lem2forerr}
	Consider a set of $m$ nonzero vectors $\{v_i : i \in [m]\}$. Then for any $u > 0$, with probability at least $1 - e^{-Cu}$,
	$$\max_{k \in [m]}\left|\left<\tilde{\epsilon}, \frac{v_k}{||v_k||}\right>\right| \leq \sigma\left\|{H^j}^T\right\|_{\rm OP}\sqrt{2\log(2m) + 2u},$$
	where $C$ is some positive constant.
\end{lemma}
\begin{proof}
	Fix any $k \in [m]$. Since $\tilde{\epsilon} = H^j\epsilon_{I_j^c}$, we can write
	\begin{align*}
	\left|\left<\tilde{\epsilon}, \frac{v_k}{||v_k||}\right>\right| &= \left|\left<H^j\epsilon_{I_j^c}, \frac{v_k}{||v_k||}\right>\right|\\
	&=\left|\epsilon_{I_j^c}^T{H^j}^T\frac{v_k}{||v_k||}\right|\\
	&=\left|\epsilon_{I_j^c}^T\frac{{H^j}^Tv_k}{\left\|{H^j}^Tv_k\right\|}\right|\frac{\left\|{H^j}^Tv_k\right\|}{\|v_k\|}\\
	&\leq \left|\left<\epsilon_{I_j^c}, \frac{{H^j}^Tv_k}{\left\|{H^j}^Tv_k\right\|}\right>\right|\left\|{H^j}^T\right\|_{\rm OP}
	\end{align*}
	Therefore,
	\begin{align*}
	&P\left(\max_{k \in [m]}\left|\left<\tilde{\epsilon}, \frac{v_k}{||v_k||}\right>\right| \leq \sigma\left\|{H^j}^T\right\|_{\rm OP}\sqrt{2\log(2m) + 2u}\right)\\ 
	%&\geq P\left(\max_{k \in [m]}\left|\left<\epsilon_{I_j^c}, \frac{{H^j}^Tv_k}{\left\|{H^j}^Tv_k\right\|}\right>\right|\left\|{H^j}^T\right\|_{\rm OP} \leq \eta(u)\right)\\
	&\geq P\left(\max_{k \in [m]}\left|\left<\epsilon_{I_j^c}, \frac{{H^j}^Tv_k}{\left\|{H^j}^Tv_k\right\|}\right>\right| \leq \sigma\sqrt{2\log(2m) + 2u}\right)\\
	&\geq 1 - e^{-Cu}.
	\end{align*}
	The last result follows from the standard maximal inequality of $m$ subgaussian random variables.
\end{proof}

\subsubsection{Properties of $H^j$}	Since Lemma \ref{lem:boundonerr} involves both $\left\|H^j\right\|_{\rm OP}$ and $\left\|{H^j}^T\right\|_{\rm OP}$, first we are going to describe how the matrix $H^j$ looks like and then provide upper bounds of those quantities. Note that, the method of extrapolation is: for any index $i \in I_j$, we always extrapolate from its right neighborhood as long as it contains $r$ elements, that is, $i + r \leq n$. Otherwise, we extrapolate from its left neighborhood.	As an example, let us consider $j = 1$ and we show how $H^1$ looks like in the following.
\begin{enumerate}
	\item{Case 1: $n = Kn_0 = (r+1)n_0$.} In this case we never have to extrapolate from the left. Thus, $H^1$ is an $(r+1)n_0 \times rn_0$ block matrix:
	$$H^1 = 
	\begin{bmatrix}
	A & 0 &\dots &0\\
	0 & A &\dots &0\\
	&\vdots\\
	0 & 0 &\dots &A
	\end{bmatrix},$$
	where $A$ is the $(r+1) \times r$ matrix: 
	$$A := 
	\begin{bmatrix}
	{r \choose 1} & -{r \choose 2} &+{r \choose 3} \dots &(-1)^{r+1}\\
	1 &0 &0 \dots &0\\
	0 &1 &0 \dots &0\\
	0 &0 &1 \dots &0\\
	&\vdots\\
	0 &0 &0 \dots &1\\
	\end{bmatrix},$$
	\item{Case 2: $n = Kn_0 + d_*$,} for some $1 \leq d_* \leq K-1$. In this case, we extrapolate from the right only when $i = Kn_0 + 1$. Therefore, $|I_1| = n_0 + 1$ and $H^1$ is an $n \times ((K - 1)n_0 + (d_* - 1))$ block matrix:
	$$H^1 = \begin{bmatrix}
	A & 0 &\dots &0 &0 &0\\
	0 & A &\dots &0 &0 &0\\
	&\vdots\\
	0 & 0 &\dots &A &0 &0\\
	0 & 0 &\dots &0 &B &0\\
	0 & 0 &\dots &0 &0 &I_{d_* - 1}\\
	\end{bmatrix},$$
	where $B$ is the $(r+2) \times r$ matrix:
	$$B := 
	\begin{bmatrix}
	{r \choose 1} & -{r \choose 2} &+{r \choose 3} \dots &(-1)^{r+1}\\
	1 &0 &0 \dots &0\\
	0 &1 &0 \dots &0\\
	0 &0 &1 \dots &0\\
	&\vdots\\
	0 &0 &0 \dots &1\\
	(-1)^{r+1} \dots &+{r \choose 3} & -{r \choose 2} &{r \choose 1}
	\end{bmatrix},$$
\end{enumerate}

In general, it is not difficult to observe that for any $j \in [K]$, the matrix $H^j$ is a block matrix of the form:
$$H^j = 
\begin{bmatrix}
H_1^j & 0 &\dots &0\\
0 & H_2^j &\dots &0\\
&\vdots\\
0 & 0 &\dots &H_m^j
\end{bmatrix},$$ for some $m \geq 1$ and $H_i^j \in \{A, B, I_d : 1 \leq d \leq K\}$, $i \in [m]$. 

\begin{comment}
Let $$\tilde{\theta}(I_j^c) := E\left[\tilde{y}(I_j^c)\right],\ \text{and}\ \ \tilde{\epsilon}(I_j^c) := \tilde{y}(I_j^c) - \tilde{\theta}(I_j^c).$$
Note that, in step (4) of the proposed estimator, we basically implement the traditional trend filtering on the data vector $\tilde{y}(I_j^c)$ and apply the results from ~\cite{ortelli2019prediction}. However, the current setup is different than that in ~\cite{ortelli2019prediction} in three aspects: (i) the true mean is not $\theta^*$ anymore, rather it is $\tilde{\theta}(I_j^c)$, (ii) the entries of the error vector $\tilde{\epsilon}(I_j^c)$ are dependent, and (iii) the original error vector $\epsilon$ is subgaussian. As a result, Proposition \ref{prop:slowratetv} and Proposition \ref{prop:fastratetv} are restatements of Theorem 2.2 in ~\cite{ortelli2019prediction} with some necessary adjustments. Nevertheless, we provide the proof in detail for sake of completeness.
\end{comment}

Next, in Lemma \ref{ABopnorm} we provide upper bounds of the quantities $||H^j||_{\rm OP}$ and $||{H^j}^{T}||_{\rm OP}$. Before that, we prove the next lemma which is going to be useful in proving Lemma \ref{ABopnorm}. 
\begin{lemma}\label{lem:opnormofblock}
	If a block matrix $H \in \R^{p \times q}$ is of the form
	$$H = 
	\begin{bmatrix}
	H_1 & 0 &\dots &0\\
	0 & H_2 &\dots &0\\
	&\vdots\\
	0 & 0 &\dots &H_m
	\end{bmatrix},$$ where $H_i \in \R^{p_i \times q_i}, i \in [m]$, $\sum_i p_i = p$ and $\sum_{i} q_i = q$, then
	$$\|H\|_{\rm OP} \leq \max_{i \in [m]} \|H_i\|_{\rm OP}.$$
\end{lemma}
\begin{proof}
	For any $x = (x_1^T, x_2^T \dots x_m^T)$ such that $x_i \in \R^{q_i}$ for all $i \in [m]$, we have
	$$Hx = 
	\begin{bmatrix}
	H_1x_1\\
	H_2x_2\\
	\vdots\\
	H_mx_m
	\end{bmatrix},$$
	implying 
	\begin{align*}
	\|Hx\|^2 = \sum_{i = 1}^m\|H_ix_i\|^2 \leq \sum_{i = 1}^m \|H_i\|_{\rm OP}^2\|x_i\|^2 \leq \max_{i \in [m]} \|H_i\|_{\rm OP}^2 \sum_{i = 1}^m\|x_i\|^2 = \max_{i \in [m]} \|H_i\|_{\rm OP}^2 \|x\|^2.
	\end{align*}
\end{proof}
Recall that,
for any $j \in [K]$, the matrix $H^j$ is a block matrix of the form of $H$ in Lemma ~\ref{lem:opnormofblock} for some $m \geq 1$ with $H_i \in \{A, B, I_d : 1 \leq d \leq K\}$, $i \in [m]$. 
Thus, in order to have an upper bound of $\|H^j\|_{\rm OP}$ and $\|{H^j}^T\|_{\rm OP}$ it suffices to have upper bounds of $\|A\|_{\rm OP}$, $\|A^T\|_{\rm OP}$, $\|B\|_{\rm OP}$ and $\|B^T\|_{\rm OP}$.
\begin{lemma}\label{ABopnorm}
	The following results hold.
	\begin{enumerate}
		\item $\|A\|_{\rm OP} \leq \left\{1 + {r \choose 1}^2 + {r \choose 2}^2 + \dots + {r \choose r}^2\right\} =: c_r$,
		\item $\|A^T\|_{\rm OP} \leq \left\{r + {r \choose 1}^2 + {r \choose 2}^2 + \dots + {r \choose r}^2\right\} =: c'_r$,
		\item $\|B\|_{\rm OP} \leq 2c_r$,
		\item $\|B^T\|_{\rm OP} \leq 2c'_r$.
	\end{enumerate}
	Thus, $\|H^j\|_{\rm OP}, \|{H^j}^T\|_{\rm OP} \leq 2c'_r$.
\end{lemma}
\begin{proof}
	We will only prove the first result since the rest follow in a similar way. Fix $x \in \R^r$ such that $||x|| = 1$. Then using Cauchy Schwarz inequality we have,
	\begin{align*}
	||Ax||^2 &= \left\{{r \choose 1}x_1 - {r \choose 2}x_2 + \dots + (-1)^{r+1}x_r\right\}^2 + x_2^2 + \dots + x_r^2\\
	&\leq \left\{{r \choose 1}^2 + {r \choose 2}^2 + \dots + 1\right\}(x_1^2 + \dots + x_r^2) + (x_1^2 + x_2^2 + \dots + x_r^2)\\
	&\leq \left\{1 + {r \choose 1}^2 + {r \choose 2}^2 + \dots + {r \choose r}^2\right\}||x||^2 = c_r. 
	\end{align*}
	Using Lemma ~\ref{lem:opnormofblock}, we have
	$$\|H^j\|_{\rm OP} \leq \max\{\|A\|_{\rm OP}, \|B\|_{\rm OP}, 1\} \leq 2c_r, \; \text{and} \; \|{H^j}^T\|_{\rm OP} \leq \max\{\|A^T\|_{\rm OP}, \|B^T\|_{\rm OP}, 1\} \leq 2c'_r.$$ 
	Since $c_r < c'_r$, one can bound both $\|H^j\|_{\rm OP}$ and $\|{H^j}^T\|_{\rm OP}$ by $2c'_r$.
\end{proof}

\section{Proofs for Singular Value Thresholding}\label{sec:svtproofs}
\begin{proof}[Proof of Theorem~\ref{thm:cvsvt}]
As in all of our applications of Theorem~\ref{thm:main}, the main task for us is to bound the prediction errors for $I = I_1,I_2$,
$$\min_{\lambda \in \Lambda} \|\hat{\theta}^{(\lambda,I)}_{I^c} - \theta^*_{I_c}\|^2 \leq \min_{\lambda \in \Lambda} \|\hat{\theta}^{(\lambda,I)} - \theta^*\|^2.$$
Therefore, our task is to bound $\|\hat{\theta}^{(\lambda,I)} - \theta^*\|^2$ for an appropriate choice of $\lambda \in \Lambda.$

Recall that $\hat{\theta}^{(\lambda,I)}$ is the SVT estimator applied to the matrix $\tilde{y}(I).$
Now we make a few observations. 
\begin{enumerate}
	\item $\tilde{y}(I)$ is unbiased for $\theta^*.$ This is because
	\begin{equation*}
	\E \tilde{y}(I)_{ij} = 2 \E [y_{ij} W_{ij}] = 2 \E y_{ij} \E W_{ij} = \E y_{ij} = \theta^*. 
	\end{equation*}

	\item Define $\tilde{\epsilon}(I) = \tilde{y}(I) - \theta^*$. Observe that we can write
	
	\begin{equation*}
	\tilde{\epsilon}(I)_{ij} = \begin{cases}
	\theta^*_{ij} + 2 \epsilon_{ij} \:\:&\text{with probability} \:\:1/2\\
	-\theta^*_{ij} \:\:&\text{with probability}\:\: 1/2.
	\end{cases}
	\end{equation*}

	This implies that the matrix $\tilde{\epsilon}$ has independent subgaussian entries with subgaussian norm at most $c \left(|\theta^*|_{\infty} + \sigma\right)$ for an absolute constant $c > 0.$ This conclusion follows from Lemma~\ref{lem:subg} (stated and proved after this proof.)
	\end{enumerate}

Based on the above two observations, we can think of $\hat{\theta}^{(\lambda,I)}$ as a usual SVT estimator under the model $\tilde{y}(I) = \theta^* + \tilde{\epsilon}(I).$ The matrix  $\tilde{\epsilon}(I)$ plays the role of the error matrix here which has  mean zero and has independent subgaussian entries with maximum subgaussian norm $c \left(|\theta^*|_{\infty} + \sigma\right).$ We can now apply Theorem~\ref{thm:svt} to conclude that as long as $\lambda \geq C \left(|\theta^*|_{\infty} + \sigma\right) \sqrt{n}$ for an appropriate absolute constant $C > 0$, the same bound as in the right hand side of Theorem~\ref{thm:svt} also holds for $\|\hat{\theta}^{(\lambda,I)} - \theta^*\|^2$. This finishes the proof. 
\end{proof}

\begin{lemma}\label{lem:subg}
	Suppose $\epsilon$ is a mean $0$ subgaussian random variable with subgaussian norm at most $\sigma.$ For a given $K > 0$, consider the random variable 
	\begin{equation*}
	X = \begin{cases}
	K + 2 \epsilon \:\:&\text{with probability} \:\:1/2\\
	-K \:\:&\text{with probability}\:\: 1/2.
	\end{cases}
	\end{equation*}
	Then $X$ is also a mean $0$ subgaussian random variable with subgaussian norm at most $C K \sigma$ where $C > 0$ is an absolute constant. 
	\end{lemma}

\begin{proof}
	We will use the following characterization of a subgaussian random variable, see Proposition $2.5.2$ in~\cite{vershynin2018high}.
	
	A random variable $U$ is subgaussian $(\sigma)$ if there exists an absolute constant $C > 0$ such that for all integers $p \geq 1$, we have 
	$$\left(\E |U|^{p}\right)^{1/p} \leq C \sigma \sqrt{p}.$$

	We will now show that $X$ satisfies the above characterization. Fix any integer $p \geq 1.$ We can write
	\begin{align*}
	&\left(\E |X|^{p}\right)^{1/p} = \left(\frac{1}{2} \E |K + 2 \epsilon|^{p} + \frac{1}{2} K^p\right)^{1/p} \leq  \left(\frac{1}{2} \E |K + 2 \epsilon|^{p}\right)^{1/p} + \left(\frac{1}{2} K^p\right)^{1/p} \leq \\& \left(\E |K + 2 \epsilon|^{p}\right)^{1/p} + K \leq K + \left(\E |2 \epsilon|^{p}\right)^{1/p} + K \leq 2K + C \sigma \sqrt{p} \leq \max\{2,C\} \left(K + \sigma\right) \sqrt{p}.
	\end{align*}
	In the above display, the first inequality follows by using the elementary inequality $(a + b)^{1/p} \leq a^{1/p} + b^{1/p}$ for any $a,b > 0$, second inequality uses the triangle inequality for the $\ell_p$ norm of a random variable, the fourth inequality uses the subgaussian characterization for $\epsilon$ and the last inequality uses the fact that $p \geq 1.$
\end{proof}

\section{Lasso}\label{sec:lasso}
\subsection{\textbf{Background and Related Literature}}
The lasso, proposed by~\cite{tibshirani1996regression} is one of the most popular tools for high dimensional regression. By now, there is a vast literature on analyzing the mean squared error of lasso. The typical statement of the results say that if the tuning parameter $\lambda$ is chosen appropriately depending on some problem parameters (which are typically unknown); then a certain MSE bound holds. However, in practice, the tuning parameter is often chosen using cross validation. The literature giving rigorous theoretical analysis of cross validated lasso is far thinner. As far as we are aware, the first few papers undertaking theoretical analysis of cross validated lasso are~\cite{lecue2012oracle},~\cite{homrighausen2014leave},~\cite{homrighausen2013lasso},~\cite{homrighausen2017risk},~\cite{miolane2018distribution}. Two papers which contain the state of the art theoretical results on cross validated lasso are the papers ~\cite{chatterjee2015prediction},~\cite{chetverikov2020cross}. The paper~\cite{chatterjee2015prediction} is the object of inspiration for the current article. They analyzed a two fold cross validated version of the constrained or primal lasso proposed in~\cite{tibshirani1996regression}. Their result gives the analogue of the so-called slow rate for Lasso (e.g, see Theorem $2.15$ in~\cite{rigollet2015high}) in the fixed design setup. On the other hand,  the paper~\cite{chetverikov2020cross} analyzes a related but different cross validated Lasso estimator and their main result gives an analogue of the fast rate for Lasso (e.g, see Theorem $2.18$ in~\cite{rigollet2015high}) under random design with certain assumptions on the distribution of the covariates and the noise variables.

To the best of our knowledge, a single cross validated lasso estimator which attains both the slow rate and the fast rate in the fixed design setup has not yet been proposed in the literature. Our goal here is to demonstrate that designing such a cross validated lasso is possible. We consider a two fold cross validated version $\hat{\beta}_{cvlasso}$ of the penalized lasso and prove two results. The first result, Theorem~\ref{thm:lassoslow} gives the so-called slow rate under essentially no assumptions on the design matrix. This extends the result of ~\cite{chatterjee2015prediction} to cross validated penalized lasso. Our second result gives the fast rate for the same estimator $\hat{\beta}_{cvlasso}$ under a standard incoherence condition on the design matrix $X.$ Thus, we are able to ensure that qualitatively both the results of~\cite{chatterjee2015prediction} and~\cite{chetverikov2020cross} hold for our $\hat{\beta}_{cvlasso}$ estimator. We now describe the $\hat{\beta}_{cvlasso}$ estimator precisely. We consider a well specified linear model $y = X \beta^* + \epsilon$ where $X$ is a fixed $n \times p$ design matrix and $\epsilon \in \R^n$ is an error vector consisting of i.i.d mean $0$ subgaussian entries. We denote the standard lasso estimator with tuning parameter $\lambda$ by $\hat{\beta}^{(\lambda)} $, defined as follows:
\begin{equation*}
\hat{\beta}^{(\lambda)} := \argmin_{\beta \in \R^p} \left[\sum_{i = 1}^{n} (y_i - x_i^{t} \beta)^2 + \lambda \|\beta\|_1\right],
\end{equation*}
where $x_i^{t}$ is the $i$th row of the design matrix $X.$

%{\color{red} describe model and notations}

%\textbf{Describe briefly the results in these two papers and give the motivation for this paper.}

\subsection{\textbf{Description of the CVLASSO estimator}}

\begin{enumerate}
	\item Set $K = 2$.
	
	\item Divide $[n]$ into $I_1,I_2$ randomly as follows. Each entry $i \in [n]$ belongs to $I_1$ or $I_2$ with probability $1/2$ independently of other entries.
	
	\item For $j \in \{1,2\}$, define 
	\begin{equation*}
	\hat{\beta}^{(\lambda,I_j)} := \argmin_{\beta \in \R^p} \left[\sum_{i \in I_j} (y_i - x_i^{t} \beta)^2 + \lambda \|\beta\|_1\right].
	\end{equation*}
	%where $x_i^{t}$ is the $i$th row of the design matrix $X.$ %Now define $\hat{\theta}^{(\lambda,I_j)} = X \hat{\beta}^{(\lambda,I_j)}.$

	\item  Consider a finite grid of possible values of the tuning parameter $\lambda$, namely $\Lambda = \{1,2,2^2,2^3,\dots,2^{N^*}\}$ where $N^*$ is chosen by the user.  For any $j \in \{1,2\}$, define $\hat{\lambda}_j$ to be the candidate in $\Lambda$ for which the prediction error is the minimum, that is,
	$$\hat{\lambda}_j := \argmin_{\lambda \in \Lambda} \sum_{i \in I_j} (y_i - x_i^{t} \hat{\beta}^{(\lambda,I_j^c)})^2.$$

	%\left|\left|y_{I_j} - \hat{\theta}^{(\lambda, I_j^c)}_{I_j}\right|\right|^2.$$
    Note that $I_1^{c} = I_2$ and vice-versa.

	\item Now define an intermediate estimator $\tilde{\theta} \in \R^n$ such that for any $j \in \{1,2\}$, if $i \in I_j$ then
	\begin{equation*}\label{eq:cvlassodefn}
	\tilde{\theta}_{i} = x_i^{t} \hat{\beta}^{(\hat{\lambda}_j, I_j^c)}.
	\end{equation*}

	\item Define $$\hat{\lambda} := \argmin_{\lambda \in \Lambda} \|X \hat{\beta}^{(\lambda)} - \tilde{\theta}\|^2$$ %where $\Lambda = \{1, 2, 2^2, 2^3,\dots, 2^{N^*}\}.$ 
	Finally, our estimator (CVLASSO) is defined to be $$\hat{\beta}_{cvlasso} = \hat{\beta}^{(\hat{\lambda})}.$$
	
	\end{enumerate}

\subsection{\textbf{Main Results}}

\begin{theorem}\label{thm:lassoslow}[Slow Rate]
Suppose $M > 0$ is a number such that the design matrix $X$ satisfies $$ \max_{j \in [p]} \frac{1}{n} \sum_{i = 1}^{n} X_{ij}^{4} \leq M.$$ Fix any $0 < \delta < 1.$ Suppose we take our  grid $\Lambda = \{1,2,2^2,\dots,2^{N^*}\}$ such that
\begin{equation*}
2^{N^*} \geq 4 \sqrt{2 \sigma^2 M^{1/2} n \log p} + \sqrt{2 \sigma^2 M^{1/2} n \log 1/\delta}.
\end{equation*}
 Then we have the following bound with probability atleast $1 - 7 \delta$ for an appropriate absolute constant $C > 0$,
%\begin{align*}
%\sqrt{SSE(\hat{\theta}_{cvlasso},\theta^*)} \leq 4 \sqrt{32} |\beta^*|_1 (M n \log p)^{1/4} + 4 (2 \log 1/\delta)^{1/4} + 8\sqrt{2} \sigma \sqrt{\log 2^{K^*}} + 8 \sqrt{2 \sigma^2 \log 4/\delta}
%\end{align*}
\begin{align*}
\E \frac{1}{n} \|X \hat{\beta}_{cvlasso} - X \beta^*\|^2 \leq \frac{C}{n} \Big[\|\beta^*\|_1^2 \sqrt{M n \log p} + \|\beta^*\|_1 \sigma &\sqrt{M n \log p/\delta} + \sigma^2 N^*\\
&+  \sigma^2 \log 1/\delta\Big].
\end{align*}
\end{theorem}

%\begin{remark}
	%A reasonable value of $C$ in the above theorem can be given. This will be clear from the proof. 
%	\end{remark}

\begin{remark}
	The above result is qualitatively similar to the result in Theorem $2.1$ in~\cite{chatterjee2015prediction}. The main difference is that while that result is about a cross validated version of the constrained lasso, our result is about the corresponding cross validated version of the penalized lasso. There are certain advantages of using the penalized form of Lasso instead of the constrained form as mentioned in Remark~\ref{rem:comp}. 
\end{remark}

\begin{remark}
	The bound in Theorem~\ref{thm:lassoslow} basically says that the MSE of $\hat{\theta}_{cvlasso}$ scales like $O(\|\beta^*\|_1^2 \sqrt{\frac{\log p}{n}})$ with high probability if $M$ and $\sigma$ are bounded away from $\infty.$ Note that, this result holds essentially without any assumptions on the design matrix. As mentioned in~\cite{chatterjee2015prediction}, this MSE scaling agrees with the persistency condition for lasso (under random design) defined in~\cite{greenshtein2004persistence} which says that if  $\|\beta^*\|_1 = o\big((\frac{n}{\log p})^{1/4}\big)$ then persistency holds. 
\end{remark}

\begin{theorem}\label{thm:lassofast}[Fast Rate]
	Suppose that the design matrix $X$ satisfies an incoherence condition
	$$\left|\frac{X^t X}{n} -  I_{p \times p}\right|_{\infty} \leq \frac{1}{64 k}$$
	where $k = \|\beta^*\|_0.$ 
	Suppose $M > 0$ is a number such that the design matrix $X$ satisfies $$ \max_{j \in [p]} \frac{1}{n} \sum_{i = 1}^{n} X_{ij}^{4} \leq M.$$
    Assume that the sample size $n$ is large enough so that 
	$$ 2 \log p < \frac{n}{2^{14} k^2 M^2}.$$ 
	
	Fix any $0 < \delta < 1.$ Suppose we take our  grid $\Lambda = \{1,2,2^2,\dots,2^{N^*}\}$ such that
	\begin{equation*}
	2^{N^*} \geq 4 \sqrt{2 \sigma^2 M^{1/2} n \log p} + \sqrt{2 \sigma^2 M^{1/2} n \log 1/\delta}.
	\end{equation*}
	Then we have the following bound with probability atleast $1 - 6 \delta - \exp\left(-\frac{n}{2^{14} k^2 M^2}\right)$ for an appropriate absolute constant $C > 0$,
	%\begin{align*}
	%\sqrt{SSE(\hat{\theta}_{cvlasso},\theta^*)} \leq 4 \sqrt{32} |\beta^*|_1 (M n \log p)^{1/4} + 4 (2 \log 1/\delta)^{1/4} + 8\sqrt{2} \sigma \sqrt{\log 2^{K^*}} + 8 \sqrt{2 \sigma^2 \log 4/\delta}
	%\end{align*}
	\begin{align*}
	\E \frac{1}{n} \|X \hat{\beta}_{cvlasso} - X \beta^*\|^2 \leq \frac{C}{n} \Bigg[\frac{k^2}{\sqrt{n}} M( \log (p/\delta))^{3/2} + \sigma^2 k \sqrt{M} \log (p/\delta) &+ \sigma^2 N^*\\
	&+  \sigma^2 \log 1/\delta\Bigg].
	\end{align*}
\end{theorem}

\begin{remark}
Assuming that $M$ and $\sigma$ are terms bounded away from $\infty$, the first term (inside the brackets) in the bound given in Theorem~\ref{thm:lassofast} which scales like $O\left(\frac{k^2}{\sqrt{n}} (\log p)^{3/2}\right)$ is dominated by the second term $O(k \log p)$ as long as $k < O(\sqrt{n/\log p})$ which we can readily check is the interesting regime where we can expect fast rates. This is because as soon as $k =  O(\sqrt{n/\log p})$, the fast SSE rate $O(k \log p) = O(\sqrt{n \log p})$ which matches the slow rate. Therefore, when $k \geq  O(\sqrt{n/\log p})$, one should use the slow rate result in Theorem~\ref{thm:lassoslow}. To summarize, the above result in Theorem~\ref{thm:lassofast} is useful in the sparse regime when $k < O(\sqrt{n/\log p})$ in which case the MSE of $\hat{\theta}_{cvlasso}$ scales like the fast rate $O(\frac{k \log p}{n})$ with high probability.
\end{remark}

\begin{remark}
	In Theorem~\ref{thm:lassofast}, we need a slightly stronger (by a factor of $2$) incoherence condition than the standard one assumed for penalized Lasso in the literature. For example, Theorem $2.18$ in~\cite{rigollet2015high} assumes that $\left|\frac{X^t X}{n} -  I_{p \times p}\right|_{\infty} \leq \frac{1}{32 k}.$ It is well known that weaker (than incoherence) assumptions on the design matrix $X$ such as the restricted isometry property are also sufficient to ensure fast rates for the lasso. Such results are likely to be true for our cross validated estimator $\hat{\theta}_{cvlasso}$ as well. However, we leave this for future research and just consider the incoherence condition because of two reasons. Firstly, it seems to be the simplest sufficient condition for fast rates available in the literature and is actually checkable in practice in contrast to the restricted isometry type properties which are computationally intractable to check. Secondly, our goal here is to simply demonstrate that both the slow rate and the fast rate are attainable for a single cross validated lasso estimator. Thus, we prefer to sacrifice some generality in exchange to demonstrate a phenomenon under a simpler sufficent condition. 
	\end{remark}

\begin{remark}
	For both of our theorems, the grid $\Lambda = \{1,2,2^2,\dots,2^{N^*}\}$ needs to satisfy that $2^{N^*} \geq C \sigma \sqrt{M^{1/2} n \log p}$ for a specified constant $C.$ This is a very mild condition to ensure in practice. The parameter $\sigma$ is the only unknown term and as explained before, even a gross over estimate can be plugged in without any serious consequences since our grid grows exponentially. Thus, the number of grid points would be $O(\log n \log \log p)$ which means we would need to solve the lasso optimization problem $O(\log n \log \log p)$ times to compute $\hat{\theta}_{cvlasso}.$ In other words, the computational complexity of $\hat{\theta}_{cvlasso}$ would be $O(\log n \log \log p)$ times the computational complexity of computing a single instance of the penalized lasso estimator. 
	\end{remark}

\begin{remark}\label{rem:comp}
	A potential advantage of the cross validated penalized lasso over constrained lasso is as follows. It is known that under certain nonsingularity conditions on the design matrix $X$, to attain fast rates for the usual constrained lasso, it is sufficient that the tuning parameter is chosen to be exactly equal to $\|\beta^*\|_1$; see Theorem $2.1$ in~\cite{chatterjee2014new}. Clearly, this is hard to achieve in practice. It is not known to what extent is this result robust to the choice of this tuning parameter. For example, Theorem $2.1$ in~\cite{chatterjee2014new} further indicates that if the tuning parameter is chosen to be $\|\beta^*\|_1 \pm 1$ then it is not possible for the constrained lasso to attain fast rates. In contrast, the penalized lasso seems to be more robust with respect to the choice of its tuning parameter. For example, if one sets $\lambda$ to be twice the ideal choice of $\lambda$ known to achieve fast rates for penalized lasso (e.g, see Theorem $2.18$ in~\cite{rigollet2015high}), then the risk at most doubles and hence the rate of convergence remains the same. 
	
	The upshot of this is that we can afford to have a grid of $\lambda$ growing exponentially and still attain fast rates for our estimator $\hat{\theta}_{cvlasso}$. This has significant computational advantages as this means the cardinality of our grid $\Lambda$ is only growing like $O(\log n \log \log p)$ which means we have to solve the lasso optimization problem at most $O(\log n \log \log p)$ times to compute $\hat{\theta}_{cvlasso}$. In contrast, it is likely that the $\Lambda$ grid needs to be much finer in resolution (with cardinality growing like $n^{\alpha}$ for some $\alpha > 0$) for the cross validated constrained Lasso proposed in~\cite{chatterjee2015prediction} to attain fast rates. %This is essentially the same reason why we considered Penalized Trend Filtering instead of its constrained counterpart in this paper. 

	%{\color{red} this is repetition. make this shorter?}

  	\end{remark}

  \begin{remark}
  	Like in~\cite{chatterjee2015prediction} we have proposed a $2$ fold cross validated version of Lasso. However, if it is so desired, one can easily construct a similar $K$ fold version as should be clear from our general framework and the description of $\hat{\theta}_{cvlasso}$. Similar risk bounds as in Theorem~\ref{thm:lassoslow} and Theorem~\ref{thm:lassofast} would hold for the $K$ fold version as well. 
  	\end{remark}

\section{Proofs for Lasso}\label{sec:appendixlasso}

%{\color{red} I will now hand over this to Anamitra to read every line very carefully. Look out for notational inconsistencies, typos and mistakes in constants. I have tried to keep the constants till the last level of proving the theorems. This is not terribly important but still let s keep it this way. If there is a way to simplify further lets discuss.}

\begin{proof}[Proof of Theorem~\ref{thm:lassoslow}]
	The main step to prove Theorem~\ref{thm:lassoslow} is the following proposition (proved after this proof). 
	
	\begin{proposition}\label{prop:lassoslow}
		Suppose $M > 0$ is a number such that the design matrix $X$ satisfies $$ \max_{j \in [p]} \frac{1}{n} \sum_{i = 1}^{n} X_{ij}^{4} \leq M.$$ Fix any $0 < \delta < 1.$ Set $$\lambda \geq 4 \left[\sqrt{2 \sigma^2 M^{1/2} n \log p} + \sqrt{2 \sigma^2 M^{1/2} n \log 1/\delta}\right].$$ Then for any $j \in \{1, 2\},$ we have the following bound with probability atleast $1 - 3 \delta$,
		\begin{align*}
		SSE(\hat{\theta}^{(\lambda, I_j^c)}_{I_j}, \theta^*_{I_j}) \leq 16 \|\beta^*\|_1^2 \left(2 \sqrt{\log p}\:n^{1/2} M^{1/2} + \sqrt{2 \log 1/\delta}\right) + \frac{3}{2} \|\beta^*\|_1 \lambda
		\end{align*}
		%where $$A = 32 |\beta^*|_1^2 (2 \sqrt{\log p}\:n^{1/2} M^{1/2} + \sqrt{2 \log 1/\delta})$$ and $$B = 3 |\beta^*|_1 \lambda.$$ 
	\end{proposition}

	Now, by our choice of $N^*$, there exists a choice of $\lambda \in \Lambda$ such that 
	$$4 \sqrt{2 \sigma^2 M^{1/2} n \log p} + \sqrt{2 \sigma^2 M^{1/2} n \log 1/\delta}  \leq \lambda \leq 8 \sqrt{2 \sigma^2 M^{1/2} n \log p} + \sqrt{2 \sigma^2 M^{1/2} n \log 1/\delta}.$$

	For this choice of $\lambda$ we can invoke Proposition~\ref{prop:lassoslow} to obtain for any $j \in \{1, 2\},$ the following bound with probability at least $1 - 3 \delta$,
	\begin{align*}
	SSE(\hat{\theta}^{(\lambda_j, I_j^c)}_{I_j}, \theta^*_{I_j}) \leq C \left[\|\beta^*\|_1^2 \sqrt{M n \log p} + \|\beta^*\|_1 \sigma \sqrt{M n \log p/\delta}\right].
	\end{align*}

	Now, plugging in the above bound in to Theorem~\ref{thm:main} finishes the proof.
\end{proof}

\begin{proof}[Proof of Proposition~\ref{prop:lassoslow}]

	Let $L = 4 \|\beta^*\|_1.$ Let us define the event
	
	\begin{equation*}
	A = \left\{\sup_{\Delta \in \R^p: \|\Delta\|_1 \leq L}  \sum_{i = 1}^{n} \eta_i (x_i^{t} \Delta)^2 \leq L^2 \left(2 \sqrt{\log p}\:n^{1/2} M^{1/2} + \sqrt{2 \log 1/\delta}\right)\right\}
	\end{equation*}
	where $\eta_i$ are i.i.d Rademacher random variables independent of the data vector $y.$
	By Lemma~\ref{lem:supell1}, $P(A^{c}) \leq 2 \delta.$

	Also define the event 
	\begin{equation*}
	B = \{\lambda \geq 4 |X^{t} \epsilon|_{\infty}\}
	\end{equation*}
	By our choice of $\lambda$ and lemma~\ref{lem:lassomax}, we also have $P(B^{c}) \leq \delta.$

	We will now work inside the event $A \cap B.$ Since $B$ holds, we know by Lemma~\ref{lem:restricted1} that 
	\begin{equation}\label{eq:l1constraint1}
	\|\hat{\beta}^{(\lambda,I_j^{c})} - \beta^*\| \leq 4 \|\beta^*\|_1 = L.
	\end{equation}
	
	Now we can write
	\begin{align*}
	&SSE(\hat{\theta}^{(\lambda, I_j^c)}_{I_j}, \theta^*_{I_j}) = \sum_{i \in I_j}  (x_i^{t} \hat{\beta}^{(\lambda, I_j^{c})} - x_i^{t} \beta^*)^2 = \\& \underbrace{\sum_{i \in I_j}  (x_i^{t} \hat{\beta}^{(\lambda, I_j^{c})} - x_i^{t} \beta^*)^2 - \sum_{i \in I_j^{c}}  (x_i^{t} \hat{\beta}^{(\lambda, I_j^{c})} - x_i^{t} \beta^*)^2}_{T_1} +  \underbrace{\sum_{i \in I_j^{c}}  (x_i^{t} \hat{\beta}^{(\lambda, I_j^{c})} - x_i^{t} \beta^*)^2}_{T_2}.
		\end{align*}

	Now note that $T_2$ is the usual squared error in a lasso problem where the design matrix is $X_{I_j^{c}}.$ %Therefore, this can be handled using standard existing methods; we outline these steps later. 
	Therefore, the main task now is to bound $T_1.$ For this, we can write the pointwise inequality
	\begin{align*}
	T_1 &\leq  \sup_{\Delta \in \R^p: \|\Delta\|_1 \leq L} \sum_{i \in I_j}  (x_i^{t} \Delta)^2 - \sum_{i \in I_j^{c}}  (x_i^{t} \Delta)^2 =  \sup_{\Delta \in \R^p: \|\Delta\|_1 \leq L}  \sum_{i = 1}^{n} \eta_i (x_i^{t} \Delta)^2 \\&\leq L^2 \left(2 \sqrt{\log p}\:n^{1/2} M^{1/2} + \sqrt{2 \log 1/\delta}\right)
	\end{align*}
	where the first inequality is because of~\eqref{eq:l1constraint1} and the last inequality is because the event $A$ holds. 
	
	%where $L$ is any deterministic number satisfying $\|\hat{\beta}^{(\lambda, I_j^{c})} - \beta^*\| \leq L$ and $\eta_i$ are i.i.d Rademacher random variables independent of the data vector $y.$ The last equality follows because of the way $I_1,I_2$ is chosen.
	
	By Lemma~\ref{lem:restricted1} again, we can also conclude that 
	$$T_2 \leq \frac{3}{2} \lambda \|\beta^*\|_{1}.$$ 
	
	The last two displays imply that when the events $A,B$ hold then 
	
	$$SSE(\hat{\theta}^{(\lambda_j, I_j^c)}_{I_j}, \theta^*_{I_j}) \leq L^2 \left(2 \sqrt{\log p}\:n^{1/2} M^{1/2} + \sqrt{2 \log 1/\delta}\right) + \frac{3}{2} \lambda \|\beta^*\|_{1}.$$

	%{\color{red} Why is $\sigma$ absent in the first term above? Check carefully tomorrow. I really want to not track constants and also keep lambda in the risk bound. So appropriately modify the writing.}

\end{proof}

\begin{proof}
	The main step to prove Theorem~\ref{thm:lassofast} is the following proposition (proved after this proof). 
	\begin{proposition}\label{prop:lassofast}
		Suppose $M > 0$ is a number such that the design matrix $X$ satisfies $$ \max_{j \in [p]} \frac{1}{n} \sum_{i = 1}^{n} X_{ij}^{4} \leq M.$$ Fix any $0 < \delta < 1.$ Set $$\lambda \geq 4 \left[\sqrt{2 \sigma^2 M^{1/2} n \log p} + \sqrt{2 \sigma^2 M^{1/2} n \log 1/\delta}\right].$$ Suppose that the design matrix $X$ satisfies an incoherence condition
		$$\left|\frac{X^t X}{n} -  I_{p \times p}\right|_{\infty} \leq \frac{1}{64 k}$$
		where $k = \|\beta^*\|_0.$ Also, we assume that the sample size $n$ is large enough so that 
		$$ 2 \log p < \frac{n}{2^{14} k^2 M^2}.$$

		Then for any $j \in \{1, 2\}$, we have the following bound with probability atleast $1 - 3 \delta - \exp\left(\frac{n}{2^{14} k^2 M^2}\right)$,
		\begin{align*}
		SSE(\hat{\theta}^{(\lambda_j, I_j^c)}_{I_j}, \theta^*_{I_j}) \leq 2^{14} \lambda^2 \frac{k^2}{n^2} \left(2 \sqrt{\log p}\:n^{1/2} M^{1/2} + \sqrt{2 \log 1/\delta}\right) + 64 \lambda^2 \frac{k}{n}.
		\end{align*}

	\end{proposition}

	Now, by our choice of $N^*$, there exists a choice of $\lambda \in \Lambda$ such that 
	\begin{align*}
	4 \left[\sqrt{2 \sigma^2 M^{1/2} n \log p} + \sqrt{2 \sigma^2 M^{1/2} n \log 1/\delta}\right]  \leq \lambda \leq 8 \Bigg[&\sqrt{2 \sigma^2 M^{1/2} n \log p}\\
	&+ \sqrt{2 \sigma^2 M^{1/2} n \log 1/\delta}\Bigg].
	\end{align*}

	For this choice of $\lambda$ we can invoke Proposition~\ref{prop:lassofast} to obtain for any $j \in \{1, 2\},$ the following bound with probability atleast $1 - 3 \delta$,
	\begin{align*}
	SSE(\hat{\theta}^{(\lambda_j, I_j^c)}_{I_j}, \theta^*_{I_j}) &\leq C \left[\frac{k^2}{n^2} \sqrt{M} n \log (p/\delta) \sqrt{\log p/\delta}\:n^{1/2} M^{1/2}\right] + C \left[\sigma \frac{k}{n} \sqrt{M^{1/2} n \log (p/\delta)}\right] \\& = C \left[\frac{k^2}{\sqrt{n}} M( \log (p/\delta))^{3/2}\right] + C \left[\sigma^2 k \sqrt{M} \log (p/\delta)\right].
	\end{align*}

	Now, plugging in the above bound in to Theorem~\ref{thm:main} finishes the proof.
	%{\color{red} Finish the statement of this theorem.}
	
\end{proof}

\begin{proof}[Proof of Proposition~\ref{prop:lassofast}]
	Fix a $j \in [2].$ Let $L = 128 \lambda \frac{k}{n}.$ Similar to the proof of Proposition~\ref{prop:lassoslow}, let us define the events
	\begin{equation*}
	A = \left\{\sup_{\Delta \in \R^p: \|\Delta\|_1 \leq L}  \sum_{i = 1}^{n} \eta_i (x_i^{t} \Delta)^2 \leq L^2 (2 \sqrt{\log p}\:n^{1/2} M^{1/2} + \sqrt{2 \log 1/\delta})\right\}
	\end{equation*}
	and 
	\begin{equation*}
	B = \{\lambda \geq 4 |X^{t} \epsilon|_{\infty}\}
	\end{equation*}
	where $\eta_i$ are i.i.d Rademacher random variables. By Lemma~\ref{lem:supell1}, $P(A^{c}) \leq 2 \delta$ and by our choice of $\lambda$ and Lemma~\ref{lem:lassomax}, we also have $P(B^{c}) \leq \delta.$

	%Also define the event 
	%\begin{equation*}
	%B = \{\lambda \geq 4 |X^{t} \epsilon|_{\infty}\}
	%\end{equation*}
	%By our choice of $\lambda$ and lemma~\ref{lem:lassomax}, we also have $P(B^{c}) \leq \delta.$

	We also define the event 
	$$C = \left\{\left|\frac{X_{I}^t X_{I}}{n} - \frac{X^t X}{2n}\right|_{\infty} \leq 1/128 k\right\}.$$
	
	By Lemma~\ref{lem:incoh} and the fact that $X$ satisfies the stated incoherence condition in this proposition, an application of triangle inequality lets us conclude that $P(C^c) \leq \exp(-\frac{n}{2^{14} k^2 M^2}).$

	We will now work inside the event $A \cap B \cap C.$ Since $B$ holds, we know by Lemma~\ref{lem:restricted2} that 
	\begin{equation}\label{eq:l1constraint2}
	\|\hat{\beta}^{(\lambda,I_j^{c})} - \beta^*\|_1 \leq L.
	\end{equation}
	
	Now we can again write
	\begin{align*}
	&SSE(\hat{\theta}^{(\lambda, I_j^c)}_{I_j}, \theta^*_{I_j}) = \sum_{i \in I_j}  (x_i^{t} \hat{\beta}^{(\lambda, I_j^{c})} - x_i^{t} \beta^*)^2 = \\& \underbrace{\sum_{i \in I_j}  (x_i^{t} \hat{\beta}^{(\lambda, I_j^{c})} - x_i^{t} \beta^*)^2 - \sum_{i \in I_j^{c}}  (x_i^{t} \hat{\beta}^{(\lambda, I_j^{c})} - x_i^{t} \beta^*)^2}_{T_1} +  \underbrace{\sum_{i \in I_j^{c}}  (x_i^{t} \hat{\beta}^{(\lambda, I_j^{c})} - x_i^{t} \beta^*)^2}_{T_2}.
		\end{align*}

	%	Now note that $T_2$ is the usual squared error in a lasso problem where the design matrix is $X_{I_j^{c}}.$ %Therefore, this can be handled using standard existing methods; we outline these steps later. 
	To bound $T_1$ we can write 
	
	\begin{align*}
	&T_1 \leq  \sup_{\Delta \in \R^p: \|\Delta\|_1 \leq L} \sum_{i \in I_j}  (x_i^{t} \Delta)^2 - \sum_{i \in I_j^{c}}  (x_i^{t} \Delta)^2 =  \sup_{\Delta \in \R^p: \|\Delta\|_1 \leq L}  \sum_{i = 1}^{n} \eta_i (x_i^{t} \Delta)^2 \leq \\& L^2 (2 \sqrt{\log p}\:n^{1/2} M^{1/2} + \sqrt{2 \log 1/\delta})
	\end{align*}
	where the first inequality is because of~\eqref{eq:l1constraint2} and the last inequality is because the event $A$ holds. 
	
	%where $L$ is any deterministic number satisfying $\|\hat{\beta}^{(\lambda, I_j^{c})} - \beta^*\| \leq L$ and $\eta_i$ are i.i.d Rademacher random variables independent of the data vector $y.$ The last equality follows because of the way $I_1,I_2$ is chosen.
	
	By Lemma~\ref{lem:restricted2} again and because the event $C$ holds, we can also conclude that 
	$$T_2 \leq 64 \lambda^2 \frac{k}{n}.$$ 
	
	The last two displays imply that when the events $A,B,C$ hold then 
	
	$$SSE(\hat{\theta}^{(\lambda_j, I_j^c)}_{I_j}, \theta^*_{I_j}) \leq 2^{14} \lambda^2 \frac{k^2}{n^2} \left(2 \sqrt{\log p}\:n^{1/2} M^{1/2} + \sqrt{2 \log 1/\delta}\right) + 64 \lambda^2 \frac{k}{n}.$$
	
\end{proof}

\subsection{Auxiliary Lemmas}

\begin{comment}

The following standard result appears as Lemma $5.2$ in~\cite{van2014probability}
\begin{lemma}[Finite Maxima of Subgaussians]\label{lem:subg}
	If $X_1,\dots,X_n$ are finitely many random variables each of which is $\sigma^2$ subgaussian then we have for all $x \geq 0,$
	\begin{equation*}
	P(\max_{i \in [n]} X_i \geq \sqrt{2 \sigma^2 \log n} + x) \leq \exp(-\frac{x^2}{2 \sigma^2}).
	\end{equation*}
	\end{lemma}

{\color{red} The above one is a repetition. Delete this theorem changing the citations properly.}
\end{comment}

\begin{lemma}\label{lem:supell1}
	Let $X$ be a $n \times p$ matrix. Let $M = \max_{1 \leq j \leq p} \left[\frac{1}{n} \sum_{i = 1}^{n} x_{ij}^4\right].$ Also for $i \in [n]$, let $\eta_i$ be i.i.d rademacher random variables. Then we have the following inequality:
	\begin{equation*}
	P\left(\sup_{\Delta \in \R^p: \|\Delta\|_1 \leq L}  \sum_{i = 1}^{n} \eta_i (x_i^{t} \Delta)^2 \leq L^2 (2 \sqrt{\log p}\:n^{1/2} M^{1/2} + x)\right) \geq 1 - 2 \exp\left(-\frac{x^2}{2}\right).
	\end{equation*}
	%Here, $M = \max_{1 \leq j \leq p} [\frac{1}{n} \sum_{i = 1}^{n} x_{ij}^4].$
	
\end{lemma}

\begin{proof}
	We can expand the squares and write
	\begin{align*}
	\sup_{\Delta \in \R^p: \|\Delta\|_1 \leq L}  \sum_{i = 1}^{n} \eta_i (x_i^{t} \Delta)^2 &= \sup_{\Delta \in \R^p: \|\Delta\|_1 \leq L} \sum_{j,k = 1}^{p} \Delta_j \Delta_k \left(\sum_{i = 1}^{n} \eta_i x_{ij} x_{ik}\right)\\ &\leq L^2 \max_{1 \leq j,k \leq p} \left|\sum_{i = 1}^{n} \eta_i x_{ij} x_{ik}\right|.
	\end{align*}
	By Lemma \ref{lem:max_subG}, we have for all $x \geq 0$,
	$$P\left(\max_{1 \leq j,k \leq p} \left|\sum_{i = 1}^{n} \eta_i x_{ij} x_{ik}\right| \leq 2 \sqrt{\log p} \left[\max_{1 \leq j,k \leq p} \sum_{i = 1}^{n} x_{ij}^2 x_{ik}^2\right]^{1/2}\: + x\right) \leq 2 \exp\left(-\frac{x^2}{2}\right).$$
	
	Now note that by the Cauchy Schwarz inequality $$\max_{1 \leq j,k \leq p} \left(\sum_{i = 1}^{n} x_{ij}^2 x_{ik}^2\right)^{1/2} \leq n^{1/2} M^{1/2}.$$ This finishes the proof. 
\end{proof}

\begin{lemma}\label{lem:lassomax}
	Let $X$ be a $n \times p$ matrix with $M = \max_{1 \leq j \leq p} \left[\frac{1}{n} \sum_{i = 1}^{n} x_{ij}^4\right].$ Let $\epsilon \in \R^n$ be a random vector of i.i.d mean $0$ subgaussian random variables with subgaussian norm $\sigma$. Then we have for any $\delta > 0$,
	\begin{equation*}
	P\left(|X^t \epsilon|_{\infty} \geq \sqrt{2 \sigma^2 M^{1/2} n \log p} + \sqrt{2 \sigma^2 M^{1/2} n \log 1/\delta}\right) \leq \delta.
	\end{equation*}
	
	\end{lemma}

\begin{proof}
	Observe that for each $i \in [p]$ we have $X[,i]^{t} \epsilon$ is subgaussian with subgaussian norm atmost $\sigma \|X[,i]\| \leq \sigma M^{1/4} n^{1/2}$ by Cauchy Schwarz inequality. Therefore, applying lemma~\ref{lem:max_subG} we have for all $x \geq 0$,
	
	$$P\left(|X^t \epsilon|_{\infty} \geq \sqrt{2 \sigma^2 M^{1/2} n \log p} + x\right) \leq \exp\left(-\frac{x^2}{2 \sigma^2 M^{1/2} n}\right).$$ 
	Setting $\delta$ to be the right hand side in the above display finishes the proof. 
\end{proof}

\begin{lemma}\label{lem:restricted1}
	If $\lambda$ is chosen such that $\lambda \geq 4 |X^{t} \epsilon|_{\infty}$ then for any $j \in \{1,2\}$ we have two pointwise inequalities. Firstly, 
	\begin{equation*}
	\|\hat{\beta}^{(\lambda, I_j)} - \beta^*\|_1 \leq 4 \|\beta^*\|_1.
	\end{equation*}
	Secondly, we also have 
	\begin{equation*}
	\sum_{i \in I_j} \left(x_i^{t} \beta^* - x_i^{t} \hat{\beta}^{(\lambda, I_j)}\right)^2 \leq \frac{3}{2} \lambda \|\beta^*\|_1.
	\end{equation*}
\end{lemma}

\begin{proof}
	By definition, we have the basic inequality:
	\begin{equation*}
	\sum_{i \in I_j} \left(y_i - x_i^{t} \hat{\beta}^{(\lambda, I_j)}\right)^2 + \lambda \left\|\hat{\beta}^{(\lambda, I_j)}\right\|_1 \leq \sum_{i \in I_j} (y_i - x_i^{t} \beta^*)^2 + \lambda \|\beta^*\|_1
	\end{equation*}
	Writing $y = X \beta^* + \epsilon$ and expanding the squares in above we obtain
	\begin{align} \notag
	\sum_{i \in I_j} \left(x_i^{t} \beta^* - x_i^{t} \hat{\beta}^{(\lambda, I_j)}\right)^2 + \frac{\lambda}{2} \left\|\hat{\beta}^{(\lambda, I_j)}\right\|_1 \leq 2 \sum_{i \in I_j} (x_i^{t} \hat{\beta}^{(\lambda, I_j)} - x_i^t \beta^*) \epsilon_i \:\: -\; &\frac{\lambda}{2} \left\|\hat{\beta}^{(\lambda, I_j)}\right\|_1\\ \label{eq:basiclasso}
	 &+ \lambda \|\beta^*\|_1.
	\end{align}

	Now using Holder's inequality and the given choice of $\lambda$ we further obtain
	\begin{align*}
	&\sum_{i \in I_j} \left(x_i^{t} \beta^* - x_i^{t} \hat{\beta}^{(\lambda, I_j)}\right)^2 + \frac{\lambda}{2} \left\|\hat{\beta}^{(\lambda, I_j)}\right\|_1\\ &\leq 2 |X_{I_j}^t \epsilon|_{\infty} \|\hat{\beta}^{(\lambda, I_j)}\|_1 + 2 |X_{I_j}^t \epsilon|_{\infty} \|\beta^*\|_1  - \frac{\lambda}{2} \|\hat{\beta}^{(\lambda, I_j)}\|_1 + \lambda \|\beta^*\|_1 \leq \frac{3}{2} \lambda \|\beta^*\|_1.
	\end{align*}

	This further implies
	\begin{equation*}
	\|\hat{\beta}^{(\lambda, I_j)}\|_1 \leq 3 \|\beta^*\|_1. 
	\end{equation*}
	We can now apply triangle inequality above to get $\|\hat{\beta}^{(\lambda, I_j)} - \beta^*\|_1 \leq 4 \|\beta^*\|_1.$
	We can also conclude from the second last display that 
	\begin{equation*}
	\sum_{i \in I_j} \left(x_i^{t} \beta^* - x_i^{t} \hat{\beta}^{(\lambda, I_j)}\right)^2 \leq \frac{3}{2} \lambda \|\beta^*\|_1.
	\end{equation*}
	%The above along with an application of triangle inequality finishes the proof. 
\end{proof}

\begin{lemma}\label{lem:restricted2}
	Fix any $j \in \{1, 2\}$. Suppose $X_{I_j}$ satisfies $$\left|\frac{X_{I_j}^t X_{I_j}}{n} - \frac{1}{2} I_{p \times p}\right|_{\infty} \leq \frac{1}{64 k}.$$
	
	  Also suppose that $\lambda$ is chosen such that $\lambda \geq 4 |X^{t} \epsilon|_{\infty}.$ Then we have two pointwise inequalities. Firstly,
	\begin{equation*}
		\|X_{I_j} (\hat{\beta}^{(\lambda, I_j)} - \beta^*)\| \leq 8 \lambda \sqrt{\frac{k}{n}}.
		\end{equation*}
		
		Secondly,
		\begin{equation*}
		\|\hat{\beta}^{(\lambda, I_j)} - \beta^*\|_1 \leq 128 \lambda \frac{k}{n}.
		\end{equation*}

		%+ \frac{\lambda}{2} |\hat{\beta}^{(\lambda, I_j)} - \beta^*|_1 \leq \sqrt{2} \lambda \sqrt{\frac{k}{n}} \|X_{I_j} (\hat{\beta}^{(\lambda, I_j)} - \beta^*)\|_2
%	\end{equation*}
\end{lemma}

\begin{proof}
	Let $S = \{j \in [p]: \beta_j \neq 0\}.$ An equivalent way to write~\eqref{eq:basiclasso} is as follows:
	\begin{equation*}
	\sum_{i \in I_j} \left(x_i^{t} \beta^* - x_i^{t} \hat{\beta}^{(\lambda, I_j)}\right)^2 \leq 2 \sum_{i \in I_j} \left(x_i^{t} \hat{\beta}^{(\lambda, I_j)} - x_i^t \beta^*\right) \epsilon_i \:\: - \lambda \|\hat{\beta}^{(\lambda, I_j)}\|_1 + \lambda \|\beta^*\|_1.
	\end{equation*}

	Now we can write 
	\begin{align*}
	2 \sum_{i \in I_j} \left(x_i^{t} \hat{\beta}^{(\lambda, I_j)} - x_i^t \beta^*\right) \epsilon_i \leq 2 |X_{I_j}^t \epsilon|_{\infty} \|\hat{\beta}^{(\lambda, I_j)} - \beta^*\|_1 \leq \frac{\lambda}{2} \|\hat{\beta}^{(\lambda, I_j)} - \beta^*\|_1,
	\end{align*}
	where we have used Holder's inequality and the fact that $\lambda \geq 4 |X^{t} \epsilon|_{\infty}$. 
	
	The last two displays imply that 
	\begin{align*}
	&\|X_{I_j} ( \hat{\beta}^{(\lambda, I_j)} - \beta^*)\|^2 + \frac{\lambda}{2} \|\hat{\beta}^{(\lambda, I_j)} - \beta^*\|_1 \leq \lambda \|\hat{\beta}^{(\lambda, I_j)} - \beta^*\|_1 + \lambda \|\beta^*\|_1 - \lambda \|\hat{\beta}^{(\lambda, I_j)}\|_1 =\\&  \lambda \|\hat{\beta}^{(\lambda, I_j)}_{S} - \beta^*_{S}\|_1 + \lambda \|\beta^*_{S}\|_1 - \lambda \|\hat{\beta}^{(\lambda, I_j)}_{S}\|_1 \leq 2 \lambda \|\hat{\beta}^{(\lambda, I_j)}_{S} - \beta^*_{S}\|_1.
	\end{align*}

	In particular, the above display implies that the vector $\hat{\Delta} = \hat{\beta}^{(\lambda, I_j)} - \beta^*$ lies in the set $$C = \{v \in \R^p: \|v\|_{1} \leq 4 \|v_{S}\|_1\}.$$ This is commonly referred to as the \textit{cone condition} in the Lasso literature. 
	
	%Now by Lemma $2.17$ in~\cite{rigollet2015high} {\color{red} state this} we can further write
	We can further write
	\begin{align*}
	\|X_{I_j} \hat{\Delta}\|^2 + \frac{\lambda}{2} \|\hat{\Delta}\|_1 \leq 2 \lambda \|\hat{\Delta}_{S}\|_1 \leq 2 \lambda \sqrt{k} \|\hat{\Delta}_{S}\| \leq 2 \lambda \sqrt{k} \|\hat{\Delta}\| \leq 2 \lambda \sqrt{k} \frac{4}{\sqrt{n}} \|X_{I_j} \hat{\Delta}\|,
	\end{align*}
	where in the second inequality we have used Cauchy Schwarz inequality and in the final inequality we have used Lemma~\ref{lem:restev}. The above display implies the conclusions of this lemma.

	%{\color{red} Correct the constants here now.}
	
	%$2.17$ in~\cite{rigollet2015high}, stated in blah because $X_{I_j}$ satisfies $Inc(k).$ The last display implies the conclusions of the lemma and finishes the proof. 
	%{\color{red} State Rigollet Hutter's lemma about incoherence and restricted eigenvalues}
	
	%Now, let us choose $\lambda = C \sigma \sqrt{\log p} \sqrt{n}.$ Under this choice of $\lambda$, we can conclude that $\|X_{I_j} \hat{\Delta}\| \leq \sqrt{k \log p}.$ This further implies that $$\frac{\lambda}{2} |\hat{\Delta}|_1  \leq k \log p$$. Dividing by $\lambda/2$ both sides and plugging in the value of $\lambda$ finishes the proof. 

	%Using Holder's inequality and then using the choice of $\lambda$ in~\eqref{eq:lambdachoice} we obtain
	%\begin{equation*}
	%\sum_{i \in I_j} (y_i - x_i^{t} \hat{\beta}^{(\lambda, I_j)})^2 \leq 2 \sum_{i \in I_j} (x_i^{t} \hat{\beta}^{(\lambda, I_j)} - x_i^t \beta^*) \epsilon_i \:\: - \lambda |\hat{\beta}^{(\lambda, I_j)}|_1 + \lambda |\beta^*|_1.
	%\end{equation*}
	
\end{proof}

The following lemma is a version of a standard result (e.g, see Lemma $2.17$ in~\cite{rigollet2015high}) in the lasso literature which says that if the design matrix $X$ satisfies an incoherence condition and a vector satisfies a certain cone condition then the matrix $X$ satisfies a restricted eigenvalue condition. 
\begin{lemma}\label{lem:restev}
	Fix any $j \in \{1, 2\}$. Suppose $X_{I_j}$ satisfies $$\left|\frac{X_{I_j}^t X_{I_j}}{n} - \frac{1}{2} I_{p \times p}\right|_{\infty} \leq \frac{1}{64 k}.$$ Let $S \subset [p]$ such that $|S| = k.$
	Then for any $v \in \R^p$ satisfying $\|v\|_{1} \leq 4 \|v_{S}\|_1$ we have 
	\begin{equation*}
	\|v\| \leq \frac{4}{\sqrt{n}} \|X_{I_j} v\|
	\end{equation*}
	\end{lemma}

\begin{proof}
	Let us denote $W = X_{I}$ within this proof. We can write
	\begin{align*}
	&\frac{1}{n} \|Wv\|^2 = \sum_{i = 1}^{p} v_i^2 \frac{(W^t W)_{ii}}{n}  + \sum_{1 \leq i \neq j \leq p} v_i v_j \frac{(W^t W)_{ij}}{n} \geq \\& \left(\frac{1}{2} - \frac{1}{64 k}\right) \|v\|^2 - \frac{1}{64 k} \|v\|_1^2 \geq \left(\frac{1}{2} - \frac{1}{64 k}\right) \|v\|^2 -  \frac{1}{64 k} 16 \|v_{S}\|_1^2 \geq \\& \left(\frac{1}{2} - \frac{1}{64 k}\right) \|v\|^2 -  \frac{1}{64 k} 16 k \|v_{S}\|^2 \geq  \left(\frac{1}{2} - \frac{1}{64 k}\right) \|v\|^2 -  \frac{1}{4} \|v\|^2 \geq \left(\frac{1}{2} - \frac{1}{64k} - \frac{1}{4}\right) \|v\|^2.
	\end{align*}
	In the second inequality we used the fact that $\|v\|_{1} \leq 4 \|v_{S}\|_1$, in the next inequality we used the Cauchy Schwarz inequality and in the last inequality we used the fact that $k \geq 1.$ 
\end{proof}

\begin{lemma}\label{lem:incoh}
	%If the design matrix $X$ satisfies $Inc(2k)$ then for any $j \in [2],$ the matrix $X_{I_j}$ satisfies $Inc(k)$ with probability atleast $1 - p^2 \exp(-\frac{2n}{(128)^2 k^2 M^2}).$
	
	Let $I = I_j$ for $j \in \{1, 2\}.$ If $\; 2 \log p < \frac{n}{2^{14} k^2 M^2}$ then 
		\begin{align*}
	P\left(\left|\frac{X_{I}^t X_{I}}{n} - \frac{X^t X}{2n}\right|_{\infty} \geq 1/128 k\right) \leq \exp\left(-\frac{n}{2^{14} k^2 M^2}\right). 
	\end{align*}
\end{lemma}

\begin{proof}
	Note that the $(i,j)$th element of $\frac{X_{I}^t X_{I}}{n}$ can be written as follows:
	$$ \frac{1}{n} \sum_{l \in I} X_{li} X_{lj} =  \frac{1}{n} \sum_{l = 1}^{n} \mathrm{1}(l \in I) X_{li} X_{lj}.$$ 
	
	Therefore, by the standard tail bound for a subgaussian random variable we can write
	\begin{align*}
	P\left(\frac{1}{n} \sum_{l = 1}^{n} [\mathrm{1}(l \in I) - 1/2] X_{li} X_{lj} \geq t\right) \leq \exp\left(-2 n t^2/ M^2\right).
	\end{align*}

	Now by a union bound argument, 
	\begin{align*}
	P\left(\left|\frac{X_{I}^t X_{I}}{n} - \frac{X^t X}{2n}\right|_{\infty} \geq t\right) &=  P\left(\max_{i,j \in [p]} \frac{1}{n} \sum_{l = 1}^{n} [\mathrm{1}(l \in I) - 1/2] X_{li} X_{lj} \geq t\right)\\
	&\leq p^2 \exp(-2 n t^2/ M^2).
	\end{align*}
	
	By setting $t = \frac{1}{128 k}$ we finish the proof. %have 
	%\begin{align*}
	%P(|2 \frac{X_{I}^t X_{I}}{n} - \frac{X^t X}{n}|_{\infty} \geq 1/64 k) \leq \exp(-\frac{2n}{(128)^2 k^2 M^2} + 2 \log p) 
	%\end{align*}
	
	%{\color{red} Now decide how to simplify this writing.}
\end{proof}

\bibliographystyle{chicago}
\bibliography{references}

\def\noopsort#1{}
\begin{thebibliography}{}

\bibitem[\protect\citeauthoryear{Arnold, Tibshirani, Arnold, and
  ByteCompile}{Arnold et~al.}{2022}]{arnold2020package}
Arnold, T.~B., R.~J. Tibshirani, M.~T. Arnold, and T.~ByteCompile (2022).
\newblock Package genlasso.
\newblock {\em Statistics\/}~{\em 39\/}(3), 1335--1371.

\bibitem[\protect\citeauthoryear{Barron, Birg{\'e}, and Massart}{Barron
  et~al.}{1999}]{BarronBirgeMassart}
Barron, A., L.~Birg{\'e}, and P.~Massart (1999).
\newblock Risk bounds for model selection via penalization.
\newblock {\em Probab. Theory Related Fields\/}~{\em 113\/}(3), 301--413.

\bibitem[\protect\citeauthoryear{Blanchard, Sch{\"a}fer, Rozenholc, and
  M{\"u}ller}{Blanchard et~al.}{2007}]{blanchard2007optimal}
Blanchard, G., C.~Sch{\"a}fer, Y.~Rozenholc, and K.-R. M{\"u}ller (2007).
\newblock Optimal dyadic decision trees.
\newblock {\em Machine Learning\/}~{\em 66\/}(2-3), 209--241.

\bibitem[\protect\citeauthoryear{Boucheron, Lugosi, and Massart}{Boucheron
  et~al.}{2013}]{boucheron2013concentration}
Boucheron, S., G.~Lugosi, and P.~Massart (2013).
\newblock {\em Concentration inequalities: A nonasymptotic theory of
  independence}.
\newblock Oxford University Press.

\bibitem[\protect\citeauthoryear{Cai, Cand{\`e}s, and Shen}{Cai
  et~al.}{2010}]{cai2010singular}
Cai, J.-F., E.~J. Cand{\`e}s, and Z.~Shen (2010).
\newblock A singular value thresholding algorithm for matrix completion.
\newblock {\em SIAM Journal on optimization\/}~{\em 20\/}(4), 1956--1982.

\bibitem[\protect\citeauthoryear{Chatterjee}{Chatterjee}{2014}]{chatterjee2014new}
Chatterjee, S. (2014).
\newblock A new perspective on least squares under convex constraint.
\newblock {\em The Annals of Statistics\/}~{\em 42\/}(6), 2340--2381.

\bibitem[\protect\citeauthoryear{Chatterjee}{Chatterjee}{2015}]{chatterjee2015matrix}
Chatterjee, S. (2015).
\newblock Matrix estimation by universal singular value thresholding.
\newblock {\em Annals of Statistics\/}~{\em 43\/}(1), 177--214.

\bibitem[\protect\citeauthoryear{Chatterjee and Goswami}{Chatterjee and
  Goswami}{2019a}]{chatterjee2019adaptive}
Chatterjee, S. and S.~Goswami (2019a).
\newblock Adaptive estimation of multivariate piecewise polynomials and bounded
  variation functions by optimal decision trees.
\newblock {\em To appear in Annals of Statistics\/}.

\bibitem[\protect\citeauthoryear{Chatterjee and Goswami}{Chatterjee and
  Goswami}{2019b}]{chatterjee2019new}
Chatterjee, S. and S.~Goswami (2019b).
\newblock New risk bounds for 2d total variation denoising.
\newblock {\em arXiv preprint arXiv:1902.01215\/}.

\bibitem[\protect\citeauthoryear{Chatterjee and Jafarov}{Chatterjee and
  Jafarov}{2015}]{chatterjee2015prediction}
Chatterjee, S. and J.~Jafarov (2015).
\newblock Prediction error of cross-validated lasso.
\newblock {\em arXiv preprint arXiv:1502.06291\/}.

\bibitem[\protect\citeauthoryear{Chatterjee and Lafferty}{Chatterjee and
  Lafferty}{2019}]{chatterjeel2019adaptive}
Chatterjee, S. and J.~Lafferty (2019).
\newblock Adaptive risk bounds in unimodal regression.
\newblock {\em Bernoulli\/}~{\em 25\/}(1), 1--25.

\bibitem[\protect\citeauthoryear{Chatterjee and Mukherjee}{Chatterjee and
  Mukherjee}{2019}]{chatterjee2019estimation}
Chatterjee, S. and S.~Mukherjee (2019).
\newblock Estimation in tournaments and graphs under monotonicity constraints.
\newblock {\em IEEE Transactions on Information Theory\/}~{\em 65\/}(6),
  3525--3539.

\bibitem[\protect\citeauthoryear{Chetverikov, Liao, and
  Chernozhukov}{Chetverikov et~al.}{2020}]{chetverikov2020cross}
Chetverikov, D., Z.~Liao, and V.~Chernozhukov (2020).
\newblock On cross-validated lasso in high dimensions.
\newblock {\em Annal. Stat.(Forthcoming)\/}~{\em 40}.

\bibitem[\protect\citeauthoryear{Donoho and Gavish}{Donoho and
  Gavish}{2014}]{donoho2014minimax}
Donoho, D. and M.~Gavish (2014).
\newblock Minimax risk of matrix denoising by singular value thresholding.
\newblock {\em Annals of Statistics\/}~{\em 42\/}(6), 2413--2440.

\bibitem[\protect\citeauthoryear{Donoho}{Donoho}{1997}]{donoho1997cart}
Donoho, D.~L. (1997).
\newblock {CART} and best-ortho-basis: a connection.
\newblock {\em The Annals of Statistics\/}~{\em 25\/}(5), 1870--1911.

\bibitem[\protect\citeauthoryear{Fang, Guntuboyina, and Sen}{Fang
  et~al.}{2021}]{fang2021multivariate}
Fang, B., A.~Guntuboyina, and B.~Sen (2021).
\newblock Multivariate extensions of isotonic regression and total variation
  denoising via entire monotonicity and hardy--krause variation.
\newblock {\em The Annals of Statistics\/}~{\em 49\/}(2), 769--792.

\bibitem[\protect\citeauthoryear{Golub, Heath, and Wahba}{Golub
  et~al.}{1979}]{golub1979generalized}
Golub, G.~H., M.~Heath, and G.~Wahba (1979).
\newblock Generalized cross-validation as a method for choosing a good ridge
  parameter.
\newblock {\em Technometrics\/}~{\em 21\/}(2), 215--223.

\bibitem[\protect\citeauthoryear{Greenshtein, Ritov, et~al.}{Greenshtein
  et~al.}{2004}]{greenshtein2004persistence}
Greenshtein, E., Y.~Ritov, et~al. (2004).
\newblock Persistence in high-dimensional linear predictor selection and the
  virtue of overparametrization.
\newblock {\em Bernoulli\/}~{\em 10\/}(6), 971--988.

\bibitem[\protect\citeauthoryear{Guntuboyina, Lieu, Chatterjee, and
  Sen}{Guntuboyina et~al.}{2020}]{guntuboyina2020adaptive}
Guntuboyina, A., D.~Lieu, S.~Chatterjee, and B.~Sen (2020).
\newblock Adaptive risk bounds in univariate total variation denoising and
  trend filtering.
\newblock {\em The Annals of Statistics\/}~{\em 48\/}(1), 205--229.

\bibitem[\protect\citeauthoryear{Hernan and Chatterjee}{Hernan and
  Chatterjee}{2021}]{hernan2021risk}
Hernan, O. and S.~Chatterjee (2021).
\newblock Risk bounds for quantile trend filtering.
\newblock {\em Biometrika\/}.

\bibitem[\protect\citeauthoryear{Hoefling}{Hoefling}{2010}]{hoefling2010path}
Hoefling, H. (2010).
\newblock A path algorithm for the fused lasso signal approximator.
\newblock {\em Journal of Computational and Graphical Statistics\/}~{\em
  19\/}(4), 984--1006.

\bibitem[\protect\citeauthoryear{Homrighausen and McDonald}{Homrighausen and
  McDonald}{2013}]{homrighausen2013lasso}
Homrighausen, D. and D.~McDonald (2013).
\newblock The lasso, persistence, and cross-validation.
\newblock In {\em International Conference on Machine Learning}, pp.\
  1031--1039. PMLR.

\bibitem[\protect\citeauthoryear{Homrighausen and McDonald}{Homrighausen and
  McDonald}{2014}]{homrighausen2014leave}
Homrighausen, D. and D.~J. McDonald (2014).
\newblock Leave-one-out cross-validation is risk consistent for lasso.
\newblock {\em Machine learning\/}~{\em 97\/}(1), 65--78.

\bibitem[\protect\citeauthoryear{Homrighausen and McDonald}{Homrighausen and
  McDonald}{2017}]{homrighausen2017risk}
Homrighausen, D. and D.~J. McDonald (2017).
\newblock Risk consistency of cross-validation with lasso-type procedures.
\newblock {\em Statistica Sinica\/}, 1017--1036.

\bibitem[\protect\citeauthoryear{H{\"u}tter and Rigollet}{H{\"u}tter and
  Rigollet}{2016}]{hutter2016optimal}
H{\"u}tter, J.-C. and P.~Rigollet (2016).
\newblock Optimal rates for total variation denoising.
\newblock In {\em Conference on Learning Theory}, pp.\  1115--1146.

\bibitem[\protect\citeauthoryear{Ki, Fang, and Guntuboyina}{Ki
  et~al.}{2021}]{ki2021mars}
Ki, D., B.~Fang, and A.~Guntuboyina (2021).
\newblock Mars via lasso.
\newblock {\em arXiv preprint arXiv:2111.11694\/}.

\bibitem[\protect\citeauthoryear{Kim, Koh, Boyd, and Gorinevsky}{Kim
  et~al.}{2009}]{kim2009ell_1}
Kim, S.-J., K.~Koh, S.~Boyd, and D.~Gorinevsky (2009).
\newblock {$\ell_1$} trend filtering.
\newblock {\em SIAM Rev.\/}~{\em 51\/}(2), 339--360.

\bibitem[\protect\citeauthoryear{Koltchinskii}{Koltchinskii}{2011}]{koltchinskii2011oracle}
Koltchinskii, V. (2011).
\newblock {\em Oracle Inequalities in Empirical Risk Minimization and Sparse
  Recovery Problems: Ecole d’Et{\'e} de Probabilit{\'e}s de Saint-Flour
  XXXVIII-2008}, Volume 2033.
\newblock Springer Science \& Business Media.

\bibitem[\protect\citeauthoryear{Lecu{\'e}, Mitchell, et~al.}{Lecu{\'e}
  et~al.}{2012}]{lecue2012oracle}
Lecu{\'e}, G., C.~Mitchell, et~al. (2012).
\newblock Oracle inequalities for cross-validation type procedures.
\newblock {\em Electronic Journal of Statistics\/}~{\em 6}, 1803--1837.

\bibitem[\protect\citeauthoryear{Li and Barron}{Li and Barron}{2000}]{LiBarron}
Li, J.~Q. and A.~R. Barron (2000).
\newblock Mixture density estimation.
\newblock In {\em Advances in Neural Information Processing Systems},
  Volume~12, pp.\  279--285. Cambridge, MA: MIT Press.

\bibitem[\protect\citeauthoryear{Li and Racine}{Li and
  Racine}{2004}]{li2004cross}
Li, Q. and J.~Racine (2004).
\newblock Cross-validated local linear nonparametric regression.
\newblock {\em Statistica Sinica\/}, 485--512.

\bibitem[\protect\citeauthoryear{Mammen and van~de Geer}{Mammen and van~de
  Geer}{1997}]{mammen1997locally}
Mammen, E. and S.~van~de Geer (1997).
\newblock Locally adaptive regression splines.
\newblock {\em The Annals of Statistics\/}~{\em 25\/}(1), 387--413.

\bibitem[\protect\citeauthoryear{Miolane and Montanari}{Miolane and
  Montanari}{2018}]{miolane2018distribution}
Miolane, L. and A.~Montanari (2018).
\newblock The distribution of the lasso: Uniform control over sparse balls and
  adaptive parameter tuning.
\newblock {\em arXiv preprint arXiv:1811.01212\/}.

\bibitem[\protect\citeauthoryear{Nemirovski}{Nemirovski}{2000}]{nemirovski2000topics}
Nemirovski, A. (2000).
\newblock Topics in non-parametric statistics.
\newblock {\em Lectures on probability theory and statistics (Saint-Flour,
  1998)\/}~{\em 1738}, 85--277.

\bibitem[\protect\citeauthoryear{Nowak, Mitra, and Willett}{Nowak
  et~al.}{2004}]{nowak2004estimating}
Nowak, R., U.~Mitra, and R.~Willett (2004).
\newblock Estimating inhomogeneous fields using wireless sensor networks.
\newblock {\em IEEE Journal on Selected Areas in Communications\/}~{\em
  22\/}(6), 999--1006.

\bibitem[\protect\citeauthoryear{Ortelli and van~de Geer}{Ortelli and van~de
  Geer}{2019}]{ortelli2019prediction}
Ortelli, F. and S.~van~de Geer (2019).
\newblock Prediction bounds for (higher order) total variation regularized
  least squares.
\newblock {\em arXiv preprint arXiv:1904.10871\/}.

\bibitem[\protect\citeauthoryear{Ortelli and van~de Geer}{Ortelli and van~de
  Geer}{2020}]{ortelli2020adaptive}
Ortelli, F. and S.~van~de Geer (2020).
\newblock Adaptive rates for total variation image denoising.
\newblock {\em Journal of Machine Learning Research\/}~{\em 21}, 247.

\bibitem[\protect\citeauthoryear{Ortelli and van~de Geer}{Ortelli and van~de
  Geer}{2021}]{ortelli2021oracle}
Ortelli, F. and S.~van~de Geer (2021).
\newblock Oracle inequalities for square root analysis estimators with
  application to total variation penalties.
\newblock {\em Information and Inference: A Journal of the IMA\/}~{\em
  10\/}(2), 483--514.

\bibitem[\protect\citeauthoryear{Padilla and Chatterjee}{Padilla and
  Chatterjee}{2021}]{padilla2021quantile}
Padilla, O. H.~M. and S.~Chatterjee (2021).
\newblock Quantile regression by dyadic cart.
\newblock {\em arXiv preprint arXiv:2110.08665\/}.

\bibitem[\protect\citeauthoryear{Politsch, Cisewski-Kehe, Croft, and
  Wasserman}{Politsch et~al.}{2020}]{politsch2020trend}
Politsch, C.~A., J.~Cisewski-Kehe, R.~A. Croft, and L.~Wasserman (2020).
\newblock Trend filtering--i. a modern statistical tool for time-domain
  astronomy and astronomical spectroscopy.
\newblock {\em Monthly Notices of the Royal Astronomical Society\/}~{\em
  492\/}(3), 4005--4018.

\bibitem[\protect\citeauthoryear{Rigollet and H{\"u}tter}{Rigollet and
  H{\"u}tter}{2015}]{rigollet2015high}
Rigollet, P. and J.-C. H{\"u}tter (2015).
\newblock High dimensional statistics.
\newblock {\em Lecture notes for course 18S997\/}.

\bibitem[\protect\citeauthoryear{Rudin, Osher, and Fatemi}{Rudin
  et~al.}{1992}]{rudin1992nonlinear}
Rudin, L.~I., S.~Osher, and E.~Fatemi (1992).
\newblock Nonlinear total variation based noise removal algorithms.
\newblock {\em Physica D: Nonlinear Phenomena\/}~{\em 60\/}(1), 259--268.

\bibitem[\protect\citeauthoryear{Sadhanala, Wang, and Tibshirani}{Sadhanala
  et~al.}{2016}]{sadhanala2016total}
Sadhanala, V., Y.-X. Wang, and R.~J. Tibshirani (2016).
\newblock Total variation classes beyond 1d: Minimax rates, and the limitations
  of linear smoothers.
\newblock In {\em Advances in Neural Information Processing Systems}, pp.\
  3513--3521.

\bibitem[\protect\citeauthoryear{Scott and Nowak}{Scott and
  Nowak}{2006}]{scott2006minimax}
Scott, C. and R.~D. Nowak (2006).
\newblock Minimax-optimal classification with dyadic decision trees.
\newblock {\em IEEE transactions on information theory\/}~{\em 52\/}(4),
  1335--1353.

\bibitem[\protect\citeauthoryear{Shah, Balakrishnan, Guntuboyina, and
  Wainwright}{Shah et~al.}{2016}]{shah2016stochastically}
Shah, N., S.~Balakrishnan, A.~Guntuboyina, and M.~Wainwright (2016).
\newblock Stochastically transitive models for pairwise comparisons:
  Statistical and computational issues.
\newblock In {\em International Conference on Machine Learning}, pp.\  11--20.

\bibitem[\protect\citeauthoryear{Shao}{Shao}{1993}]{shao1993linear}
Shao, J. (1993).
\newblock Linear model selection by cross-validation.
\newblock {\em Journal of the American statistical Association\/}~{\em
  88\/}(422), 486--494.

\bibitem[\protect\citeauthoryear{Tibshirani}{Tibshirani}{1996}]{tibshirani1996regression}
Tibshirani, R. (1996).
\newblock Regression shrinkage and selection via the lasso.
\newblock {\em Journal of the Royal Statistical Society. Series B
  (Methodological)\/}, 267--288.

\bibitem[\protect\citeauthoryear{Tibshirani}{Tibshirani}{2014}]{tibshirani2014adaptive}
Tibshirani, R.~J. (2014).
\newblock Adaptive piecewise polynomial estimation via trend filtering.
\newblock {\em The Annals of Statistics\/}~{\em 42\/}(1), 285--323.

\bibitem[\protect\citeauthoryear{Tibshirani}{Tibshirani}{2020}]{tibshirani2020divided}
Tibshirani, R.~J. (2020).
\newblock Divided differences, falling factorials, and discrete splines:
  Another look at trend filtering and related problems.
\newblock {\em arXiv preprint arXiv:2003.03886\/}.

\bibitem[\protect\citeauthoryear{Tibshirani and Taylor}{Tibshirani and
  Taylor}{2011}]{tibshirani2011solution}
Tibshirani, R.~J. and J.~Taylor (2011).
\newblock The solution path of the generalized lasso.
\newblock {\em The annals of statistics\/}~{\em 39\/}(3), 1335--1371.

\bibitem[\protect\citeauthoryear{van~de Geer and Ortelli}{van~de Geer and
  Ortelli}{2019}]{van2019prediction}
van~de Geer, S. and F.~Ortelli (2019).
\newblock Prediction bounds for (higher order) total variation regularized
  least squares.
\newblock {\em arXiv preprint arXiv:1904.10871\/}.

\bibitem[\protect\citeauthoryear{van Handel}{van
  Handel}{2014}]{van2014probability}
van Handel, R. (2014).
\newblock Probability in high dimension.
\newblock Technical report, PRINCETON UNIV NJ.

\bibitem[\protect\citeauthoryear{Vershynin}{Vershynin}{2018}]{vershynin2018high}
Vershynin, R. (2018).
\newblock {\em High-dimensional probability: An introduction with applications
  in data science}, Volume~47.
\newblock Cambridge university press.

\bibitem[\protect\citeauthoryear{Wainwright}{Wainwright}{2019}]{wainwright2019high}
Wainwright, M.~J. (2019).
\newblock {\em High-dimensional statistics: A non-asymptotic viewpoint},
  Volume~48.
\newblock Cambridge University Press.

\bibitem[\protect\citeauthoryear{Wang, Sharpnack, Smola, and Tibshirani}{Wang
  et~al.}{2016}]{wang2016trend}
Wang, Y.-X., J.~Sharpnack, A.~Smola, and R.~J. Tibshirani (2016).
\newblock Trend filtering on graphs.
\newblock {\em Journal of Machine Learning Research\/}~{\em 17\/}(105), 1--41.

\bibitem[\protect\citeauthoryear{Wang, Smola, and Tibshirani}{Wang
  et~al.}{2014}]{wang2014falling}
Wang, Y.-X., A.~J. Smola, and R.~J. Tibshirani (2014).
\newblock The falling factorial basis and its statistical applications.
\newblock In {\em ICML}, pp.\  730--738.

\bibitem[\protect\citeauthoryear{Willett and Nowak}{Willett and
  Nowak}{2007}]{willett2007multiscale}
Willett, R.~M. and R.~D. Nowak (2007).
\newblock Multiscale poisson intensity and density estimation.
\newblock {\em IEEE Transactions on Information Theory\/}~{\em 53\/}(9),
  3171--3187.

\bibitem[\protect\citeauthoryear{Wong}{Wong}{1983}]{wong1983consistency}
Wong, W.~H. (1983).
\newblock On the consistency of cross-validation in kernel nonparametric
  regression.
\newblock {\em The Annals of Statistics\/}~{\em 11\/}(4), 1136--1141.

\end{thebibliography}
\def\noopsort#1{}

\end{document}